\documentclass[10pt,a4paper,reqno]{amsart}
\usepackage{color}
\usepackage{verbatim}
\usepackage{mathtools}
\usepackage{amstext}
\usepackage{amsthm}
\usepackage{amssymb,mathrsfs}
\usepackage{kotex}
\usepackage[unicode=true,bookmarks=false,breaklinks=false,pdfborder={0 0 1},colorlinks=false]{hyperref}
\hypersetup{
 colorlinks,linkcolor={red!60!black},citecolor={green!60!black},urlcolor={blue!60!black}}

\makeatletter
\numberwithin{equation}{section}
\numberwithin{figure}{section}
\theoremstyle{plain}

\usepackage{lineno}
\usepackage{footnote}
\usepackage{enumitem}
\usepackage{setspace}
\usepackage{tikz}
\usepackage{cite}
\makesavenoteenv{tblr}

\newtheorem{theorem}{Theorem}[section]
\newtheorem{proposition}[theorem]{Proposition}
\newtheorem{lemma}[theorem]{Lemma}

\newtheorem{corollary}[theorem]{Corollary}

\theoremstyle{remark}
\newtheorem{remark}[theorem]{\protect\remarkname}

\providecommand{\remarkname}{Remark}

\global\long\def\Im{\mathrm{Im}}%
\global\long\def\Re{\mathrm{Re}}%
\global\long\def\norm#1{\|#1\|}%
\global\long\def\tilde#1{\widetilde{#1}}%
\global\long\def\bbR{\mathbf{\mathbb{R}}}%
\global\long\def\bbC{\mathbf{\mathbb{C}}}%
\global\long\def\bbN{\mathbf{\mathbb{N}}}%
\global\long\def\calH{\mathcal{H}}%
\global\long\def\gmm{\gamma}%

\newcommand{\Ci}{\mathrm{Ci}}
\newcommand{\Si}{\mathrm{Si}}
\DeclareMathOperator{\Log}{Log}

\makeatother

\providecommand{\remarkname}{Remark}

\theoremstyle{definition}
\newtheorem{definition}[theorem]{Definition}

\title[Blow-up dynamics for the mass-critical Half-wave equation]{Finite-time blow-up for the mass-critical half-wave equation with negative energy}
\author{Jeongheon Park}
\date{}


\begin{document}

\begin{abstract}
    We study the one-dimensional focusing mass-critical half-wave equation
    \[
        i\partial_tu=|D|u-|u|^2u.
    \]
    For even initial data with negative energy and mass slightly above the ground-state mass, we prove finite-time blow-up and obtain the upper bound
    \[
        \norm{u(t)}_{\dot H^{1/2}}
        \lesssim
        \frac{|\log(T-t)|^{1/4}}{\sqrt{T-t}}
        \qquad\text{as }t\uparrow T.
    \]
    This is the first finite-time blow-up result for negative-energy solutions to the mass-critical half-wave equation in the near-ground-state regime.
    
    The proof uses the modulation analysis developed by Merle--Rapha\"el \cite{MerleRaphael2005AnnMath}. For the half-wave equation, the key missing ingredient is a coercivity estimate for a nonlocal quadratic form generated by the scaling direction. Unlike the local NLS, the half-wave equation does not admit the ODE methods used to establish the corresponding coercivity estimate. Instead, we prove the coercivity estimate by an analytic reduction followed by a rigorous computer-assisted proof.
    
    The spectral analysis part reduces the coercivity problem to a finite collection of spectral inequalities by combining constrained Morse index arguments with the Birman--Schwinger principle. The computer-assisted part certifies these inequalities by interval arithmetic using a validated approximation of the ground state obtained by compactifying the real line. Together with the modulation analysis, this coercivity theorem gives the finite-time blow-up result.
\end{abstract}

\maketitle

\tableofcontents

\setcounter{tocdepth}{1}

\section{Introduction}
We consider the \emph{mass-critical one-dimensional focusing half-wave equation}
\begin{equation}\tag{HW}\label{eq:halfwave}
    \begin{cases}
        i\partial_t u = |D|u - |u|^2u, \qquad (t,x)\in \bbR\times\bbR,\\
        u(0,x)=u_0(x),
    \end{cases}
\end{equation}
where $|D|=|\partial_x|$ is the Fourier multiplier operator with symbol $|\xi|$. This equation appears in several physical settings, including continuum limits of lattice models \cite{KLS2013CMP}, wave turbulence models \cite{Physics1997JNS,Physics2001PhysD}, and gravitational collapse \cite{Physics2007CPAM,FrohlichLenzmann2007CPAM}.

The half-wave equation \eqref{eq:halfwave} may be viewed as a close relative of the mass-critical nonlinear Schr\"odinger equation, as the two models share the same $L^2$-critical scaling. The blow-up dynamics of the mass-critical NLS have been studied extensively, aided by its rich structure, including pseudo-conformal, Galilean, and virial identities. By contrast, the blow-up theory for \eqref{eq:halfwave} remains largely underdeveloped. The half-wave equation lacks Lorentz, pseudo-conformal, and Galilean symmetries, while the nonlocality of the operator $|D|$ further obstructs the direct use of standard NLS techniques. In particular, the Glassey virial argument does not convert negative energy into finite-time blow-up for \eqref{eq:halfwave}. Consequently, even some of the most basic questions concerning its blow-up dynamics remain open.

The present paper proves that near-soliton solutions with negative energy blow up in finite time under even symmetry. We also obtain an upper bound on the blow-up rate. This is analogous to NLS results developed by Merle and Rapha\"el \cite{MerleRaphael2005AnnMath,MerleRaphael2003GAFA,Raphael2005MathAnnalen,MerleRaphael2004Invent,MerleRaphael2006JAMS,MerleRaphael2005CMP}. The dynamical part follows the modulation analysis and local virial argument of Merle--Rapha\"el. By contrast, establishing the spectral property for \eqref{eq:halfwave} requires substantial input from spectral theory together with a computer-assisted proof.

\subsection{Symmetries and solitons}
Equation \eqref{eq:halfwave} is invariant under the $L^2$-critical scaling
\[
    u_\lambda(t,x)=\lambda^{\frac12}u(\lambda t,\lambda x), \qquad \lambda>0,
\]
which preserves the $L^2$-norm. It also enjoys space-time translation and phase rotation symmetries:
\begin{equation*}
    u(t,x)\mapsto u(t+t_0,x+x_0), \qquad
    u(t,x)\mapsto e^{i\gamma}u(t,x).
\end{equation*}
The corresponding conserved quantities are mass, momentum, and energy:
\begin{equation}\label{eq:conservation law}
    \begin{aligned}
        M(u)&=\int_{\bbR}|u(t,x)|^2\,dx,\\
        P(u)&=\int_{\bbR} i\overline{u}(t,x)\,\partial_xu(t,x)\,dx,\\
        E(u)&=\frac12\int_{\bbR}||D|^{1/2}u(t,x)|^2\,dx -\frac14\int_{\bbR}|u(t,x)|^4\,dx.
    \end{aligned}
\end{equation}
We call \eqref{eq:halfwave} mass-critical, or $L^2$-critical, because the scaling symmetry preserves the mass.

The equation \eqref{eq:halfwave} admits the stationary solution $e^{it}Q$, which represents a solitary wave, or soliton. Here, the profile $Q$ solves
\begin{equation*}
    |D|Q+Q-Q^3=0.
\end{equation*}
The positive even solution \(Q\), which is an optimizer of the sharp Gagliardo--Nirenberg inequality associated with \eqref{eq:halfwave}, is called the ground state and plays a pivotal role in the dynamics of \eqref{eq:halfwave}. Frank and Lenzmann \cite{FrankLenzmann2013Acta} proved the uniqueness of the ground state $Q\in H^{1/2}(\bbR)$, up to translation, and also established the nondegeneracy of the associated linearized operator. In particular, the sharp inequality yields global existence below the mass $M(Q)$, so that $M(Q)$ appears as the natural threshold mass in the dynamics of \eqref{eq:halfwave}.

It is noteworthy that \eqref{eq:halfwave} also admits traveling waves below the threshold. For $|v|<1$, there are solutions of the form $e^{it}Q_v(x-vt)$, where
\begin{equation*}
    |D|Q_v+Q_v-|Q_v|^2Q_v+iv\partial_xQ_v=0.
\end{equation*}
Krieger, Lenzmann, and Rapha\"el \cite{KLR2013ARMAhalfwave} showed that $Q_v\to Q$ as $v\to0$, while $M(Q_v)\to0$ as $|v|\to1$. Thus, unlike for the mass-critical NLS, solitary waves exist below the ground-state mass $M(Q)$; indeed, moving solitary waves persist with arbitrarily small mass.

\subsection{Previous results}
We briefly recall some previous results for \eqref{eq:halfwave}. The dynamical theory is still far from complete, but there are a few results. Krieger, Lenzmann, and Rapha\"el \cite{KLR2013ARMAhalfwave} constructed minimal mass ($M(u)=M(Q)$) blow-up solutions in one dimension, with blow-up rate $\norm{|D|^{1/2}u(t)}_{L^2}\sim |t|^{-1}$ as $t\to 0$.
See also \cite{Georgiev2Dhalfwaveblowup,Georgiev3Dhalfwaveblowup} for related constructions in higher dimensions.

A subthreshold two-soliton solution with transient turbulent behavior was constructed in \cite{GLPR2018annPDE}: its high Sobolev norms grow and then enter a long-time saturation phase. Slightly above the threshold, Bourgain--Wang type solutions with the same blow-up rate as the minimal blow-up solution were constructed by Kim, Kwon, and the author \cite{KimKwonPark2025arXiv}, and they also established their instability.

Other aspects of the half-wave equation have also been studied. We mention weak stability of multi-soliton sums \cite{Li2024WeakStability}, multi-bubble blow-up \cite{CaoSuZhang2022}, inhomogeneous mass-critical blow-up \cite{Li2022Inhomogeneous}, non-dispersive traveling waves and the failure of small-data scattering \cite{BellazziniGeorgievVisciglia2018}, and ill-posedness results for the cubic half-wave equation and related fractional NLS models \cite{ChoffrutPocovnicu2018}. Related negative-energy blow-up alternatives for $L^2$-critical fractional NLS models were obtained in \cite{BHL2016JFA,Dinh2019Nonlinearity}. In the one-dimensional $L^2$-critical fractional NLS, Lan \cite{Lan2022IMRN} proved finite-time blow-up and described the near-soliton dynamics for $\beta$ sufficiently close to $2$, using a perturbative argument to transfer the spectral property from the classical NLS. For the half-wave maps equation, see \cite{GL2018LMP,GL2025sigma,GL2026arXiv} and references therein.

\subsection{Main results}
We now state the dynamical result of the paper, which gives finite-time
blow-up in the negative-energy regime for \eqref{eq:halfwave}.
\begin{theorem}\label{thm:main}
    There exist constants $\alpha^*>0$ and $C^*>0$ with the following property.
    Let $u_0\in H^{1/2}_{\mathrm e}(\mathbb{R})$ and assume
    \begin{equation}\label{eq:main-assumptions}
        \begin{gathered}
            0<\alpha_0=\alpha(u_0):=\norm{u_0}_{L^2}^2-\norm{Q}_{L^2}^2<\alpha^*,\quad 
            E_0:=E(u_0)<0.
        \end{gathered}
    \end{equation}
    Let $u(t)$ be the corresponding solution to \eqref{eq:halfwave} on its maximal forward lifespan $[0,T)$. Then $T<+\infty$, and $u(t)$ blows up at time $T$. Moreover, as $t \to T$,
    \begin{equation}\label{eq:main-upper-rate}
        \norm{u(t)}_{\dot H^{1/2}}
        \leq C^*\frac{|\log(T-t)|^{1/4}}{\sqrt{T-t}}.
    \end{equation}
\end{theorem}
Theorem~\ref{thm:main} gives the first finite-time blow-up result for negative-energy solutions to \eqref{eq:halfwave} in the near-soliton regime. In view of \eqref{eq:main-upper-rate}, the blow-up solutions described here are of a different nature from the previously known finite-time blow-up solutions in \cite{KLR2013ARMAhalfwave,KimKwonPark2025arXiv}. The bound \eqref{eq:main-upper-rate} is not expected to be sharp. In the subsequent work by Kim, Kwon, and the author \cite{KimKwonPark2026arXiv}, it is improved to the log-log upper bound
\begin{equation*}
    \norm{u(t)}_{\dot H^{1/2}}
    \lesssim
    \left(\frac{\log|\log(T-t)|}{T-t}\right)^{1/2}.
\end{equation*}
Building on the foundation of the present paper, \cite{KimKwonPark2026arXiv} further refines the approximate profile and modulation analysis to obtain the log-log upper bound.

We now describe the spectral ingredient used in the modulation analysis. 
After decomposing the solution into the modulated $Q$, with radiation denoted by $\epsilon$, we use a local virial identity adapted to the scaling direction. Let
\begin{equation*}
	\Lambda=\frac12+y\partial_y,
	\qquad
	L_Qf=|D|f+f-Q^2f-2\Re\{Q^2 f\}.
\end{equation*}
Here $\Lambda$ is the generator of the $L^2$-critical scaling, and $L_Q$ is the $\bbR$-linear linearized operator around $Q$. In the local virial computation \eqref{eq:proof strategy; motivation for H}, the leading quadratic terms in the radiation are expressed through the commutator form
\begin{equation*}
	\mathbf H(\epsilon,\epsilon)=([L_Q,\Lambda]\epsilon,\epsilon)_r .
\end{equation*}
We use the weighted Sobolev norm
\begin{equation*}
    \norm{f}_{\mathcal H^{1/2}}^2\coloneqq
    \norm{f}_{\dot H^{1/2}}^2+\norm{\langle y\rangle^{-2}f}_{L^2}^2.
\end{equation*}
The weighted $L^2$ component reflects the profile weight $|yQQ'|\lesssim\langle y\rangle^{-4}$ appearing in the local virial computation. The following theorem gives the coercivity of $\mathbf H$ in the form needed for the proof of Theorem~\ref{thm:main}.

\begin{theorem}[Spectral structure of the bilinear form $\mathbf H$]\label{thm:Spectral structure of the bilinear form H}
    There exist universal constants $\delta_0, C>0$ such that for all $\epsilon\in H^{1/2}_{\mathrm e}$,
    \begin{equation}\label{positivity of H}
        \mathbf H(\epsilon,\epsilon) \geq \delta_0 \|\epsilon\|_{\calH^{1/2}}^2-C((\epsilon,Q)_r^2+(\epsilon,\Lambda Q)_r^2+(\epsilon, i\Lambda Q)_r^2+ (\epsilon,i\Lambda^2 Q)_r^2).
    \end{equation}
\end{theorem}
In Part~\ref{part:blow-up-dynamics}, we assume the coercivity estimate \eqref{positivity of H} and use it to prove Theorem~\ref{thm:main}. This estimate controls the radiation term in the local virial argument, while the finite-dimensional terms on the right-hand side are handled through the modulation analysis. The remainder of the paper, beginning with Part~\ref{part:spectral-reduction}, is devoted to proving Theorem~\ref{thm:Spectral structure of the bilinear form H}.

\begin{remark}[Method and novelties] 
    The proof of Theorem~\ref{thm:main} follows the modulation analysis and
    local virial method developed by Merle and Rapha\"el for the
    mass-critical NLS \cite{MerleRaphael2005AnnMath,MerleRaphael2003GAFA,Raphael2005MathAnnalen,MerleRaphael2004Invent,MerleRaphael2006JAMS,MerleRaphael2005CMP}. In the half-wave setting, the local virial identity leads to the quadratic form \(\mathbf H\). The new input is the coercivity estimate in Theorem~\ref{thm:Spectral structure of the bilinear form H} and its proof is the main contribution of the paper. 
    We prove this theorem by reducing the coercivity problem to a computable finite-dimensional problem through a spectral analytic reduction and then verifying the reduced problem by a computer-assisted proof.

    \emph{Reduction of the spectral problem.}
    Unlike the mass-critical NLS, the half-wave equation involves the nonlocal operator $|D|$, so the standard ODE arguments used in the NLS setting are no longer available. Instead, we reduce the coercivity of $\mathbf H$ to a finite collection of explicit inequalities by combining the Birman--Schwinger principle with the symmetric resolvent identity. A further difficulty is that the symbol $|\xi|$ is singular at the zero frequency, preventing a direct application of this strategy. To overcome this issue, we introduce a positive regularization parameter $\mu$, which removes the zero-frequency singularity while preserving the relevant spectral information.

    \emph{Approximation of $Q$.}
    The finite-dimensional inequalities obtained above still depend on the ground state $Q$, so a sufficiently accurate explicit approximation is required for rigorous verification. We construct such an approximation by combining a Newton--Kantorovich argument with the uniqueness theorem due to Frank--Lenzmann--Silvestre \cite{FrankLenzmannSilvestre2016CPAM}, which yields a quantitative bound for $\mathfrak g-Q$. A major additional difficulty is that, since $|D|$ is nonlocal, approximating $Q$ requires controlling its algebraically decaying tail. To overcome this difficulty, we construct the approximate profile $\mathfrak g$ by \emph{compactifying the real line} and factoring out the algebraic decay. This compactification procedure is one of the main new ingredients of the paper.
\end{remark}

\begin{remark}[Comparison with mass-critical NLS]
    We recall related results for the $L^2$-critical nonlinear Schr\"odinger equation
    \begin{equation}\label{eq:NLS}\tag{NLS}
        i \partial_t u + \Delta u + |u|^\frac{4}{d} u = 0,\quad  (t,x) \in I \times \bbR^d.
    \end{equation}
    The negative-energy blow-up dynamics near the soliton for \eqref{eq:NLS} are stable and governed by the log-log blow-up rate $\|u(t)\|_{\dot H^1}^{-1}\sim \sqrt{2\pi(T-t)/\log|\log(T-t)|}$.
    These blow-up dynamics in the near-soliton regime were developed in the celebrated works of Merle and Rapha\"el  \cite{MerleRaphael2005AnnMath,MerleRaphael2003GAFA,Raphael2005MathAnnalen,MerleRaphael2004Invent,MerleRaphael2006JAMS,MerleRaphael2005CMP} by a local virial argument. More recently, Kwak--Kwon \cite{KwakKwon2026arXiv} proved the aforementioned log-log lower bound for \textit{general} radial solutions. 

    The local virial estimate for \eqref{eq:NLS} requires a spectral property for two Schr\"odinger-type operators. For \(d=1\), Merle--Rapha\"el \cite{MerleRaphael2005AnnMath} use the explicit formula for \(Q\) to prove the spectral property. For \(d=2,3,4\), Fibich--Merle--Rapha\"el \cite{FibichMerleRaphel2006PhisicaD} give a computer-assisted proof. They use the Sturm oscillation theorem and an ODE shooting method. For \(d=5,\cdots,10\), and in the radial case for \(d=11,12\), Yang--Roudenko--Zhao \cite{YangRoudenkoZhao2018Nonlinearity} give a computer-assisted proof using similar methods. These ODE-based arguments are not available for the half-wave equation, since the corresponding linearized operators involve the nonlocal operator $|D|$. In particular, the spectral problem cannot be reduced directly to a Sturm oscillation problem.
\end{remark}

\begin{remark}[Outside even symmetry] It is possible to extend our result to general solutions. 
    Since $Q$ is even, $[L_Q,\Lambda]$ preserves parity. Writing $\epsilon=\epsilon^e+\epsilon^o$ for the even and odd parts of $\epsilon$,
    \[
        \mathbf H(\epsilon,\epsilon)=\mathbf H(\epsilon^e,\epsilon^e)+\mathbf H(\epsilon^o,\epsilon^o).
    \]
    Theorem~\ref{thm:Spectral structure of the bilinear form H} gives the coercivity estimate for the even part.

    Let \(G_1\) be the odd function satisfying \(L_Q[iG_1]=-iQ'\). For the odd part, the expected estimate is
    \begin{equation}\label{eq:odd coercivity intro}
        \mathbf H(\epsilon^o,\epsilon^o)
        \geq
        \delta_{\mathrm o}
        \norm{\epsilon^o}_{\mathcal H^{1/2}}^2
        -
        C_{\mathrm o}
        (\epsilon,G_1)_r^2.
    \end{equation}

   Without the evenness assumption, the modulation analysis requires the translation \(x(t)\) and the velocity parameter \(v(t)\). If \eqref{eq:odd coercivity intro} were proved, the modulation argument could be adapted to \(\beta=1\) using a \(v\)-dependent approximate blow-up profile, as in \cite[Sections~2--4]{Lan2022IMRN}. This would remove the evenness assumption from Theorem~\ref{thm:main}.
\end{remark}

\subsection{Strategy of the proof}
\mbox{}\par\nobreak\medskip

As explained above, the proof is divided into two parts. In Part~\ref{part:blow-up-dynamics}, we assume Theorem~\ref{thm:Spectral structure of the bilinear form H} and prove Theorem~\ref{thm:main} through a modulation analysis. This part follows the overall roadmap developed by Merle and Rapha\"el for the mass-critical NLS \cite{MerleRaphael2005AnnMath}. We therefore focus below on Parts~\ref{part:spectral-reduction}--\ref{part:interval-verification}, namely the proof of Theorem~\ref{thm:Spectral structure of the bilinear form H}, where the main new ingredients of the present paper enter.

\medskip
\noindent\textbf{1. Blow-up dynamics.}
Under the near-soliton and negative-energy conditions \eqref{eq:main-assumptions}, the tube stability near \(Q\) allows us to write 
\begin{equation}\label{eq:proof strategy solution decomposition}
    u(t,x)
    =
    \frac{1}{\lambda(t)^{1/2}}
    \bigl(Q+\epsilon(t)\bigr)
    \left(\frac{x}{\lambda(t)}\right)
    e^{-i\gamma(t)}.
\end{equation}
Standard local well-posedness theory implies that, for any finite-time blow-up solution, the scaling parameter \(\lambda(t)\) tends to zero as \(t\) approaches the maximal lifespan. The blow-up analysis reduces to tracking the dynamics of \(\lambda(t)\).
Hence, we use the scale-adapted Hamiltonian pairing \((\partial_t u,i\Lambda u)_r\). By the scaling virial identity and conservation of energy,
\begin{equation}\label{eq:virial intro}
    (\partial_t u,i\Lambda u)_r=-E_0.
\end{equation}
Introducing the renormalized variables
\[
    \tfrac{ds}{dt}=\tfrac{1}{\lambda},
    \qquad
    y=\tfrac{x}{\lambda}
\]
with the decomposition \eqref{eq:proof strategy solution decomposition} and \eqref{eq:virial intro} yields
\begin{equation}\label{eq:proof strategy; motivation for H}
    \partial_s(\epsilon, i\Lambda Q)_r = \frac{1}{2}\mathbf H(\epsilon,\epsilon) + \lambda|E_0| - \frac{\lambda_s}{\lambda}(\epsilon,i\Lambda^2Q)_r -\tilde{\gmm}_s(\epsilon,\Lambda Q)_r + \text{h.o.t.},
\end{equation}
where $\mathbf{H}(\epsilon,\epsilon) \coloneqq ([L_Q,\Lambda]\epsilon,\epsilon)_r$. We fix \(\lambda(t)\) and \(\gamma(t)\) by imposing
\begin{equation}\label{eq:proof strategy orthogonal conditions}
    (\epsilon,i\Lambda^2Q)_r=0,
    \qquad
    (\epsilon,\Lambda Q)_r=0
\end{equation}
so that the two modulation terms in the identity \eqref{eq:proof strategy; motivation for H}  vanish.

Assuming Theorem~\ref{thm:Spectral structure of the bilinear form H}, the
identity \eqref{eq:proof strategy; motivation for H}, together with the
orthogonality conditions \eqref{eq:proof strategy orthogonal conditions}, gives
the monotonicity estimate
\[
    \partial_s(\epsilon,i\Lambda Q)_r
    \gtrsim
    \norm{\epsilon}_{\mathcal H^{1/2}}^2
    +
    \lambda |E_0|
    -
    C\left(
        (\epsilon,Q)_r^2
        +
        (\epsilon,i\Lambda Q)_r^2
    \right).
\]
The remaining negative directions are then removed by correcting the
monotonicity estimate. The near-soliton conservation laws control the
\((\epsilon,Q)_r\)-direction, while the \((\epsilon,i\Lambda Q)_r\)-direction is
absorbed by a correction along the generalized kernel direction \(R_1^\perp\) satisfying
\(L_Q[iR_1^\perp]=i\Lambda Q\). This produces a corrected scaling functional \(\tilde b\) with
\[
    \tilde b_s+C\tilde b^4
    \gtrsim
    \lambda |E_0|
    +
    \norm{\epsilon}_{\mathcal H^{1/2}}^2,
    \qquad
    \tilde b\sim -(\log\lambda)_s .
\]
Integrating this differential inequality and using the almost monotonicity of
\(\lambda\) gives finite-time blow-up and the upper bound
\eqref{eq:main-upper-rate}.

\medskip
\noindent\textbf{2. Coercivity of \(\mathbf H\).}
We first consider the Morse index problem. By applying the Birman--Schwinger principle and the symmetric resolvent identity, we reduce the determination of the Morse index to two concrete tasks: establishing suitable bounds for a compact operator and determining the sign of the determinant of an associated finite-dimensional matrix. This reduction forms Part~\ref{part:spectral-reduction}. 

Writing \(\epsilon=\epsilon_1+i\epsilon_2 \in H^{1/2}_{\mathrm e}\), the commutator form \(\mathbf H\) splits into two scalar quadratic forms.
\[
    \mathbf H(\epsilon,\epsilon)
    =
    \bigl((|D|+6V_Q)\epsilon_1,\epsilon_1\bigr)_r
    +
    \bigl((|D|+2V_Q)\epsilon_2,\epsilon_2\bigr)_r,
    \qquad
    V_Q\coloneqq yQ'Q.
\]
Let
\begin{align*}
    X_1&\coloneqq\{\epsilon_1 \in H^{1/2}_e : (\epsilon_1,Q)_r = (\epsilon_1,\Lambda Q)_r = 0\},\\
    X_2&\coloneqq \{\epsilon_2 \in H^{1/2}_e : (\epsilon_2,\Lambda Q)_r = (\epsilon_2,\Lambda^2 Q)_r = 0\}.
\end{align*}
Since \(V_Q\leq0\), it is enough to prove the stronger statement obtained by replacing the coefficient \(6\) and \(2\) with slightly larger constants \(c_1>6\) and \(c_2>2\). Indeed, if the operators \(|D|+c_jV_Q\) have Morse index zero on \(X_j\), then the desired coercivity \eqref{positivity of H} follows.

By the constrained index formula, the Morse index of the restriction to $X_j$ is given by the difference between the full Morse index of $|D|+c_jV_Q$ and the inertia of the $2\times2$ matrix associated with the constraints.

We then use the Birman--Schwinger formulation to count the full Morse index and, through the symmetric resolvent identity, to analyze the inertia of the constraint matrix. However, the zero-frequency singularity of the multiplier \(|D|^{-1}\) obstructs a direct analysis of both the Birman--Schwinger eigenvalue problem and the constraint matrix. We therefore introduce a regularized family of operators \(\mathcal L_{c}^{(\mu)}\) and the associated Birman--Schwinger operator \(K_Q^{(\mu)}\) for \(\mu>0\):
\[
    \mathcal L_c^{(\mu)}
    \coloneqq
    |D|+\mu+cV_Q,\qquad 
    K_Q^{(\mu)}
    \coloneqq
    |V_Q|^{1/2}(|D|+\mu)^{-1}|V_Q|^{1/2}
    \quad\text{on }L^2_{\mathrm e}.
\]
For each fixed \(\mu>0\), the operator \(K_Q^{(\mu)}\) is compact, nonnegative, and self-adjoint on \(L^2_{\mathrm e}\). The Birman--Schwinger principle gives
\[
    n_-^{\mathrm e}\bigl(\mathcal L_c^{(\mu)}\bigr)
    =
    n_+\bigl(c^{-1};K_Q^{(\mu)}\bigr),
    \qquad
    \dim\ker\mathcal L_c^{(\mu)}
    =
    \dim\ker\bigl(I-cK_Q^{(\mu)}\bigr).
\]
As \(\mu\to0^+\), the free resolvent \((|D|+\mu)^{-1}\) contains a logarithmically divergent
constant term, which appears in the Birman--Schwinger operator as a rank-one divergent term:
\[
    K_Q^{(\mu)}
    =
    \ell_\mu |V_Q|^{1/2}\otimes |V_Q|^{1/2}
    +
    T_Q
    +
    E_Q^{(\mu)},
    \qquad
    \ell_\mu\sim\log(\mu^{-1}),
    \qquad
    \norm{E_Q^{(\mu)}}_{L^2_{\mathrm e}\to L^2_{\mathrm e}}\to0.
\]
Consequently, writing \(\Phi_Q\coloneqq \tfrac{|V_Q|^{1/2}}{\norm{|V_Q|^{1/2}}_{L^2}}\), we obtain
\[
    \lambda_1(K_Q^{(\mu)})\to+\infty,
    \qquad
    \lambda_2(K_Q^{(\mu)})
    \to
    \norm{
        P_{\Phi_Q^\perp}T_QP_{\Phi_Q^\perp}
    }_{L^2_{\mathrm e}\to L^2_{\mathrm e}},
\]
where $P_{\Phi_Q^\perp}$ denotes the orthogonal projection onto \(\operatorname{span}\{\Phi_Q\}^\perp\).

Hence, once we show the projected logarithmic bound
\begin{equation*}
    \norm{
        P_{\Phi_Q^\perp}T_QP_{\Phi_Q^\perp}
    }_{L^2_{\mathrm e}\to L^2_{\mathrm e}}
    <
    \min\{c_1^{-1},c_2^{-1}\},
\end{equation*}
we deduce that, for $j=1,2$, the full Morse index is $1$ and \(\mathcal L_{c_j}^{(\mu)}\) is invertible for all sufficiently small \(\mu>0\).

Next, we consider the constraint matrix. Writing 
\[
    \Xi_1=\{Q,\Lambda Q\},
    \qquad
    \Xi_2=\{\Lambda Q,\Lambda^2Q\},
    \qquad
    \Xi_j=\{\Psi_1^{(j)},\Psi_2^{(j)}\},
\]
the matrix is given by
\[
    \mathbf M_j^{(\mu)}
    =
    \left(
        \bigl(
            (\mathcal L_{c_j}^{(\mu)})^{-1}\Psi_m^{(j)},
            \Psi_n^{(j)}
        \bigr)_r
    \right)_{1\leq m,n\leq2}.
\]
The symmetric resolvent identity separates the logarithmically divergent
rank-one terms in these entries. These divergent terms cancel and hence we obtain a finite matrix \(\mathbf M_{0,j}\) such that
\[
    \mathbf M_j^{(\mu)}
    =
    \mathbf M_{0,j}
    +
    o_{\mu\to0^+}^{\mathrm{ent}}(1).
\]
Therefore, Theorem~\ref{thm:Spectral structure of the bilinear form H} is reduced to proving, for $j=1,2$,
\begin{equation}\label{reduced by part 2 intro}
    \norm{
        P_{\Phi_Q^\perp}T_QP_{\Phi_Q^\perp}
    }_{L^2_{\mathrm e}\to L^2_{\mathrm e}}
    <
    \min\{c_1^{-1},c_2^{-1}\},
    \quad
    \det\mathbf M_{0,j}<0,
    \quad c_1>6,\ c_2>2.
\end{equation}

Second, to verify \eqref{reduced by part 2 intro} by a computer-assisted proof, we first construct a good approximation \(\mathfrak g\) of $Q$ and establish quantitative bounds for \(\mathfrak g - Q\). This forms Part~\ref{part:ground-state-approximation}.

The construction of \(\mathfrak g\) is not conventional, rather adapted to the nonlocal nature of the problem. Due to the nonlocality of $|D|$, we compactify the domain $\mathbb R$ via the change of variables $y=\tan \theta$ and incorporate the algebraic decay of $Q$ into the Fourier ansatz. More precisely, we set 
\[
    \mathfrak g(\tan\theta)
    =
    \cos^2\theta
    \sum_{k=0}^{200}U_k\cos(2k\theta),
    \qquad
    U_k\in\mathbb Q.
\]

We will find $\{U_k\}$ via a numerical Newton iteration so that the profile error $|D|\mathfrak g +\mathfrak g -\mathfrak g^3$ is small. Then, we use the Newton--Kantorovich theorem and the uniqueness result given by Frank--Lenzmann--Silvestre \cite{FrankLenzmannSilvestre2016CPAM} to prove that the approximation error $Q-\mathfrak g$ is sufficiently small.

Indeed, we apply a Newton--Kantorovich argument to \(\mathcal F(f)=|D|f+f-f^3\). A quantitative inverse bound for \(\mathcal F'(\mathfrak g)\) yields an exact even solution  \(\widetilde Q\) in an explicit \(H^1\)-neighborhood of \(\mathfrak g\). We prove that \(\widetilde Q\) is positive and that \(\mathcal F'(\widetilde Q)\) has Morse index one. The uniqueness theorem due to Frank--Lenzmann--Silvestre \cite{FrankLenzmannSilvestre2016CPAM} then gives \(\widetilde Q=Q\), yielding a rigorous bound for \(Q-\mathfrak g\).

Third, in Parts~\ref{part:spectral-estimates}--\ref{part:interval-verification}, we reduce the quantities appearing in \eqref{reduced by part 2 intro} to finite-dimensional problems with explicit error bounds in terms of $\mathfrak g$. Since $\mathfrak g$ is given explicitly, the resulting quantities can be rigorously evaluated by interval arithmetic with all rounding errors enclosed. Then, we prove \eqref{reduced by part 2 intro} using quantitative bounds for \(\mathfrak g-Q\) obtained from Part~\ref{part:ground-state-approximation}.

\subsection{Organization of the paper}
The paper is organized as follows. Section~2 collects notation and preliminary facts. In Part~\ref{part:blow-up-dynamics}, we prove Theorem~\ref{thm:main} assuming Theorem~\ref{thm:Spectral structure of the bilinear form H}. More precisely, we carry out the modulation analysis and
conclude finite-time blow-up together with the logarithmic upper bound on the
blow-up rate.

In Parts~\ref{part:spectral-reduction}--\ref{part:interval-verification} we prove Theorem~\ref{thm:Spectral structure of the bilinear form H}. In Part~\ref{part:spectral-reduction}, we reduce the coercivity theorem to \eqref{reduced by part 2 intro}, using a constrained Morse index formula and a Birman--Schwinger renormalization at zero frequency. In Part~\ref{part:ground-state-approximation}, we construct the compactified approximate profile \(\mathfrak g\), validate its closeness to the exact ground state \(Q\) by a Newton--Kantorovich argument, and derive the perturbation estimates used later. In Part~\ref{part:spectral-estimates}, we carry out finite-dimensional spectral reductions, including
the projected logarithmic Birman--Schwinger norm bound and the limiting constraint-matrix estimates, and obtain the corresponding estimates for \(Q\). Part~\ref{part:interval-verification} records the interval-arithmetic certificates and outward-rounded numerical bounds used in the proof.

The appendices contain the auxiliary analytic estimates, compactification
identities, quadrature bounds, and finite-dimensional approximation tools used
in Parts~\ref{part:ground-state-approximation}--\ref{part:interval-verification}.

\vspace{5bp}
\noindent\textbf{Acknowledgements.}
The author sincerely thanks Taegyu Kim and Soonsik Kwon for helpful discussions and valuable comments. 
This work was carried out in parallel with a companion paper by Taegyu Kim, Soonsik Kwon, and the author \cite{KimKwonPark2026arXiv}.
The author was partially supported by the National Research Foundation of Korea, RS-2019-NR040050, RS-2022-NR069873, and RS-2024-00333393.

\section{Notation and preliminaries}

\subsection{Notation and function spaces}

For nonnegative quantities \(A\) and \(B\), we write
\(A\lesssim B\) if $A\leq CB$ for some constant \(C>0\).  If the constant is allowed to depend on a parameter \(a\), we write \(A\lesssim_a B\).  We write $A\sim B$ when \(A\lesssim B\) and \(B\lesssim A\). We use the notation \(\mathcal O(B)\) with the same quantitative meaning as \(\lesssim B\).

If \(X\) is a normed space and \(a_\mu\geq0\), then
\[
    f_\mu=\mathcal O_X(a_\mu)
\]
means
\[
    \norm{f_\mu}_X\lesssim a_\mu.
\]
For finite-dimensional matrix-valued quantities, the notation
\[
    \mathbf A_\mu
    =
    \mathcal O_{\mathrm{ent}}(a_\mu)
\]
means that every entry of \(\mathbf A_\mu\) is
\(\mathcal O(a_\mu)\).  Similarly,
\[
    \mathbf A_\mu
    =
    o_{\mu\to0^+}^{\mathrm{ent}}(1)
\]
means entrywise convergence to zero as \(\mu\to0^+\).  All logarithms are
natural logarithms.

Let \(\Omega\subset\mathbb R\) be measurable and let
\(\mathbb K\in\{\mathbb R,\mathbb C\}\).  We use the standard Lebesgue spaces
\(L^p(\Omega;\mathbb K)\), \(1\leq p\leq\infty\), with norms
\[
    \norm{f}_{L^p(\Omega)}
    \coloneqq
    \left(
        \int_\Omega |f(x)|^p\,dx
    \right)^{1/p}
\]
for \(1\leq p<\infty\), with the usual modification for \(p=\infty\).
We omit \(\Omega\) and \(\mathbb K\) when they are clear from the context.

For \(f_1,f_2\in L^2(\Omega;\mathbb C)\), we set
\[
    (f_1,f_2)_{L^2(\Omega)}
    \coloneqq
    \int_\Omega
        f_1(x)\overline{f_2(x)}
    \,dx
\]
and
\[
    (f_1,f_2)_{r,\Omega}
    \coloneqq
    \operatorname{Re}
    (f_1,f_2)_{L^2(\Omega)}.
\]
When \(\Omega=\mathbb R\), we simply write
\[
    (f_1,f_2)_{L^2},
    \qquad
    (f_1,f_2)_r.
\]

We use the following non-unitary Fourier transform convention:
\[
    \widehat f(\xi)
    \coloneqq
    \int_{\mathbb R}
        e^{-ix\xi}f(x)
    \,dx,
    \qquad
    f(x)
    =
    \frac{1}{2\pi}
    \int_{\mathbb R}
        e^{ix\xi}\widehat f(\xi)
    \,d\xi.
\]
With this convention, Plancherel's identity reads
\[
    (f_1,f_2)_{L^2}
    =
    \frac{1}{2\pi}
    \int_{\mathbb R}
        \widehat f_1(\xi)
        \overline{\widehat f_2(\xi)}
    \,d\xi,
    \qquad
    \norm{f}_{L^2}^2
    =
    \frac{1}{2\pi}
    \norm{\widehat f}_{L^2}^2.
\]
We also use
\[
    \widehat{\partial_xf}(\xi)
    =
    i\xi\widehat f(\xi),
    \qquad
    \widehat{xf}(\xi)
    =
    i\partial_\xi\widehat f(\xi),
\]
and
\[
    \widehat{f_1*f_2}
    =
    \widehat f_1\,\widehat f_2,
    \qquad
    \widehat{f_1f_2}
    =
    \frac{1}{2\pi}
    \widehat f_1*\widehat f_2.
\]

For \(s\in\mathbb R\), we write \(|D|^s\) and
\(\langle\nabla\rangle^s\) for the Fourier multipliers defined by
\[
    \widehat{|D|^sf}(\xi)
    =
    |\xi|^s\widehat f(\xi),
    \qquad
    \widehat{\langle\nabla\rangle^sf}(\xi)
    =
    \langle\xi\rangle^s\widehat f(\xi),
\]
where
\[
    \langle\xi\rangle
    \coloneqq
    (1+|\xi|^2)^{1/2}.
\]
The homogeneous and inhomogeneous Sobolev norms on \(\mathbb R\) are
\[
    \norm{f}_{\dot H^s}
    \coloneqq
    \norm{|D|^sf}_{L^2},
    \qquad
    \norm{f}_{H^s}
    \coloneqq
    \norm{\langle\nabla\rangle^sf}_{L^2}.
\]

For an integer \(k\geq0\) and an interval or the whole line
\(\Omega\subset\mathbb R\), we also use
\[
    \norm{f}_{W^{k,2}(\Omega)}^2
    \coloneqq
    \sum_{\ell=0}^{k}
        \norm{\partial_x^\ell f}_{L^2(\Omega)}^2.
\]
On \(\mathbb R\), the spaces \(H^k\) and \(W^{k,2}\) coincide and their norms
are equivalent.  We nevertheless distinguish the Bessel-potential norm
\(\norm{\cdot}_{H^k}\) from the derivative norm
\(\norm{\cdot}_{W^{k,2}}\).

If \(\Omega\) is invariant under \(x\mapsto-x\), we define
\[
    L^2_{\mathrm e}(\Omega;\mathbb K)
    \coloneqq
    \left\{
        f\in L^2(\Omega;\mathbb K):
        f(-x)=f(x)
        \ \text{for a.e. }x\in\Omega
    \right\},
\]
and
\[
    L^2_{\mathrm o}(\Omega;\mathbb K)
    \coloneqq
    \left\{
        f\in L^2(\Omega;\mathbb K):
        f(-x)=-f(x)
        \ \text{for a.e. }x\in\Omega
    \right\}.
\]
The spaces
\[
    H^s_{\mathrm e},
    \qquad
    H^s_{\mathrm o},
    \qquad
    W^{k,2}_{\mathrm e},
    \qquad
    W^{k,2}_{\mathrm o}
\]
are defined similarly.  Unless otherwise stated, scalar spectral operators
are considered on real-valued spaces.

For \(\mathbf x=(x_1,\dots,x_N)\in\mathbb R^N\), we set
\[
    \norm{\mathbf x}_2
    \coloneqq
    \left(
        \sum_{j=1}^{N}|x_j|^2
    \right)^{1/2},
    \qquad
    (\mathbf x,\mathbf y)_{\mathbb R^N}
    \coloneqq
    \sum_{j=1}^{N}x_jy_j.
\]
For \(\mathbf A\in\mathbb R^{M\times N}\), the notation
\[
    \norm{\mathbf A}_2
    \coloneqq
    \sup_{\mathbf x\neq0}
    \frac{
        \norm{\mathbf A\mathbf x}_2
    }{
        \norm{\mathbf x}_2
    }
\]
denotes the induced Euclidean operator norm. For a real matrix \(\mathbf A=(A_{ij})\), we set
\[
    |\mathbf A|_{\mathrm{entry}}
    \coloneqq
    (|A_{ij}|)_{ij}.
\]
Inequalities between nonnegative matrices are understood componentwise.

We write
\[
    \langle x\rangle
    \coloneqq
    (1+x^2)^{1/2},
\]
and denote the generator of the \(L^2\)-scaling by
\[
    \Lambda f
    \coloneqq
    \left.
        \frac{d}{d\lambda}
    \right|_{\lambda=1}
        \lambda^{1/2}f(\lambda\cdot)
    =
    \frac12f+x\partial_xf.
\]
The same notation is used when the spatial variable is denoted by \(y\).
We define
\[
    \mathcal H^{1/2}
    \coloneqq
    \left\{
        f\in\dot H^{1/2}(\mathbb R):
        \langle x\rangle^{-2}f\in L^2(\mathbb R)
    \right\},
\]
with norm
\[
    \norm{f}_{\mathcal H^{1/2}}^2
    \coloneqq
    \norm{f}_{\dot H^{1/2}}^2
    +
    \norm{\langle x\rangle^{-2}f}_{L^2}^2.
\]

For a measurable set \(E\), we write \(\mathbf 1_E\) for its characteristic
function.  If \(X\) is a normed space, we write
\[
    B_X(f,r)
    \coloneqq
    \left\{
        v\in X:
        \norm{v-f}_X<r
    \right\}.
\]

\subsection{Operator and spectral notation}

Let \(X\) and \(Y\) be normed spaces and let \(A:X\to Y\) be bounded.  We
write
\[
    \norm{A}_{X\to Y}
    \coloneqq
    \sup_{f\in X\setminus\{0\}}
    \frac{
        \norm{Af}_Y
    }{
        \norm{f}_X
    }.
\]
When the underlying \(L^2\)-spaces are clear, we occasionally use
\[
    \norm{A}_{2\to2}
    \coloneqq
    \norm{A}_{L^2\to L^2}.
\]

For a Hilbert--Schmidt operator \(A\), we denote its Hilbert--Schmidt norm by
\(\norm{A}_{\mathrm{HS}}\).  In particular, if \(A\) is an integral operator
on \(L^2(\Omega)\) with kernel \(K_A\), then
\[
    \norm{A}_{\mathrm{HS}}^2
    =
    \int_\Omega
    \int_\Omega
        |K_A(x,y)|^2
    \,dx\,dy.
\]

We denote by \(A^*\) the Hilbert-space adjoint of \(A\), by
\(\mathbf A^\top\) the transpose of a real matrix, and by \(I_X\) the identity
operator on \(X\).  When the underlying space is clear, we simply write \(I\).

For vectors \(f_1,f_2\) in a real or complex Hilbert space, the rank-one
operator \(f_1\otimes f_2\) is defined by
\[
    (f_1\otimes f_2)v
    \coloneqq
    (v,f_2)f_1,
\]
where the relevant Hilbert-space pairing is understood.  If \(X\) is a closed
subspace, \(P_X\) denotes the orthogonal projection onto \(X\).  In particular,
if \(\Phi\) is normalized, then
\[
    P_{\Phi^\perp}
    =
    I-\Phi\otimes\Phi.
\]

For a closed operator \(A\), we denote its spectrum, essential spectrum,
point spectrum, resolvent set, and kernel by
\[
    \sigma(A),
    \qquad
    \sigma_{\mathrm{ess}}(A),
    \qquad
    \sigma_p(A),
    \qquad
    \rho(A),
    \qquad
    \ker A,
\]
respectively. For a bounded operator \(A\), we denote its spectral radius by
\[
    r(A)
    \coloneqq
    \sup\left\{
        |\lambda|:
        \lambda\in\sigma(A)
    \right\}.
\]
We use the commutator convention
\[
    [A,B]
    \coloneqq
    AB-BA.
\]
Operator inequalities between self-adjoint operators are understood in the
sense of quadratic forms.

For a compact nonnegative self-adjoint operator \(K\), we write
\[
    \lambda_1(K)
    \geq
    \lambda_2(K)
    \geq
    \cdots
    \geq0
\]
for its eigenvalues, counted with multiplicity and arranged in non-increasing
order.  For \(a>0\), we set
\[
    n_+(a;K)
    \coloneqq
    \#\left\{
        k:
        \lambda_k(K)>a
    \right\}.
\]

For a self-adjoint operator, a real symmetric matrix, or a symmetric
quadratic form, \(n_-\) denotes the number of strictly negative directions,
counted with multiplicity, whenever this number is finite.  We similarly write
\(n_{\leq0}\) for the number of non-positive directions.

If a self-adjoint operator \(A\) preserves parity, we define
\[
    n_-^{\mathrm e}(A)
    \coloneqq
    n_-\left(
        A\big|_{L^2_{\mathrm e}}
    \right),
    \qquad
    n_-^{\mathrm o}(A)
    \coloneqq
    n_-\left(
        A\big|_{L^2_{\mathrm o}}
    \right).
\]
When the full Morse index is finite,
\[
    n_-(A)
    =
    n_-^{\mathrm e}(A)
    +
    n_-^{\mathrm o}(A).
\]

For a real symmetric matrix \(\mathbf A\), we write
\[
    \lambda_{\max}(\mathbf A),
    \qquad
    \lambda_{\min}(\mathbf A)
\]
for its largest and smallest eigenvalues.

\subsection{Standard operator and spectral tools}

We record the standard perturbative facts used below.

\begin{lemma}[Neumann-series criterion]\label{lem:Neumann-series-criterion}
    Let \(X\) be a Banach space and let \(B:X\to X\) satisfy
    \[
        \norm{B}_{X\to X}<1.
    \]
    Then \(I-B\) is invertible and
    \[
        (I-B)^{-1}
        =
        \sum_{k=0}^{\infty}B^k,
        \qquad
        \norm{(I-B)^{-1}}_{X\to X}
        \leq
        \frac{1}{
            1-\norm{B}_{X\to X}
        }.
    \]
    More generally, if \(A:X\to X\) is invertible and
    \(E:X\to X\) satisfies
    \[
        \norm{A^{-1}E}_{X\to X}<1,
    \]
    then \(A+E\) is invertible and
    \[
        \norm{(A+E)^{-1}}_{X\to X}
        \leq
        \frac{
            \norm{A^{-1}}_{X\to X}
        }{
            1-\norm{A^{-1}E}_{X\to X}
        }.
    \]
\end{lemma}

\begin{lemma}[Weyl perturbation inequality]
\label{lem:Weyl-perturbation-inequality}
    Let \(A\) and \(B\) be compact nonnegative self-adjoint operators on a
    Hilbert space.  Then
    \[
        \left|
            \lambda_k(A)-\lambda_k(B)
        \right|
        \leq
        \norm{A-B},
        \qquad
        k\geq1.
    \]
    In particular, if
    \(\mathbf A,\mathbf B\in\mathbb R^{N\times N}\) are symmetric, then
    \[
        \left|
            \lambda_k(\mathbf A)-\lambda_k(\mathbf B)
        \right|
        \leq
        \norm{\mathbf A-\mathbf B}_2,
        \qquad
        1\leq k\leq N.
    \]
\end{lemma}

\begin{theorem}[Weyl's essential-spectrum theorem]
    Let \(A\) be self-adjoint and let \(B\) be relatively compact with
    respect to \(A\).  Then
    \[
        \sigma_{\mathrm{ess}}(A+B)
        =
        \sigma_{\mathrm{ess}}(A).
    \]
\end{theorem}

A bounded operator \(K\) on a real \(L^2\)-space is called positivity
improving if \(f\geq0\) and \(f\not\equiv0\), then \(Kf>0\) almost everywhere.

\begin{theorem}[Jentzsch's theorem]
\label{thm:Jentzsch-theorem}
    Let \(K\) be a compact, nonnegative, self-adjoint, and positivity
    improving operator on a real \(L^2\)-space.  Then
    \(\lambda_1(K)>0\) is simple, and the corresponding normalized
    eigenfunction may be chosen strictly positive almost everywhere.  In
    particular,
    \[
        \lambda_1(K)>\lambda_2(K).
    \]
\end{theorem}

\subsection{Basic analytic tools}

We first record the weighted estimate used in Part~\ref{part:blow-up-dynamics}.  For a locally
integrable function \(f\) and a bounded interval \(J\subset\mathbb R\), set
\begin{equation}\label{eq:average-f}
    f_J
    \coloneqq
    \frac{1}{|J|}
    \int_J f(x)\,dx.
\end{equation}
Within the proof below, we use the seminorm
\[
    \norm{f}_{\mathrm{BMO}}
    \coloneqq
    \sup_{J\subset\mathbb R}
    \frac{1}{|J|}
    \int_J
        |f(x)-f_J|
    \,dx,
\]
where the supremum is taken over bounded intervals.  The classical embedding
\[
    \dot H^{1/2}(\mathbb R)
    \hookrightarrow
    \mathrm{BMO}(\mathbb R)
\]
gives
\begin{equation}\label{eq:H12-BMO}
    \norm{f}_{\mathrm{BMO}}
    \lesssim
    \norm{f}_{\dot H^{1/2}(\mathbb R)}.
\end{equation}
See, for instance, \cite[Chapter~VI]{Steinharmonic-book}.

\begin{lemma}[Fractional weighted Hardy inequality]
    For every \(f\in\mathcal H^{1/2}\), we have
    \begin{equation}\label{eq:Hardy inequ}
        \norm{\langle x\rangle^{-1}f}_{L^2}
        \lesssim
        \norm{f}_{\mathcal H^{1/2}}.
    \end{equation}
\end{lemma}

\begin{proof}
    By \eqref{eq:H12-BMO}, it is enough to prove
    \[
        \norm{\langle x\rangle^{-1}f}_{L^2}
        \lesssim
        \norm{f}_{\mathrm{BMO}}
        +
        \norm{\langle x\rangle^{-2}f}_{L^2}.
    \]
    Let \(J_k \coloneqq (-2^k,2^k)\) and \(m_k \coloneqq f_{J_k}\) for \(k\geq0\).
    Since \(\langle x\rangle^{-2}\simeq1\) on \(J_0\),
    \[
        |m_0|
        \lesssim
        \norm{\langle x\rangle^{-2}f}_{L^2}.
    \]
    For \(k\geq1\),
    \[
        \begin{aligned}
            |m_k-m_{k-1}|
            \leq
            \frac{1}{|J_{k-1}|}
            \int_{J_{k-1}}
                |f-m_k|
            \,dx
            \lesssim
            \norm{f}_{\mathrm{BMO}},
    \end{aligned}
    \]
    and hence
    \[
        |m_k|
        \lesssim
        \norm{\langle x\rangle^{-2}f}_{L^2}
           +
        k\norm{f}_{\mathrm{BMO}}.
    \]
    By the John--Nirenberg inequality,
    \[
        \int_{J_k}
            |f-m_k|^2
        \,dx
        \lesssim
        |J_k|
        \norm{f}_{\mathrm{BMO}}^2
        \lesssim
        2^k
        \norm{f}_{\mathrm{BMO}}^2.
    \]
    Now, set
    \[
        A_0
        \coloneqq
        J_0,
        \qquad
        A_k
        \coloneqq
        \left\{
            2^{k-1}\leq|x|<2^k
        \right\},
        \qquad
        k\geq1.
    \]
    Then
    \[
        \int_{A_0}
            \langle x\rangle^{-2}|f|^2
        \,dx
        \lesssim
        \norm{\langle x\rangle^{-2}f}_{L^2}^2,
    \]
    while, for \(k\geq1\),
    \[
    \begin{aligned}
        \int_{A_k}
            \langle x\rangle^{-2}|f|^2
        \,dx
        &\lesssim
        2^{-2k}
        \int_{J_k}
            |f-m_k|^2
        \,dx
        +
        2^{-2k}|J_k||m_k|^2
        \\
        &\lesssim
        2^{-k}
        \norm{f}_{\mathrm{BMO}}^2
        +
        2^{-k}
        \left(
            \norm{\langle x\rangle^{-2}f}_{L^2}
            +
            k\norm{f}_{\mathrm{BMO}}
        \right)^2.
    \end{aligned}
    \]
    Summing in \(k\) yields
    \[
        \norm{\langle x\rangle^{-1}f}_{L^2}^2
        \lesssim
        \norm{\langle x\rangle^{-2}f}_{L^2}^2
        +
        \norm{f}_{\mathrm{BMO}}^2,
    \]
    and \eqref{eq:Hardy inequ} follows from \eqref{eq:H12-BMO}.
\end{proof}

We next recall the following fractional commutator estimate.

\begin{lemma}[Fractional commutator estimate
{\protect\cite[Appendix~E]{KLR2013ARMAhalfwave}}]
    Let \(0\leq\alpha\leq1\).  Assume that
    \(f_1,f_2:\mathbb R\to\mathbb C\) satisfy
    \[
        |D|^\alpha f_1\in L^2(\mathbb R),
        \qquad
        \widehat f_2\in L^1(\mathbb R).
    \]
    Then
    \begin{equation}\label{eq:frac-comm}
        \norm{
            |D|^\alpha(f_1f_2)
            -
            f_1|D|^\alpha f_2
        }_{L^2}
        \leq
        C_\alpha
        \norm{|D|^\alpha f_1}_{L^2}
        \norm{\widehat f_2}_{L^1},
    \end{equation}
    where \(C_\alpha>0\) depends only on \(\alpha\).
\end{lemma}

We also use the following interval Poincar\'e inequality.

\begin{lemma}[Wirtinger--Poincar\'e inequality]
\label{lem:wirtinger}
    Let \(J\subset\mathbb R\) be an interval of length \(|J|=\ell\), and let
    \(f\in H^1(J)\).  Then
    \begin{equation*}
        \norm{f-f_J}_{L^2(J)}
        \leq
        \frac{\ell}{\pi}
        \norm{f'}_{L^2(J)},
    \end{equation*}
    where \(f_J\) is defined in \eqref{eq:average-f}.
\end{lemma}

Finally, we record the elementary Sobolev estimates with explicit constants.

\begin{lemma}[Sobolev and product estimates]
    For every \(f\in H^1(\bbR)\),
    \begin{equation}\label{eq:H1-Linfty}
        \norm{f}_{L^\infty}
        \leq
        \frac1{\sqrt2}
        \norm{f}_{H^1}.
    \end{equation}
    Moreover, for \(f_1,f_2\in H^1(\bbR)\),
    \begin{equation}\label{eq:H1-product}
        \norm{f_1f_2}_{H^1}
        \leq
        \sqrt{\frac52}
        \norm{f_1}_{H^1}
        \norm{f_2}_{H^1}.
    \end{equation}
\end{lemma}

\begin{proof}
    Fourier inversion, Cauchy--Schwarz, and Plancherel's identity give
    \[
        \begin{aligned}
            |f(y)|
            \leq
            \frac1{2\pi}
            \left(
                \int_{\bbR}
                    \frac{d\xi}{1+\xi^2}
            \right)^{1/2}
            \left(
                \int_{\bbR}
                    (1+\xi^2)
                    |\widehat f(\xi)|^2
                \,d\xi
            \right)^{1/2}
            =
            \frac1{\sqrt2}
            \norm{f}_{H^1}.
        \end{aligned}
    \]
    Hence
    \[
        \norm{f_1f_2}_{L^2}
        \leq
        \frac1{\sqrt2}
        \norm{f_1}_{H^1}
        \norm{f_2}_{H^1},
        \qquad
        \norm{(f_1f_2)'}_{L^2}
        \leq
        \sqrt2
        \norm{f_1}_{H^1}
        \norm{f_2}_{H^1}.
    \]
    Squaring and adding gives \eqref{eq:H1-product}.
\end{proof}

\subsection{Local Cauchy theory}

We recall the local theory for \eqref{eq:halfwave}.

\begin{theorem}[Local Cauchy theory
{\protect\cite[Theorem~D.1]{KLR2013ARMAhalfwave}}]
\label{thm:Cauchy-theory}
    Let \(s\geq\frac12\).  For every
    \(u_0\in H^s(\mathbb R)\), there exists a unique maximal forward solution
    \[
        u\in C([0,T);H^s(\mathbb R))
    \]
    to \eqref{eq:halfwave}, where
    \(0<T=T(u_0)\leq+\infty\).  Moreover:
    \begin{enumerate}[label=\textup{(\roman*)}]
        \item
        if \(T<+\infty\), then \(\norm{u(t)}_{H^{1/2}} \to +\infty\) as \(t\uparrow T\).

        \item
        if \(s>\frac12\), then the flow map is locally Lipschitz on bounded
        subsets of \(H^s(\mathbb R)\).
    \end{enumerate}
\end{theorem}

The uniqueness of the flow implies that even initial data generate even
solutions.  The mass and energy in \eqref{eq:conservation law} are conserved
on the maximal lifespan.

\subsection{The ground state and linearized facts}

We fix the unique positive even solution \(Q\) of
\begin{equation*}
    |D|Q+Q-Q^3=0.
\end{equation*}
The existence, uniqueness up to translation, and nondegeneracy of \(Q\) were
proved in \cite{FrankLenzmann2013Acta}.  The corresponding solitary wave is
\(e^{it}Q\), and
\[
    E(Q)=0.
\]
The ground state \(Q\) is the optimizer of the sharp Gagliardo--Nirenberg
inequality
\begin{equation}\label{eq:sharp-Gagliardo-Nirenberg}
    \norm{f}_{L^4}^4
    \leq
    \frac{2}{
        \norm{Q}_{L^2}^2
    }
    \norm{f}_{L^2}^2
    \norm{f}_{\dot H^{1/2}}^2,
    \qquad
    f\in H^{1/2}(\mathbb R).
\end{equation}
Consequently,
\begin{equation}\label{eq:energy-lower-bound}
    E(f)
    \geq
    \frac12
    \norm{f}_{\dot H^{1/2}}^2
    \left(
        1-
        \frac{
            \norm{f}_{L^2}^2
        }{
            \norm{Q}_{L^2}^2
        }
    \right).
\end{equation}
In particular, solutions with initial mass strictly below \(M(Q)\) are global.

For a complex-valued function \(f\), we define the real-linear operator
\[
    L_Q[f]
    \coloneqq
    |D|f+f-Q^2f-2Q^2\operatorname{Re}f.
\]
The operator \(L_Q\) is self-adjoint on the real Hilbert space
\(L^2(\mathbb R;\mathbb C)\), with domain
\(H^1(\mathbb R;\mathbb C)\), with respect to the real pairing
\((\cdot,\cdot)_r\). If \(f=f_1+if_2\), where \(f_1\) and \(f_2\) are real-valued, then
\[
    L_Q[f]
    =
    \bigl(|D|+1-3Q^2\bigr)f_1
    +
    i\bigl(|D|+1-Q^2\bigr)f_2.
\]

The phase and scaling identities, together with the first generalized kernel
relation, are
\begin{equation}\label{eq:even-kernel-relations}
    L_Q[iQ]=0,
    \qquad
    L_Q[\Lambda Q]=-Q,
    \qquad
    L_Q[iR_1]=i\Lambda Q.
\end{equation}
Here \(R_1\) denotes the real-valued even profile denoted by \(S_1\) in
\cite[Proposition~4.1]{KLR2013ARMAhalfwave}.  We set
\begin{equation}\label{eq:R_1-positivity, e_1 definition}
    e_1
    \coloneqq
    \bigl(i^{-1}L_Q[iR_1],R_1\bigr)_r
    =
    (\Lambda Q,R_1)_r
    >
    0.
\end{equation}
Moreover,
\begin{equation}\label{eq:R_1-decay}
    |R_1(y)|
    +
    |\Lambda R_1(y)|
    +
    |\Lambda^2R_1(y)|
    \lesssim
    \langle y\rangle^{-2}.
\end{equation}

We recall the spectral structure of \(L_Q\) proved in
\cite{FrankLenzmann2013Acta}.  The operator has one negative direction in the
even sector and no negative direction in the odd sector:
\begin{equation*}
    n_-^{\mathrm e}(L_Q)=1,
    \qquad
    n_-^{\mathrm o}(L_Q)=0.
\end{equation*}
Its parity kernels are
\begin{equation*}
    \ker\left(
        L_Q\big|_{L^{2}_{\mathrm e}}
    \right)
    =
    \operatorname{span}_{\mathbb R}\{iQ\},
    \qquad
    \ker\left(
        L_Q\big|_{L^{2}_{\mathrm o}}
    \right)
    =
    \operatorname{span}_{\mathbb R}\{Q'\}.
\end{equation*}

Let \(\phi\in H^{1/2}_{\mathrm e}(\mathbb R;\mathbb R)\) be the positive
\(L^2\)-normalized negative eigenfunction:
\[
    L_Q[\phi]
    =
    e_0\phi,
    \qquad
    e_0<0,
    \qquad
    \norm{\phi}_{L^2}=1.
\]

\begin{lemma}[Even-sector coercivity of \(L_Q\),
{\protect\cite{FrankLenzmann2013Acta}}]
\label{lem:Coercivity 1}
    There exists a universal constant \(\kappa_0>0\) such that, for every
    \(f\in H^{1/2}_{\mathrm e}(\mathbb R;\mathbb C)\),
    \begin{equation*}
        (L_Q[f],f)_r
        \geq
        \kappa_0
        \norm{f}_{H^{1/2}}^2
        -
        \frac{1}{\kappa_0}
        \left[
            (f,\phi)_r^2
            +
            (f,iQ)_r^2
        \right].
    \end{equation*}
\end{lemma}

Finally, the half-wave operator satisfies
\begin{equation}\label{eq:commutator-formula}
    [|D|,\Lambda]
    =
    |D|.
\end{equation}
Consequently, if \(f=f_1+if_2\), with \(f_1,f_2\) real-valued, then
\begin{equation*}
    [L_Q,\Lambda]f
    =
    \bigl(
        |D|+6yQ'(y)Q(y)
    \bigr)f_1
    +
    i
    \bigl(
        |D|+2yQ'(y)Q(y)
    \bigr)f_2.
\end{equation*}

\part{Blow-up dynamics}\label{part:blow-up-dynamics}

In Part~\ref{part:blow-up-dynamics}, we prove Theorem~\ref{thm:main} using Theorem~\ref{thm:Spectral structure of the bilinear form H}. The latter theorem is proved in Parts~\ref{part:spectral-reduction}--\ref{part:interval-verification}.

\section{Modulation analysis and blow-up dynamics}
Let \(u_0\in H^{1/2}_{\mathrm e}(\bbR)\) satisfy
\begin{equation}\label{eq:condition1}
    0 < \alpha_0 \coloneqq\norm{u_0}_{L^2}^2- \norm{Q}_{L^2}^2<\alpha^*,\qquad
    E_0\coloneqq E(u_0)<0,
\end{equation}
where \(\alpha^*>0\) will be fixed below. Let
\[
    u\in C\bigl([0,T);H^{1/2}_{\mathrm e}(\bbR)\bigr)
\]
be the corresponding maximal forward solution to \eqref{eq:halfwave}.

\subsection{Modulation decomposition near the ground state}
\begin{lemma}\label{lem:orthogonal decomposition of the solution u}
    Let $u$ be a solution to \eqref{eq:halfwave} satisfying \eqref{eq:condition1}. There exists $\alpha^*>0$ such that, if $\alpha_0<\alpha^*$, then there are continuous parameters $\lambda:[0,T)\to\bbR_+$ and $\gmm:[0,T)\to\bbR$ such that
    \begin{equation}\label{eq:definition of epsilon(t)}
        \epsilon(t,y) := e^{i\gmm(t)}\lambda^{1/2}(t)u(t,\lambda(t)y) - Q(y)
    \end{equation}
    satisfies the following properties:
    \begin{equation}\label{eq:orthogonal condition on ep}
        (\epsilon, \Lambda Q)_r = (\epsilon, i\Lambda^2 Q)_r = 0,
    \end{equation}
    and
    \begin{equation}\label{eq:smallenss of ep and lambda esitmate}
        \left|1-\lambda(t)\frac{\norm{u(t)}_{\dot{H}^{1/2}}^2}{\norm{Q}_{\dot{H}^{1/2}}^2}\right| + \norm{\epsilon(t)}_{H^{1/2}} \leq \delta(\alpha_0),\quad \text{where} \ \delta(\alpha_0) \to 0 \text{ as } \alpha_0 \to 0,
    \end{equation}
    for all $t \in [0,T)$.
\end{lemma}

\begin{lemma}
    There exists $\alpha^*>0$ so that for $\alpha_0\in (0,\alpha^*)$, we have
    \begin{equation}\label{eq:refined bound of ep(t)}
        \norm{\epsilon(t)}_{H^{1/2}} \lesssim \sqrt{\alpha_0},\quad t \in [0,T).
    \end{equation}
\end{lemma}
\begin{proof}
    From mass and energy conservation laws, we have
    \begin{equation}\label{eq:identity from mass and energy conservation law}
        \begin{aligned}
            &\norm{\epsilon}_{L^2}^2 + 2(\epsilon,Q)_r = \alpha_0,\\
            &\norm{\epsilon}_{\dot H^{1/2}}^2- 2(\epsilon,Q)_r - \int \left[2Q^2|\epsilon|^2+\Re\{\epsilon\}^2Q^2\right] = -2\lambda|E_0| + \frac{1}{2}\int \mathcal{N}_{\geq 3}[\epsilon],
        \end{aligned}
    \end{equation}
    where $\mathcal{N}_{\geq 3}[\epsilon]$ is defined by
    \begin{equation}\label{eq:definition of N_geq 3}
        \mathcal{N}_{\geq 3}[\epsilon] \coloneqq 4\Re\{Q\epsilon\}|\epsilon|^2 + |\epsilon|^4.
    \end{equation}
    Hence, we have
    \begin{equation}\label{eq: L_Q bilinear form estimate 1}
        (L_Q[\epsilon],\epsilon)_r \leq \alpha_0 + C\norm{\epsilon}_{H^{1/2}}\norm{\epsilon}_{L^2}^2.
    \end{equation}
    Let
    \[
        \epsilon^\perp
        \coloneqq
        \epsilon-a_1\Lambda Q-a_2iQ,
    \]
    where \(a_1 \coloneqq \tfrac{(\epsilon,\phi)_r}{(\phi,\Lambda Q)_r}\), \(a_2 \coloneqq
    \frac{(\epsilon,iQ)_r}{\norm{Q}_{L^2}^2}\) and \(\phi\) is the positive \(L^2\)-normalized eigenfunction introduced in Lemma~\ref{lem:Coercivity 1}. The constants $a_1$ and $a_2$ are well defined because $(\phi,\Lambda Q)_r=-e_0^{-1}(\phi,Q)_r>0$, as $Q,\phi>0$. Furthermore, the definitions of $a_1$ and $a_2$, together with \eqref{eq:orthogonal condition on ep}, give $a_1 =-\tfrac{(\epsilon^\perp,\Lambda Q)_r}{\norm{\Lambda Q}_{L^2}^2}$, $a_2=\tfrac{(\epsilon^\perp,i\Lambda^2 Q)_r}{\norm{\Lambda Q}_{L^2}^2}$. Therefore, we obtain
    \begin{equation}\label{eq:equivalence of L^2 norm of ep and tilde ep}
        (\epsilon^\perp, \phi)_r = (\epsilon^\perp,iQ)_r = 0,\quad \norm{\epsilon^\perp}_{L^2} \sim \norm{\epsilon}_{L^2}
    \end{equation}
    Since \(\epsilon^\perp\) is an even function, using Lemma~\ref{lem:Coercivity 1}, we have
    \begin{equation}\label{eq:L_Q corecivity of tilde ep}
        \norm{\epsilon^\perp}_{H^{1/2}}^2 \lesssim(L_Q[\epsilon^\perp],\epsilon^\perp)_r = (L_Q[\epsilon],\epsilon)_r + 2a_1(\epsilon,Q)_r.
    \end{equation}
    In the last equality, we use $L_Q[\Lambda Q]=-Q$. From \eqref{eq:identity from mass and energy conservation law} we have $|(\epsilon,Q)_r|\lesssim \alpha_0 + \norm{\epsilon}_{L^2}^2$. Hence \eqref{eq:L_Q corecivity of tilde ep} and \eqref{eq:equivalence of L^2 norm of ep and tilde ep} give
    \[
        \norm{\epsilon}_{L^2}^2 \lesssim \norm{\epsilon^\perp}_{L^2}^2 \lesssim (L_Q[\epsilon^\perp],\epsilon^\perp)_r \lesssim \alpha_0 + \norm{\epsilon}_{H^{1/2}}\norm{\epsilon}_{L^2}^2.
    \]
    In the last inequality, we use \eqref{eq: L_Q bilinear form estimate 1}. From \eqref{eq:smallenss of ep and lambda esitmate}, possibly reducing $\alpha^*>0$ we obtain
    \[
        \norm{\epsilon}_{L^2}\lesssim \sqrt{\alpha_0}.
    \]
    Substituting this estimate into \eqref{eq: L_Q bilinear form estimate 1} gives 
    \[
        \norm{\epsilon}_{H^{1/2}}^2 = (L_Q[\epsilon],\epsilon)_r + \int Q^2|\epsilon|^2 + 2Q^2\Re\{\epsilon\}^2 \lesssim \alpha_0,
    \]
    which proves \eqref{eq:refined bound of ep(t)}.
\end{proof}

\begin{remark}
    The modulation orthogonality conditions
    \eqref{eq:orthogonal condition on ep} do not by themselves imply
    coercivity of \(L_Q\). Indeed,
    \[
        (Q,\Lambda Q)_r=(Q,i\Lambda^2Q)_r=0,
    \]
    whereas
    \[
        (L_Q[Q],Q)_r=-2\int_{\bbR}Q^4\,dy<0.
    \]
    Thus the constrained space still contains the negative direction \(Q\).
\end{remark}

We now introduce the renormalized time variable
\[
    s(t)\coloneqq \int_0^t\frac{dt'}{\lambda(t')},
    \qquad
    \frac{ds}{dt}=\frac1{\lambda(t)}.
\]
We claim that \(s(t)\to+\infty\) as \(t\uparrow T\). Suppose first that
\(T=+\infty\). By mass conservation and \eqref{eq:energy-lower-bound},
\[
    E_0
    \geq
    -\frac{\alpha_0}{2\norm{Q}_{L^2}^2}
    \norm{u(t)}_{\dot H^{1/2}}^2.
\]
Since \(E_0<0\), this gives
\(\norm{u(t)}_{\dot H^{1/2}}^2
\geq 2\norm{Q}_{L^2}^2|E_0|/\alpha_0\). Combining this estimate with
\eqref{eq:smallenss of ep and lambda esitmate}, and reducing
\(\alpha^*>0\) if necessary, we obtain
\(\lambda(t)\lesssim_{u_0}1\). Hence
\(s(t)\gtrsim_{u_0}t\to+\infty\).

Suppose instead that \(T<+\infty\). By the trivial self-similar estimate in
Lemma~\ref{lem:trivial lower bound}, one has
\(\lambda(t)\lesssim T-t\). Therefore,
\[
    s(t)
    \gtrsim
    \int_0^t\frac{dt'}{T-t'}
    \longrightarrow+\infty
    \qquad\text{as }t\uparrow T.
\]
Thus the renormalized time variable \(s\) ranges over \([0,+\infty)\). From now on, we work in the renormalized variables
\begin{equation}\label{eq:renomalised variable}
    s=s(t),\qquad y\coloneqq \frac{x}{\lambda(t)}.
\end{equation}

We first record the following estimates.
\begin{lemma}
    Let $k \in \mathbb{N}$. Then
    \begin{equation}\label{eq:Lambda^k Q control}
        \norm{\langle y \rangle^{2}\Lambda^k Q}_{L^\infty} \lesssim_k 1.
    \end{equation}
    Moreover, the following estimates hold:
    \begin{enumerate}
        \item Control of linear terms:
        \begin{equation}\label{eq:linear control}
            |(\epsilon, \langle y \rangle^{-2})_r|
            \lesssim \norm{\epsilon}_{\mathcal{H}^{1/2}}.
        \end{equation}

        \item Control of second order terms:
        \begin{equation}\label{eq:quadratic control}
            |(\mathrm{NL}_{\geq 2}(\epsilon),\langle y \rangle^{-2})_r|
            \lesssim \norm{\epsilon}_{\mathcal{H}^{1/2}}^2,
        \end{equation}
        where $\mathrm{NL}_{\geq 2}(\epsilon)$ is defined in \eqref{eq:definition of quadratic term}.

        \item Control of higher order terms:
        \begin{equation}\label{eq:cubic control}
            \int |\mathcal{N}_{\geq 3}[\epsilon]|
            + |(|\epsilon|^2\epsilon, \langle y \rangle^{-2})_r|
            \lesssim \sqrt{\alpha_0}\norm{\epsilon}_{\mathcal{H}^{1/2}}^2,
        \end{equation}
        where $\mathcal{N}_{\geq 3}[\epsilon]$ are defined in \eqref{eq:definition of N_geq 3}.
    \end{enumerate}
\end{lemma}

\begin{proof}
    First, by \eqref{eq:commutator-formula}, for all $k \in \mathbb{N}$ we have
    \[
        |D|(\Lambda^k Q) = [|D|,\Lambda]\Lambda^{k-1} Q + \Lambda |D|(\Lambda^{k-1}Q), 
    \]
    and, by induction, we obtain
    \begin{equation}\label{D(Lambda^kQ)}
        |D|(\Lambda^k Q) = (1+\Lambda)^k |D|Q.
    \end{equation}

    We claim that $L_Q[\Lambda^k Q] = \Pi(Q,\Lambda Q,\cdots,\Lambda^{k-1} Q)$, where $\Pi(Q,\Lambda Q,\cdots,\Lambda^{k-1} Q)$ denotes a finite linear combination of products of elements in $\{Q,\Lambda Q,\cdots,\Lambda^{k-1} Q\}$. For $k=1$, \eqref{eq:even-kernel-relations} yields $L_Q[\Lambda Q] = -Q$. Now assume that the claim holds for $k = m$. Then, using \eqref{eq:commutator-formula} and \eqref{D(Lambda^kQ)}, we obtain
    \begin{equation}\label{induction for k}
        \begin{aligned}
            L_Q[\Lambda^{m+1}Q] &= [L_Q,\Lambda] \Lambda^m Q + \Lambda L_Q[\Lambda^m Q] \\
            &= (|D|+6yQ'Q)\Lambda^m Q + \Lambda L_Q[\Lambda^m Q] \\
            & = (1+\Lambda)^m |D|Q + 6(\Lambda Q - \tfrac{1}{2}Q)\Lambda^mQ + \Lambda L_Q[\Lambda^m Q] \\
            & = \Pi(Q,\Lambda Q,\cdots,\Lambda^m Q).
        \end{aligned}
    \end{equation}
    The last equality follows from \(|D|Q = - Q + Q^3\). Next, since $|\Lambda Q(y)| \lesssim \langle y \rangle^{-2}$ and $(\Lambda^k Q,\nabla Q)_r = 0$ for all $k\in \mathbb{N}$, we obtain \eqref{eq:Lambda^k Q control} by applying \cite[Lemma A.1]{KLR2013ARMAhalfwave}.

    From the Cauchy--Schwarz inequality and \eqref{eq:Hardy inequ}, we have
    \[
        |(\epsilon,\langle y\rangle^{-2})_r| \lesssim \norm{\langle y\rangle^{-1}\epsilon}_{L^2} \lesssim \norm{\epsilon}_{\mathcal{H}^{1/2}},
    \]
    which yields \eqref{eq:linear control}. Finally, the Cauchy--Schwarz inequality gives \eqref{eq:quadratic control}, while the Gagliardo--Nirenberg inequality and \eqref{eq:smallenss of ep and lambda esitmate} give \eqref{eq:cubic control}.
\end{proof}

\begin{proposition}
    There exists $\alpha^*>0$ such that, for $\alpha_0 \in (0,\alpha^*)$, one can choose modulation parameters $\{\lambda(s),\gmm(s)\}$ of class $C^1$ with the following properties:
    \begin{enumerate}
    \item Equation for $\epsilon$: The error term $\epsilon$ satisfies
    \begin{equation}\label{eq:flow of ep(s)}
        \partial_s \epsilon + iL_Q[\epsilon] - i\mathrm{NL}_{\geq 2}(\epsilon)
        = \frac{\lambda_s}{\lambda}\Lambda(Q+\epsilon) - i\tilde{\gmm}_s (Q+\epsilon),
    \end{equation}
    where we set \(\tilde{\gmm}_s \coloneqq -1-\gmm_s\),
    and the nonlinear term $\mathrm{NL}_{\geq 2}(\epsilon)$ is defined by
    \begin{equation}\label{eq:definition of quadratic term}
        \mathrm{NL}_{\geq 2}(\epsilon)
        \coloneqq 2\Re\{Q\epsilon\}\,\epsilon + |\epsilon|^2 Q + |\epsilon|^2 \epsilon.
    \end{equation}

    \item Energy-type control and modulation bounds: We have
    \begin{equation}\label{eq:energy type control}
        \bigl|\lambda(s) E_0 + (\epsilon,Q)_r \bigr|
        \lesssim \|\epsilon\|_{\mathcal{H}^{1/2}}^2,
    \end{equation}
    as well as the \textit{a priori} bounds on the modulation parameters
    \begin{equation}\label{eq:a priori estimate for the modulation parameters}
        \left|\frac{\lambda_s}{\lambda}\right|
        + |\tilde{\gmm}_s|
        \lesssim \|\epsilon\|_{\mathcal{H}^{1/2}}.
    \end{equation}
    \end{enumerate}
\end{proposition}

\begin{proof}
    Substituting the decomposition defined by \eqref{eq:definition of epsilon(t)} into \eqref{eq:halfwave}, and rewriting the resulting equation in the renormalized variables
    \eqref{eq:renomalised variable}, we obtain \eqref{eq:flow of ep(s)}. The decay $Q(y)\lesssim\langle y\rangle^{-2}$, together with \eqref{eq:quadratic control}, \eqref{eq:cubic control}, \eqref{eq:identity from mass and energy conservation law} and energy conservation, gives \eqref{eq:energy type control}.

    Next, taking the inner product of \eqref{eq:flow of ep(s)} with $\Lambda Q$, we obtain
    \begin{equation*}
        \frac{\lambda_s}{\lambda}\left(\norm{\Lambda Q}_{L^2}^2 - (\epsilon,\Lambda^2Q)_r\right)
        = (iL_Q[\epsilon],\Lambda Q)_r - (i\mathrm{NL}_{\geq 2}(\epsilon),\Lambda Q)_r + \tilde{\gmm}_s(i\epsilon,\Lambda Q)_r.
    \end{equation*}
    Since $(L_Q[\epsilon],\Lambda Q)_r = -(\epsilon, Q)_r$, by \eqref{eq:Lambda^k Q control} and \eqref{eq:linear control} we have
    \begin{align*}
        |(iL_Q[\epsilon],\Lambda Q)_r|
        = |(\epsilon, iQ)_r - 2(\epsilon, i Q^2\Lambda Q)_r|
        \lesssim \norm{\epsilon}_{\mathcal{H}^{1/2}}.
    \end{align*}
    Moreover, using \eqref{eq:linear control} and \eqref{eq:quadratic control}, and possibly reducing $\alpha^*$, we obtain
    \begin{equation}\label{eq:lambda control1}
        \left|\frac{\lambda_s}{\lambda}\right| \lesssim (1+|\tilde{\gmm}_s|)\norm{\epsilon}_{\mathcal{H}^{1/2}}.
    \end{equation}

    On the other hand, taking the inner product of \eqref{eq:flow of ep(s)} with $i\Lambda^2 Q$, we obtain 
    \begin{equation*}
        \tilde{\gmm}_s\left(\norm{\Lambda Q}_{L^2}^2 - (\epsilon, \Lambda^2 Q)_r\right)
        = (L_Q[\epsilon],\Lambda^2 Q)_r - (\mathrm{NL}_{\geq 2}(\epsilon),\Lambda^2 Q)_r + \frac{\lambda_s}{\lambda}(\epsilon,i\Lambda^3 Q)_r.
    \end{equation*}
    From \eqref{eq:Lambda^k Q control}, we have $|\Lambda^2 Q|(y) \lesssim \langle y \rangle^{-2}$. Furthermore, we have $|L_Q[\Lambda^2 Q]|(y) \lesssim \langle y \rangle^{-2}$ from \eqref{induction for k} for $k=2$. Thus, we obtain
    \begin{equation*}
        |(L_Q[\epsilon],\Lambda^2 Q)_r| \lesssim \norm{\epsilon}_{\mathcal{H}^{1/2}}.
    \end{equation*}
    Moreover, using \eqref{eq:Lambda^k Q control}, \eqref{eq:linear control}, and \eqref{eq:quadratic control}, and possibly reducing $\alpha^*>0$, we obtain 
    \begin{equation}\label{eq:gmm control1}
        |\tilde{\gmm}_s| \lesssim \left(1+\left|\frac{\lambda_s}{\lambda}\right|\right)\norm{\epsilon}_{\mathcal{H}^{1/2}}.
    \end{equation}
    Estimates~\eqref{eq:lambda control1} and \eqref{eq:gmm control1} give the \textit{a priori} bound \eqref{eq:a priori estimate for the modulation parameters}.
\end{proof}

\subsection{Quadratic form $\mathbf H$ and local virial estimates}
For the \(L^2\)-critical NLS, the virial identity together with the
negative-energy condition gives finite-time blow-up through a convexity
argument. For the half-wave equation, no comparable virial convexity
argument is available.

The evolution of the scale is associated with the quantity
\[
    (\partial_tu,i\Lambda u)_r.
\]
Here \(\Lambda\) is the generator of the \(L^2\)-scaling, and the real
inner product extracts the component of the flow in the scaling direction.
The half-wave equation has the Hamiltonian form
\(\partial_tu=-i\nabla E(u)\), and the energy is conserved. Hence, for
smooth and sufficiently decaying solutions,
\[
    (\partial_tu,i\Lambda u)_r = -(\nabla E(u),\Lambda u)_r = -\left.\frac{d}{da}E(u^{(a)})
    \right|_{a=1}=-E(u),
\]
where \( u^{(a)}(t,x)\coloneqq a^{1/2}u(t,ax)\). This formal computation is the starting point for the local virial identity used in the modulation analysis due to Merle--Rapha\"el
\cite{MerleRaphael2005AnnMath}. In the variables \((s,y)\) defined in \eqref{eq:renomalised variable}, expansion around the ground state \(Q\) expresses this identity through \(\partial_s(\epsilon,i\Lambda Q)_r\). Its leading quadratic term is the form \(\mathbf H\) defined in \eqref{eq:definition of H}, and the coercivity of \(\mathbf H\) in Theorem~\ref{thm:Spectral structure of the bilinear form H} gives \eqref{eq:raw-monotonicity}.

\begin{proposition}[Local virial identity]
    Let \(u\) be a solution to \eqref{eq:halfwave} whose initial datum
    satisfies \eqref{eq:condition1}. Then
    \begin{equation}\label{eq:lambda_s law}
        \partial_s(\epsilon,i\Lambda Q)_r
        =
        \frac12\mathbf H(\epsilon,\epsilon)
        +
        \lambda|E_0|
        +
        \mathcal O\left(
            \sqrt{\alpha_0}
            \norm{\epsilon}_{\mathcal H^{1/2}}^2
        \right),
    \end{equation}
    where
    \begin{equation}\label{eq:definition of H}
        \mathbf H(\epsilon,\epsilon)
        \coloneqq
        ([L_Q,\Lambda]\epsilon,\epsilon)_r.
    \end{equation}
    Moreover, after possibly reducing \(\alpha^*>0\), there exists a
    universal constant \(\delta_1>0\) such that
    \begin{equation}\label{eq:raw-monotonicity}
        \partial_s(\epsilon,i\Lambda Q)_r
        \geq
        \frac{\delta_1}{2}
        \norm{\epsilon}_{\mathcal H^{1/2}}^2
        +
        \lambda|E_0|
        -
        \frac{2}{\delta_1}
        \left[
            (\epsilon,Q)_r^2
            +
            (\epsilon,i\Lambda Q)_r^2
        \right].
    \end{equation}
\end{proposition}

\begin{proof}
    Taking the real inner product of \eqref{eq:flow of ep(s)} with
    \(i\Lambda Q\), and using \eqref{eq:orthogonal condition on ep}, we obtain
    \begin{equation}\label{eq:inner product with iLambda Q}
        \begin{aligned}
            \partial_s(\epsilon,i\Lambda Q)_r
            &=
            -(L_Q[\epsilon],\Lambda Q)_r
            -
            \frac{\lambda_s}{\lambda}
            (\epsilon,i\Lambda^2Q)_r
            -
            \tilde{\gmm}_s(\epsilon,\Lambda Q)_r
            +
            (\mathrm{NL}_{\geq 2}(\epsilon),\Lambda Q)_r
            \\
            &=
            (\epsilon,Q)_r
            +
            (\mathrm{NL}_{\geq 2}(\epsilon),\Lambda Q)_r.
        \end{aligned}
    \end{equation}
    Here, we used \(L_Q[\Lambda Q]=-Q\) and the orthogonality conditions
    in \eqref{eq:orthogonal condition on ep}.

    On the other hand, the energy identity in
    \eqref{eq:identity from mass and energy conservation law} gives
    \begin{equation}\label{eq:expansion of (epsilon,Q)_r}
        \begin{aligned}
            (\epsilon,Q)_r
            &=
            \frac12\norm{\epsilon}_{\dot H^{1/2}}^2
            -
            \frac12
            \int_{\bbR}
            \left[
                Q^2|\epsilon|^2
                +
                2Q^2\Re\{\epsilon\}^2
            \right]
            +
            \lambda|E_0|
            -
            \frac14
            \int_{\bbR}
            \mathcal N_{\geq 3}[\epsilon].
        \end{aligned}
    \end{equation}
    Moreover, we have
    \begin{equation*}
        \begin{aligned}
            \frac12\mathbf H(\epsilon,\epsilon)
            &=
            \frac12\norm{\epsilon}_{\dot H^{1/2}}^2
            -
            \frac12
            \int_{\bbR}
            \left[
                Q^2|\epsilon|^2
                +
                2Q^2\Re\{\epsilon\}^2
            \right]
            +
            (\mathrm{NL}_2(\epsilon),\Lambda Q)_r,
        \end{aligned}
    \end{equation*}
    where \(\mathrm{NL}_2(\epsilon) \coloneqq \mathrm{NL}_{\geq 2}(\epsilon) - |\epsilon|^2\epsilon\). Substituting \eqref{eq:expansion of (epsilon,Q)_r} into \eqref{eq:inner product with iLambda Q} yields
    \begin{equation}\label{eq:d/ds(ep,iLambda Q) sim H}
        \partial_s(\epsilon,i\Lambda Q)_r
        =
        \frac12\mathbf H(\epsilon,\epsilon)
        +
        \lambda|E_0|
        -
        \frac14\int_{\bbR}\mathcal N_{\geq 3}[\epsilon]
        +
        (|\epsilon|^2\epsilon,\Lambda Q)_r.
    \end{equation}
    Now, applying \eqref{eq:cubic control} to the last two terms in
    \eqref{eq:d/ds(ep,iLambda Q) sim H} yields \eqref{eq:lambda_s law}.

    The identity \eqref{eq:lambda_s law}, Theorem~\ref{thm:Spectral structure of the bilinear form H}, and the orthogonality conditions \eqref{eq:orthogonal condition on ep} give 
    \eqref{eq:raw-monotonicity} after possibly reducing \(\alpha^*>0\).
\end{proof}

The term \((\epsilon,Q)_r^2\) in
\eqref{eq:raw-monotonicity} is removed by differentiating \((\epsilon,Q)_r(\epsilon,i\Lambda Q)_r\). This gives the following local virial estimate.

\begin{lemma}[Local virial estimate]
    Possibly after reducing \(\alpha^*>0\), there exists a universal constant \(\delta_1>0\) such that, for every \(\alpha_0\in (0,\alpha^*)\), the following hold.

    \smallskip

    \noindent
    \textnormal{(i)} For all \(s \geq 0\),
    \begin{equation}\label{eq:corrected-repulsivity}
        \partial_s\left[\left(1+\frac{2}{\delta_1}(\epsilon,Q)_r\right)(\epsilon,i\Lambda Q)_r\right]
        \geq \frac{\delta_1}{4}\norm{\epsilon}_{\mathcal H^{1/2}}^2
        + \lambda|E_0|
        - \frac{4}{\delta_1}(\epsilon,i\Lambda Q)_r^2.
    \end{equation}

    \noindent
    \textnormal{(ii)} For all \(s_2\geq s_1 \geq 0\),
    \begin{equation}\label{eq:corrected-repulsivity-integrated}
        \begin{aligned}
            &\left[\left(1+\frac{2}{\delta_1}(\epsilon,Q)_r\right)(\epsilon,i\Lambda Q)_r\right]_{s_1}^{s_2} \\
            &\qquad\geq
            \frac{\delta_1}{4}\int_{s_1}^{s_2}\norm{\epsilon}_{\mathcal H^{1/2}}^2\,ds
            + |E_0|\int_{s_1}^{s_2}\lambda(s)\,ds
            - \frac{4}{\delta_1}\int_{s_1}^{s_2}(\epsilon,i\Lambda Q)_r^2\,ds.
        \end{aligned}
    \end{equation}
\end{lemma}

\begin{proof}
    Taking the inner product of \eqref{eq:flow of ep(s)} with \(Q\), and using \eqref{eq:quadratic control}, we obtain
    \begin{equation}\label{eq:epQ-law}
        \partial_s(\epsilon,Q)_r
        = -\frac{\lambda_s}{\lambda}(\epsilon,\Lambda Q)_r
        + \tilde{\gmm}_s(\epsilon,iQ)_r
        + \mathcal{O}\bigl(\norm{\epsilon}_{\mathcal H^{1/2}}^2\bigr).
    \end{equation}
    The identities \eqref{eq:inner product with iLambda Q} and \eqref{eq:epQ-law}, together with the estimates \eqref{eq:a priori estimate for the modulation parameters},
    \eqref{eq:linear control}, and \eqref{eq:refined bound of ep(t)}, give
    \begin{equation}\label{eq:product-derivative-estimate}
        \left|
        \partial_s\left[(\epsilon,Q)_r(\epsilon,i\Lambda Q)_r\right]
        - (\epsilon,Q)_r^2
        \right|
        \lesssim \sqrt{\alpha_0}\norm{\epsilon}_{\mathcal H^{1/2}}^2.
    \end{equation}

    We now differentiate the corrected quantity:
    \[
        \partial_s\left[\left(1+\frac{2}{\delta_1}(\epsilon,Q)_r\right)(\epsilon,i\Lambda Q)_r\right]
        =
        \partial_s(\epsilon,i\Lambda Q)_r
        + \frac{2}{\delta_1}\partial_s\Bigl\{(\epsilon,Q)_r(\epsilon,i\Lambda Q)_r\Bigr\}.
    \]
    Combining \eqref{eq:raw-monotonicity} and \eqref{eq:product-derivative-estimate}, we obtain
    \begin{align*}
        \partial_s\left[\left(1+\frac{2}{\delta_1}(\epsilon,Q)_r\right)(\epsilon,i\Lambda Q)_r\right] 
        \geq
        \frac{\delta_1}{2}\norm{\epsilon}_{\mathcal H^{1/2}}^2
        + \lambda|E_0|
        - \frac{2}{\delta_1}(\epsilon,i\Lambda Q)_r^2
        - C\sqrt{\alpha_0}\norm{\epsilon}_{\mathcal H^{1/2}}^2.
    \end{align*}
    After possibly reducing \(\alpha^*>0\), the last term is absorbed into
    \(\frac{\delta_1}{2}\norm{\epsilon}_{\mathcal H^{1/2}}^2\), and we arrive at
    \eqref{eq:corrected-repulsivity}. Integrating \eqref{eq:corrected-repulsivity} from
    \(s_1\) to \(s_2\) gives \eqref{eq:corrected-repulsivity-integrated}.
\end{proof}

\subsection{Almost monotonicity of the scaling parameter}
In this subsection, we use the monotonicity estimate \eqref{eq:corrected-repulsivity} to analyze the sign of \((\epsilon,i\Lambda Q)_r\) and relate its time integral to the variation of \(\log\lambda\). We then derive the almost monotonicity of the scaling parameter.

\begin{proposition}[Almost monotonicity of the scaling parameter]
    Possibly reducing \(\alpha^*>0\), for every
    \(\alpha_0\in(0,\alpha^*)\), there exists a unique
    \(s_0\in[0,\infty)\) such that the following properties hold.

    \begin{enumerate}
        \item
        We have
        \begin{equation}\label{eq:unique vanishing point}
            \begin{cases}
                (\epsilon,i\Lambda Q)_r(s)<0,
                & 0\leq s<s_0,\\
                (\epsilon,i\Lambda Q)_r(s_0)\geq0,
                & s=s_0,\\
                (\epsilon,i\Lambda Q)_r(s)>0,
                & s>s_0.
            \end{cases}
        \end{equation}
        Moreover, if \(s_0>0\), then
        \[
            (\epsilon,i\Lambda Q)_r(s_0)=0.
        \]

        \item
        For every \(s_2\geq s_1\geq s_0\), we have
        \begin{equation}\label{eq:lambda_s/lambda sim (ep,iLambda Q)}
            \frac{1}{2}
            \int_{s_1}^{s_2}(\epsilon,i\Lambda Q)_r\,ds
            -
            C\sqrt{\alpha_0}
            \leq
            -e_1\log\left(\frac{\lambda(s_2)}{\lambda(s_1)}\right)
            \leq
            \frac{3}{2}
            \int_{s_1}^{s_2}(\epsilon,i\Lambda Q)_r\,ds
            +
            C\sqrt{\alpha_0},
        \end{equation}
        and
        \begin{equation}\label{eq:almost monotonicity of lambda}
            \lambda(s_2)<2\lambda(s_1).
        \end{equation}
    \end{enumerate}
\end{proposition}

\begin{proof}
    Taking the real inner product of \eqref{eq:flow of ep(s)} with $R_1$, and using \eqref{eq:R_1-positivity, e_1 definition}, \eqref{eq:R_1-decay}, \eqref{eq:linear control}, \eqref{eq:quadratic control}, and \eqref{eq:a priori estimate for the modulation parameters}, gives
    \begin{equation}\label{eq:innerproduct with S_1}
        \frac{d}{ds}(\epsilon,R_1)_r = (\epsilon,i\Lambda Q)_r + e_1\frac{\lambda_s}{\lambda} + \mathcal{O}(\norm{\epsilon}_{\mathcal{H}^{1/2}}^2).
    \end{equation}
    Integrating \eqref{eq:innerproduct with S_1} over $[s_1,s_2]$ and using \eqref{eq:refined bound of ep(t)} gives
    \begin{equation}\label{eq:monotonicity ineq 1}
        \left|\int_{s_1}^{s_2} (\epsilon,i\Lambda Q)_r+e_1\log\left(\frac{\lambda(s_2)}{\lambda(s_1)}\right)\right| \lesssim \sqrt{\alpha_0} + \int_{s_1}^{s_2}\norm{\epsilon}_{\mathcal{H}^{1/2}}^2.
    \end{equation}
    Possibly reducing $\alpha^*>0$, for $\alpha_0 \in (0,\alpha^*)$ so that 
    \begin{equation}\label{eq:smallness of alpha^* 1}
        \left|\frac{2}{\delta_1}(\epsilon,Q)_r\right| < \frac{1}{2},
    \end{equation}
    from \eqref{eq:corrected-repulsivity-integrated} and \eqref{eq:refined bound of ep(t)} we obtain
    \begin{equation}\label{eq:monotonicity ineq 2}
        \begin{aligned}
            \int_{s_1}^{s_2}\norm{\epsilon}_{\mathcal{H}^{1/2}}^2(s)\,ds &\leq \frac{16}{\delta_1^2}\int_{s_1}^{s_2} (\epsilon,i\Lambda Q)_r^2(s)\,ds - \frac{4|E_0|}{\delta_1}\int_{s_1}^{s_2} \lambda(s)\,ds 
            \\
            &\quad + \frac{12}{\delta_1}\max_{s=s_1,\,s_2}|(\epsilon,i\Lambda Q)_r(s)| \\
            & \leq C(\delta_1) \left(\int_{s_1}^{s_2}(\epsilon,i\Lambda Q)_r^2(s) \,ds + \sqrt{\alpha_0}\right),
        \end{aligned}
    \end{equation}
    where $C(\delta_1)>0$ depends only on $\delta_1$.
    From \eqref{eq:monotonicity ineq 1} and \eqref{eq:monotonicity ineq 2} and possibly reducing $\alpha^*>0$ for $\alpha_0 \in (0,\alpha^*)$ we obtain
    \begin{equation}\label{eq:monotonicity ineq 3}
        \left|\int_{s_1}^{s_2} (\epsilon,i\Lambda Q)_r(s)\,ds+e_1\log\left(\frac{\lambda(s_2)}{\lambda(s_1)}\right)\right| \leq C\sqrt{\alpha_0} + \frac{1}{2}\int_{s_1}^{s_2}|(\epsilon,i\Lambda Q)_r(s)|\,ds,
    \end{equation}
    for some constant $C>0$.

        We now prove \eqref{eq:unique vanishing point}. We first claim that if
    \(s'\in[0,\infty)\) satisfies
    \[
        (\epsilon,i\Lambda Q)_r(s')=0,
    \]
    then
    \begin{equation}\label{eq:d_s(ep,iLambda Q)>0}
        \frac{d}{ds}(\epsilon,i\Lambda Q)_r(s')>0.
    \end{equation}
    Indeed, evaluating \eqref{eq:corrected-repulsivity} at \(s=s'\), we obtain
    \[
        \left(
            1+\frac{2}{\delta_1}(\epsilon,Q)_r(s')
        \right)
        \frac{d}{ds}(\epsilon,i\Lambda Q)_r(s')
        \geq
        \frac{\delta_1}{4}
        \norm{\epsilon(s')}_{\mathcal H^{1/2}}^2
        +
        \lambda(s')|E_0|
        >0.
    \]
    Since \eqref{eq:smallness of alpha^* 1} gives
    \[
        1+\frac{2}{\delta_1}(\epsilon,Q)_r(s')>\frac12,
    \]
    the claim follows. Since \((\epsilon,i\Lambda Q)_r\) is of class \(C^1\), 
    \eqref{eq:d_s(ep,iLambda Q)>0} shows that every zero is a strict crossing
    from negative to positive. If there were two distinct zeros, continuity
    would force an intermediate crossing from positive to negative, which is
    impossible. Hence, \((\epsilon,i\Lambda Q)_r\) has at most one zero on
    \([0,\infty)\).

    Suppose first that
    \[
        (\epsilon,i\Lambda Q)_r(0)\geq0.
    \]
    We set \(s_0=0\). By \eqref{eq:d_s(ep,iLambda Q)>0} and the fact that
    there is at most one zero, we have
    \[
        (\epsilon,i\Lambda Q)_r(s)>0,
        \qquad s>0.
    \]

    It remains to consider the case
    \[
        (\epsilon,i\Lambda Q)_r(0)<0.
    \]
    We claim that \((\epsilon,i\Lambda Q)_r\) must vanish at some positive
    time. Suppose, towards a contradiction, that
    \begin{equation}\label{eq:sign condition <0}
        (\epsilon,i\Lambda Q)_r(s)<0,
        \qquad s\geq0.
    \end{equation}
    Applying \eqref{eq:monotonicity ineq 3} with \(s_1=0\) and \(s_2=s\),
    we obtain
    \[
        \frac12
        \int_0^s
            |(\epsilon,i\Lambda Q)_r(s')|
        \,ds'
        -
        C\sqrt{\alpha_0}
        \leq
        e_1\log\left(\frac{\lambda(s)}{\lambda(0)}\right)
        \leq
        \frac32
        \int_0^s
            |(\epsilon,i\Lambda Q)_r(s')|
        \,ds'
        +
        C\sqrt{\alpha_0}.
    \]
    If
    \[
        \int_0^\infty
            |(\epsilon,i\Lambda Q)_r(s)|
        \,ds
        =
        \infty,
    \]
    then \(\lambda(s)\to\infty\) as \(s\to\infty\). By
    \eqref{eq:smallenss of ep and lambda esitmate}, this would imply
    \(\norm{u(t)}_{\dot H^{1/2}}\to0\). Moreover, the Gagliardo--Nirenberg inequality and mass conservation give
    \[
        \norm{u(t)}_{L^4}^4
        \lesssim
        \norm{u_0}_{L^2}^2
        \norm{u(t)}_{\dot H^{1/2}}^2
        \to0.
    \]
    Hence \(E(u(t))\to0\), contradicting the energy conservation \(E(u(t))=E_0<0\). Therefore,
    \begin{equation}\label{eq:(ep,iLambda Q)_r in L^1}
        \int_0^\infty
            |(\epsilon,i\Lambda Q)_r(s)|
        \,ds
        <
        \infty.
    \end{equation}
    Combining \eqref{eq:(ep,iLambda Q)_r in L^1} with
    \eqref{eq:monotonicity ineq 3}, we obtain constants
    \(\lambda_1,\lambda_2>0\) such that
    \begin{equation}\label{eq:bouondedness of lambda(s)}
        0<\lambda_1\leq\lambda(s)\leq\lambda_2<\infty,
        \qquad s\geq0.
    \end{equation}

    Moreover, \eqref{eq:inner product with iLambda Q} and
    \eqref{eq:refined bound of ep(t)} imply that
    \(\partial_s(\epsilon,i\Lambda Q)_r\) is uniformly bounded. Hence
    \((\epsilon,i\Lambda Q)_r\) is uniformly continuous on
    \([0,\infty)\). Together with \eqref{eq:(ep,iLambda Q)_r in L^1}, this yields
    \begin{equation}\label{eq:asymptotic behavior of (ep,iLambda Q)}
        (\epsilon,i\Lambda Q)_r(s)\to0
        \qquad\text{as }s\to\infty.
    \end{equation}
    On the other hand, by \eqref{eq:corrected-repulsivity},
    \[
        \begin{aligned}
            \partial_s
            \left[
                \left(
                    1+\frac{2}{\delta_1}(\epsilon,Q)_r
                \right)
                (\epsilon,i\Lambda Q)_r
            \right]
            &\geq
            \lambda(s)|E_0|
            -
            \frac{4}{\delta_1}
            (\epsilon,i\Lambda Q)_r^2.
        \end{aligned}
    \]
    In view of \eqref{eq:bouondedness of lambda(s)} and
    \eqref{eq:asymptotic behavior of (ep,iLambda Q)}, there exists
    \(s_*>0\) such that, for every \(s\geq s_*\),
    \[
        \partial_s
        \left[
            \left(
                1+\frac{2}{\delta_1}(\epsilon,Q)_r
            \right)
            (\epsilon,i\Lambda Q)_r
        \right]
        \geq
        \frac12\lambda_1|E_0|
        >0.
    \]
    This contradicts \eqref{eq:sign condition <0},
    \eqref{eq:smallness of alpha^* 1}, and
    \eqref{eq:asymptotic behavior of (ep,iLambda Q)}, since the quantity
    inside the derivative is negative and converges to zero. Therefore,
    \((\epsilon,i\Lambda Q)_r\) has a zero \(s_0>0\). The uniqueness of the
    zero and \eqref{eq:d_s(ep,iLambda Q)>0} yield
    \[
        (\epsilon,i\Lambda Q)_r(s)<0
        \quad\text{for }0\leq s<s_0,
        \qquad
        (\epsilon,i\Lambda Q)_r(s)>0
        \quad\text{for }s>s_0.
    \]
    This proves \eqref{eq:unique vanishing point}.
    For \(s_2\geq s_1\geq s_0\), we have
    \[
        (\epsilon,i\Lambda Q)_r(s)\geq0,
        \qquad s\in[s_1,s_2].
    \]
    Hence \eqref{eq:monotonicity ineq 3} gives
    \eqref{eq:lambda_s/lambda sim (ep,iLambda Q)}.
    Finally, we prove \eqref{eq:almost monotonicity of lambda}. Suppose, towards
    a contradiction, that there exist \(s_0\leq s_1'<s_2'\) such that
    \[
        \lambda(s_2')\geq2\lambda(s_1').
    \]
    By \eqref{eq:lambda_s/lambda sim (ep,iLambda Q)},
    \[
        \begin{aligned}
            e_1\log2
            \leq
            e_1\log\left(\frac{\lambda(s_2')}{\lambda(s_1')}\right)
            \leq
            C\sqrt{\alpha_0}
            -
            \frac12
            \int_{s_1'}^{s_2'}
                (\epsilon,i\Lambda Q)_r(s)
            \,ds
            \leq
            C\sqrt{\alpha_0}.
        \end{aligned}
    \]
    After possibly reducing \(\alpha^*>0\) so that
    \(C\sqrt{\alpha^*}<e_1\log2\), this is a contradiction. Therefore,
    \eqref{eq:almost monotonicity of lambda} holds.
\end{proof}

By time-translation invariance, we may assume that $s_0=0$.

\subsection{Finite-time blow-up and the logarithmic upper bound}
In this subsection, we complete the proof of the main theorem. We first prove that
\[
    \lambda(s)\to0
    \qquad
    \text{as }s\to\infty,
\]
which shows that the solution concentrates at arbitrarily small scales. We then introduce a refined correction in the generalized kernel direction \(R_1\) and derive an improved differential inequality for the scaling parameter. Together with the almost monotonicity of \(\lambda\), this gives finite-time blow-up and the logarithmic upper bound on the blow-up rate.

\begin{lemma}\label{lem:finite or infinite time blow up}
    Under \eqref{eq:condition1}, we have
    \begin{equation}\label{eq:lambda goes to zero as s to infty}
        \lim_{s\to \infty} \lambda(s) = 0.
    \end{equation}
\end{lemma}

\begin{proof}
    Suppose that \eqref{eq:lambda goes to zero as s to infty} fails.
    Then there exist \(s_n\to\infty\) and \(\lambda_0>0\) such that
    \(\lambda(s_n)\geq\lambda_0\). By \eqref{eq:almost monotonicity of lambda},
    \(\lambda(s)\geq \tfrac{\lambda_0}{2}\) for all \(s\geq0\).

    Combining \eqref{eq:energy-lower-bound} and \eqref{eq:smallenss of ep and lambda esitmate}, we obtain \(\lambda(s)\lesssim_{u_0}1\). Since \((\epsilon,i\Lambda Q)_r\geq0\),
    \eqref{eq:lambda_s/lambda sim (ep,iLambda Q)} implies
    \[
        \int_0^\infty(\epsilon,i\Lambda Q)_r\,ds<\infty,
        \qquad
        \int_0^\infty(\epsilon,i\Lambda Q)_r^2\,ds<\infty,
    \]
    where the second bound follows from \eqref{eq:refined bound of ep(t)}.
    Integrating \eqref{eq:corrected-repulsivity} over \([0,S]\) gives
    \[
        \left[
            \left(
                1+\frac{2}{\delta_1}(\epsilon,Q)_r
            \right)
            (\epsilon,i\Lambda Q)_r
        \right]_0^S
        \geq
        \frac{\lambda_0|E_0|}{2}S
        -
        \frac4{\delta_1}
        \int_0^\infty(\epsilon,i\Lambda Q)_r^2\,ds.
    \]
    The left-hand side is bounded by \eqref{eq:refined bound of ep(t)}, whereas the right-hand side tends to \(+\infty\). This contradiction proves \eqref{eq:lambda goes to zero as s to infty}.
\end{proof}

The preceding estimates do not yet give the blow-up rate bound in Theorem~\ref{thm:main}. We therefore refine the monotonicity argument using a modification of the generalized kernel
element \(R_1\) that satisfies the orthogonality condition \eqref{eq:orthogonal condition on ep}.

We introduce the modified first corrector
\begin{equation*}
    R_1^\perp
    \coloneqq
    R_1
    +
    \frac{
        (R_1,\Lambda^2Q)_r
    }{
        \norm{\Lambda Q}_{L^2}^2
    }Q.
\end{equation*}
Since \((Q,\Lambda^2Q)_r = -\norm{\Lambda Q}_{L^2}^2\), \((Q,\Lambda Q)_r=0\), and \(L_Q[iQ]=0\), the profile \(R_1^\perp\) satisfies
\begin{equation}\label{eq:R1-perp-identities}
    (R_1^\perp,\Lambda^2Q)_r=0,
    \qquad
    (R_1^\perp,\Lambda Q)_r=e_1,
    \qquad
    L_Q[iR_1^\perp]=i\Lambda Q.
\end{equation}
We define the scaling coefficient, the corrected perturbation, and the
corrected scaling coefficient by
\begin{equation}\label{eq:def-b-tilde-b-and-refined-epsilon}
    b
    \coloneqq
    \frac{(\epsilon,i\Lambda Q)_r}{e_1},
    \qquad
    \tilde{\epsilon}
    \coloneqq
    \epsilon-ibR_1^\perp,
    \qquad
    \tilde b
    \coloneqq
    e_1\left(
        1-
        \frac{(\epsilon,\Lambda R_1^\perp)_r}{e_1}
    \right)b.
\end{equation}

We now state the refined monotonicity estimate satisfied by
\(\tilde{\epsilon}\).

\begin{proposition}[Refined local virial estimate]\label{prop:refined coercivity w.r.t b}
    There exist universal constants \(C,\delta_2>0\) such that, for every
    \(\alpha_0\in(0,\alpha^*)\), there exists \(\tilde{s}_2'>0\) for which
    \begin{equation}\label{eq:refined coercivity w.r.t b}
        \tilde b_s
        +
        Cb^4
        \geq
        \frac12\lambda|E_0|
        +
        \frac{\delta_2}{2}
        \norm{\tilde{\epsilon}}_{\mathcal H^{1/2}}^2
    \end{equation}
    for every \(s\geq\tilde{s}_2'\).
\end{proposition}

\begin{proof}
    We first record the identities
    \[
        [L_Q,\Lambda](iR_1^\perp)
        =
        L_Q[i\Lambda R_1^\perp]
        -
        i\Lambda^2Q,
        \qquad 
        \mathbf H(iR_1^\perp,iR_1^\perp)=0.
    \]
    Indeed, using \(L_Q[iR_1^\perp]=i\Lambda Q\), the self-adjointness of
    \(L_Q\), and \((R_1^\perp,\Lambda^2Q)_r=0\), we obtain
    \[
        \begin{aligned}
            \mathbf H(iR_1^\perp,iR_1^\perp)
            =
            \bigl(
                L_Q[i\Lambda R_1^\perp]-i\Lambda^2Q,
                iR_1^\perp
            \bigr)_r
            =
            (\Lambda R_1^\perp,\Lambda Q)_r
            -
            (\Lambda^2Q,R_1^\perp)_r
            =
            0.
        \end{aligned}
    \]
    Moreover, by \eqref{eq:R1-perp-identities}, \eqref{eq:def-b-tilde-b-and-refined-epsilon}, and the orthogonality conditions in \eqref{eq:orthogonal condition on ep},
    \[
        (\tilde{\epsilon},\Lambda Q)_r
        =
        (\tilde{\epsilon},i\Lambda Q)_r
        =
        (\tilde{\epsilon},i\Lambda^2Q)_r
        =
        0,
        \qquad
        (\tilde{\epsilon},Q)_r=(\epsilon,Q)_r.
    \]
    Expanding the commutator form around \(\tilde{\epsilon}\), we therefore
    obtain
    \begin{equation}\label{eq:H-epsilon-refined-expansion}
        \mathbf H(\epsilon,\epsilon)
        =
        \mathbf H(\tilde{\epsilon},\tilde{\epsilon})
        +
        2b
        \bigl(
            \epsilon,
            L_Q[i\Lambda R_1^\perp]
        \bigr)_r.
    \end{equation}
    Substituting \eqref{eq:H-epsilon-refined-expansion} into
    \eqref{eq:d/ds(ep,iLambda Q) sim H} gives
    \begin{equation}\label{eq:refined-scaling-identity}
        \begin{aligned}
            \partial_s(\epsilon,i\Lambda Q)_r
            =
            \lambda|E_0|
            +
            \frac12
            \mathbf H(\tilde{\epsilon},\tilde{\epsilon})
            +
            b
            \bigl(
                \epsilon,
                L_Q[i\Lambda R_1^\perp]
            \bigr)_r
            -
            \frac14
            \int_{\bbR}
                \mathcal N_{\geq3}[\epsilon]
            +
            (|\epsilon|^2\epsilon,\Lambda Q)_r.
        \end{aligned}
    \end{equation}
    By Theorem~\ref{thm:Spectral structure of the bilinear form H}, after
    decreasing \(\delta_2>0\) if necessary, we have
    \begin{equation}\label{eq:H-refined-coercivity}
        \frac12
        \mathbf H(\tilde{\epsilon},\tilde{\epsilon})
        \geq
        2\delta_2
        \norm{\tilde{\epsilon}}_{\mathcal H^{1/2}}^2
        -
        C(\epsilon,Q)_r^2.
    \end{equation}
    Since \(\epsilon = \tilde{\epsilon} + ibR_1^\perp\), and \(Q\), \(\Lambda Q\), and \(R_1^\perp\) are real-valued, every cubic term contains at least one factor of \(\tilde\epsilon\), while the contribution containing only \(ibR_1^\perp\) is of order \(b^4\). The decay estimates for \(Q\) and \(R_1^\perp\), together with \eqref{eq:linear control}--\eqref{eq:cubic control}, therefore give
    \[
        \begin{aligned}
            \left|
                -
                \frac14
                \int_{\bbR}
                    \mathcal N_{\geq3}[\epsilon]
                +
                (|\epsilon|^2\epsilon,\Lambda Q)_r
            \right|
            &\lesssim
            \norm{\tilde{\epsilon}}_{\mathcal H^{1/2}}
            \left(
                \norm{\tilde{\epsilon}}_{\mathcal H^{1/2}}
                +
                |b|
            \right)^2
            +
            \left(
                \norm{\tilde{\epsilon}}_{\mathcal H^{1/2}}
                +
                |b|
            \right)^4
            \\
            &\leq
            \frac{\delta_2}{4}
            \norm{\tilde{\epsilon}}_{\mathcal H^{1/2}}^2
            +
            Cb^4,
        \end{aligned}
    \]
    after possibly reducing \(\alpha^*>0\). We next control the remaining \(Q\)-direction. By \eqref{eq:energy type control},
    \[
        |(\epsilon,Q)_r|
        \lesssim
        \lambda|E_0|
        +
        \norm{\epsilon}_{\mathcal H^{1/2}}^2.
    \]
    Since \(\norm{\epsilon}_{\mathcal H^{1/2}} \lesssim \norm{\tilde{\epsilon}}_{\mathcal H^{1/2}} + |b|\) and \(\lambda(s)\to0\) as \(s\to\infty\), for every sufficiently small
    fixed constant \(\eta>0\), there exists \(\tilde{s}'>0\) such that
    \begin{equation}\label{eq:epsilon-Q-refined-bound}
        |(\epsilon,Q)_r|^2
        \leq
        \eta
        \left(
            \lambda|E_0|
            +
            \norm{\tilde{\epsilon}}_{\mathcal H^{1/2}}^2
        \right)
        +
        C_\eta b^4
    \end{equation}
    for every \(s\geq\tilde{s}'\). From \eqref{eq:refined-scaling-identity}, \eqref{eq:H-refined-coercivity}, and \eqref{eq:epsilon-Q-refined-bound}, after choosing \(\eta>0\) sufficiently small, we obtain
    \begin{equation}\label{refined monotonicity ineq 1}
        \partial_s(\epsilon,i\Lambda Q)_r
        +
        Cb^4
        \geq
        \frac34\lambda|E_0|
        +
        \delta_2
        \norm{\tilde{\epsilon}}_{\mathcal H^{1/2}}^2
        +
        b
        \bigl(
            \epsilon,
            L_Q[i\Lambda R_1^\perp]
        \bigr)_r
    \end{equation}
    for all sufficiently large \(s\).

    We now rewrite the last term as a time derivative. Taking the real inner
    product of \eqref{eq:flow of ep(s)} with \(\Lambda R_1^\perp\), and using
    \((R_1^\perp,\Lambda^2Q)_r=0\), we obtain
    \begin{equation}\label{refined monotonicity ineq 2}
        \begin{aligned}
            \bigl(
                \epsilon,
                L_Q[i\Lambda R_1^\perp]
            \bigr)_r
            =
            \partial_s
            (\epsilon,\Lambda R_1^\perp)_r
            +
            \frac{\lambda_s}{\lambda}
            (\epsilon,\Lambda^2R_1^\perp)_r
            -
            \tilde{\gmm}_s
            (\epsilon,i\Lambda R_1^\perp)_r
            +
            \bigl(
                \mathrm{NL}_{\geq2}(\epsilon),
                i\Lambda R_1^\perp
            \bigr)_r.
        \end{aligned}
    \end{equation}
    From the definitions of $b$ and $\tilde b$ in \eqref{eq:def-b-tilde-b-and-refined-epsilon} the product rule yields
    \[
        \begin{aligned}
            \tilde b_s
            =
            \partial_s(\epsilon,i\Lambda Q)_r
            -
            b
            \partial_s(\epsilon,\Lambda R_1^\perp)_r
            -
            \frac{(\epsilon,\Lambda R_1^\perp)_r}{e_1}
            \partial_s(\epsilon,i\Lambda Q)_r.
        \end{aligned}
    \]
    The product-rule identity above, together with \eqref{refined monotonicity ineq 1}, \eqref{refined monotonicity ineq 2}, and \eqref{eq:inner product with iLambda Q}, gives
    \begin{equation}\label{eq:corrected-refined-identity}
        \tilde b_s + Cb^4 
        \geq
        \frac34\lambda|E_0| + \delta_2 \norm{\tilde{\epsilon}}_{\mathcal H^{1/2}}^2 + 
        \mathcal{E}_b,
    \end{equation}
    where
    \begin{equation*}
        \begin{aligned}
            \mathcal E_b
            \coloneqq{}&
            b
            \left[
                \frac{\lambda_s}{\lambda}
                (\epsilon,\Lambda^2R_1^\perp)_r
                -
                \tilde{\gmm}_s
                (\epsilon,i\Lambda R_1^\perp)_r
                +
                \bigl(
                    \mathrm{NL}_{\geq2}(\epsilon),
                    i\Lambda R_1^\perp
                \bigr)_r
            \right]
            \\
            &-
            \frac{
                (\epsilon,\Lambda R_1^\perp)_r
            }{
                e_1
            }
            \left[
                (\epsilon,Q)_r
                +
                \bigl(
                    \mathrm{NL}_{\geq2}(\epsilon),
                    \Lambda Q
                \bigr)_r
            \right].
        \end{aligned}
    \end{equation*}
    Using \(\epsilon=\tilde{\epsilon}+ibR_1^\perp\) and decay of \(R_1^\perp\), together with
    \eqref{eq:linear control}, \eqref{eq:quadratic control}, and
    \eqref{eq:a priori estimate for the modulation parameters}, we obtain
    \[
        \begin{aligned}
            |\mathcal E_b|
            \lesssim
            \sqrt{\alpha_0}
            \norm{\tilde{\epsilon}}_{\mathcal H^{1/2}}^2
            +
            b^2\norm{\tilde{\epsilon}}_{\mathcal H^{1/2}}
            +
            b^4
            +
            \norm{\tilde{\epsilon}}_{\mathcal H^{1/2}}
            |(\epsilon,Q)_r|.
        \end{aligned}
    \]
    Hence, by \eqref{eq:epsilon-Q-refined-bound} and Young's inequality, after
    possibly reducing \(\alpha^*>0\) and increasing \(\tilde{s}_2'>0\),
    \begin{equation}\label{eq:corrected-refined-error-bound}
        |\mathcal E_b|
        \leq
        Cb^4
        +
        \frac14\lambda|E_0|
        +
        \frac{\delta_2}{2}
        \norm{\tilde{\epsilon}}_{\mathcal H^{1/2}}^2,
    \end{equation}
    for every \(s\geq\tilde{s}_2'\). Substituting \eqref{eq:corrected-refined-error-bound} into \eqref{eq:corrected-refined-identity} yields
    \eqref{eq:refined coercivity w.r.t b}.
\end{proof}

\begin{lemma}
    There exist universal constants \(B>0\) and \(s_3>0\) such that
    \begin{equation}\label{eq:refined differential ineq}
        \lambda(s)
        \leq
        \exp\left(-\frac{B}{b(s)^2}\right),
        \qquad
        s\geq s_3.
    \end{equation}
\end{lemma}

\begin{proof}
    Recall that \(b(s)>0\) for \(s>0\). By the refined bound on
    \(\epsilon\), after possibly reducing \(\alpha^*>0\), we have
    \begin{equation}\label{eq:tilde-b-comparison}
        \frac12 b(s)
        \leq
        \widetilde b(s)
        \leq
        \frac32 b(s).
    \end{equation}
    Proposition~\ref{prop:refined coercivity w.r.t b} therefore gives
    \[
        \widetilde b_s+C\widetilde b^4\geq0
    \]
    for all sufficiently large \(s\). Since \(\widetilde b>0\),
    \[
        \frac{d}{ds}\bigl(\widetilde b^{-3}\bigr)
        =
        -3\widetilde b^{-4}\widetilde b_s
        \leq C.
    \]
    Integrating this inequality and using \eqref{eq:tilde-b-comparison}, we
    obtain
    \begin{equation}\label{refined law for (ep,iLambda Q)}
        b(s)\gtrsim s^{-1/3}
    \end{equation}
    for all sufficiently large \(s\). On the other hand, since \((\epsilon,i\Lambda Q)_r=e_1b\), \eqref{eq:lambda_s/lambda sim (ep,iLambda Q)} gives
    \[
        -\log\lambda(s)
        \gtrsim
        \int_{s_3}^{s}b(s')\,ds'
        -C
        \gtrsim
        s^{2/3}
    \]
    for all sufficiently large \(s\). Together with
    \eqref{refined law for (ep,iLambda Q)}, this yields
    \[
        \frac1{b(s)^2}
        \lesssim
        s^{2/3}
        \lesssim
        -\log\lambda(s).
    \]
    After decreasing \(B>0\) if necessary, we obtain
    \eqref{eq:refined differential ineq}.
\end{proof}

\begin{proof}[Proof of Theorem~\ref{thm:main}]
    Let \(t_{n_0}\) be the first time such that
    \[
        \lambda(t_{n_0})=2^{-n_0}.
    \]
    Once \(t_n\) is defined, let \(t_{n+1}\) be the first time after
    \(t_n\) such that
    \[
        \lambda(t_{n+1})=2^{-(n+1)}.
    \]
    Set \(s_n\coloneqq s(t_n)\). These times exist since
    \(\lambda(s)\to0\) as \(s\to\infty\) from Lemma~\ref{lem:finite or infinite time blow up}. Now, by the definition of \(t_{n+1}\) and
    \eqref{eq:almost monotonicity of lambda}, for every \(s\in[s_n,s_{n+1}]\),
    \begin{equation}\label{eq:dyadic-scale-comparison}
        2^{-(n+1)}
        \leq
        \lambda(s)
        <
        2^{-(n-1)}.
    \end{equation}
    Since \((\epsilon,i\Lambda Q)_r=e_1b\),
    \eqref{eq:lambda_s/lambda sim (ep,iLambda Q)} gives
    \[
        \int_{s_n}^{s_{n+1}} b(s)\,ds
        \lesssim 1.
    \]
    On the other hand, \eqref{eq:refined differential ineq} implies
    \[
        b(s)
        \gtrsim
        \frac{1}{\sqrt{|\log\lambda(s)|}}
    \]
    for all sufficiently large \(s\). Hence, using \(\tfrac{ds}{dt}=\tfrac{1}{\lambda}\),
    \[
        1
        \gtrsim
        \int_{s_n}^{s_{n+1}}b(s)\,ds
        \gtrsim
        \int_{t_n}^{t_{n+1}}
            \frac{dt}{
                \lambda(t)\sqrt{|\log\lambda(t)|}
            }.
    \]
    By \eqref{eq:dyadic-scale-comparison},
    \[
        t_{n+1}-t_n
        \lesssim
        2^{-n}\sqrt{n}.
    \]
    Since \(\sum_{n\geq n_0}2^{-n}\sqrt n<\infty\), the sequence
    \(\{t_n\}\) converges to a finite limit. Moreover,
    \eqref{eq:smallenss of ep and lambda esitmate} gives
    \[
        \norm{u(t_n)}_{\dot H^{1/2}}^2
        \sim
        \lambda(t_n)^{-1}
        =
        2^n
        \longrightarrow
        \infty.
    \]
    The blow-up alternative therefore shows that this finite limit is the
    maximal forward time \(T\). Consequently, for \(t\in[t_n,t_{n+1}]\), \eqref{eq:dyadic-scale-comparison} gives
    \begin{equation}\label{eq:dyadic-time-tail}
        T-t \leq T-t_n
        \lesssim
        \sum_{k\geq n}2^{-k}\sqrt{k}
        \lesssim
        2^{-n}\sqrt n
        \lesssim 
        \lambda(t)\sqrt{|\log\lambda(t)|}.
    \end{equation}
    It remains to express the logarithm in terms of \(T-t\). Since \(\sqrt{|\log\lambda(t)|}\leq\lambda(t)^{-1/2}\), \eqref{eq:dyadic-time-tail} gives
    \(\lambda(t)\gtrsim(T-t)^2\), and hence \(|\log\lambda(t)|\lesssim|\log(T-t)|\). Therefore,
    \[
        T-t
        \lesssim
        \lambda(t)\sqrt{|\log\lambda(t)|}
        \lesssim
        \lambda(t)\sqrt{|\log(T-t)|}.
    \]
    Finally, \eqref{eq:smallenss of ep and lambda esitmate} yields
    \[
        \norm{u(t)}_{\dot H^{1/2}}
        \lesssim
        \lambda(t)^{-1/2}
        \lesssim
        \frac{
            |\log(T-t)|^{1/4}
        }{
            \sqrt{T-t}
        }.
    \]
    This proves Theorem~\ref{thm:main}.
\end{proof}

\part{Spectral reduction of the coercivity problem}\label{part:spectral-reduction}
In this part, we reduce Theorem~\ref{thm:Spectral structure of the bilinear form H} to Proposition~\ref{prop:reduced-BS-bound-Q} and Proposition~\ref{prop:certified-limiting-matrix-conditions}, whose proofs are completed in Part~\ref{part:spectral-estimates}, using the profile estimates from Part~\ref{part:ground-state-approximation} and the interval bounds from Part~\ref{part:interval-verification}

Writing \(\epsilon=\epsilon_1+i\epsilon_2\), we decompose \(\mathbf H\) into two scalar quadratic forms associated with the nonlocal operators
\[
    |D|+6yQ'Q,
    \qquad
    |D|+2yQ'Q.
\]
In Section~\ref{sec:coercivity-constrained-nonnegativity}, we reduce \eqref{positivity of H} to the non-negativity of \(\mathcal L_{c_j}\coloneqq|D|+c_jyQ'Q\) on \(X_j\), for constants \(c_1>6\), \(c_2>2\). In Section~\ref{sec:regularized-constrained-index}, we regularize \(\mathcal L_{c_j}\) to \(\mathcal L_{c_j}^{(\mu)}=|D|+\mu+c_jV_Q\) and apply the constrained index formula of Theorem~\ref{thm:constrained-index-formula} to further reduce this non-negativity to the Morse index and constraint matrix conditions of Proposition~\ref{prop:shifted-index-determinant-conditions}. Sections~\ref{sec:BS-principle} and~\ref{sec:limiting-constraint-matrices} reduce these conditions to Proposition~\ref{prop:reduced-BS-bound-Q} and Proposition~\ref{prop:certified-limiting-matrix-conditions}, respectively, via the Birman--Schwinger principle and the symmetric resolvent identity.

Throughout this part, all operators are considered on the real even subspace, and orthogonality conditions and Morse indices are understood with respect to the inner product \((\cdot,\cdot)_r\). We set
\[
    V_Q(y)\coloneqq yQ'(y)Q(y).
\]
We also define the two ordered constraint families by 
\[
    \Xi_1\coloneqq\{Q,\Lambda Q\},
    \qquad
    \Xi_2\coloneqq\{\Lambda Q,\Lambda^2Q\}.
\]
Writing
\[
    \Xi_j=\{\Psi_1^{(j)},\Psi_2^{(j)}\},
    \qquad
    j=1,2,
\]
we define the corresponding constrained spaces by
\[
    X_j
    \coloneqq
    \left\{
        f\in H^{1/2}_{\mathrm e}:
        (f,\Psi_m^{(j)})_r=0,\quad m=1,2
    \right\}.
\]

\section{From coercivity to constrained non-negativity}\label{sec:coercivity-constrained-nonnegativity}
In this section, we reduce the coercivity estimate \eqref{positivity of H}
to the non-negativity of \(|D|+c_jV_Q\) on \(X_j\), \(j=1,2\), for suitable
constants \(c_1>6\) and \(c_2>2\).

\begin{proposition}[Constrained non-negativity with larger constants]\label{prop:Non-negativity of D+cyQ'Q}
    There exist constants $c_1>6$ and $c_2>2$ such that the following holds.
    Let $\epsilon=\epsilon_1+i\epsilon_2\in H^{1/2}_{\mathrm e}$ satisfy
    \begin{equation*}
        \epsilon_j \in X_j, \quad \text{for} \quad j=1,2.
    \end{equation*}
    Then, for each $j\in \{1,2\}$ we have
    \[
        (\mathcal L_{c_j}\epsilon_j,\epsilon_j)_r\geq 0, \quad \text{where} \quad \mathcal L_{c_j}\coloneqq |D|+c_jV_Q.
    \]
\end{proposition}

We now show that Proposition~\ref{prop:Non-negativity of D+cyQ'Q} implies Theorem~\ref{thm:Spectral structure of the bilinear form H}.

\begin{proof}[Proof of Theorem~\ref{thm:Spectral structure of the bilinear form H} assuming Proposition~\ref{prop:Non-negativity of D+cyQ'Q}]
    We first prove the coercive estimate under the constraints \(\epsilon_j\in X_j\), \(j=1,2\). From \cite[Proposition 1.1]{FrankLenzmann2013Acta} we have $V_Q=yQ'(y)Q(y)\leq 0$ for all $y\in \bbR$. Since $\mathbf H(\epsilon,\epsilon) = ([L_Q,\Lambda]\epsilon,\epsilon)_r$, choose $\delta\in\left(\max\{\tfrac{6}{c_1},\tfrac{2}{c_2}\},1\right)$ so that
    \begin{equation*}
        \begin{aligned}
            \mathbf H(\epsilon,\epsilon) &= (|D|\epsilon_1+6V_Q\epsilon_1,\epsilon_1)_r + (|D|\epsilon_2+2V_Q\epsilon_2,\epsilon_2)_r \\
            &\geq \delta \left[(|D|\epsilon_1 + c_1V_Q\epsilon_1,\epsilon_1)_r + (|D|\epsilon_2 + c_2V_Q\epsilon_2,\epsilon_2)_r\right] \\
            &\quad + (1-\delta)\norm{\epsilon}_{\dot H^{1/2}}^2 + \min\{|c_1\delta-6|,|c_2\delta-2|\}\int |\epsilon|^2|V_Q|.
        \end{aligned}
    \end{equation*}
    Given $c_1, c_2$, by choosing $\delta>0$, there is a universal constant $\tilde \delta_0>0$ so that 
    \begin{equation*}
        \mathbf H(\epsilon,\epsilon) \geq \tilde \delta_0\cdot \left(\norm{\epsilon}_{\dot H^{1/2}}^2 + \int |\epsilon(y)|^2|V_Q(y)|\,dy\right).
    \end{equation*}
    Therefore, to prove Theorem~\ref{thm:Spectral structure of the bilinear form H} it suffices to show that
    \begin{equation}\label{eq: coercivity to non-negativity 1}
        \int |\epsilon(y)|^2\langle y \rangle^{-4}\,dy \lesssim  \norm{\epsilon}_{\dot H^{1/2}}^2 + \int |\epsilon(y)|^2|V_Q(y)|\,dy.
    \end{equation}
    To prove \eqref{eq: coercivity to non-negativity 1} we use a compactness argument.
    Let $\psi$ be a smooth function with $\mathrm{supp}\,\psi \subset [-2,2]$ such that
    $\psi(y) = 1$ for $y \in [-1,1]$ and $\int \psi = 0$. We decompose
    $\epsilon = \psi \epsilon + (1-\psi)\epsilon$. Since $|V_Q(y)|\sim\langle y \rangle^{-4}$ as $|y| \to \infty$ from Lemma~\ref{lem:tail-lower-yQpQ}, it suffices to show that
    \[
        \int |\psi\epsilon(y)|^2\langle y \rangle^{-4} \,dy
        \lesssim \text{(right-hand side of } \eqref{eq: coercivity to non-negativity 1}\text{)}.
    \]

    We now argue by contradiction. Suppose that there exists a sequence 
    $\{\epsilon^{(n)}\}_n \subset H^{1/2}$ such that
    $\norm{\psi\epsilon^{(n)}}_{L^2(\langle y \rangle^{-4}dy)} = 1$, but
    \[
        \norm{\epsilon^{(n)}}_{\dot{H}^{1/2}}^2 
        + \int |\epsilon^{(n)}(y)|^2 |V_Q(y)| \,dy \to 0 
        \quad\text{as } n \to \infty.
    \]
    Note that $\|\psi\epsilon^{(n)}\|_{L^2} \geq 1$. Moreover, by \eqref{eq:frac-comm} we have
    \begin{align*}
        \norm{\psi\epsilon^{(n)}}_{\dot{H}^{1/2}}^2 
        &\lesssim \norm{|D|^{1/2}(\psi\epsilon^{(n)}) - \epsilon^{(n)} |D|^{1/2}\psi}_{L^2}^2 
           + \norm{\epsilon^{(n)} |D|^{1/2}\psi}_{L^2}^2 \\
        &\lesssim \norm{\epsilon^{(n)}}_{\dot{H}^{1/2}}^2 
           + \int |\epsilon^{(n)}(y)|^2 \langle y \rangle^{-5} \,dy.
    \end{align*}
    In the last inequality, we used \(|(|D|^{1/2}\psi)(y)|\lesssim \langle y\rangle^{-5/2}\),
    which follows from \(\int_{\bbR}\psi=0\). Since \(\psi=1\) on \([-1,1]\) and \(|V_Q(y)|\gtrsim\langle y\rangle^{-4}\) outside this interval, \(\langle y\rangle^{-5} \lesssim|\psi(y)|^2\langle y\rangle^{-4}+|V_Q(y)|\). Hence \(\{\psi\epsilon^{(n)}\}_n\) is bounded in \(H^{1/2}(\bbR)\). After passing to a subsequence, \(\psi\epsilon^{(n)} \to  \psi\epsilon^{(\infty)}\) strongly in \(L^2(\bbR)\) and almost everywhere.

    Since \(\psi\) has compact support and \(\langle y\rangle^{-4}\) is bounded,
    \[
        \norm{
            \psi\epsilon^{(\infty)}
        }_{L^2(\langle y\rangle^{-4}dy)}
        =
        1.
    \]
    On the other hand, Fatou's lemma gives
    \[
        \int_{\bbR}
            |V_Q|
            |\psi\epsilon^{(\infty)}|^2
        \,dy
        \leq
        \norm{\psi}_{L^\infty}^2\liminf_{n\to\infty}
        \int_{\bbR}
            |V_Q|
            |\epsilon^{(n)}|^2
        \,dy
        =
        0.
    \]
    Since \(|V_Q|>0\) almost everywhere, this gives
    \(\psi\epsilon^{(\infty)}=0\), a contradiction.
    We have therefore proved that, whenever
    \(\epsilon_j\in X_j\) for \(j=1,2\),
    \begin{equation}\label{eq:constrained-coercivity-H}
        \mathbf H(\epsilon,\epsilon)
        \geq
        \widetilde\delta_0
        \norm{\epsilon}_{\mathcal H^{1/2}}^2
    \end{equation}
    for some universal constant \(\widetilde\delta_0>0\).
    It remains to remove the orthogonality conditions. Let
    \(\epsilon=\epsilon_1+i\epsilon_2\in H^{1/2}_{\mathrm e}\) be arbitrary,
    and let
    \[
        \epsilon_j^X
        \coloneqq
        P_{X_j}\epsilon_j,
        \qquad
        j=1,2,
        \qquad
        \epsilon^X
        \coloneqq
        \epsilon_1^X+i\epsilon_2^X,
    \]
    where \(P_{X_j}\) denotes the \(L^2\)-orthogonal projection onto \(X_j\).
    Since each constraint family \(\Xi_j\) is linearly independent, the
    corresponding Gram matrix is invertible, and hence
    \[
        \norm{\epsilon-\epsilon^X}_{\mathcal H^{1/2}}^2
        \lesssim
        (\epsilon,Q)_r^2
        +
        (\epsilon,\Lambda Q)_r^2
        +
        (\epsilon,i\Lambda Q)_r^2
        +
        (\epsilon,i\Lambda^2Q)_r^2.
    \]
    Applying \eqref{eq:constrained-coercivity-H} to \(\epsilon^X\), and using
    the continuity of \(\mathbf H\) on \(H^{1/2}_{\mathrm e}\), we obtain
    \[
        \mathbf H(\epsilon,\epsilon)
        \geq
        \frac{\widetilde\delta_0}{2}
        \norm{\epsilon^X}_{\mathcal H^{1/2}}^2
        -
        C
        \left[
            (\epsilon,Q)_r^2
            +
            (\epsilon,\Lambda Q)_r^2
            +
            (\epsilon,i\Lambda Q)_r^2
            +
            (\epsilon,i\Lambda^2Q)_r^2
        \right].
    \]
    Finally,
    \[
        \norm{\epsilon}_{\mathcal H^{1/2}}^2
        \lesssim
        \norm{\epsilon^X}_{\mathcal H^{1/2}}^2
        +
        \norm{\epsilon-\epsilon^X}_{\mathcal H^{1/2}}^2,
    \]
    and therefore, after decreasing \(\delta_0>0\) if necessary, we obtain \eqref{positivity of H}.
\end{proof}

\section{Reduction to regularized constrained Morse indices}\label{sec:regularized-constrained-index}
In this section, we reduce the non-negativity of \(\mathcal L_{c_j}\) on \(X_j\) to a constrained Morse index problem for the \(\mu\)-regularized operators \(\mathcal L_{c_j}^{(\mu)}\coloneqq|D|+\mu+c_jV_Q\), using the abstract constrained index formula of Theorem~\ref{thm:constrained-index-formula}. This reduces Proposition~\ref{prop:Non-negativity of D+cyQ'Q} to the full Morse index and constraint matrix conditions of Proposition~\ref{prop:shifted-index-determinant-conditions} for the regularized operators.

\begin{theorem}[Constrained index formula {\cite[Theorem~1.5]{TranZhou23}}]
\label{thm:constrained-index-formula}
    Let $\mathscr H$ be a real Hilbert space and let
    $S:\mathscr H\times\mathscr H\to\mathbb R$ be a continuous symmetric bilinear
    form with $n_-(S)<\infty$. Let
    $\phi_1,\dots,\phi_k\in\mathscr H^*$ be linearly independent, and assume that
    there exist $u_1,\dots,u_k\in\mathscr H$ such that
    \[
        S(u_m,v)=\phi_m(v)
        \qquad
        \text{for all }v\in\mathscr H,\quad m=1,\dots,k.
    \]
    Let
    \[
        \mathscr X\coloneqq\bigcap_{m=1}^k\ker\phi_m,
        \qquad
        \mathbf M\coloneqq
        \bigl(S(u_m,u_n)\bigr)_{1\le m,n\le k}.
    \]
    Then
    \[
        n_-\bigl(S|_{\mathscr X}\bigr)
        =
        n_-(S)-n_{\le0}(\mathbf M),
    \]
    where $n_{\le0}(\mathbf M)$ denotes the number of non-positive eigenvalues of
    $\mathbf M$, counted with multiplicity. In particular, if $\mathbf M$ is
    non-singular, then
    \[
        n_-\bigl(S|_{\mathscr X}\bigr)
        =
        n_-(S)-n_-(\mathbf M).
    \]
\end{theorem}

The goal is to prove that \(\mathcal L_{c_j}\) is nonnegative on \(X_j\) in terms of the quadratic form. Equivalently, the Morse index of \(\mathcal L_{c_j}\) on \(X_j\) is zero. By Theorem~\ref{thm:constrained-index-formula}, this constrained Morse index is determined by the full Morse index of \(\mathcal L_{c_j}\) and the inertia of the associated \(2\times2\) constraint matrix.

At the bottom of the essential spectrum, the free resolvent has the singular Fourier multiplier \(|\xi|^{-1}\). Hence the Birman--Schwinger formulation cannot be applied directly at \(\mu=0\). We therefore introduce a positive parameter \(\mu>0\) and consider the regularized family of operators 
\begin{equation*}
    \mathcal L_c^{(\mu)}
    \coloneqq
    |D|+\mu+cV_Q.
\end{equation*}
We regard \(\mathcal L_c^{(\mu)}\) as the self-adjoint operator on \(L^2_{\mathrm e}\) associated with the closed quadratic form on \(H^{1/2}_{\mathrm e}\) given by
\[
    f\mapsto
    \norm{|D|^{1/2}f}_{L^2}^2
    +
    \mu\norm{f}_{L^2}^2
    +
    c\int_{\mathbb R}V_Q(y)f(y)^2\,dy.
\]
For each fixed \(\mu>0\), \(0\) can belong to the spectrum only as an eigenvalue, and \(\ker\mathcal L_c^{(\mu)}=\{0\}\) implies \(0\in\rho(\mathcal L_c^{(\mu)})\). Once the regularized quadratic forms are shown to be nonnegative on \(X_j\) for sufficiently small \(\mu>0\), letting \(\mu\to0^+\) gives the nonnegativity of \(\mathcal L_{c_j}\) on \(X_j\). The following proposition records the Morse index and constraint matrix conditions required for the regularized family.

\begin{proposition}[Index and determinant conditions for regularized operator]
\label{prop:shifted-index-determinant-conditions}
    There exist constants \(c_1>6\), \(c_2>2\), and \(\mu_0>0\) such that, for
    each \(j=1,2\) and every \(0<\mu<\mu_0\), the operator $\mathcal L_{c_j}^{(\mu)}$
    satisfies
    \begin{equation}\label{eq:invertibility and morse index 1}
        0\in\rho(\mathcal L_{c_j}^{(\mu)}),
        \qquad
        n_-^{\mathrm{e}}(\mathcal L_{c_j}^{(\mu)})=1.
    \end{equation}
    Moreover, if \(\Xi_j=\{\Psi_1^{(j)},\Psi_2^{(j)}\}\), then the constraint matrix
    \[
        \mathbf M_j^{(\mu)}
        \coloneqq
        \left(
            \bigl(
                (\mathcal L_{c_j}^{(\mu)})^{-1}\Psi_m^{(j)},
                \Psi_n^{(j)}
            \bigr)_r
        \right)_{1\le m,n\le2}
    \]
    satisfies
    \begin{equation}\label{eq: constraint matrix det<0}
        \det \mathbf M_j^{(\mu)}<0.
    \end{equation}
\end{proposition}

We now explain why Proposition~\ref{prop:shifted-index-determinant-conditions}
implies Proposition~\ref{prop:Non-negativity of D+cyQ'Q}.

\begin{proof}[Proof of Proposition~\ref{prop:Non-negativity of D+cyQ'Q}
assuming Proposition~\ref{prop:shifted-index-determinant-conditions}]
    Let $c_1>6$, $c_2>2$, and $\mu_0>0$ be given by
    Proposition~\ref{prop:shifted-index-determinant-conditions}. Fix
    $j\in\{1,2\}$ and $0<\mu<\mu_0$.
    
    We apply Theorem~\ref{thm:constrained-index-formula} to the Hilbert space $H^{1/2}_{\mathrm e}(\mathbb R;\mathbb R)$, the symmetric bilinear form $S^{(\mu)}(f,g) \coloneqq (\mathcal L_{c_j}^{(\mu)}f,g)_r$, and the constraint functionals
    \[
        \phi_{m}(f)\coloneqq (f,\Psi_m^{(j)})_r,
        \qquad m=1,2.
    \]
    The two functionals \(\phi_1,\phi_2\) are linearly independent, since the
    two functions in \(\Xi_j\) are linearly independent in \(L^2_{\mathrm e}\).
    From \eqref{eq:invertibility and morse index 1} we have $0\in\rho(\mathcal L_{c_j}^{(\mu)})$ and hence the vectors
    \[
        u_m^{(\mu)}
        \coloneqq
        (\mathcal L_{c_j}^{(\mu)})^{-1}\Psi_m^{(j)}
    \]
    are well-defined and the constrained space is exactly $X_j$. Moreover, the matrix satisfies
    \[
        \bigl(S^{(\mu)}(u_m^{(\mu)},u_n^{(\mu)})\bigr)_{m,n=1}^2 = \mathbf M^{(\mu)}_j.
    \]
    By Proposition~\ref{prop:shifted-index-determinant-conditions},
    \[
        n_-^{\mathrm{e}}(\mathcal L_{c_j}^{(\mu)})=1,
        \qquad
        \det \mathbf M_j^{(\mu)}<0.
    \]
    Since $\mathbf M_j^{(\mu)}$ is a real symmetric $2\times2$ matrix, the determinant
    condition implies
    \[
        n_-(\mathbf M_j^{(\mu)})=1.
    \]
    Thus Theorem~\ref{thm:constrained-index-formula} gives
    \[
        n_-\left(S^{(\mu)}|_{X_j}\right)
        =
        n_-^{\mathrm{e}}(\mathcal L_{c_j}^{(\mu)})
        -
        n_-(\mathbf M_j^{(\mu)})
        =
        0.
    \]
    Hence, for each $\mu \in (0,\mu_0)$ we have
    \[
        (\mathcal L_{c_j}^{(\mu)}f,f)_r =
        (\mathcal L_{c_j} f,f)_r+\mu\norm{f}_{L^2}^2 \geq 0
        \qquad
        \text{for every }f\in X_j.
    \]
    Letting $\mu\to0^+$ implies Proposition~\ref{prop:Non-negativity of D+cyQ'Q}.
\end{proof}

It remains to prove Proposition~\ref{prop:shifted-index-determinant-conditions}.
Sections~\ref{sec:BS-principle} and \ref{sec:limiting-constraint-matrices} establish respectively the conditions in \eqref{eq:invertibility and morse index 1} and \eqref{eq: constraint matrix det<0}.

\section{Renormalization and Birman--Schwinger reduction}\label{sec:BS-principle}
In this section, we prove the Morse index and invertibility conditions in \eqref{eq:invertibility and morse index 1}. Directly counting the negative eigenvalues of \(\mathcal L_{c_j}^{(\mu)}\) is difficult for every sufficiently small \(\mu>0\). We therefore use the Birman--Schwinger principle. This principle converts the problem into an eigenvalue counting problem for a compact nonnegative operator \(K_Q^{(\mu)}\). The Morse index is then determined by the number of eigenvalues of \(K_Q^{(\mu)}\) exceeding \(c^{-1}\). Related estimates for the number of negative eigenvalues of fractional Schr\"odinger operators at the critical exponent are given in \cite{BreteauxFaupinGrasselli2025JST}. These estimates, however, do not give the profile-dependent strict bound required below.

As \(\mu\to0^+\), the zero-frequency singularity of \((|D|+\mu)^{-1}\) produces a logarithmically divergent rank-one contribution to \(K_Q^{(\mu)}\). Consequently, \(K_Q^{(\mu)}\) does not converge in operator norm. Separating this contribution isolates the first eigendirection, along which the largest eigenvalue diverges, while the second eigenvalue
converges to the norm of the projected logarithmic operator. The Morse index and invertibility conditions thus reduce to a strict norm bound for this projected operator.

The next two propositions carry out this reduction. The first establishes the Birman--Schwinger correspondence and the basic spectral properties of \(K_Q^{(\mu)}\) for fixed \(\mu>0\). The second isolates the logarithmically divergent rank-one contribution as \(\mu\to0^+\) and identifies the operator governing the remaining spectrum.

\begin{proposition}[Birman--Schwinger principle and basic spectral properties]\label{prop:BS-principle}
    Let \(c>0\) and \(\mu>0\), and define
    \begin{equation}\label{eq:def-BS-operator-Q}
        K_Q^{(\mu)}
        \coloneqq
        |V_Q|^{1/2}(|D|+\mu)^{-1}|V_Q|^{1/2}
        \qquad\text{on }L^2_{\mathrm e}.
    \end{equation}
    Then the following hold.
    \begin{enumerate}
        \item The essential spectrum of \(\mathcal L_c^{(\mu)}\) satisfies
        \begin{equation}\label{eq:BS-essential-spectrum}
            \sigma_{\mathrm{ess}}\bigl(\mathcal L_c^{(\mu)}\bigr)
            =
            [\mu,\infty).
        \end{equation}
        In particular, \(n_-^{\mathrm{e}}(\mathcal L_c^{(\mu)})<\infty\), and
        \begin{equation}\label{eq:BS-zero-resolvent-kernel}
            0\in\rho\bigl(\mathcal L_c^{(\mu)}\bigr)
            \quad \text{if and only if} \quad 
            \ker\mathcal L_c^{(\mu)}=\{0\}.
        \end{equation}

        \item The operator \(K_Q^{(\mu)}\) is compact, nonnegative,
        self-adjoint, and positivity improving, in the sense that if
        \(0\leq v\in L^2_{\mathrm e}\) and \(v\not\equiv0\), then
        \begin{equation}\label{eq:BS-positivity-improving}
            K_Q^{(\mu)}v>0
            \qquad\text{a.e. on }\mathbb R.
        \end{equation}

        \item One has the kernel identity
        \begin{equation}\label{eq:BS-kernel-identity}
            \dim\ker\mathcal L_c^{(\mu)}
            =
            \dim\ker\bigl(I-cK_Q^{(\mu)}\bigr),
        \end{equation}
        and the index identity
        \begin{equation}\label{eq:BS-index-identity}
            n_-^{\mathrm{e}}\bigl(\mathcal L_c^{(\mu)}\bigr)
            =
            n_+\bigl(c^{-1};K_Q^{(\mu)}\bigr).
        \end{equation}

        \item If \(0<\mu_1<\mu_2\), then
        \begin{equation}\label{eq:BS-operator-monotonicity}
            K_Q^{(\mu_1)}\geq K_Q^{(\mu_2)}
        \end{equation}
        in the operator order. Consequently, for each \(k\geq1\),
        \begin{equation}\label{eq:BS-eigenvalue-monotonicity}
            \lambda_k\bigl(K_Q^{(\mu_1)}\bigr)
            \geq
            \lambda_k\bigl(K_Q^{(\mu_2)}\bigr).
        \end{equation}
    \end{enumerate}
\end{proposition}

\begin{proof}
    Since \(V_Q\leq0\), we have
    \[
        \mathcal L_c^{(\mu)}
        =
        |D|+\mu-c|V_Q|.
    \]
    Moreover, since \(V_Q\in L^\infty\) and \(V_Q(y)\to0\) as
    \(|y|\to\infty\), multiplication by \(V_Q\) is relatively compact with
    respect to \(|D|+\mu\). Hence Weyl's theorem gives
    \[
        \sigma_{\mathrm{ess}}\bigl(\mathcal L_c^{(\mu)}\bigr)
        =
        \sigma_{\mathrm{ess}}(|D|+\mu)
        =
        [\mu,\infty),
    \]
    which proves \eqref{eq:BS-essential-spectrum}. In particular, the spectrum
    below \(0\) consists only of finitely many
    \(L^2_{\mathrm e}\)-eigenvalues, counted with finite multiplicity. By
    the min--max principle, \(n_-^{\mathrm{e}}(\mathcal L_c^{(\mu)})\) is the number of
    these negative eigenvalues, counted with multiplicity. Since \(0\) lies
    below the essential spectrum, \(0\in\sigma(\mathcal L_c^{(\mu)})\) can
    occur only as an \(L^2_{\mathrm e}\)-eigenvalue. Therefore
    \eqref{eq:BS-zero-resolvent-kernel} follows.

    We next prove the stated properties of \(K_Q^{(\mu)}\). The operator
    \((|D|+\mu)^{-1}:L^2_{\mathrm e}\to H^1_{\mathrm e}\) is bounded,
    and multiplication by \(|V_Q|^{1/2}\) is compact from
    \(H^1_{\mathrm e}\) to \(L^2_{\mathrm e}\). Thus
    \(K_Q^{(\mu)}\), defined in \eqref{eq:def-BS-operator-Q}, is compact on
    \(L^2_{\mathrm e}\). It is self-adjoint by symmetry. Its nonnegativity follows from  \eqref{eq:def-BS-operator-Q}.

    We now prove positivity improving. Let \(G^{(\mu)}\) be the integral kernel
    of \((|D|+\mu)^{-1}\). By \eqref{eq:K-rescaled} and
    \eqref{eq:G-poisson-real},
    \[
        G^{(\mu)}(y)>0
        \qquad
        \text{for }y\neq0.
    \]
    Hence \((|D|+\mu)^{-1}\) maps every nonzero nonnegative function to a
    strictly positive function. Since \(|V_Q|^{1/2}>0\) a.e., if
    \(0\leq v\in L^2_{\mathrm e}\) and \(v\not\equiv0\), then
    \[
        K_Q^{(\mu)}v
        =
        |V_Q|^{1/2}(|D|+\mu)^{-1}|V_Q|^{1/2}v
        >0
        \qquad\text{a.e.}
    \]
    This proves \eqref{eq:BS-positivity-improving}.

    We prove \eqref{eq:BS-kernel-identity}. If
    \(f\in\ker\mathcal L_c^{(\mu)}\), then
    \((|D|+\mu)f=c|V_Q|f\). Setting
    \(v=|V_Q|^{1/2}f\), we obtain
    \(v=cK_Q^{(\mu)}v\), hence
    \(v\in\ker(I-cK_Q^{(\mu)})\). Conversely, if
    \(v\in\ker(I-cK_Q^{(\mu)})\), set
    \[
        f
        =
        c(|D|+\mu)^{-1}|V_Q|^{1/2}v.
    \]
    Then \(|V_Q|^{1/2}f=cK_Q^{(\mu)}v=v\), and therefore
    \[
        (|D|+\mu)f
        =
        c|V_Q|^{1/2}v
        =
        c|V_Q|f.
    \]
    Thus \(f\in\ker\mathcal L_c^{(\mu)}\). These two
    constructions are inverse to each other, which proves
    \eqref{eq:BS-kernel-identity}.

    It remains to prove \eqref{eq:BS-index-identity}. The quadratic form of
    \(\mathcal L_c^{(\mu)}\) is
    \[
        f\mapsto
        \norm{(|D|+\mu)^{1/2}f}_{L^2}^2
        -
        c\norm{|V_Q|^{1/2}f}_{L^2}^2,
        \qquad
        f\in H^{1/2}_{\mathrm e}.
    \]
    Under the isomorphism
    \(f\mapsto (|D|+\mu)^{1/2}f\) from \(H^{1/2}_{\mathrm e}\) onto
    \(L^2_{\mathrm e}\), this form becomes
    \[
        g\mapsto
        \norm{g}_{L^2}^2
        -
        c\Bigl(
            (|D|+\mu)^{-1/2}|V_Q|(|D|+\mu)^{-1/2}g,
            g
        \Bigr)_r .
    \]
    Hence the min--max principle gives
    \[
        n_-^{\mathrm{e}}\bigl(\mathcal L_c^{(\mu)}\bigr)
        =
        n_+\left(
            c^{-1};
            (|D|+\mu)^{-1/2}|V_Q|(|D|+\mu)^{-1/2}
        \right).
    \]
    The nonzero spectra of
    \[
        (|D|+\mu)^{-1/2}|V_Q|^{1/2}
        \,|V_Q|^{1/2}(|D|+\mu)^{-1/2}
    \]
    and
    \[
        |V_Q|^{1/2}(|D|+\mu)^{-1/2}
        (|D|+\mu)^{-1/2}|V_Q|^{1/2}
        =
        K_Q^{(\mu)}
    \]
    agree, counted with multiplicity. Therefore
    \[
        n_-^{\mathrm{e}}\bigl(\mathcal L_c^{(\mu)}\bigr)
        =
        n_+\bigl(c^{-1};K_Q^{(\mu)}\bigr),
    \]
    which proves \eqref{eq:BS-index-identity}.

    Finally, let \(0<\mu_1<\mu_2\). Since
    \[
        (|D|+\mu_1)^{-1}\geq (|D|+\mu_2)^{-1}
    \]
    in the operator order, multiplying on both sides by \(|V_Q|^{1/2}\) gives \eqref{eq:BS-operator-monotonicity}. The eigenvalue monotonicity \eqref{eq:BS-eigenvalue-monotonicity} follows from the min--max principle.
\end{proof}

By Proposition~\ref{prop:BS-principle}, \(K_Q^{(\mu)}\) is compact, nonnegative, self-adjoint, and positivity improving. Hence Theorem~\ref{thm:Jentzsch-theorem} implies that \(\lambda_1(K_Q^{(\mu)})\) is simple and that its normalized eigenfunction may be chosen strictly positive almost everywhere.

\begin{proposition}[Rank-one decomposition and the reduced operator]\label{prop:rank-one-decomposition}
    Let \(K_Q^{(\mu)}\) be defined by \eqref{eq:def-BS-operator-Q}. Then the
    following hold.
    \begin{enumerate}
        \item For each \(\mu>0\), the operator \(K_Q^{(\mu)}\) admits the decomposition: for \(v\in L^2\),
    \begin{equation}\label{decomposition of K_g(mu)}
        \begin{aligned}
            K_Q^{(\mu)}[v] &= \frac{1}{\pi}|V_Q|^{1/2}\Bigl(\log\frac{1}{\mu}-\gamma_E\Bigr)(|V_Q|^{1/2},v)_r + T_Q[v] + E_Q^{(\mu)}[v],
        \end{aligned}
    \end{equation}
    where \(T_Q\) and \(E_Q^{(\mu)}\) are defined by
    \begin{equation*}
        \begin{aligned}
            T_Q[v]&\coloneqq -\frac{1}{\pi}|V_Q|^{1/2}\bigl(\log|\cdot| \ast (|V_Q|^{1/2}v)\bigr),\\
            E_Q^{(\mu)}[v]&\coloneqq \frac{1}{\pi}|V_Q|^{1/2}\bigl(r_{\mu}\ast(|V_Q|^{1/2}v)\bigr).
        \end{aligned}
    \end{equation*}
    Here, \(\gamma_E>0\) denotes the Euler--Mascheroni constant, and \(r_\mu\) is given by
    \begin{equation*}
        r_{\mu}(y):=\int_0^{\infty} (1-\cos(\xi y))\,\frac{\mu}{\xi(\xi+\mu)}\,d\xi .
    \end{equation*}
        \item For each \(\mu>0\), the operators \(T_Q\) and \(E_Q^{(\mu)}:L^2\to L^2\) are compact and symmetric. Furthermore, for each \(\mu>0\) we have
        \begin{equation}\label{estimate of error}
            \norm{E_Q^{(\mu)}}_{L^2\to L^2}\;\lesssim_Q \mu.
        \end{equation}
        \item  Let \(\Phi_Q\coloneqq \dfrac{\sqrt{|V_Q|}}{\|\sqrt{|V_Q|}\|_{L^2}}\), and \(P_{\Phi_Q^\perp}\) be a projection operator defined by
        \begin{equation*}
            P_{\Phi_Q^\perp}[v] \coloneqq v-(v,\Phi_Q)_r \Phi_Q. 
        \end{equation*}
        Set
        \begin{equation}\label{eq:definition of nu_mu}
            \nu_\mu
            \coloneqq
            \frac{1}{\pi}
            \Bigl(\log\frac{1}{\mu}-\gamma_E\Bigr)
            \norm{|V_Q|^{1/2}}_{L^2}^2.
        \end{equation}
        Then
        \begin{equation}\label{eq:lambda1-asymptotic-KQ}
            \lambda_1(K_Q^{(\mu)})
            =
            \nu_\mu+O(1)
            \to+\infty
            \qquad
            \text{as }\mu\to0^+.
        \end{equation}
        Let \(\varphi_{1,Q}^{(\mu)}\) denote the positive even \(L^2\)-normalized eigenfunction
        corresponding to \(\lambda_1(K_Q^{(\mu)})\). Then
        \begin{equation}\label{asymptotic for e_g(mu)}
            \lim_{\mu\to 0+}\norm{\varphi_{1,Q}^{(\mu)}- \Phi_Q}_{L^2}=0,
        \end{equation}
        and
        \begin{equation}\label{asymptotic for theta_2(mu)}
            \lim_{\mu\to 0+}\lambda_2(K_Q^{(\mu)})=\norm{P_{\Phi_Q^{\perp}} T_Q P_{\Phi_Q^{\perp}}}_{L^2_{\mathrm e}\to L^2_{\mathrm e}}.
        \end{equation}
    \end{enumerate}
\end{proposition}

\begin{proof}
    We work on \(L^2_{\mathrm e}\), which is invariant under all the operators below since \(V_Q\) is even.

    \textbf{Step 1: decomposition of \(K_Q^{(\mu)}\).}
    By Lemma~\ref{lem:Kmu-smallmu}, for every \(y\neq z\) and \(\mu>0\),
    \[
        G^{(\mu)}(y-z)
        =
        \frac{1}{\pi}\Bigl(\log\frac{1}{\mu}-\gamma_E\Bigr)
        -\frac{1}{\pi}\log|y-z|
        +\frac{1}{\pi}r_\mu(y-z).
    \]
    Since
    \[
        K_Q^{(\mu)}[v](y)
        =
        |V_Q(y)|^{1/2}\int_\bbR G^{(\mu)}(y-z)|V_Q(z)|^{1/2}v(z)\,dz,
    \]
    substituting the above identity gives
    \[
        K_Q^{(\mu)}[v]
        =
        \frac{1}{\pi}|V_Q|^{1/2}
        \Bigl(\log\frac{1}{\mu}-\gamma_E\Bigr)(|V_Q|^{1/2},v)_r
        +T_Q[v]
        +E_Q^{(\mu)}[v].
    \]
    This proves \eqref{decomposition of K_g(mu)}.

    \textbf{Step 2: compactness and the error bound.}
    Since \(|V_Q(y)|\lesssim \langle y\rangle^{-4}\), we have
    \[
        (|V_Q|\ast |V_Q|)(y)\lesssim_Q \langle y\rangle^{-4}.
    \]
    Therefore
    \[
        \iint_{\bbR^2}
        |V_Q(y)|\,(\log|y-z|)^2\,|V_Q(z)|\,dy\,dz
        =
        \int_\bbR (|V_Q|\ast |V_Q|)(y)(\log|y|)^2\,dy
        <\infty,
    \]
    and hence \(T_Q\) is Hilbert--Schmidt, hence compact. Symmetry is immediate from the symmetry of the kernel \(\log|y-z|\). Similarly, by \eqref{eq:rmu-global-bound},
    \[
        |r_\mu(y)|\lesssim \mu |y|.
    \]
    Hence
    \[
        \iint_{\bbR^2}
        |V_Q(y)|\,|r_\mu(y-z)|^2\,|V_Q(z)|\,dy\,dz
        \lesssim
        \mu^2
        \int_\bbR (|V_Q|\ast |V_Q|)(y)|y|^2\,dy
        <\infty,
    \]
    and thus \(E_Q^{(\mu)}\) is also Hilbert--Schmidt and symmetric. In particular,
    \[
        \norm{E_Q^{(\mu)}}_{L^2\to L^2}
        \leq
        \norm{E_Q^{(\mu)}}_{\mathrm{HS}}
        \lesssim_Q \mu,
    \]
    which proves \eqref{estimate of error}.

    \textbf{Step 3: asymptotics of the first eigendirection.}
    From the definition of $\nu_\mu$ in \eqref{eq:definition of nu_mu}, we have
    \[
        K_Q^{(\mu)}
        =
        \nu_\mu\Phi_Q\otimes \Phi_Q
        +T_Q
        +E_Q^{(\mu)}.
    \]
    Since \(T_Q\) is bounded and \(\norm{E_Q^{(\mu)}}_{L^2\to L^2}\lesssim 1\), evaluating the Rayleigh quotient at \(\Phi_Q\) and taking the supremum over unit vectors gives
    \begin{equation}\label{eq: nu_mu sim lambda_1(K_p^mu)}
        \nu_\mu-C
        \leq
        \lambda_1(K_Q^{(\mu)})
        \leq
        \nu_\mu+C
    \end{equation}
    for some \(C>0\) independent of \(\mu\). In particular,
    \begin{equation}\label{eq:first eigenvalue blow up}
        \lambda_1(K_Q^{(\mu)})=\nu_\mu+O(1)
        \qquad\text{and}\qquad
        \lambda_1(K_Q^{(\mu)})\to\infty
        \quad\text{as }\mu\to 0^+,
    \end{equation}
    and this yields \eqref{eq:lambda1-asymptotic-KQ}.
    Now, we decompose \(\varphi_{1,Q}^{(\mu)}\) by
    \begin{equation}\label{eq:decomposition of e_p(mu)}
        \varphi_{1,Q}^{(\mu)}=(\varphi_{1,Q}^{(\mu)},\Phi_Q)_r \Phi_Q+h_\mu,
        \qquad
        h_\mu\in \Phi_Q^\perp.
    \end{equation}
    Since \(\norm{\varphi_{1,Q}^{(\mu)}}_{L^2}=1\), we have
    \[
        (\varphi_{1,Q}^{(\mu)},\Phi_Q)_r^2+\norm{h_\mu}_{L^2}^2=1.
    \]
    Applying \(P_{\Phi_Q^\perp}\) to both sides of \(K_Q^{(\mu)}[\varphi_{1,Q}^{(\mu)}] =\lambda_1(K_Q^{(\mu)})\varphi_{1,Q}^{(\mu)}\) and substituting the decomposition \(\varphi_{1,Q}^{(\mu)}\) to \eqref{eq:decomposition of e_p(mu)} gives
    \[
        \Bigl(
            \lambda_1(K_Q^{(\mu)})I
            -
            P_{\Phi_Q^\perp} (T_Q+E_Q^{(\mu)}) P_{\Phi_Q^\perp}
        \Bigr)[h_\mu]
        =
        (\varphi_{1,Q}^{(\mu)},\Phi_Q)_r P_{\Phi_Q^\perp} (T_Q+E_Q^{(\mu)})[\Phi_Q].
    \]
    Since \(\lambda_1(K_Q^{(\mu)})\to\infty\) and \(T_Q+E_Q^{(\mu)}\) remains uniformly bounded,
    Lemma~\ref{lem:Neumann-series-criterion} shows that the operator on the left is invertible on \(\operatorname{span}\{\Phi_Q\}^\perp\) for all sufficiently small \(\mu\), and
    \[
        \norm{h_\mu}_{L^2}
        \lesssim
        \frac{|(\varphi_{1,Q}^{(\mu)},\Phi_Q)_r|}{\lambda_1(K_Q^{(\mu)})}
        \leq \frac{1}{\lambda_1(K_Q^{(\mu)})}.
    \]
    In the last inequality, we use \(|(\varphi_{1,Q}^{(\mu)},\Phi_Q)_r| \leq 1\). Using \eqref{eq:first eigenvalue blow up} we obtain \(\norm{h_\mu}_{L^2}\to 0\), and hence \((\varphi_{1,Q}^{(\mu)},\Phi_Q)_r\to 1\). Therefore, using \eqref{eq:decomposition of e_p(mu)} we conclude
    \[
        \norm{\varphi_{1,Q}^{(\mu)}- (\varphi_{1,Q}^{(\mu)},\Phi_Q)_r\Phi_Q}_{L^2}
        =
        \norm{h_\mu}_{L^2}
        \to 0,
    \]
    which proves \eqref{asymptotic for e_g(mu)}.

    \textbf{Step 4: reduction of the second eigenvalue.}
    Since \(K_Q^{(\mu)}\) is compact, self-adjoint, and nonnegative on \(L^2_{\mathrm e}\), we have
    \[
        \lambda_2(K_Q^{(\mu)})
        =
        \norm{
            P_{(\varphi_{1,Q}^{(\mu)})^\perp} 
            K_Q^{(\mu)} 
            P_{(\varphi_{1,Q}^{(\mu)})^\perp}
        }_{L^2_{\mathrm e}\to L^2_{\mathrm e}}.
    \]
    Also, since \(\varphi_{1,Q}^{(\mu)}\to \Phi_Q\) in \(L^2\),
    \[
        \norm{
            P_{(\varphi_{1,Q}^{(\mu)})^\perp}
            -
            P_{\Phi_Q^\perp}
        }_{L^2\to L^2}
        \to 0.
    \]
    Next, using the decomposition of \(K_Q^{(\mu)}\),
    \[
        P_{(\varphi_{1,Q}^{(\mu)})^\perp} 
        K_Q^{(\mu)} 
        P_{(\varphi_{1,Q}^{(\mu)})^\perp}
        =
        \nu_\mu
        P_{(\varphi_{1,Q}^{(\mu)})^\perp} 
        (\Phi_Q\otimes \Phi_Q) 
        P_{(\varphi_{1,Q}^{(\mu)})^\perp} 
        +
        P_{(\varphi_{1,Q}^{(\mu)})^\perp}(T_Q + E_Q^{(\mu)})
        P_{(\varphi_{1,Q}^{(\mu)})^\perp}.
    \]
    The contribution containing \(E_Q^{(\mu)}\) tends to zero in operator norm
    by \eqref{estimate of error}. For the rank-one part,
    \[
        \norm{
            P_{(\varphi_{1,Q}^{(\mu)})^\perp} 
            (\Phi_Q\otimes \Phi_Q) 
            P_{(\varphi_{1,Q}^{(\mu)})^\perp}
        }_{L^2\to L^2}
        =
        \norm{P_{(\varphi_{1,Q}^{(\mu)})^\perp}[\Phi_Q]}_{L^2}^2
        =
        1-((\varphi_{1,Q}^{(\mu)}),\Phi_Q)_r^2
        =
        \norm{h_\mu}_{L^2}^2.
    \]
    Since \(\nu_\mu\sim \lambda_1(K_Q^{(\mu)})\) and
    \(\norm{h_\mu}_{L^2}\lesssim \lambda_1(K_Q^{(\mu)})^{-1}\), we obtain
    \[
        \nu_\mu
        \norm{
            P_{(\varphi_{1,Q}^{(\mu)})^\perp} 
            (\Phi_Q\otimes \Phi_Q) 
            P_{(\varphi_{1,Q}^{(\mu)})^\perp}
        }_{L^2\to L^2}
        \to 0.
    \]
    Finally,
    \[
        \norm{
            P_{(\varphi_{1,Q}^{(\mu)})^\perp} (T_Q + E^{(\mu)}_Q) P_{(\varphi_{1,Q}^{(\mu)})^\perp}
            -
            P_{\Phi_Q^\perp} T_Q  P_{\Phi_Q^\perp}
        }_{L^2\to L^2}
        \to 0
    \]
    because \(E_Q^{(\mu)} \to 0\) as \(\mu \to 0^+\), \(T_Q\) is bounded and
    \(P_{(\varphi_{1,Q}^{(\mu)})^\perp}\to P_{\Phi_Q^\perp}\) in operator norm. Hence
    \[
        \norm{
            P_{(\varphi_{1,Q}^{(\mu)})^\perp} 
            K_Q^{(\mu)}  
            P_{(\varphi_{1,Q}^{(\mu)})^\perp}
            -
            P_{\Phi_Q^\perp} 
            T_Q  
            P_{\Phi_Q^\perp}
        }_{L^2_{\mathrm e}\to L^2_{\mathrm e}}
        \to 0.
    \]
    Taking operator norms gives
    \[
        \lim_{\mu\to 0+}\lambda_2(K_Q^{(\mu)})
        =
        \norm{
            P_{\Phi_Q^\perp}  
            T_Q  
            P_{\Phi_Q^\perp}
        }_{L^2_{\mathrm e}\to L^2_{\mathrm e}},
    \]
    which is \eqref{asymptotic for theta_2(mu)}.
\end{proof}

\begin{corollary}[Projected logarithmic criterion]\label{cor:projected-log-criterion}
    Let \(c>0\). Assume that
    \begin{equation}\label{eq:projected-log-criterion-assumption}
        \norm{
            P_{\Phi_Q^\perp}T_QP_{\Phi_Q^\perp}
        }_{L^2_{\mathrm e}\to L^2_{\mathrm e}}
        < c^{-1}.
    \end{equation}
    Then there exists \(\mu_c>0\) such that, for every \(0<\mu<\mu_c\),
    \begin{equation*}
        n_-^{\mathrm{e}}\bigl(\mathcal L_c^{(\mu)}\bigr)=1,
        \qquad
        0\in\rho\bigl(\mathcal L_c^{(\mu)}\bigr).
    \end{equation*}
\end{corollary}

\begin{proof}
    By Proposition~\ref{prop:rank-one-decomposition},
    \[
        \lambda_1(K_Q^{(\mu)})\to+\infty\quad \text{and} \quad
        \lambda_2(K_Q^{(\mu)})
        \nearrow
        \norm{
            P_{\Phi_Q^\perp}T_QP_{\Phi_Q^\perp}
        }_{L^2_{\mathrm e}\to L^2_{\mathrm e}} \quad\text{as}\quad \mu \to 0^+.
    \]
    The convergence is monotone, in view of \eqref{eq:BS-eigenvalue-monotonicity}.
    Hence \eqref{eq:projected-log-criterion-assumption} implies that, for all
    sufficiently small \(\mu>0\),
    \[
        \lambda_1(K_Q^{(\mu)})>c^{-1}>
        \lambda_2(K_Q^{(\mu)}).
    \]
    Therefore
    \[
        n_+\bigl(c^{-1};K_Q^{(\mu)}\bigr)=1,
        \qquad
        c^{-1}\notin\sigma(K_Q^{(\mu)}).
    \]
    By Proposition~\ref{prop:BS-principle}, this gives
    \[
        n_-^{\mathrm{e}}\bigl(\mathcal L_c^{(\mu)}\bigr)=1 \quad \text{and} \quad
        \ker\mathcal L_c^{(\mu)}=\{0\}.
    \]
    Using \eqref{eq:BS-zero-resolvent-kernel}, we obtain
    \(0\in\rho(\mathcal L_c^{(\mu)})\).
\end{proof}

By Corollary~\ref{cor:projected-log-criterion}, \eqref{eq:invertibility and morse index 1} follows once we prove the following quantitative bound.

\begin{proposition}[Reduced Birman--Schwinger bound]
\label{prop:reduced-BS-bound-Q}
    We have
    \begin{equation*}
        \norm{
            P_{\Phi_Q^\perp}T_QP_{\Phi_Q^\perp}
        }_{L^2_{\mathrm e}\to L^2_{\mathrm e}}
        <0.1571.
    \end{equation*}
\end{proposition}

The proof of Proposition~\ref{prop:reduced-BS-bound-Q} is completed in
Part~\ref{part:spectral-estimates}, using the profile estimates from
Part~\ref{part:ground-state-approximation} and the interval bounds from
Part~\ref{part:interval-verification}. We record its consequence below.

\begin{corollary}[Morse index and invertibility for the regularized family]
\label{cor:shifted-morse-index-invertibility}
    Let
    \[
        c_1=6.1,
        \qquad
        c_2=2.1.
    \]
    Then there exists \(\mu_{\mathrm{BS}}>0\) such that, for \(j=1,2\) and
    \(0<\mu<\mu_{\mathrm{BS}}\),
    \[
        n_-^{\mathrm{e}}\bigl(\mathcal L_{c_j}^{(\mu)}\bigr)=1,
        \qquad
        0\in\rho\bigl(\mathcal L_{c_j}^{(\mu)}\bigr).
    \]
\end{corollary}

\begin{proof}
    By Proposition~\ref{prop:reduced-BS-bound-Q},
    \[
        \norm{
            P_{\Phi_Q^\perp}T_QP_{\Phi_Q^\perp}
        }_{L^2_{\mathrm e}\to L^2_{\mathrm e}}
        <0.1571.
    \]
    Since
    \[
        0.1571<\frac1{6.1},
        \qquad
        0.1571<\frac1{2.1},
    \]
    Corollary~\ref{cor:projected-log-criterion} applies with \(c=c_1\) and
    \(c=c_2\). Taking
    \[
        \mu_{\mathrm{BS}}
        \coloneqq
        \min\{\mu_{c_1},\mu_{c_2}\}
    \]
    completes the proof.
\end{proof}

\section{Limiting constraint matrices and completion of the spectral reduction}\label{sec:limiting-constraint-matrices}
The preceding section established the invertibility of \(\mathcal L_{c_j}^{(\mu)}\) and \(n_-^{\mathrm e}(\mathcal L_{c_j}^{(\mu)})=1\), for \(c_1=6.1\) and \(c_2=2.1\); see Corollary~\ref{cor:shifted-morse-index-invertibility}. It remains to prove the determinant condition \eqref{eq: constraint matrix det<0} for the constraint matrices \(\mathbf M_j^{(\mu)}\).

As \(\mu\to0^+\), the resolvents \((\mathcal L_{c_j}^{(\mu)})^{-1}\) do not converge as bounded operators on \(L^2_{\mathrm e}\), since zero is the bottom of the essential spectrum of \(\mathcal L_{c_j}\). Our key observation is that the divergence is confined to a single rank-one term, which cancels between the free resolvent and the Birman--Schwinger correction.

We first use the symmetric resolvent identity to separate the free resolvent contribution from the Birman--Schwinger correction.

\begin{lemma}[Symmetric resolvent identity]
    For \(j=1,2\) and \(0<\mu<\mu_0\), we have
    \begin{equation}\label{eq:woodbury-BS-resolvent}
        (\mathcal L_{c_j}^{(\mu)})^{-1}
        =
        G^{(\mu)}
        +
        c_jG^{(\mu)} |V_Q|^{1/2}
        (I-c_jK_Q^{(\mu)})^{-1}
        |V_Q|^{1/2}G^{(\mu)}.
    \end{equation}
\end{lemma}

\begin{proof}
    The invertibility of $\mathcal L_{c_j}^{(\mu)}$ for $\mu\in(0,\mu_0)$ follows from Corollary~\ref{cor:shifted-morse-index-invertibility}.
    The identity follows by multiplying the right-hand side of
    \eqref{eq:woodbury-BS-resolvent} by \(\mathcal L_{c_j}^{(\mu)}\), using
    \[
        I+c_jK_Q^{(\mu)}(I-c_jK_Q^{(\mu)})^{-1}
        =
        (I-c_jK_Q^{(\mu)})^{-1}.
    \]
\end{proof}

For the rest of this subsection, fix \(j\in\{1,2\}\), and write
\[
    \Xi_j=\{\Psi_1^{(j)},\Psi_2^{(j)}\},
    \qquad
    c\coloneqq c_j.
\]
Set
\begin{equation}\label{eq:equation for S_0j}
    \ell_\mu
    \coloneqq
    \frac1\pi\left(\log\frac1\mu-\gamma_E\right),
    \qquad
    S_{\mu,j}
    \coloneqq
    I-c_jT_Q-c_jE_Q^{(\mu)},
    \qquad
    S_{0,j}\coloneqq I-c_jT_Q.
\end{equation}
For \(m=1,2\), define
\[
    (\Psi_m^{(j)})^{\log}
    \coloneqq
    -\frac1\pi\log|\cdot|*\Psi_m^{(j)},
    \qquad
    \psi_m^{(j)}
    \coloneqq
    |V_Q|^{1/2}(\Psi_m^{(j)})^{\log}.
\]
Define
\[
    (\mathbf m_j)_m
    \coloneqq
    \int_{\mathbb R}\Psi_m^{(j)}(y)\,dy.
\]
Whenever \(S_{\mu,j}\) is invertible, define
\[
    \mathfrak s_j(\mu)
    \coloneqq
    \bigl(
        |V_Q|^{1/2},
        S_{\mu,j}^{-1}|V_Q|^{1/2}
    \bigr)_r,
\]
\[
    (\mathbf a_j(\mu))_m
    \coloneqq
    \bigl(
        \psi_m^{(j)},
        S_{\mu,j}^{-1}|V_Q|^{1/2}
    \bigr)_r,
\]
\[
    (\mathbf R_j(\mu))_{mn}
    \coloneqq
    \bigl(
        \psi_m^{(j)},
        S_{\mu,j}^{-1}\psi_n^{(j)}
    \bigr)_r,
\]
and
\[
    (\mathbf C_{\log,j})_{mn}
    \coloneqq
    \bigl(
        (\Psi_m^{(j)})^{\log},
        \Psi_n^{(j)}
    \bigr)_r.
\]

\begin{proposition}[Explicit formula for the limiting constraint matrices]
\label{prop:constraint-entry-formula}
    Assume that \(S_{0,j}\) is invertible and that
    \[
        \mathfrak s_j(0^+)
        \coloneqq
        \bigl(
            |V_Q|^{1/2},
            S_{0,j}^{-1}|V_Q|^{1/2}
        \bigr)_r
        \neq0.
    \]
    Define
    \begin{equation}\label{eq:limiting constraint matrix}
        \mathbf M_{0,j}
        \coloneqq
        \mathbf C_{\log,j}
        +
        c_j\mathbf R_j(0^+)
        -
        \frac{1}{\mathfrak s_j(0^+)}
        \left(
            \sqrt{c_j}\mathbf a_j(0^+)
            +
            \frac{1}{\sqrt{c_j}}\mathbf m_j
        \right)
        \left(
            \sqrt{c_j}\mathbf a_j(0^+)
            +
            \frac{1}{\sqrt{c_j}}\mathbf m_j
        \right)^\top .
    \end{equation}
    Then, as \(\mu\to0^+\),
    \[
        \mathbf M_j^{(\mu)}
        =
        \mathbf M_{0,j}
        +
        o_{\mu\to0^+}^{\mathrm{ent}}(1).
    \]
    In particular,
    \[
        \det \mathbf M_j^{(\mu)}
        \to
        \det \mathbf M_{0,j}.
    \]
\end{proposition}

\begin{proof}
    For notational convenience, we fix \(j\in\{1,2\}\) throughout the proof.
    Write
    \[
        \mathbf M^{(\mu)}\coloneqq \mathbf M_j^{(\mu)}.
    \]
    From Proposition~\ref{prop:rank-one-decomposition}, we have
    \[
        S_{\mu,j}=S_{0,j}-c_jE_Q^{(\mu)},
        \qquad
        \norm{E_Q^{(\mu)}}_{L^2_{\mathrm e}\to L^2_{\mathrm e}}
        \to0
        \quad\text{as }\mu\to0^+.
    \]
    Hence, if \(S_{0,j}\) is invertible, then \(S_{\mu,j}\) is invertible for all
    sufficiently small \(\mu>0\), and
    \begin{equation*}
        \norm{S_{\mu,j}^{-1}-S_{0,j}^{-1}}_{L^2_{\mathrm e}\to L^2_{\mathrm e}}
        \to0
        \qquad\text{as }\mu\to0^+.
    \end{equation*}
    Therefore
    \begin{equation}\label{eq:s-a-R-to-limits-j}
        \mathfrak s_j(\mu)=\mathfrak s_j(0^+)+o(1),
        \qquad
        \mathbf a_j(\mu)=\mathbf a_j(0^+)+o_{\mu\to0^+}^{\mathrm{ent}}(1),
        \qquad
        \mathbf R_j(\mu)=\mathbf R_j(0^+)+o_{\mu\to0^+}^{\mathrm{ent}}(1).
    \end{equation}
    Since \(\mathfrak s_j(0^+)\neq0\), we also have
    \begin{equation}\label{eq:one-over-s-limit-j}
        \frac1{\mathfrak s_j(\mu)}
        =
        \frac1{\mathfrak s_j(0^+)}
        +
        o(1).
    \end{equation}

    We decompose
    \[
        \mathbf M^{(\mu)}
        =
        \mathbf M_{\mathrm{free},j}^{(\mu)}
        +
        \mathbf M_{\mathrm{corr},j}^{(\mu)}
    \]
    according to the symmetric resolvent identity
    \eqref{eq:woodbury-BS-resolvent}. Thus
    \[
        \bigl(\mathbf M_{\mathrm{free},j}^{(\mu)}\bigr)_{mn}
        \coloneqq
        \bigl(
            G^{(\mu)}\Psi_m^{(j)},
            \Psi_n^{(j)}
        \bigr)_r,
    \]
    and
    \[
        \bigl(\mathbf M_{\mathrm{corr},j}^{(\mu)}\bigr)_{mn}
        \coloneqq
        c_j
        \bigl(
            (I-c_jK_Q^{(\mu)})^{-1}\psi_{m,\mu}^{(j)},
            \psi_{n,\mu}^{(j)}
        \bigr)_r,
    \]
    where, for this proof only,
    \[
        \psi_{m,\mu}^{(j)}
        \coloneqq
        |V_Q|^{1/2}G^{(\mu)}\Psi_m^{(j)}.
    \]

    We first expand the free part. By the small-\(\mu\) decomposition of
    \(G^{(\mu)}\) in \eqref{eq:Kmu-smallmu}, for \(m=1,2\),
    \[
        G^{(\mu)}\Psi_m^{(j)}
        =
        \ell_\mu(\mathbf m_j)_m
        +
        (\Psi_m^{(j)})^{\log}
        +
        \Psi_{m,\mu}^{(j),r},
        \qquad
        \Psi_{m,\mu}^{(j),r}
        \coloneqq
        \frac1\pi r_\mu*\Psi_m^{(j)}.
    \]
    Since \(\Psi_m^{(j)}\in\{Q,\Lambda Q,\Lambda^2Q\}\), we have
    \(|\Psi_m^{(j)}(y)|\lesssim\langle y\rangle^{-2}\). Using the bound on
    \(r_\mu\) in \eqref{eq:rmu-global-bound}, we obtain
    \begin{equation}\label{eq:mulogmu bounds}
        \bigl(
            \Psi_{m,\mu}^{(j),r},
            \Psi_n^{(j)}
        \bigr)_r
        =
        \mathcal O\bigl(\mu\log(1/\mu)\bigr),
        \qquad
        \norm{|V_Q|^{1/2}\Psi_{m,\mu}^{(j),r}}_{L^2}
        =
        \mathcal O\bigl(\mu\log(1/\mu)\bigr).
    \end{equation}
    Indeed, for \(0<\mu\leq\frac12\), splitting at \(|y|=\mu^{-1}\) and using the
    linear and logarithmic bounds for \(r_\mu\), together with
    \(|\Psi_m^{(j)}(y)|\lesssim\langle y\rangle^{-2}\), gives
    \[
        \int_{\mathbb R}r_\mu(y)|\Psi_m^{(j)}(y)|\,dy
        \lesssim
        \mu\log(1/\mu).
    \]
    Since \(r_\mu(x-y)\lesssim r_\mu(x)+r_\mu(y)\) and
    \(\||V_Q|^{1/2}r_\mu\|_{L^2}\lesssim\mu\), the two estimates \eqref{eq:mulogmu bounds} follow.
    Hence
    \begin{equation}\label{eq:free-matrix-expansion-j}
        \mathbf M_{\mathrm{free},j}^{(\mu)}
        =
        \ell_\mu\mathbf m_j\mathbf m_j^\top
        +
        \mathbf C_{\log,j}
        +
        \mathcal O_{\mathrm{ent}}\bigl(\mu\log(1/\mu)\bigr).
    \end{equation}
    Multiplying the above expansion by \(|V_Q|^{1/2}\), we also get
    \begin{equation}\label{eq:psi-mu-expansion-j}
        \psi_{m,\mu}^{(j)}
        =
        \ell_\mu(\mathbf m_j)_m|V_Q|^{1/2}
        +
        \psi_m^{(j)}
        +
        \mathcal O_{L^2}\bigl(\mu\log(1/\mu)\bigr).
    \end{equation}

    Next, the rank-one decomposition of \(K_Q^{(\mu)}\) gives
    \[
        I-c_jK_Q^{(\mu)}
        =
        S_{\mu,j}
        -
        c_j\ell_\mu
        \bigl(
            |V_Q|^{1/2}
            \otimes
            |V_Q|^{1/2}
        \bigr).
    \]
    Since \(I-c_jK_Q^{(\mu)}\) and \(S_{\mu,j}\) are invertible for all
    sufficiently small \(\mu>0\), the rank-one symmetric resolvent identity yields
    \begin{equation}\label{eq:rank-one-Woodbury-j}
        (I-c_jK_Q^{(\mu)})^{-1}
        =
        S_{\mu,j}^{-1}
        +
        \frac{c_j\ell_\mu}{1-c_j\ell_\mu\mathfrak s_j(\mu)}
        \bigl(S_{\mu,j}^{-1}|V_Q|^{1/2}\bigr)
        \otimes
        \bigl(S_{\mu,j}^{-1}|V_Q|^{1/2}\bigr).
    \end{equation}

    Define, only within this proof,
    \[
        \widetilde{\mathbf R}_j^{(\mu)}
        \coloneqq
        \left(
            \bigl(
                \psi_{m,\mu}^{(j)},
                S_{\mu,j}^{-1}\psi_{n,\mu}^{(j)}
            \bigr)_r
        \right)_{m,n=1}^2,
    \]
    and
    \[
        \widetilde{\mathbf a}_j^{(\mu)}
        \coloneqq
        \left(
            \bigl(
                \psi_{m,\mu}^{(j)},
                S_{\mu,j}^{-1}|V_Q|^{1/2}
            \bigr)_r
        \right)_{m=1}^2.
    \]
    Then \eqref{eq:rank-one-Woodbury-j} gives
    \begin{equation}\label{eq:Mcorr-pre-j}
        \mathbf M_{\mathrm{corr},j}^{(\mu)}
        =
        c_j\widetilde{\mathbf R}_j^{(\mu)}
        +
        \frac{c_j^2\ell_\mu}{1-c_j\ell_\mu\mathfrak s_j(\mu)}
        \widetilde{\mathbf a}_j^{(\mu)}
        \bigl(\widetilde{\mathbf a}_j^{(\mu)}\bigr)^\top .
    \end{equation}
    Using \eqref{eq:psi-mu-expansion-j} and the boundedness of
    \(S_{\mu,j}^{-1}\), we obtain
    \begin{equation}\label{eq:Rtilde-expansion-j}
        \widetilde{\mathbf R}_j^{(\mu)}
        =
        \ell_\mu^2\mathfrak s_j(\mu)\mathbf m_j\mathbf m_j^\top
        +
        \ell_\mu
        \bigl(
            \mathbf m_j\mathbf a_j(\mu)^\top
            +
            \mathbf a_j(\mu)\mathbf m_j^\top
        \bigr)
        +
        \mathbf R_j(\mu)
        +
        \mathcal O_{\mathrm{ent}}\bigl(\mu(\log(1/\mu))^2\bigr),
    \end{equation}
    and
    \begin{equation}\label{eq:atilde-expansion-j}
        \widetilde{\mathbf a}_j^{(\mu)}
        =
        \ell_\mu\mathfrak s_j(\mu)\mathbf m_j
        +
        \mathbf a_j(\mu)
        +
        \mathcal O\bigl(\mu\log(1/\mu)\bigr).
    \end{equation}
    Substituting \eqref{eq:Rtilde-expansion-j} and
    \eqref{eq:atilde-expansion-j} into \eqref{eq:Mcorr-pre-j}, we obtain
    \begin{equation}\label{eq:Mcorr-rearranged-j}
    \begin{aligned}
        \mathbf M_{\mathrm{corr},j}^{(\mu)}
        &=
        \frac{
            c_j\ell_\mu^2\mathfrak s_j(\mu)
        }{
            1-c_j\ell_\mu\mathfrak s_j(\mu)
        }
        \mathbf m_j\mathbf m_j^\top
        \\
        &\quad+
        \frac{
            c_j\ell_\mu
        }{
            1-c_j\ell_\mu\mathfrak s_j(\mu)
        }
        \bigl(
            \mathbf m_j\mathbf a_j(\mu)^\top
            +
            \mathbf a_j(\mu)\mathbf m_j^\top
        \bigr)
        \\
        &\quad+
        c_j\mathbf R_j(\mu)
        +
        \frac{
            c_j^2\ell_\mu
        }{
            1-c_j\ell_\mu\mathfrak s_j(\mu)
        }
        \mathbf a_j(\mu)\mathbf a_j(\mu)^\top
        +
        \mathcal O_{\mathrm{ent}}\bigl(\mu(\log(1/\mu))^2\bigr).
    \end{aligned}
    \end{equation}
    Since \(\ell_\mu\to+\infty\), \(\mathfrak s_j(\mu)\to\mathfrak s_j(0^+)\),
    and \(\mathfrak s_j(0^+)\neq0\), we have
    \begin{equation}\label{eq:coefficient limits}
        \begin{aligned}
            \frac{
                c_j\ell_\mu^2\mathfrak s_j(\mu)
            }{
                1-c_j\ell_\mu\mathfrak s_j(\mu)
            }
            &=
            -\ell_\mu
            -
            \frac{1}{c_j\mathfrak s_j(0^+)}
            +
            o(1),
            \\
            \frac{
                c_j\ell_\mu
            }{
                1-c_j\ell_\mu\mathfrak s_j(\mu)
            }
            &=
            -
            \frac{1}{\mathfrak s_j(0^+)}
            +
            o(1),
            \\
            \frac{
                c_j^2\ell_\mu
            }{
                1-c_j\ell_\mu\mathfrak s_j(\mu)
            }
            &=
            -
            \frac{c_j}{\mathfrak s_j(0^+)}
            +
            o(1).
        \end{aligned}
    \end{equation}
    Combining these coefficient limits \eqref{eq:coefficient limits} with
    \eqref{eq:s-a-R-to-limits-j}, \eqref{eq:one-over-s-limit-j}, and
    \eqref{eq:Mcorr-rearranged-j}, we obtain
    \begin{equation}\label{eq:Mcorr-final-j}
    \begin{aligned}
        \mathbf M_{\mathrm{corr},j}^{(\mu)}
        &=
        -\ell_\mu\mathbf m_j\mathbf m_j^\top
        +
        c_j\mathbf R_j(0^+)
        \\
        &\quad
        -
        \frac{1}{\mathfrak s_j(0^+)}
        \left(
            \sqrt{c_j}\mathbf a_j(0^+)
            +
            \frac{1}{\sqrt{c_j}}\mathbf m_j
        \right)
        \left(
            \sqrt{c_j}\mathbf a_j(0^+)
            +
            \frac{1}{\sqrt{c_j}}\mathbf m_j
        \right)^\top
        +
        o_{\mu\to0^+}^{\mathrm{ent}}(1).
    \end{aligned}
    \end{equation}
    Finally, combining \eqref{eq:free-matrix-expansion-j} and
    \eqref{eq:Mcorr-final-j}, the divergent terms
    \(\ell_\mu\mathbf m_j\mathbf m_j^\top\) cancel. Therefore
    \[
        \mathbf M_j^{(\mu)}
        =
        \mathbf M_{0,j}
        +
        o_{\mu\to0^+}^{\mathrm{ent}}(1),
    \]
    which proves the proposition.
\end{proof}

It remains to prove that \(S_{0,j}\) is invertible and \(\mathfrak s_j(0^+)\neq0\), and that \(\det\mathbf M_{0,j}<0\).

\begin{proposition}[Certified limiting matrix conditions]
\label{prop:certified-limiting-matrix-conditions}
    Let \(c_1=6.1\) and \(c_2=2.1\). For \(j=1,2\), the operator
    \[
        S_{0,j}=I-c_jT_Q
    \]
    is invertible on \(L^2_{\mathrm e}\), and
    \begin{equation*}
        \mathfrak s_j(0^+)
        =
        \bigl(
            |V_Q|^{1/2},
            S_{0,j}^{-1}|V_Q|^{1/2}
        \bigr)_r
        \neq0.
    \end{equation*}
    Moreover, the limiting constraint matrices satisfy
    \begin{equation*}
        \det \mathbf M_{0,1}<0,
        \qquad
        \det \mathbf M_{0,2}<0.
    \end{equation*}
\end{proposition}

Proposition~\ref{prop:certified-limiting-matrix-conditions} is proved at the
end of Part~\ref{part:spectral-estimates}, using the profile estimates from
Part~\ref{part:ground-state-approximation} and the interval bounds from
Part~\ref{part:interval-verification}.

\begin{corollary}[Finite-\(\mu\) determinant signs]
\label{cor:finite-mu-determinant-signs}
    After possibly decreasing \(\mu_0>0\), for every \(j=1,2\) and
    \(0<\mu<\mu_0\), we have
    \[
        \det \mathbf M_j^{(\mu)}<0.
    \]
    In particular,
    \[
        n_-(\mathbf M_j^{(\mu)})=1.
    \]
\end{corollary}

\begin{proof}
    Propositions~\ref{prop:constraint-entry-formula} and
    \ref{prop:certified-limiting-matrix-conditions} give
    \[
        \mathbf M_j^{(\mu)}
        =
        \mathbf M_{0,j}
        +
        o_{\mu\to0^+}^{\mathrm{ent}}(1),
        \qquad
        \det \mathbf M_{0,j}<0.
    \]
    Hence \(\det\mathbf M_j^{(\mu)}<0\) for all sufficiently small
    \(\mu>0\). Since \(\mathbf M_j^{(\mu)}\) is a real symmetric
    \(2\times2\) matrix, it has exactly one negative eigenvalue.
\end{proof}

\begin{proof}[Proof of Proposition~\ref{prop:shifted-index-determinant-conditions}]
    We take \(c_1=6.1\) and \(c_2=2.1\). By Corollary~\ref{cor:shifted-morse-index-invertibility}, \(0\in\rho(\mathcal L_{c_j}^{(\mu)})\) and \(n_-^{\mathrm e}(\mathcal L_{c_j}^{(\mu)})=1\). By Corollary~\ref{cor:finite-mu-determinant-signs}, after possibly decreasing \(\mu_0>0\), \(\det \mathbf M_j^{(\mu)}<0\).
\end{proof}

This reduces Theorem~\ref{thm:Spectral structure of the bilinear form H} to Propositions~\ref{prop:reduced-BS-bound-Q} and \ref{prop:certified-limiting-matrix-conditions}, proved in Parts~\ref{part:ground-state-approximation}--\ref{part:interval-verification}.

\part{Approximation of the ground state}\label{part:ground-state-approximation}
Propositions~\ref{prop:reduced-BS-bound-Q} and \ref{prop:certified-limiting-matrix-conditions} depend quantitatively on the ground state \(Q\), for which no closed formula is available. In Part~\ref{part:ground-state-approximation}, we construct an explicit approximate profile \(\mathfrak g\) and prove the
quantitative estimates for \(Q-\mathfrak g\) used in Parts~\ref{part:spectral-estimates} and~\ref{part:interval-verification}.

In Section~\ref{sec:construction-approximate-profile}, we define \(\mathfrak g\) and bound its residual in the ground-state equation. In Section~\ref{sec:NK-validation-profile}, we apply the Newton--Kantorovich theorem to obtain an exact even solution near \(\mathfrak g\), and then
identify this solution with \(Q\). In Section~\ref{sec:higher-order-profile-estimates}, we
derive the higher-order and weighted estimates required in Section~\ref{sec:transfer from g to Q}. Finally, in Section~\ref{sec:inverse-bound-g}, we prove the inverse estimate for
\(\mathcal F'(\mathfrak g)\) used in the Newton--Kantorovich argument.

\section{Construction of the approximate profile}\label{sec:construction-approximate-profile}

In this section, we construct an approximate profile \(\mathfrak g\) for the ground state \(Q\). The nonlocality of \(|D|\) makes it difficult to construct an approximate profile directly on \(\bbR\). Indeed, the construction near the origin depends on the algebraic decay of the profile. We therefore compactify the domain
\[
    I\coloneqq\left(-\frac{\pi}{2},\frac{\pi}{2}\right),
    \qquad
    y=\tan\theta,
    \qquad
    \theta\in I,
\]
and define the compactified profile \(Q_{\mathrm{cp}}\) by
\begin{equation}\label{eq:def-Qcp}
    Q(\tan\theta)
    =
    \cos^2\theta\,Q_{\mathrm{cp}}(\theta),
    \qquad
    \theta\in I.
\end{equation}
Since \(Q(y)\) has algebraic decay of order \(\langle y\rangle^{-2}\), the function \(Q(\tan\theta)\) vanishes at the endpoints of \(I\) with the factor \(\cos^2\theta\). Indeed, the change of variables gives \(\langle y\rangle^{-2} = (1+\tan^2\theta)^{-1} = \cos^2\theta\). Thus the factor \(\cos^2\theta\) extracts the prescribed algebraic decay, leaving \(Q_{\mathrm{cp}}\) to be approximated on the bounded interval \(I\).

Moreover, the conformal covariance of \(|D|\) converts the nonlocal operator on \(\bbR\) into \(\cos^2\theta\,|\partial_\theta|\) on \(I\). Here \(|\partial_\theta|\) denotes the Fourier multiplier on \(\pi\)-periodic functions with symbol \(2|n|\). We rewrite the ground-state equation in terms of \(Q_{\mathrm{cp}}\).

\begin{lemma}[Compactified ground-state equation]
    The compactified profile \(Q_{\mathrm{cp}}\) defined by \eqref{eq:def-Qcp}
    satisfies
    \begin{equation}\label{eq:compactified-ground-state-equation}
        \bigl(|\partial_\theta|(\cos^2\theta\,Q_{\mathrm{cp}})\bigr)(\theta)
        +Q_{\mathrm{cp}}(\theta)
        =
        \cos^4\theta\,Q_{\mathrm{cp}}(\theta)^3,
        \qquad \theta\in I,
    \end{equation}
    where \(|\partial_\theta|\) is the Fourier multiplier on \(\pi\)-periodic
    functions with symbol \(2|n|\).
\end{lemma}

\begin{proof}
    We use the conformal covariance identity
    \[
        (|D|Q)(\tan\theta)
        =
        \cos^2\theta\,
        |\partial_\theta|\bigl(Q(\tan\cdot)\bigr)(\theta).
    \]
    Evaluating the soliton equation at \(y=\tan\theta\) and using
    \eqref{eq:def-Qcp}, we obtain
    \[
        0
        =
        \cos^2\theta\,
        |\partial_\theta|(\cos^2\theta\,Q_{\mathrm{cp}})
        +\cos^2\theta\,Q_{\mathrm{cp}}
        -\cos^6\theta\,Q_{\mathrm{cp}}^3.
    \]
    Dividing by \(\cos^2\theta\) gives
    \eqref{eq:compactified-ground-state-equation}.
\end{proof}

We now approximate \(Q_{\mathrm{cp}}\) by a finite even cosine polynomial to make profile error small.

\begin{definition}[The approximate profile \(\mathfrak g\)]
\label{def:approximate-profile-g}
    Let \(\{U_k\}_{k=0}^{200}\subset\mathbb{Q}\) be the coefficients
    recorded in
    \[
        \texttt{profile\_g\_U\_cos\_J200.csv}.
    \]
    Define
    \begin{equation}\label{eq:def-U-profile}
        U(\theta)
        \coloneqq
        \sum_{k=0}^{200}U_k\cos(2k\theta),
        \qquad \theta\in I,
    \end{equation}
    and
    \begin{equation*}
        \mathfrak g(y)
        \coloneqq
        \frac{1}{1+y^2}U(\arctan y).
    \end{equation*}
    Equivalently,
    \begin{equation}\label{eq:def-g-compactified}
        \mathfrak g(\tan\theta)
        =
        \cos^2\theta\,U(\theta),
        \qquad \theta\in I.
    \end{equation}
\end{definition}

We now estimate the residual of \(\mathfrak g\) in the ground-state equation. Define
\begin{equation*}
    \mathcal F:H^1_{\mathrm e}\to L^2_{\mathrm e},
    \qquad
    \mathcal F(f)\coloneqq |D|f+f-f^3.
\end{equation*}
The change of variables \(y=\tan\theta\) gives
\begin{equation*}
    \norm{\mathcal F(\mathfrak g)}_{L^2(\bbR)}^2
    =
    \int_I
        \cos^2\theta
        \left|
            |\partial_\theta|(\cos^2\theta\,U)
            +
            U
            -
            \cos^4\theta\,U^3
        \right|^2
    \,d\theta .
\end{equation*}

Since \(U\) is a finite even cosine polynomial, the function
\[
    \cos\theta
    \left(
        |\partial_\theta|(\cos^2\theta\,U)
        +U
        -\cos^4\theta\,U^3
    \right)
\]
has a finite cosine expansion. Its \(L^2(I)\)-norm is calculated by Parseval's identity. The higher Sobolev norms of \(\mathfrak g\) and \(\mathcal F(\mathfrak g)\), and the \(L^2\)-norms of \(\Lambda^a\mathfrak g\) and \(\Lambda^a\mathcal F(\mathfrak g)\) for \(a=1,2\), are calculated in the same way. The corresponding identities on \(I\) are recorded in
Appendix~\ref{app:compactification-profile-calculus}, and the resulting bounds are stated in
Proposition~\ref{prop:part5-g-profile-residual-bounds}. We denote by \(\delta_{\mathrm{res}}^{\mathfrak g}\) the resulting outward-rounded upper bound for \(\norm{\mathcal F(\mathfrak g)}_{L^2(\bbR)}\).

\begin{proposition}[Certified approximate profile]
\label{prop:certified-approximate-profile}
    The profile \(\mathfrak g\) in
    Definition~\ref{def:approximate-profile-g} belongs to \(C^\infty(\bbR)\),
    is strictly positive on \(\bbR\), and satisfies
    \begin{equation*}
        \norm{\mathcal F(\mathfrak g)}_{L^2}
        \leq
        \delta_{\mathrm{res}}^{\mathfrak g}.
    \end{equation*}
    Moreover,
    \begin{equation*}
        \delta_{\mathrm{res}}^{\mathfrak g}
        \leq
        2.51\times10^{-14}.
    \end{equation*}
\end{proposition}

\begin{proof}
    The explicit formula in Definition~\ref{def:approximate-profile-g} gives \(\mathfrak g\in C^\infty(\bbR)\). The positivity and residual bounds follow from Proposition~\ref{prop:part5-g-profile-residual-bounds}.
\end{proof}

\section{Newton--Kantorovich validation of the approximate profile}\label{sec:NK-validation-profile}
In the preceding section, we constructed a positive even profile \(\mathfrak g\) with a small residual in the ground-state equation. We now prove a quantitative \(H^1\)-bound for \(Q-\mathfrak g\).

In Subsection~\ref{sec:NK-exact-solution}, the Newton--Kantorovich theorem gives an exact even solution \(\widetilde Q\) in an explicit \(H^1\)-neighborhood of \(\mathfrak g\). The required inverse estimate for \(\mathcal F'(\mathfrak g)\) is proved in Section~\ref{sec:inverse-bound-g}. In Subsection~\ref{sec:identification-ground-state}, we prove that \(\widetilde Q\) is positive and that \(\mathcal F'(\widetilde Q)\) has Morse index one. The uniqueness theorem
due to Frank--Lenzmann--Silvestre \cite{FrankLenzmannSilvestre2016CPAM} then identifies \(\widetilde Q\) with \(Q\).

\begin{proposition}[Certified closeness of \(\mathfrak g\) to \(Q\)]
\label{prop:closeness-of-g-to-Q}
    Let \(Q\) be the positive even ground state fixed above. Then
    \begin{equation}\label{eq:certified-H1-closeness-Q-g}
        \norm{Q-\mathfrak g}_{H^1}
        \leq r_{\mathrm{NK}},
    \end{equation}
    where \(r_{\mathrm{NK}}\) is the Newton--Kantorovich radius defined in
    \eqref{eq:NK-radius-bound}. In particular, the certified enclosure gives
    \[
        r_{\mathrm{NK}}
        \leq 5.617\times10^{-13}.
    \]
\end{proposition}

\subsection{Newton--Kantorovich validation near \(\mathfrak g\)}
\label{sec:NK-exact-solution}

We verify the Newton--Kantorovich hypotheses at the base point \(\mathfrak g\) to prove the existence of an exact even solution \(\widetilde Q\) in a quantitative \(H^1\)-neighborhood of \(\mathfrak g\).

\begin{theorem}[Newton--Kantorovich criterion for \(\mathcal F\)]\label{thm:NK-for-F}
    Assume that
    \begin{equation}\label{eq:NK-bijectivity}
        \mathcal F'(\mathfrak g)\big|_{L^2_{\mathrm e}}:H^1_{\mathrm e}\to L^2_{\mathrm e}
        \quad\text{is bijective}.
    \end{equation}
    Suppose that there exist constants
    \(A_{\mathrm{NK}},B_{\mathrm{NK}},C_{\mathrm{NK}}>0\) such that
    \begin{equation}\label{eq:NK-step-and-inverse-bounds}
        \norm{(\mathcal F'(\mathfrak g))^{-1}\mathcal F(\mathfrak g)}_{H^1}
        \leq A_{\mathrm{NK}},
        \qquad
        \norm{(\mathcal F'(\mathfrak g))^{-1}}_
        {L^2_{\mathrm e}\to H^1_{\mathrm e}}
        \leq B_{\mathrm{NK}},
    \end{equation}
    and
    \begin{equation}\label{eq:NK-Lipschitz-bound}
        \norm{\mathcal F'(f_1)-\mathcal F'(f_2)}_
        {H^1_{\mathrm e}\to L^2_{\mathrm e}}
        \leq
        C_{\mathrm{NK}}\norm{f_1-f_2}_{H^1}
    \end{equation}
    for all
    \[
        f_1,f_2\in
        B_{H^1_{\mathrm e}}
        \left(\mathfrak g,\frac{1}{B_{\mathrm{NK}}C_{\mathrm{NK}}}\right).
    \]
    If
    \begin{equation*}
        0<
        A_{\mathrm{NK}}B_{\mathrm{NK}}C_{\mathrm{NK}}
        \leq
        \frac12,
    \end{equation*}
    then there exists
    \[
        \tilde Q\in H^1_{\mathrm e}
    \]
    such that
    \begin{equation*}
        \mathcal F(\tilde Q)=0
    \end{equation*}
    and
    \begin{equation}\label{eq:NK-radius-bound}
        \norm{\tilde Q-\mathfrak g}_{H^1}
        \leq
        r_{\mathrm{NK}}
        \coloneqq
        \frac{
            1-\sqrt{1-2A_{\mathrm{NK}}B_{\mathrm{NK}}C_{\mathrm{NK}}}
        }{
            B_{\mathrm{NK}}C_{\mathrm{NK}}
        }.
    \end{equation}
\end{theorem}

For the nonlinear map \(\mathcal F\), the linearization at \(\mathfrak g\) is
\[
    \mathcal F'(\mathfrak g)
    =
    |D|+1-3\mathfrak g^2.
\]
The residual bound required in Theorem~\ref{thm:NK-for-F} follows from Proposition~\ref{prop:certified-approximate-profile}. It remains to prove the bijectivity of \(\mathcal F'(\mathfrak g)\) on the even subspace, establish the inverse bound in \eqref{eq:NK-step-and-inverse-bounds}, and verify the local Lipschitz estimate \eqref{eq:NK-Lipschitz-bound}. The following proposition gives the bijectivity and an \(L^2\)-inverse bound.

\begin{proposition}[Certified spectral structure of $\mathcal F'(\mathfrak g)$]
\label{prop:invertibility-Fprimeg}
    We have
    \begin{equation}\label{eq:certified-L2-inverse-Fprimeg}
        n_-^{\mathrm e}\left(
            \mathcal F'(\mathfrak g)
        \right)
        =
        1,
        \qquad
        0\in
        \rho\left(
            \mathcal F'(\mathfrak g)\big|_{L^2_{\mathrm e}}
        \right),
        \qquad
        \norm{[\mathcal F'(\mathfrak g)\big|_{L^2_{\mathrm e}}]^{-1}}_
        {L^2_{\mathrm e}\to L^2_{\mathrm e}}
        \leq 2.
    \end{equation}
\end{proposition}

The proof of Proposition~\ref{prop:invertibility-Fprimeg} is deferred to
Section~\ref{sec:inverse-bound-g}. The resolvent statement in
\eqref{eq:certified-L2-inverse-Fprimeg} gives the bijectivity required in
\eqref{eq:NK-bijectivity}.

\begin{lemma}[Choice of Newton--Kantorovich constants]
    Define
    \begin{equation}\label{eq:def-NK-constants}
        \begin{aligned}
            A_{\mathrm{NK}}
            &\coloneqq
            \delta_{\mathrm{res}}^{\mathfrak g}
            \left(1+6\norm{\mathfrak g}_{L^\infty}^2\right),\\
            B_{\mathrm{NK}}
            &\coloneqq
            1+6\norm{\mathfrak g}_{L^\infty}^2,\\
            C_{\mathrm{NK}}
            &\coloneqq
            3\norm{\mathfrak g}_{L^\infty}
            +
            \sqrt{
                9\norm{\mathfrak g}_{L^\infty}^2
                +
                \frac{6}{B_{\mathrm{NK}}}
            }.
        \end{aligned}
    \end{equation}
    Then \(A_{\mathrm{NK}},B_{\mathrm{NK}},C_{\mathrm{NK}}\) satisfy
    \eqref{eq:NK-step-and-inverse-bounds} and
    \eqref{eq:NK-Lipschitz-bound}.
\end{lemma}

\begin{proof}
    Let \(f\in L^2_{\mathrm e}\) and set \(v=(\mathcal F'(\mathfrak g)\big|_{L^2_{\mathrm e}})^{-1}f\). Since
    \((|D|+1)v=f+3\mathfrak g^2v\), Proposition~\ref{prop:invertibility-Fprimeg}
    and \eqref{eq:def-NK-constants} give
    \[
        \norm{v}_{H^1}
        \le
        \norm{f}_{L^2}
        +
        3\norm{\mathfrak g}_{L^\infty}^2\norm{v}_{L^2}
        \le
        B_{\mathrm{NK}}\norm{f}_{L^2}.
    \]
    Hence
    \[
        \norm{(\mathcal F'(\mathfrak g)\big|_{L^2_{\mathrm e}})^{-1}}_
        {L^2_{\mathrm e}\to H^1_{\mathrm e}}
        \le B_{\mathrm{NK}}.
    \]
    Applying this with \(f=\mathcal F(\mathfrak g)\) and using
    \(\norm{\mathcal F(\mathfrak g)}_{L^2}\le\delta_{\mathrm{res}}^{\mathfrak g}\),
    we also get
    \[
        \norm{(\mathcal F'(\mathfrak g)\big|_{L^2_{\mathrm e}})^{-1}\mathcal F(\mathfrak g)}_{H^1}
        \le
        B_{\mathrm{NK}}\delta_{\mathrm{res}}^{\mathfrak g}
        =
        A_{\mathrm{NK}}.
    \]
    This proves \eqref{eq:NK-step-and-inverse-bounds}.

    For the Lipschitz bound, use
    \(\mathcal F'(f_1)-\mathcal F'(f_2)=-3(f_1+f_2)(f_1-f_2)\).
    On the ball appearing in \eqref{eq:NK-Lipschitz-bound}, the Sobolev
    embedding gives
    \[
        \norm{f_i}_{L^\infty}
        \le
        \norm{\mathfrak g}_{L^\infty}
        +
        \frac{1}{B_{\mathrm{NK}}C_{\mathrm{NK}}},
        \qquad i=1,2.
    \]
    Therefore
    \[
        \norm{\mathcal F'(f_1)-\mathcal F'(f_2)}_
        {H^1_{\mathrm e}\to L^2_{\mathrm e}}
        \le
        6\left(
            \norm{\mathfrak g}_{L^\infty}
            +
            \frac{1}{B_{\mathrm{NK}}C_{\mathrm{NK}}}
        \right)
        \norm{f_1-f_2}_{H^1}.
    \]
    By \eqref{eq:def-NK-constants}, \(C_{\mathrm{NK}}\) is the positive
    solution of
    \[
        C
        =
        6\norm{\mathfrak g}_{L^\infty}
        +
        \frac{6}{B_{\mathrm{NK}}C}.
    \]
    Thus the last coefficient is bounded by \(C_{\mathrm{NK}}\), proving
    \eqref{eq:NK-Lipschitz-bound}.
\end{proof}

The preceding lemma verifies the analytic hypotheses of Theorem~\ref{thm:NK-for-F}. It remains to verify the smallness condition using the certified bounds for \(\delta_{\mathrm{res}}^{\mathfrak g}\) and \(\norm{\mathfrak g}_{L^\infty}\).

\begin{corollary}[Validated exact solution near \(\mathfrak g\)]\label{cor:certified-NK-radius}
    Let \(A_{\mathrm{NK}},B_{\mathrm{NK}},C_{\mathrm{NK}}\) be defined in
    \eqref{eq:def-NK-constants}, and let \(r_{\mathrm{NK}}\) be defined in
    \eqref{eq:NK-radius-bound}.  Then
    \begin{equation*}
        A_{\mathrm{NK}}B_{\mathrm{NK}}C_{\mathrm{NK}}
        <
        \frac12,
        \qquad
        r_{\mathrm{NK}}
        \leq
        5.617\times10^{-13}.
    \end{equation*}
    Therefore, there exists \(\tilde Q\in H^1_{\mathrm e}\) such that
    \begin{equation}\label{eq:tildeQ-g-NK-radius}
        \mathcal F(\tilde Q)=0,
        \qquad
        \norm{\tilde Q-\mathfrak g}_{H^1}
        \leq
        r_{\mathrm{NK}}.
    \end{equation}
\end{corollary}

Corollary~\ref{cor:certified-NK-radius} follows from Proposition~\ref{prop:nk-Q-transfer-bounds}.

\subsection{Identification with the ground state}
\label{sec:identification-ground-state}

Let \(\widetilde Q\) be the exact even solution obtained in Section~\ref{sec:NK-exact-solution}. We prove that \(\widetilde Q\) is positive and that \(\mathcal F'(\tilde Q)\) has Morse index one. The uniqueness theorem of Frank--Lenzmann--Silvestre then gives \(\widetilde Q=Q\).

We first establish positivity and strict monotonicity of the exact solution obtained in \eqref{eq:tildeQ-g-NK-radius}.

\begin{lemma}[Positivity and monotonicity of the validated solution]
\label{lem:positivity-monotonicity-tildeQ}
    The solution \(\tilde Q\) in \eqref{eq:tildeQ-g-NK-radius} is strictly
    positive on \(\bbR\) and satisfies
    \begin{equation}\label{eq:tildeQ-strict-monotonicity}
        \tilde Q'(y)<0,
        \qquad
        y>0.
    \end{equation}
\end{lemma}

\begin{proof}
    Since \(\mathcal F(\tilde Q)=0\), we have
    \begin{equation}\label{eq:tildeQ-resolvent-representation}
        \tilde Q
        =
        (|D|+1)^{-1}(\tilde Q^3).
    \end{equation}
    Since \(\tilde Q\in H^1(\bbR)\hookrightarrow L^\infty(\bbR)\), we have \(\tilde Q^3\in H^1(\bbR)\), and hence \(\tilde Q\in H^2(\bbR)\).

    We first prove that \(\tilde Q\geq0\). Assume that \(\widetilde Q\) is negative somewhere. Since \(\widetilde Q\in H^1(\bbR)\), it is continuous and satisfies
    \(\widetilde Q(y)\to0\) as \(|y|\to\infty\). Hence, there exists \(y_{\min}\in\bbR\) such that
    \[
        m \coloneqq \widetilde Q(y_{\min})=\min_{\bbR}\widetilde Q<0.
    \]
    At the minimum point \(y_{\min}\), 
    \[
        |D|\tilde Q(y_{\min})
        =
        c_*
        \operatorname{p.v.}
        \int_{\bbR}
            \frac{
                \tilde Q(y_{\min})-\tilde Q(z)
            }{
                (y_{\min}-z)^2
            }
        \,dz
        \leq0
    \]
    for some constant \(c_*>0\). Evaluating \(\mathcal F(\tilde Q)=0\) at \(y_{\min}\), 
    \[
        m^3-m
        =
        |D|\tilde Q(y_{\min})
        \leq0.
    \]
    Since \(m<0\), we have \(m\leq-1\).

    On the other hand, Proposition~\ref{prop:certified-approximate-profile}
    gives \(\mathfrak g(y_{\min})>0\). Hence, by the Sobolev embedding and
    \eqref{eq:tildeQ-g-NK-radius},
    \[
        m
        =
        \tilde Q(y_{\min})
        \geq
        \mathfrak g(y_{\min})
        -
        \norm{\tilde Q-\mathfrak g}_{L^\infty}
        >
        -
        \frac{r_{\mathrm{NK}}}{\sqrt2}.
    \]
    Corollary~\ref{cor:certified-NK-radius} gives
    \(\tfrac{r_{\mathrm{NK}}}{\sqrt2}<1\), contradicting \(m\leq-1\).
    Thus \(\tilde Q\geq0\).

    We next show that \(\tilde Q\not\equiv0\). If not, Proposition~\ref{prop:part5-g-profile-residual-bounds} and \eqref{eq:tildeQ-g-NK-radius} imply
    \[
        \frac12
        <
        \mathfrak g(0)
        \leq
        \norm{\mathfrak g-\tilde Q}_{L^\infty}
        \leq
        \frac{r_{\mathrm{NK}}}{\sqrt2}.
    \]
    This contradicts Corollary~\ref{cor:certified-NK-radius}. Since the
    convolution kernel of \((|D|+1)^{-1}\) is strictly positive, \eqref{eq:tildeQ-resolvent-representation} yields
    \[
        \tilde Q(y)>0,
        \qquad
        y\in\bbR.
    \]

    By \cite[Proposition~3.1(ii)]{FrankLenzmannSilvestre2016CPAM} and \(\tilde Q \in H^1_{\mathrm e}\), we obtain \eqref{eq:tildeQ-strict-monotonicity}.
\end{proof}

We next determine the Morse index of the linearized operator at \(\widetilde Q\).

\begin{lemma}[Morse index of the validated solution]
\label{lem:morse-index-tildeQ}
    The linearized operator
    \[
        \mathcal F'(\tilde Q)
        =
        |D|+1-3\tilde Q^2
    \]
    satisfies
    \begin{equation}\label{eq:tildeQ-full-morse-index}
        n_-\bigl(\mathcal F'(\tilde Q)\bigr)=1.
    \end{equation}
\end{lemma}

\begin{proof}
    We prove \eqref{eq:tildeQ-full-morse-index} by considering the even and odd subspaces separately. 
    
    On \(L^2_{\mathrm e}\), \eqref{eq:tildeQ-g-NK-radius} and \eqref{eq:H1-Linfty} yield
    \begin{equation*}
        \norm{\tilde Q^2-\mathfrak g^2}_{L^\infty}
        \leq
        \frac{r_{\mathrm{NK}}}{\sqrt2}
        \left(
            2\norm{\mathfrak g}_{L^\infty}
            +
            \frac{r_{\mathrm{NK}}}{\sqrt2}
        \right).
    \end{equation*}
    Proposition~\ref{prop:invertibility-Fprimeg}, Corollary~\ref{cor:certified-NK-radius}, and the certified bound for \(\norm{\mathfrak g}_{L^\infty}\) in Proposition~\ref{prop:part5-g-profile-residual-bounds} give
    \begin{equation*}
        3
        \norm{[\mathcal F'(\mathfrak g)\big|_{L^2_{\mathrm e}}]^{-1}}_
        {L^2_{\mathrm e}\to L^2_{\mathrm e}}
        \norm{\tilde Q^2-\mathfrak g^2}_{L^\infty}
        <1.
    \end{equation*}

    Hence, for \(0\leq\tau\leq1\)
    \[
        \mathcal L_\tau \coloneqq\mathcal F'(\mathfrak g)
        -
        3\tau(\tilde Q^2-\mathfrak g^2)
        =
        \mathcal F'(\mathfrak g)
        \left[
            I
            -
            3\tau
            [\mathcal F'(\mathfrak g)\big|_{L^2_{\mathrm e}}]^{-1}
            (\tilde Q^2-\mathfrak g^2)
        \right].
    \]
    By Lemma~\ref{lem:Neumann-series-criterion}, \(\mathcal L_\tau\) is invertible on \(L^2_{\mathrm e}\) for every \(\tau\in[0,1]\). Since \(\mathcal L_\tau\) is a continuous family of self-adjoint operators, \(n_-^{\mathrm e}(\mathcal L_\tau)\) is constant
    on \([0,1]\). Therefore, Proposition~\ref{prop:invertibility-Fprimeg} yields
    \begin{equation}\label{eq:tildeQ-even-morse-index}
        n_-^{\mathrm e}\left(
            \mathcal F'(\tilde Q)
        \right)
        =
        1.
    \end{equation}

    Differentiating \(\mathcal F(\widetilde Q)=0\) with respect to \(y\), we obtain
    \begin{equation*}
        \mathcal F'(\widetilde Q)[\widetilde Q']=0.
    \end{equation*}
    By Lemma~\ref{lem:positivity-monotonicity-tildeQ}, \(-\tilde Q'\) is an odd function satisfying
    \[
        -\tilde Q'(y)>0,
        \qquad
        y>0.
    \]
    We claim that the odd restriction of \(\mathcal F'(\tilde Q)\) has no
    negative eigenvalue. Assume that \(\mathcal F'(\widetilde Q)\big|_{L^2_{\mathrm o}}\) has a negative eigenvalue. Its lowest eigenvalue is then negative and discrete, since
    the essential spectrum begins at \(1\).
    
    By \cite[Lemma~C.3]{FrankLenzmann2013Acta}, the corresponding odd
    eigenfunction \(\phi_\mathrm{o}\) can be chosen strictly positive on \((0,\infty)\).
    Then the inner product \((\phi_\mathrm{o},-\tilde Q')_r\) is strictly positive.
    This contradicts the orthogonality of eigenfunctions of a self-adjoint operator corresponding to distinct eigenvalues. Therefore, we conclude
    \begin{equation}\label{eq:tildeQ-odd-morse-index}
        n_-^{\mathrm o}\left(
            \mathcal F'(\tilde Q)
        \right)
        =
        0.
    \end{equation}
    Since \(\tilde Q\) is even, \(\mathcal F'(\tilde Q)\) preserves the even-odd decomposition
    \[
        L^2(\bbR)
        =
        L^2_{\mathrm e}(\bbR)
        \oplus
        L^2_{\mathrm o}(\bbR).
    \]
    Combining \eqref{eq:tildeQ-even-morse-index} and \eqref{eq:tildeQ-odd-morse-index}, we obtain \eqref{eq:tildeQ-full-morse-index}.
\end{proof}

\begin{proof}[Proof of Proposition~\ref{prop:closeness-of-g-to-Q}]
    By Lemma~\ref{lem:positivity-monotonicity-tildeQ},
    \(\tilde Q\) is a nontrivial positive even solution of
    \[
        |D|\tilde Q+\tilde Q-\tilde Q^3=0.
    \]
    Lemma~\ref{lem:morse-index-tildeQ} shows that \(n_-(\mathcal F'(\tilde Q))=1\). Hence \(\tilde Q\) is a ground state in the sense of \cite[Definition~3.1]{FrankLenzmannSilvestre2016CPAM}. Therefore, \cite[Theorem~4]{FrankLenzmannSilvestre2016CPAM} gives
    \[
        \tilde Q=Q.
    \]
    Substituting this identity into \eqref{eq:tildeQ-g-NK-radius} gives \eqref{eq:certified-H1-closeness-Q-g}.
\end{proof}

\section{Higher-order profile estimates and potential bounds}\label{sec:higher-order-profile-estimates}
In this section, we set
\begin{equation}\label{eq:def-e_Q}
    e_Q\coloneqq Q-\mathfrak g.
\end{equation}
The preceding section gives a quantitative \(H^1\)-estimate for \(e_Q\).
However, Section~\ref{sec:transfer from g to Q} requires higher-order and weighted \(W^{k,2}\)-bounds for \(e_Q\). In this section, we derive these quantitative bounds from the equation satisfied by \(e_Q\), the explicit form of \(\mathfrak g\), and the profile error \(\mathcal F(\mathfrak g)\).

We obtain \(W^{2,2}\)- and \(W^{3,2}\)-bounds for \(e_Q\), together with pointwise bounds for its first two derivatives. Using the resolvent representation \(Q=(|D|+1)^{-1}Q^3\), we obtain global weighted bounds for \(Q\) and its derivatives. We estimate \(\Lambda^ae_Q\), \(a=1,2\), in the \(L^2\), \(L^1\), and logarithmically weighted \(L^1\) norms required below. Finally, we obtain the square-root potential estimate used in Section~\ref{sec:transfer from g to Q}.

\begin{lemma}[Higher-order transfer from the \(H^1\) bound]
\label{lem:higher-order-transfer-from-H1}
    Let the constants
    \(\Delta_{\mathrm{cub}}^{W^{1,2}}\), \(\Delta_{|D|Q}^{W^{1,2}}\), \(\Delta_Q^{W^{2,2}}\),
    \(\Delta_{Q,0}^{L^\infty}\), \(\Delta_{Q,1}^{L^\infty}\), \(C_{\mathrm{cub},0}^{L^2}\),
    \(C_{\mathrm{cub},1}^{L^2}\), \(C_{\mathrm{cub},2}^{L^2}\),
    \(\Delta_{\mathrm{cub}}^{W^{2,2}}\), \(\Delta_{|D|Q}^{W^{2,2}}\), \(\Delta_Q^{W^{3,2}}\), and
    \(\Delta_{Q,2}^{L^\infty}\) be defined by \eqref{eq:cubic-w12-error-definition}--\eqref{eq:higher-sobolev-error-definitions}.
    Then
    \begin{equation}\label{eq:e_Q unweighted higer order bound}
        \norm{e_Q}_{W^{2,2}}\le \Delta_Q^{W^{2,2}},
        \qquad
        \norm{e_Q}_{W^{3,2}}\le \Delta_Q^{W^{3,2}},
    \end{equation}
    and
    \begin{equation}\label{eq:e_Q L^infty bound}
        \norm{e_Q}_{L^\infty}\le \Delta_{Q,0}^{L^\infty},
        \qquad
        \norm{\partial_y e_Q}_{L^\infty}\le \Delta_{Q,1}^{L^\infty},
        \qquad
        \norm{\partial_y^2 e_Q}_{L^\infty}\le \Delta_{Q,2}^{L^\infty}.
    \end{equation}
    In particular, the constants on the right-hand sides are determined
    entirely by \(r_{\mathrm{NK}}\), the explicit profile \(\mathfrak g\),
    and the residual \(r_{\mathfrak g}\).
\end{lemma}

\begin{proof}
    Since \(\mathcal F(Q)=0\) and \(r_{\mathfrak g}=\mathcal F(\mathfrak g)\),
    the difference \(e_Q=Q-\mathfrak g\) satisfies
    \[
        |D|e_Q
        =
        -e_Q+(Q^3-\mathfrak g^3)-r_{\mathfrak g}.
    \]
    We first estimate the cubic difference in \(W^{1,2}\). Since the \(H^1\)- and \(W^{1,2}\)-norms coincide, \eqref{eq:H1-Linfty} gives \(\norm{Q}_{H^1}\le \norm{\mathfrak g}_{H^1}+r_{\mathrm{NK}}\). The one-dimensional \(H^1\)-product estimate \eqref{eq:H1-product} then gives
    \[
        \norm{Q^3-\mathfrak g^3}_{W^{1,2}}
        \le
        \Delta_{\mathrm{cub}}^{W^{1,2}}.
    \]
    Here, \(\Delta_{\mathrm{cub}}^{W^{1,2}}\) is defined in \eqref{eq:cubic-w12-error-definition}, which is determined by $r_{\mathrm{NK}}$ and $\mathfrak g$. Hence, we have
    \[
        \norm{|D|e_Q}_{W^{1,2}}
        \le
        r_{\mathrm{NK}}
        +
        \Delta_{\mathrm{cub}}^{W^{1,2}}
        +
        \norm{r_{\mathfrak g}}_{W^{1,2}}
        =
        \Delta_{|D|Q}^{W^{1,2}}.
    \]
    Therefore
    \[
        \norm{e_Q}_{W^{2,2}}
        \le
        \left(r_{\mathrm{NK}}^2+(\Delta_{|D|Q}^{W^{1,2}})^2\right)^{1/2}
        =
        \Delta_Q^{W^{2,2}},
    \]
    which yields the first bound in \eqref{eq:e_Q unweighted higer order bound}.
    Expanding \(Q^3-\mathfrak g^3=3\mathfrak g^2e_Q+3\mathfrak g e_Q^2+e_Q^3\)
    and applying the product rule, the three quantities
    \(C_{\mathrm{cub},0}^{L^2}\), \(C_{\mathrm{cub},1}^{L^2}\), and
    \(C_{\mathrm{cub},2}^{L^2}\) bound $\norm{\partial_y^a (Q^3-\mathfrak g^3)}_{L^2}$ for $a=0,1,2$, respectively. Consequently, we obtain
    \[
        \norm{Q^3-\mathfrak g^3}_{W^{2,2}}
        \le
        \Delta_{\mathrm{cub}}^{W^{2,2}},
    \]
    where $\Delta_{\mathrm{cub}}^{W^{2,2}}$ is defined in \eqref{eq:cubic-w22-error-definition}. Using the equation for \(e_Q\) again,
    \[
        \norm{|D|e_Q}_{W^{2,2}}
        \le
        \Delta_Q^{W^{2,2}}
        +
        \Delta_{\mathrm{cub}}^{W^{2,2}}
        +
        \norm{r_{\mathfrak g}}_{W^{2,2}}
        =
        \Delta_{|D|Q}^{W^{2,2}}.
    \]
    Hence
    \[
        \norm{e_Q}_{W^{3,2}}
        \le
        \left((\Delta_Q^{W^{2,2}})^2+(\Delta_{|D|Q}^{W^{2,2}})^2\right)^{1/2}
        =
        \Delta_Q^{W^{3,2}},
    \]
    which yields the second bound in \eqref{eq:e_Q unweighted higer order bound}.
    Finally, applying the one-dimensional Sobolev embedding \eqref{eq:H1-Linfty}, we conclude \eqref{eq:e_Q L^infty bound}.
\end{proof}

Combining these estimates with the resolvent representation of \(Q\) and pointwise bounds for the resolvent kernel, we obtain the global weighted estimates below.

\begin{lemma}[Weighted bounds for \(Q\) and \(V_Q\)]
\label{lem:weighted-Q-V_Q-bounds}
    Fix
    \[
        R_M\geq2,
        \qquad
        R_Q\geq2R_M.
    \]  
    These are the only free parameters, since all quantities involving \(\mathfrak g\) and \(Q-\mathfrak g\) have already been fixed in Definition~\ref{def:unweighted-transfer-constants}. For this choice of \(R_Q\) and \(R_M\), the quantities \(C_{Q,a}=C_{Q,a}(R_Q,R_M)\) for \(a=0,1,2\), \(\mathsf V^{L^\infty}=\mathsf V^{L^\infty}(R_Q,R_M)\), \(\mathsf V^{L^1}=\mathsf V^{L^1}(R_Q,R_M)\) and \(\mathsf V_{\omega^2}^{L^1} = \mathsf V_{\omega^2}^{L^1}(R_Q,R_M)\) are defined by
    \eqref{eq:def-CQ-core}--\eqref{eq:def-VQ-omega-bound}. Then
    \begin{equation}\label{eq:weighted bound for Q}
        \norm{\langle y\rangle^2 Q}_{L^\infty}\le C_{Q,0},
        \qquad
        \norm{\langle y\rangle^2 yQ'}_{L^\infty}\le C_{Q,1},
        \qquad
        \norm{\langle y\rangle^2 y^2Q''}_{L^\infty}\le C_{Q,2}.
    \end{equation}
    Consequently, we have
    \begin{equation}\label{eq:V_Q- L^infyy, L^1 bound}
        \norm{V_Q}_{L^\infty}\le \mathsf V^{L^\infty},
        \qquad
        \norm{V_Q}_{L^1}\le \mathsf V^{L^1},
    \end{equation}
    and
    \begin{equation}\label{eq:V_Q sqrt bound}
        \norm{|V_Q|^{1/2}}_{L^2}\le \left(\mathsf V^{L^1}\right)^{1/2},
        \qquad
        \norm{\omega |V_Q|^{1/2}}_{L^2}\le \left(\mathsf V_{\omega^2}^{L^1}\right)^{1/2}.
    \end{equation}
    For each fixed choice of \(R_Q\) and \(R_M\), these bounds are
    obtained from the constants in Lemma~\ref{lem:higher-order-transfer-from-H1}.
\end{lemma}

\begin{proof}
    On \(|y|\le R_Q\), we write \(\partial_y^a Q = \partial_y^a\mathfrak g + \partial_y^ae_Q\) for \(a=0,1,2\). Lemma~\ref{lem:higher-order-transfer-from-H1} gives
    \begin{equation}\label{eq:Q-interior-weighted-bounds}
        \begin{aligned}
            \sup_{|y|\le R_Q}\langle y\rangle^2 |Q(y)|
            \le C_{Q,0}^{\mathrm{int}},\quad
            \sup_{|y|\le R_Q}\langle y\rangle^2 |yQ'(y)|
            \le C_{Q,1}^{\mathrm{int}},\quad
            \sup_{|y|\le R_Q}\langle y\rangle^2 |y^2Q''(y)|
            \le C_{Q,2}^{\mathrm{int}},
        \end{aligned}
    \end{equation}
    where $C_{Q,a}^{\mathrm{int}}$ are defined in \eqref{eq:def-CQ-core}.
    By the definition of \(\Delta_{\mathrm{cub}}^{L^1}\) in \eqref{eq:cubic-l1-error-definition},
    \[
        \norm{Q^3-\mathfrak g^3}_{L^1} \leq \norm{e_Q}_{L^2}(\norm{Q}_{L^\infty}\norm{Q}_{L^2}+\norm{Q}_{L^\infty}\norm{\mathfrak g}_{L^2}+\norm{\mathfrak g}_{L^2}\norm{\mathfrak g}_{L^\infty})
        \leq
        \Delta_{\mathrm{cub}}^{L^1}.
    \]
    By \(R_M\geq2\) and \(R_Q\geq2R_M\), Lemma~\ref{lem:app-exterior-weighted-Q-bounds} gives
    \begin{equation}\label{eq:Q-exterior-weighted-bounds}
        \sup_{|y|\geq R_Q}
        \langle y\rangle^2|y^a\partial_y^aQ(y)|
        \leq
        C_{Q,a}^{\mathrm{ext}},
        \qquad
        a=0,1,2.
    \end{equation}
    Combining \eqref {eq:Q-interior-weighted-bounds} and \eqref{eq:Q-exterior-weighted-bounds} gives \eqref{eq:weighted bound for Q}.
    Since \(V_Q = yQ'Q\), \eqref{eq:weighted bound for Q} gives
    \[
        |V_Q(y)|=|yQ'(y)Q(y)|
        \le C_{Q,0}C_{Q,1}\langle y\rangle^{-4}.
    \]
    Since \(\int_{\bbR}\langle y\rangle^{-4}\,dy=\pi/2\), we obtain \eqref{eq:V_Q- L^infyy, L^1 bound}. Finally,
    \[
        \norm{|V_Q|^{1/2}}_{L^2}^2
        =
        \norm{V_Q}_{L^1},
        \qquad
        \norm{\omega|V_Q|^{1/2}}_{L^2}^2
        =
        \norm{\omega^2V_Q}_{L^1},
    \]
    which gives \eqref{eq:V_Q sqrt bound}.
\end{proof}

We estimate \(\Lambda e_Q\) and \(\Lambda^2e_Q\). Combining the identities for \(\mathcal F'(f)\Lambda^a f\) with the inverse bound for \(\mathcal F'(\mathfrak g)\) gives the required \(L^2\)-estimates. The weighted pointwise bounds then give the \(L^1\) and logarithmically weighted \(L^1\) estimates.

\begin{lemma}[scaling and \(L^1\) estimates]
\label{lem:scaling-and-L1-transfer}
    Let \(e_Q\coloneqq Q-\mathfrak g\).  Fix the parameters
    \(R_Q,R_M,R_\omega,R_E\) as in
    Definitions~\ref{def:weighted-q-constants} and
    \ref{def:dilation-L1-transfer-radii}, and let
    \(\Delta_{\Lambda^aQ}^{L^2}\), \(C_{e,a}^{\mathrm{tail}}\), \(\Delta_{\omega Q}^{L^2}\),
    \(\Delta_{\Lambda^aQ}^{L^1}\), and \(\Delta_{\omega\Lambda^aQ}^{L^1}\) be defined by
    \eqref{eq:def-Delta-Q-L2}--\eqref{eq:def-Delta-omega-Lambda-a-Q-L1}.  Assume that
    \(\kappa_\Lambda>0\).  Then, for \(a=0,1,2\),
    \begin{equation}\label{eq:scaling L^2 bound}
        \norm{\Lambda^a e_Q}_{L^2}\le \Delta_{\Lambda^aQ}^{L^2}.
    \end{equation}
    Moreover, with \(\omega(y)=\log(1+|y|)\),
    \begin{equation}\label{eq:scaling L^1, wL^1 bound}
        \norm{\Lambda^a e_Q}_{L^1}\le \Delta_{\Lambda^aQ}^{L^1},
        \qquad
        \norm{\omega\Lambda^a e_Q}_{L^1}\le
        \Delta_{\omega\Lambda^aQ}^{L^1}.
    \end{equation}
\end{lemma}

\begin{proof}
    The case \(a=0\) follows from the definition
    \(\Delta_Q^{L^2}=r_{\mathrm{NK}}\) in \eqref{eq:def-Delta-Q-L2}.
    We use the scaling identities
    \[
        \begin{aligned}
            \mathcal F'(f)\Lambda f
            &=
            \Lambda\mathcal F(f)-f+\mathcal F(f),
            \\
            \mathcal F'(f)\Lambda^2 f
            &=
            \Lambda^2\mathcal F(f)
            +
            2\Lambda\mathcal F(f)
            +
            \mathcal F(f)
            -
            f
            -
            2\Lambda f
            +
            6f(\Lambda f)^2.
        \end{aligned}
    \]
    These follow from \([|D|,\Lambda]=|D|\) and the identities
    \[
        \Lambda(f^3)=3f^2\Lambda f-f^3,
        \qquad
        \Lambda^2(f^3)
        =
        3f^2\Lambda^2 f
        +
        6f(\Lambda f)^2
        -
        6f^2\Lambda f
        +
        f^3 .
    \]

    By the first identity for \(Q\) and \(\mathfrak g\), and since \(\mathcal F(Q)=0\) and \(r_{\mathfrak g}=\mathcal F(\mathfrak g)\), we obtain
    \[
        \mathcal F'(\mathfrak g)\Lambda e_Q
        =
        -e_Q
        -
        \Lambda r_{\mathfrak g}
        -
        r_{\mathfrak g}
        +
        3(Q^2-\mathfrak g^2)\Lambda Q .
    \]
    By \eqref{eq:def-Delta-Q2-Linf-kappa-Lambda},
    \[
        \norm{Q^2-\mathfrak g^2}_{L^\infty}
        \le
        \norm{e_Q}_{L^\infty}
        \left(
            \norm{Q}_{L^\infty}
            +
            \norm{\mathfrak g}_{L^\infty}
        \right)
        \le
        \Delta_{Q^2}^{L^\infty}.
    \]
    Moreover, by \eqref{eq:def-QLambda-inputs},
    \[
        \norm{\Lambda Q}_{L^2}
        \le
        \norm{\Lambda\mathfrak g}_{L^2}
        +
        \norm{\Lambda e_Q}_{L^2}.
    \]
    Therefore, the inverse bound \(\norm{[\mathcal F'(\mathfrak g)]^{-1}}_{L^2_{\mathrm e}\to L^2_{\mathrm e}} \leq K_{\mathrm{inv}}\) gives
    \[
        \norm{\Lambda e_Q}_{L^2}
        \le
        K_{\mathrm{inv}}
        \left(
            \Delta_Q^{L^2}
            +
            \norm{\Lambda r_{\mathfrak g}}_{L^2}
            +
            \norm{r_{\mathfrak g}}_{L^2}
            +
            3\Delta_{Q^2}^{L^\infty}\norm{\Lambda\mathfrak g}_{L^2}
            +
            3\Delta_{Q^2}^{L^\infty}\norm{\Lambda e_Q}_{L^2}
        \right).
    \]
    Since \(\kappa_\Lambda=1-3K_{\mathrm{inv}}\Delta_{Q^2}^{L^\infty}>0\), we obtain \eqref{eq:scaling L^2 bound} for \(a=1\); see \eqref{eq:def-Delta-Lambda-Q-L2} for the definition of \(\Delta_{\Lambda Q}^{L^2}\). By the second scaling identity for \(Q\) and \(\mathfrak g\), subtracting gives
    \[
        \begin{aligned}
        \mathcal F'(\mathfrak g)\Lambda^2 e_Q
        =
        -e_Q
        -
        2\Lambda e_Q
        -
        \Lambda^2 r_{\mathfrak g}
        -
        2\Lambda r_{\mathfrak g}
        -
        r_{\mathfrak g}+
        6\left[
            Q(\Lambda Q)^2
            -
            \mathfrak g(\Lambda\mathfrak g)^2
        \right]
        +
        3(Q^2-\mathfrak g^2)\Lambda^2 Q .
        \end{aligned}
    \]
    By \eqref{eq:def-QLambda-inputs}, \eqref{eq:def-GLambda-Linfty}, and
    \eqref{eq:def-CLambda}, the cubic difference satisfies
    \[
        \begin{aligned}
        \norm{
            Q(\Lambda Q)^2
            -
            \mathfrak g(\Lambda\mathfrak g)^2
        }_{L^2}
        &\le
        \norm{e_Q}_{L^\infty}
        \norm{\Lambda Q}_{L^\infty}
        \norm{\Lambda Q}_{L^2} 
        +
        \norm{\mathfrak g}_{L^\infty}
        \left(
            \norm{\Lambda Q}_{L^\infty}
            +
            \norm{\Lambda\mathfrak g}_{L^\infty}
        \right)
        \norm{\Lambda e_Q}_{L^2}  \\
        &\le
        \Delta_{\Lambda,\mathrm{cub}}^{L^2} .
        \end{aligned}
    \]
    Since \( \norm{\Lambda^2 Q}_{L^2}\le\norm{\Lambda^2\mathfrak g}_{L^2}+\norm{\Lambda^2 e_Q}_{L^2}\) we have
    \[
        \begin{aligned}
        \norm{\Lambda^2 e_Q}_{L^2}
        \le
        K_{\mathrm{inv}}
        \bigl(
            &\Delta_Q^{L^2}
            +
            2\Delta_{\Lambda Q}^{L^2}
            +
            \norm{\Lambda^2r_{\mathfrak g}}_{L^2}
            +
            2\norm{\Lambda r_{\mathfrak g}}_{L^2}
            +
            \norm{r_{\mathfrak g}}_{L^2}  \\
            &+
            6\Delta_{\Lambda,\mathrm{cub}}^{L^2}
            +
            3\Delta_{Q^2}^{L^\infty}\norm{\Lambda^2\mathfrak g}_{L^2}
            +
            3\Delta_{Q^2}^{L^\infty}\norm{\Lambda^2 e_Q}_{L^2}
        \bigr).
        \end{aligned}
    \]
    Using \(\kappa_\Lambda>0\) once more, we obtain \eqref{eq:scaling L^2 bound} for \(a=2\); see \eqref{eq:def-Delta-Lambda2-Q-L2} for \(\Delta_{\Lambda^2Q}^{L^2}\). It remains to estimate in \(L^1\) and \(L^1(\omega dy)\).
    
    Splitting \(\bbR\) into \(\{|y|\leq R_E\}\) and \(\{|y|>R_E\}\), the Cauchy--Schwarz inequality gives
    \begin{equation}\label{eq:interior L1 control for Lambdaa eQ}
        \norm{\Lambda^a e_Q}_{L^1(|y|\leq R_E)}
        \leq
        \sqrt{2R_E}\,
        \Delta_{\Lambda^aQ}^{L^2},
        \qquad
        a=0,1,2.
    \end{equation}
    On the exterior region, the definitions
    \eqref{eq:def-Ctail0,1,2} give
    \[
        |\Lambda^a e_Q(y)|
        \leq
        C_{e,a}^{\mathrm{tail}}\langle y\rangle^{-2},
        \qquad
        |y|>R_E,
        \qquad
        a=0,1,2.
    \]
    By Lemma~\ref{lem:app-elementary-log-tail-bounds} with \(R=R_E\), combined with \eqref{eq:app-tail-L1-bound}, \eqref{eq:app-tail-omega-L1-bound}, and \eqref{eq:def-T0-Tomega}, we obtain
    \begin{equation}\label{eq:exterior region control for F_a}
        \norm{\Lambda^a e_Q}_{L^1(|y|>R_E)}
        \leq
        \mathsf T_0(C_{e,a}^{\mathrm{tail}},R_E),
        \qquad
        \norm{\omega\Lambda^a e_Q}_{L^1(|y|>R_E)}
        \leq
        \mathsf T_\omega(C_{e,a}^{\mathrm{tail}},R_E),
    \end{equation}
    for \(a=0,1,2\). For the logarithmically weighted estimate, we first
    use \eqref{eq:def-Iomega} and \eqref{eq:def-Delta-omega-Q-L2} to obtain
    \begin{equation}\label{eq:interior region control for F_a}
        \norm{\omega e_Q}_{L^2}
        \leq
        \Delta_{\omega Q}^{L^2}.
    \end{equation}
    Hence, using \eqref{eq:interior region control for F_a} for \(a=0\),
    and \(\omega(y)\leq\log(1+R_E)\) on \(|y|\leq R_E\) for \(a=1,2\),
    we have
    \begin{equation}\label{eq:local contribution of a=0,1,2}
        \begin{aligned}
            \norm{\omega e_Q}_{L^1(|y|\leq R_E)}
            &\leq
            \sqrt{2R_E}\,\Delta_{\omega Q}^{L^2},
            \\
            \norm{\omega\Lambda^a e_Q}_{L^1(|y|\leq R_E)}
            &\leq
            \sqrt{2R_E}\log(1+R_E)\,
            \Delta_{\Lambda^aQ}^{L^2},
            \quad a=1,2.
        \end{aligned}
    \end{equation}
    Combining \eqref{eq:interior L1 control for Lambdaa eQ}, \eqref{eq:exterior region control for F_a}, and \eqref{eq:local contribution of a=0,1,2} gives \eqref{eq:scaling L^1, wL^1 bound}, with the constants as defined in \eqref{eq:def-Delta-Lambda-a-Q-L1}, \eqref{eq:def-Delta-omega-Q-L1}, and \eqref{eq:def-Delta-omega-Lambda-a-Q-L1}.
\end{proof}

\begin{lemma}[Square-root potential transfer]
\label{lem:square-root-potential-transfer}
    Recall \(V_Q=yQ'Q\) and define \(V_{\mathfrak g}\coloneqq y\mathfrak g'\mathfrak g\), and
    \[
        \delta_V\coloneqq |V_Q|^{1/2}-|V_{\mathfrak g}|^{1/2}.
    \]
    Let \(\Delta_{\sqrt V}^{L^2}\) and \(\Delta_{\omega\sqrt V}^{L^2}\) be defined by
    \eqref{eq:def-I-y2-langle4}--\eqref{eq:def-DeltaVomega}.  Then
    \begin{equation}\label{eq:difference of delta_V}
        \norm{\delta_V}_{L^2}\le \Delta_{\sqrt V}^{L^2},
        \qquad
        \norm{\omega\delta_V}_{L^2}\le \Delta_{\omega\sqrt V}^{L^2}.
    \end{equation}
\end{lemma}

\begin{proof}
    We use \(\bigl||a|^{1/2}-|b|^{1/2}\bigr|^2\le |a-b|\). Since \(V_Q-V_{\mathfrak g}
    = yQe_Q' + y\mathfrak g'e_Q\), the Cauchy--Schwarz inequality gives
    \[
        \norm{\delta_V}_{L^2}^2
        \le
        \norm{yQ}_{L^2}\norm{e_Q'}_{L^2}
        +
        \norm{y\mathfrak g'}_{L^2}\norm{e_Q}_{L^2}.
    \]
    By \eqref{eq:def-yQ-radii} and Lemma~\ref{lem:weighted-Q-V_Q-bounds},
    \[
        \norm{yQ}_{L^2}\le C_{yQ}^{L^2}.
    \]
    Also, by Proposition~\ref{prop:closeness-of-g-to-Q}, \(\norm{e_Q}_{L^2},
    \norm{e_Q'}_{L^2} \le r_{\mathrm{NK}}\).
    Hence \eqref{eq:def-DeltaV0} gives the first inequality in \eqref{eq:difference of delta_V}. For the weighted estimate, since
    \[
        \norm{\omega\delta_V}_{L^2}^2
        \le
        \norm{\omega^2 yQ}_{L^2}\norm{e_Q'}_{L^2}
        +
        \norm{\omega^2 y\mathfrak g'}_{L^2}\norm{e_Q}_{L^2},
    \]
    using \eqref{eq:def-yQ-radii}, we have
    \[
        \norm{\omega^2yQ}_{L^2}\le C_{\omega^2yQ}^{L^2}.
    \]
    Therefore \eqref{eq:def-DeltaVomega} gives the second estimate in  \eqref{eq:difference of delta_V}.
\end{proof}

These estimates are used in Section~\ref{sec:transfer from g to Q} to transfer the Birman--Schwinger bound from \(\mathfrak g\) to \(Q\).

\section{The inverse bound for the linearized operator at \(\mathfrak g\)}\label{sec:inverse-bound-g}
In this section, we prove Proposition~\ref{prop:invertibility-Fprimeg}. Subsection~\ref{subsec:NK-BS-endpoint-reduction} reduces \eqref{eq:certified-L2-inverse-Fprimeg} to the eigenvalue bounds in Lemma~\ref{lem:endpoint-BS-inequalities-NK} by the
Birman--Schwinger principle. Subsection~\ref{subsec:validated-endpoint-NK-inverse} compares
\(K_{\mathrm{NK}}^{(\mu)}\), for \(\mu\in\{\frac12,\frac32\}\), with finite-dimensional Galerkin
matrices and bounds the matrix-evaluation error. The required interval eigenvalue bounds are stated in Proposition~\ref{prop:CAP-endpoint-NK-data}.

\subsection{Birman--Schwinger reduction}\label{subsec:NK-BS-endpoint-reduction}

For \(\mu>0\), define the Birman--Schwinger operator 
\begin{equation}\label{eq:NK-BS-operator-def}
    K_{\mathrm{NK}}^{(\mu)}
    \coloneqq
    3\mathfrak g(|D|+\mu)^{-1}\mathfrak g
    \qquad\text{on }L^2_{\mathrm e}.
\end{equation}
By Proposition~\ref{prop:certified-approximate-profile}, \(\mathfrak g>0\). Since the kernel of \((|D|+\mu)^{-1}\) is strictly positive, \(K_{\mathrm{NK}}^{(\mu)}\) is compact, nonnegative,
self-adjoint, and positivity improving.

\begin{lemma}[Birman--Schwinger bounds at \(\mu=\frac12,\frac32\)]\label{lem:endpoint-BS-inequalities-NK}
    The operators \(K_{\mathrm{NK}}^{(\mu)}\) satisfy
    \begin{equation}\label{eq:NK-endpoint-eigenvalue-bounds}
        \lambda_1(K_{\mathrm{NK}}^{(1/2)})>1,
        \qquad
        \lambda_2(K_{\mathrm{NK}}^{(1/2)})<1,
        \qquad
        \lambda_1(K_{\mathrm{NK}}^{(3/2)})>1.
    \end{equation}
\end{lemma}

\begin{proof}[Proof of Proposition~\ref{prop:invertibility-Fprimeg}
assuming Lemma~\ref{lem:endpoint-BS-inequalities-NK}]
    By the same Birman--Schwinger argument as in Proposition~\ref{prop:BS-principle}, applied to \(\mathcal F'(\mathfrak g)=|D|+1-3\mathfrak g^2\), for every
    \(\mu>0\) we have
    \begin{equation}\label{eq:NK-BS-correspondence1}
        1-\mu\in\sigma_p\left( \mathcal F'(\mathfrak g)\big|_{L^2_{\mathrm e}}\right)
        \quad\text{if and only if}\quad
        1\in\sigma(K_{\mathrm{NK}}^{(\mu)}),
    \end{equation}
    and
    \begin{equation}\label{eq:NK-BS-correspondence2}
        n_-^{\mathrm e}\left(\mathcal F'(\mathfrak g)\right)=n_+\left(1;K_{\mathrm{NK}}^{(1)}\right).
    \end{equation}
    Since \(K_{\mathrm{NK}}^{(\mu)}\) decreases in \(\mu\), \eqref{eq:NK-endpoint-eigenvalue-bounds} gives, for \(\frac12\leq\mu\leq\frac32\),
    \begin{equation}\label{eq:monotonicity for K_NK}
        \lambda_1(K_{\mathrm{NK}}^{(\mu)})
        \geq
        \lambda_1( K_{\mathrm{NK}}^{(3/2)})>1,
        \qquad
        \lambda_2(K_{\mathrm{NK}}^{(\mu)})
        \leq
        \lambda_2(K_{\mathrm{NK}}^{(1/2)})
        <
        1.
    \end{equation}
    Hence \(1\notin\sigma(K_{\mathrm{NK}}^{(\mu)})\) for \(\mu\in[\frac12,\frac32]\). By \eqref{eq:NK-BS-correspondence1},
    \[
        \sigma_p
        \left(
            \mathcal F'(\mathfrak g)\big|_{L^2_{\mathrm e}}
        \right)
        \cap
        \left[
            -\frac12,
            \frac12
        \right]
        =
        \varnothing.
    \]
    Since \(\mathcal F'(\mathfrak g)\) is a compact perturbation of
    \(|D|+1\), its essential spectrum is \([1,\infty)\). Therefore,
    \[
        0
        \in
        \rho
        \left(
            \mathcal F'(\mathfrak g)\big|_{L^2_{\mathrm e}}
        \right),
        \qquad
        \norm{
            \left[
                \mathcal F'(\mathfrak g)\big|_{L^2_{\mathrm e}}
            \right]^{-1}
        }_{L^2_{\mathrm e}\to L^2_{\mathrm e}}
        \leq
        2.
    \]
    Taking \(\mu=1\) in \eqref{eq:monotonicity for K_NK} gives
    \[
        \lambda_1(K_{\mathrm{NK}}^{(1)})
        >1
        >\lambda_2(K_{\mathrm{NK}}^{(1)}).
    \]
    Thus
    \(n_+(1;K_{\mathrm{NK}}^{(1)})=1\), and \eqref{eq:NK-BS-correspondence2} yields
    \[
        n_-^{\mathrm e}
        \left(
            \mathcal F'(\mathfrak g)
        \right)
        =
        1.
    \]
\end{proof}

It remains to prove Lemma~\ref{lem:endpoint-BS-inequalities-NK}. 

\subsection{Galerkin and quadrature estimates}\label{subsec:validated-endpoint-NK-inverse}
For each endpoint \(\mu\in\{\frac12,\frac32\}\), we approximate \(K_{\mathrm{NK}}^{(\mu)}\) by an exact piecewise-constant Galerkin operator \(K_{\mathrm{NK},Y,h}^{(\mu)}\) in \eqref{eq:NK-exact-Galerkin-matrix-def}, with an explicit bound \(\delta_{\mathrm{an}}^{(\mu)}(Y,h)\) in \eqref{eq:NK-analytic-error-def} for the difference. We then approximate the matrix entries of \(\mathbf K_{\mathrm{NK},Y,h}^{(\mu)}\) in \eqref{eq:NK-evaluated-matrix-entry} by a finite quadrature evaluation, with an explicit bound \(\delta_{\mathrm{eval}}^{(\mu)}(Y,h,m,n_{\mathrm{GL}})\) in \eqref{eq:knk-evaluation-error-definition}. Lemma~\ref{lem:Weyl-perturbation-inequality} combines the two bounds into \(\varepsilon_{\mathrm{BS}}^{(\mu)}(Y,h,m,n_{\mathrm{GL}})\).

Fix \(\mu\in\{\tfrac12,\tfrac32\}\). Let \(Y>0\) and \(h>0\) with \(N\coloneqq \tfrac Yh\in\mathbb N\). The parameter values and the resulting eigenvalue and error bounds are stated in Proposition~\ref{prop:CAP-endpoint-NK-data}.

Recall \(K_{\mathrm{NK}}^{(\mu)}\) from \eqref{eq:NK-BS-operator-def}. Set
\[
    y_j\coloneqq jh,\qquad 0\leq j\leq N,
\]
and let \(\{I_j\}_{j=1}^N\) be the uniform partition of \([0,Y]\) given by
\[
    I_j=[y_{j-1},y_j)\quad (1\leq j\leq N-1),
    \qquad
    I_N=[y_{N-1},y_N].
\]
For \(1\leq j\leq N\), define the even normalized indicator
\[
    \varphi_j^{(Y,h)}(x)
    \coloneqq
    (2h)^{-1/2}\mathbf 1_{I_j}(|x|),
    \qquad x\in\mathbb R.
\]
Then \(\{\varphi_j^{(Y,h)}\}_{j=1}^N\) is an orthonormal family in
\(L^2_{\mathrm e}(\mathbb R)\). Define $\iota_{Y,h}:\mathbb R^N\to L^2_{\mathrm e}(\mathbb R)$ and $\Pi_{Y,h}$ by
\[
    \qquad
    \iota_{Y,h}\mathbf x
    \coloneqq
    \sum_{j=1}^N x_j\varphi_j^{(Y,h)}
    \quad \text{and} \quad
    \Pi_{Y,h}\coloneqq \iota_{Y,h}\iota_{Y,h}^*.
\]
Thus \(\Pi_{Y,h}\) is the orthogonal projection onto
\[
    \mathcal S_{Y,h}
    \coloneqq
    \operatorname{span}\{\varphi_j^{(Y,h)}:1\leq j\leq N\}
    \subset L^2_{\mathrm e}(\mathbb R).
\]
We define the exact Galerkin matrix by
\begin{equation}\label{eq:NK-exact-Galerkin-matrix-def}
    \mathbf K_{\mathrm{NK},Y,h}^{(\mu)}
    \coloneqq
    \iota_{Y,h}^*
    K_{\mathrm{NK}}^{(\mu)}
    \iota_{Y,h}
    \in\mathbb R^{N\times N}.
\end{equation}
Equivalently, if \(G^{(\mu)}\) denotes the integral kernel of
\((|D|+\mu)^{-1}\), then
\begin{equation}\label{eq:NK-exact-Galerkin-matrix-entry}
    \bigl(\mathbf K_{\mathrm{NK},Y,h}^{(\mu)}\bigr)_{ij}
    =
    \frac{3}{h}
    \int_{I_i}\int_{I_j}
        \mathfrak g(y)\mathfrak g(z)
        \left(
            G^{(\mu)}(y-z)+G^{(\mu)}(y+z)
        \right)
    \,dz\,dy,
    \qquad
    1\leq i,j\leq N.
\end{equation}
The corresponding finite-rank operator on \(L^2_{\mathrm e}(\mathbb R)\) is
\begin{equation*}
    K_{\mathrm{NK},Y,h}^{(\mu)}
    \coloneqq
    \iota_{Y,h}
    \mathbf K_{\mathrm{NK},Y,h}^{(\mu)}
    \iota_{Y,h}^*
    =
    \Pi_{Y,h}K_{\mathrm{NK}}^{(\mu)}\Pi_{Y,h}.
\end{equation*}
The matrix \(\mathbf K_{\mathrm{NK},Y,h}^{(\mu)}\) is symmetric, and
\[
    \norm{
        K_{\mathrm{NK},Y,h}^{(\mu)}
    }_{L^2_{\mathrm e}\to L^2_{\mathrm e}}
    =
    \norm{
        \mathbf K_{\mathrm{NK},Y,h}^{(\mu)}
    }_2 .
\]

The following lemma estimates the operator-norm difference between
\(K_{\mathrm{NK}}^{(\mu)}\) and the Galerkin operator
\(K_{\mathrm{NK},Y,h}^{(\mu)}\). The proof separates the error arising from
the spatial truncation to \([-Y,Y]\) from that introduced by the orthogonal
projection onto \(\mathcal S_{Y,h}\).

\begin{lemma}[Reduction to the exact Galerkin operator]
    For \(\mu\in\{\frac12,\frac32\}\), define
    \[
        C_{\mathrm{HS}}^{(\mu)}
        \coloneqq
        \norm{U}_{L^\infty(I)}^2
        \left(
            \frac{128}{\pi^2\mu^4}
            +
            \frac{16}{5\pi\mu}
        \right)^{1/2}
    \]
    and
    \begin{equation}\label{eq:NK-analytic-error-def}
        \delta_{\mathrm{an}}^{(\mu)}(Y,h)
        \coloneqq
        C_{\mathrm{HS}}^{(\mu)}Y^{-3}
        +
        \frac{2h}{\pi}
        \left(
            \frac{3}{\mu}
            \norm{\mathfrak g}_{L^\infty}
            \norm{\mathfrak g'}_{L^\infty}
            +
            3\norm{\mathfrak g}_{L^\infty}^2
        \right).
    \end{equation}
    Then, for $Y\geq 30$,
    \begin{equation}\label{eq:NK-analytic-Galerkin-reduction}
        \norm{
            K_{\mathrm{NK}}^{(\mu)}
            -
            K_{\mathrm{NK},Y,h}^{(\mu)}
        }_{L^2_{\mathrm e}\to L^2_{\mathrm e}}
        \leq
        \delta_{\mathrm{an}}^{(\mu)}(Y,h).
    \end{equation}
\end{lemma}

\begin{proof}
    Let \(P_Y\) denote multiplication by \(\mathbf 1_{[-Y,Y]}\) on
    \(L^2_{\mathrm e}(\mathbb R)\). Since
    \(\Pi_{Y,h}=P_Y\Pi_{Y,h}=\Pi_{Y,h}P_Y\), we have
    \[
        K_{\mathrm{NK},Y,h}^{(\mu)}
        =
        \Pi_{Y,h}K_{\mathrm{NK}}^{(\mu)}\Pi_{Y,h}
        =
        \Pi_{Y,h}P_YK_{\mathrm{NK}}^{(\mu)}P_Y\Pi_{Y,h}.
    \]
    Therefore, we estimate
    \begin{align}
        \norm{
            K_{\mathrm{NK}}^{(\mu)}
            -
            K_{\mathrm{NK},Y,h}^{(\mu)}
        }_{L^2_{\mathrm e}\to L^2_{\mathrm e}}
        &\leq
        \norm{
            K_{\mathrm{NK}}^{(\mu)}
            -
            P_YK_{\mathrm{NK}}^{(\mu)}P_Y
        }_{L^2_{\mathrm e}\to L^2_{\mathrm e}}\label{eq:NK-truncation-error-term}
        \\
        &\quad+
        \norm{
            P_YK_{\mathrm{NK}}^{(\mu)}P_Y
            -
            \Pi_{Y,h}P_YK_{\mathrm{NK}}^{(\mu)}P_Y\Pi_{Y,h}
        }_{L^2_{\mathrm e}\to L^2_{\mathrm e}}.\label{eq:NK-projection-error-term}
    \end{align}

    We first estimate the truncation term in \eqref{eq:NK-truncation-error-term}. It is enough to estimate the corresponding operator on \(L^2(\mathbb R)\). Set \(\Omega_Y\coloneqq \mathbb R^2\setminus[-Y,Y]^2\). By the Hilbert--Schmidt bound and \(\Omega_Y\subset\{|y|>Y\}\cup\{|z|>Y\}\),
    \begin{equation}\label{eq:NK-truncation-HS-bound}
        \begin{aligned}
            \norm{
                K_{\mathrm{NK}}^{(\mu)}
                -
                P_YK_{\mathrm{NK}}^{(\mu)}P_Y
            }_{L^2\to L^2}^2
            &\leq
            9
            \iint_{\Omega_Y}
                |\mathfrak g(y)|^2
                |G^{(\mu)}(y-z)|^2
                |\mathfrak g(z)|^2
            \,dz\,dy
            \\
            &\leq
            18
            \int_{|y|>Y}
                |\mathfrak g(y)|^2
                \int_{\mathbb R}
                    |G^{(\mu)}(y-z)|^2
                    |\mathfrak g(z)|^2
                \,dz\,dy .
        \end{aligned}
    \end{equation}
    For \(|y|\geq Y\geq30\), using
    \(|\mathfrak g(y)|\leq\norm{U}_{L^\infty(I)}|y|^{-2}\), the change of
    variables \(\tau=y-z\), and the split
    \(\{|\tau|\geq |y|/2\}\cup\{|\tau|<|y|/2\}\), we get
    \begin{equation}\label{eq:NK-kernel-tail-bound}
        \begin{aligned}
            \int_{\mathbb R}
                |G^{(\mu)}(y-z)|^2
                |\mathfrak g(z)|^2
            \,dz
            &\leq
            \norm{U}_{L^\infty(I)}^2\left(
            \int_{|\tau|\geq |y|/2}
                |G^{(\mu)}(\tau)|^2
            \,d\tau
            +
            \frac{16}{|y|^4}
            \int_{\mathbb R}|G^{(\mu)}(\tau)|^2\,d\tau\right)
            \\
            &\leq
            \norm{U}_{L^\infty(I)}^2
            \left(
                \frac{64}{3\pi^2\mu^4}|y|^{-3}
                +
                \frac{16}{\pi\mu}|y|^{-4}
            \right),
        \end{aligned}
    \end{equation}
    where we used the \(L^2\)-identity in \eqref{eq:Kmu-basic-bounds} and the
    tail bound obtained by integrating the pointwise estimate in
    \eqref{eq:Kmu-basic-bounds}. Substituting \eqref{eq:NK-kernel-tail-bound} into \eqref{eq:NK-truncation-HS-bound} and using \(Y\geq30\) gives
    \begin{equation}\label{eq:NK-truncation-error-term f}
        \begin{aligned}
            \norm{K_{\mathrm{NK}}^{(\mu)}-P_YK_{\mathrm{NK}}^{(\mu)}P_Y}_{L^2_{\mathrm e}\to L^2_{\mathrm e}}
            \leq
            C_{\mathrm{HS}}^{(\mu)}Y^{-3}
        \end{aligned}
    \end{equation}

    We next estimate the projection term in \eqref{eq:NK-projection-error-term}. Since
    \(P_YK_{\mathrm{NK}}^{(\mu)}P_Y\) is self-adjoint and \(\Pi_{Y,h}\) is an
    orthogonal projection, we have
    \begin{equation}\label{eq:NK-projection-error-term 1}
        \begin{aligned}
            \norm{
                P_YK_{\mathrm{NK}}^{(\mu)}P_Y
                -
                \Pi_{Y,h}P_YK_{\mathrm{NK}}^{(\mu)}P_Y\Pi_{Y,h}
            }_{L^2_{\mathrm e}\to L^2_{\mathrm e}}
            \leq
            2
            \norm{
                (I-\Pi_{Y,h})P_YK_{\mathrm{NK}}^{(\mu)}P_Y
            }_{L^2_{\mathrm e}\to L^2_{\mathrm e}}.
        \end{aligned}
    \end{equation}
    For \(v\coloneqq P_YK_{\mathrm{NK}}^{(\mu)}P_Yf\) where \(f\in L^2_{\mathrm e}(\mathbb R)\), applying Lemma~\ref{lem:wirtinger} on each \(I_j\) and its reflection, and
    using the evenness of \(v\), we obtain
    \[
    \begin{aligned}
            \norm{(I-\Pi_{Y,h})v}_{L^2(\mathbb R)}
        \leq
        \frac{h}{\pi}
        \norm{
            \partial_yK_{\mathrm{NK}}^{(\mu)}P_Yf
        }_{L^2((-Y,Y))}
        \leq
        \frac{h}{\pi}
        \norm{
            \partial_yK_{\mathrm{NK}}^{(\mu)}P_Yf
        }_{L^2(\mathbb R)}.
    \end{aligned}
    \]
    Moreover,
    \[
    \begin{aligned}
        \partial_y K_{\mathrm{NK}}^{(\mu)}P_Yf
        =
        3\mathfrak g'
        (|D|+\mu)^{-1}\bigl(\mathfrak g P_Yf\bigr)
        +
        3\mathfrak g
        \partial_y(|D|+\mu)^{-1}\bigl(\mathfrak g P_Yf\bigr).
    \end{aligned}
    \]
    Together with
    \[
        \norm{(|D|+\mu)^{-1}}_{L^2\to L^2}\leq \frac1\mu,
        \qquad
        \norm{\partial_y(|D|+\mu)^{-1}}_{L^2\to L^2}\leq 1,
    \]
    this inequalities yield
    \begin{equation}\label{eq:NK-projection-error-term 2}
        \norm{
            \partial_y K_{\mathrm{NK}}^{(\mu)}P_Y
        }_{L^2_{\mathrm e}\to L^2(\mathbb R)}
        \leq
        \frac{3}{\mu}
        \norm{\mathfrak g}_{L^\infty}
        \norm{\mathfrak g'}_{L^\infty}
        +
        3\norm{\mathfrak g}_{L^\infty}^2.
    \end{equation}
    Combining \eqref{eq:NK-truncation-error-term f},
    \eqref{eq:NK-projection-error-term 1}, and
    \eqref{eq:NK-projection-error-term 2}, we obtain
    \eqref{eq:NK-analytic-Galerkin-reduction}.
\end{proof}

Although the Galerkin matrix \(\mathbf K_{\mathrm{NK},Y,h}^{(\mu)}\) in \eqref{eq:NK-exact-Galerkin-matrix-def} is finite-dimensional, its entries in \eqref{eq:NK-exact-Galerkin-matrix-entry} are still defined by double integrals. Furthermore, the integrals involve the logarithmic singularity of \(G^{(\mu)}\).

We evaluate these entries in two regions. Near the diagonal, we separate the logarithmic part of \(G^{(\mu)}\) and integrate it exactly. Away from the diagonal, we approximate the remaining integrals by tensor-product Gauss--Legendre quadrature.

Let \(m\geq1\) be the near-diagonal parameter, and let \(n_{\mathrm{GL}}\geq2\) be the Gauss--Legendre order. Set \(n\coloneqq n_{\mathrm{GL}}\). We use the nodes and weights \(\{\chi_p,\omega_p\}_{p=1}^n\) from \eqref{eq:GL-weights}, equivalently the quadrature rule \(\mathcal Q_n\) from \eqref{eq:GL-quadrature}. By Lemma~\ref{lem:GL-analytic}, the quadrature error is explicitly bounded on the Bernstein ellipse \(E_\varrho\) from \eqref{eq:Bernstein-ellipse}. By Lemma~\ref{lem:NK-endpoint-matrix-evaluation-error}, the resulting matrix error is explicitly bounded.

For \(1\leq i\leq N\) and \(1\leq p\leq n\), set
\[
    \bar y_i\coloneqq \left(i-\frac12\right)h,
    \qquad
    y_{i,p}\coloneqq \bar y_i+\frac h2\chi_p,
    \qquad
    \omega_p^{(h)}\coloneqq \frac h2\omega_p,
    \qquad
    g_i\coloneqq \mathfrak g(\bar y_i).
\]
From Lemma~\ref{lem:Kmu-local} we decompose $G^{(\mu)}$ as 
\[
    G^{(\mu)}(d)
    =
    -\frac1\pi\log d+R_\mu(d),
    \qquad d>0.
\]
We use the logarithmic integrals \(L_d^{-}(h)\) and \(L_\ell^{+}(h)\) defined in
\eqref{eq:log-cell-diff-def} and \eqref{eq:log-cell-sum-def}. 

For \(1\leq i,j\leq N\) and \(1\leq p,q\leq n\), define
\[
    \mathcal K_{ij,pq}^{-,(\mu)}
    \coloneqq
    \begin{cases}
        R_\mu(|y_{i,p}-y_{j,q}|), & |i-j|\leq m,\\
        G^{(\mu)}(|y_{i,p}-y_{j,q}|), & |i-j|>m,
    \end{cases}
\]
and
\[
    \mathcal K_{ij,pq}^{+,(\mu)}
    \coloneqq
    \begin{cases}
        R_\mu(y_{i,p}+y_{j,q}), & i+j-1\leq m,\\
        G^{(\mu)}(y_{i,p}+y_{j,q}), & i+j-1>m.
    \end{cases}
\]
We therefore define the matrix
\[
    \mathbf K_{\mathrm{NK},Y,h,m,n_{\mathrm{GL}}}^{(\mu),\mathrm{eval}}
    \in\mathbb R^{N\times N}
\]
by
\begin{equation}\label{eq:NK-evaluated-matrix-entry}
\begin{aligned}
    \bigl(
        \mathbf K_{\mathrm{NK},Y,h,m,n_{\mathrm{GL}}}^{(\mu),\mathrm{eval}}
    \bigr)_{ij}
    &\coloneqq
    -\frac{3}{\pi h}g_ig_j
    \left[
        \mathbf 1_{\{|i-j|\leq m\}}L_{|i-j|}^{-}(h)
        +
        \mathbf 1_{\{i+j-1\leq m\}}L_{i+j-2}^{+}(h)
    \right]
    \\
    &\quad+
    \frac3h
    \sum_{p=1}^{n}\sum_{q=1}^{n}
        \omega_p^{(h)}\omega_q^{(h)}
        \mathfrak g(y_{i,p})\mathfrak g(y_{j,q})
        \left(
            \mathcal K_{ij,pq}^{-,(\mu)}
            +
            \mathcal K_{ij,pq}^{+,(\mu)}
        \right).
\end{aligned}
\end{equation}

\begin{lemma}[Endpoint matrix evaluation error]\label{lem:NK-endpoint-matrix-evaluation-error}
    Let \(\mu\in\{\frac12,\frac32\}\), \(m\ge1\),
    \(n_{\mathrm{GL}}\ge2\), \(N\in\mathbb N\), and \(h>0\) be fixed, with
    \(Y=Nh\). Set
    \[
        \varrho_m\coloneqq 2m+\sqrt{4m^2+1},
        \quad 
        a_{\mathrm{off}}^{(\mu)}
        \coloneqq
        \mu h
        \left(
            m+\frac12-\sqrt{m^2+\frac14}
        \right),
    \]
    and assume
    \[
        r_{\mathrm{near}} \coloneqq (m+1)h\le \frac{1}{\mu},
        \qquad
        \delta_{\mathrm{off}} \coloneqq hm<1.
    \]
    and define $\delta_{\mathrm{eval}}^{(\mu)}(Y,h,m,n_{\mathrm{GL}})$ by
    \begin{equation}\label{eq:knk-evaluation-error-definition}
        \begin{aligned}
            \delta_{\mathrm{eval}}^{(\mu)}(Y,h,m,n_{\mathrm{GL}})
            \coloneqq{}&
            \frac{3(2m+1)}{\pi}
            \norm{\mathfrak g}_{L^\infty}
            \norm{\mathfrak g'}_{L^\infty}
            h^2
            \left(
                1+|\log h|+\log(m+1)
            \right)
            \\
            &+
            \frac92(2m+1)L_{\mathrm{reg}}^{(\mu)}h^2
            +
            \frac{128Y}{5(\varrho_m^2-1)}
            \frac{
                \norm{\mathfrak g}_{L^\infty}
                \mathfrak G_{\delta_{\mathrm{off}}}
            }{
                \pi a_{\mathrm{off}}^{(\mu)}
            }
            \varrho_m^{-2n_{\mathrm{GL}}+2},
        \end{aligned}
    \end{equation}
    where $L_{\mathrm{reg}}^{(\mu)}$ is defined by
    \begin{equation}\label{eq:NK-local-Lipschitz-bound-def}
        L_{\mathrm{reg}}^{(\mu)}
        \coloneqq
        \norm{\mathfrak g}_{L^\infty}
        \norm{\mathfrak g'}_{L^\infty}
        \mathsf R_{\mu}^{\mathrm{loc}}(r_{\mathrm{near}})
        +
        \norm{\mathfrak g}_{L^\infty}^{2}
        \mathsf L_{\mu}^{\mathrm{loc}}(r_{\mathrm{near}}),
    \end{equation}
    and where $\mathsf R_{\mu}^{\mathrm{loc}}$ and $\mathsf L_{\mu}^{\mathrm{loc}}$ are defined in \eqref{eq:local-rmu-bound-definitions}. Also, \(\mathfrak G_{\delta_{\mathrm{off}}}\)
    is defined in \eqref{eq:g-strip-bound}. Then
    \begin{equation}\label{eq:NK-endpoint-matrix-evaluation-error}
        \norm{
            \mathbf K_{\mathrm{NK},Y,h}^{(\mu)}
            -
            \mathbf K_{\mathrm{NK},Y,h,m,n_{\mathrm{GL}}}^{(\mu),\mathrm{eval}}
        }_2
        \le
        \delta_{\mathrm{eval}}^{(\mu)}(Y,h,m,n_{\mathrm{GL}}).
    \end{equation}
\end{lemma}

\begin{proof}
    Let
    \[
        \mathbf E^{(\mu)}
        \coloneqq
        \mathbf K_{\mathrm{NK},Y,h}^{(\mu)}
        -
        \mathbf K_{\mathrm{NK},Y,h,m,n_{\mathrm{GL}}}^{(\mu),\mathrm{eval}} .
    \]
    Since both matrices are real symmetric, \(\mathbf E^{(\mu)}\) is also real symmetric. Thus, we have
    \begin{equation*}
        \norm{\mathbf E^{(\mu)}}_2
        \le
        \max_{1\le i\le N}
        \sum_{j=1}^{N}
        \left|
            \bigl(\mathbf E^{(\mu)}\bigr)_{ij}
        \right|.
    \end{equation*}
    It remains to bound the row sum uniformly in \(i\). We use the shorthand
    \[
        d_-(y,z)\coloneqq |y-z|,
        \quad
        d_+(y,z)\coloneqq y+z,
        \quad \text{and} \quad
        \nu_{ij}^-
        \coloneqq
        \mathbf 1_{\{|i-j|\le m\}},
        \quad
        \nu_{ij}^+
        \coloneqq
        \mathbf 1_{\{i+j-1\le m\}}.
    \]
    By \eqref{eq:NK-exact-Galerkin-matrix-entry},
    \eqref{eq:NK-evaluated-matrix-entry},
    \eqref{eq:log-cell-diff-id}, and \eqref{eq:log-cell-sum-id}, we may write
    \[
        \bigl(\mathbf E^{(\mu)}\bigr)_{ij}
        =
        E_{ij}^{\log}
        +
        E_{ij}^{\mathrm{reg}}
        +
        E_{ij}^{\mathrm{off}},
    \]
    where
    \begin{align*}
        E_{ij}^{\log}
        &\coloneqq
        -\frac3{\pi h}
        \int_{I_i}\int_{I_j}
        \left(
            \mathfrak g(y)\mathfrak g(z)-g_ig_j
        \right)
        \left[
            \nu_{ij}^-\log|y-z|
            +
            \nu_{ij}^+\log(y+z)
        \right]\,dz\,dy,\\
        E_{ij}^{\mathrm{reg}}
        &\coloneqq
        \frac3h
        \sum_{\sigma\in\{-,+\}}
        \nu_{ij}^{\sigma}
        \bigg[
        \int_{I_i}\int_{I_j}
            \mathfrak g(y)\mathfrak g(z)
            R_\mu(d_\sigma(y,z))\,dz\,dy
        \\
        &\qquad\qquad\qquad\qquad
        -
        \sum_{p=1}^{n_{\mathrm{GL}}}
        \sum_{q=1}^{n_{\mathrm{GL}}}
            \omega_p^{(h)}\omega_q^{(h)}
            \mathfrak g(y_{i,p})\mathfrak g(y_{j,q})
            R_\mu(d_\sigma(y_{i,p},y_{j,q}))
        \bigg],\\
        E_{ij}^{\mathrm{off}}
        &\coloneqq
        \frac3h
        \sum_{\sigma\in\{-,+\}}
        (1-\nu_{ij}^{\sigma})
        \bigg[
        \int_{I_i}\int_{I_j}
            \mathfrak g(y)\mathfrak g(z)
            G^{(\mu)}(d_\sigma(y,z))\,dz\,dy
        \\
        &\qquad\qquad\qquad\qquad\qquad
        -
        \sum_{p=1}^{n_{\mathrm{GL}}}
        \sum_{q=1}^{n_{\mathrm{GL}}}
            \omega_p^{(h)}\omega_q^{(h)}
            \mathfrak g(y_{i,p})\mathfrak g(y_{j,q})
            G^{(\mu)}(d_\sigma(y_{i,p},y_{j,q}))
        \bigg].
    \end{align*}
    It therefore suffices to bound the three row sums separately.
    
    Note that for each \(i\), we have
    \begin{equation}\label{eq:NK-near-diagonal-count}
        \sum_{j=1}^{N}
        \left(
            \nu_{ij}^-+\nu_{ij}^+
        \right)
        \le
        2m+1.
    \end{equation}

    We estimate the logarithmic near-field error. For
    \((y,z)\in I_i\times I_j\), since $|y-\bar y_i| \le \tfrac{h}{2}$ and $|z-\bar y_j| \le \tfrac{h}{2}$, the mean value theorem gives
    \begin{equation}\label{eq:NK-product-midpoint-error}
        \left|
            \mathfrak g(y)\mathfrak g(z)-g_ig_j
        \right|
        \le
        h
        \norm{\mathfrak g}_{L^\infty}
        \norm{\mathfrak g'}_{L^\infty}.
    \end{equation}
    If \(\nu_{ij}^\sigma=1\), then \(d_\sigma(y,z)\le r_{\mathrm{near}}\) on
    \(I_i\times I_j\). By Lemma~\ref{lem:exact-log-cell-integrals}, after scaling the rectangle $I_i\times I_j$ by \(h\), for $m\ge 1$ we have
    \begin{equation}\label{eq:NK-log-rectangle-absolute-bound}
        \int_{I_i}\int_{I_j}
        |\log d_\sigma(y,z)|\,dz\,dy
        \le
        h^2
        \left(
            1+|\log h|+\log(m+1)
        \right).
    \end{equation}
    Hence, using \eqref{eq:NK-product-midpoint-error},
    \eqref{eq:NK-log-rectangle-absolute-bound}, and
    \eqref{eq:NK-near-diagonal-count}, we obtain
    \begin{equation}\label{eq:NK-log-row-error}
        \sum_{j=1}^{N}|E_{ij}^{\log}|
        \le
        \frac{3(2m+1)}{\pi}
        \norm{\mathfrak g}_{L^\infty}
        \norm{\mathfrak g'}_{L^\infty}
        h^2
        \left(
            1+|\log h|+\log(m+1)
        \right).
    \end{equation}

    We next estimate the regular near-field error.  By
    Corollary~\ref{cor:local-rmu-bounds}, the assumption
    \(r_{\mathrm{near}}\le \mu^{-1}\) implies that, for \(0\le d\le
    r_{\mathrm{near}}\),
    \[
        |R_\mu(d)|
        \le
        \mathsf R_\mu^{\mathrm{loc}}(r_{\mathrm{near}}),
        \qquad
        |R_\mu'(d)|
        \le
        \mathsf L_\mu^{\mathrm{loc}}(r_{\mathrm{near}})
        \quad (d>0).
    \]
    Consequently, for every near-diagonal pair \(\nu_{ij}^\sigma=1\), the
    function
    \[
        (y,z)\mapsto
        \mathfrak g(y)\mathfrak g(z)R_\mu(d_\sigma(y,z))
    \]
    is Lipschitz on \(I_i\times I_j\), with Lipschitz constant \(L_{\mathrm{reg}}^{(\mu)}\) defined in \eqref{eq:NK-local-Lipschitz-bound-def}. More precisely, for \((y,z),(y',z')\in I_i\times I_j\),
    \begin{equation}\label{eq:NK-local-integrand-Lipschitz-bound}
        \left|
            \mathfrak g(y)\mathfrak g(z)R_\mu(d_\sigma(y,z))
            -
            \mathfrak g(y')\mathfrak g(z')R_\mu(d_\sigma(y',z'))
        \right|
        \le
        L_{\mathrm{reg}}^{(\mu)}
        \left(
            |y-y'|+|z-z'|
        \right).
    \end{equation}
    Since \(\sum_p\omega_p^{(h)}=h\), using \eqref{eq:NK-local-integrand-Lipschitz-bound} at $(y',z') = (\bar y_i,\bar y_j)$, we obtain
    \[
    \begin{aligned}
        &\left|
        \int_{I_i}\int_{I_j}
            \mathfrak g(y)\mathfrak g(z)R_\mu(d_\sigma(y,z))\,dz\,dy
        \right.\left.
        -
        \sum_{p,q=1}^{n_{\mathrm{GL}}}
            \omega_p^{(h)}\omega_q^{(h)}
            \mathfrak g(y_{i,p})\mathfrak g(y_{j,q})
            R_\mu(d_\sigma(y_{i,p},y_{j,q}))
        \right|
        \le
        \frac32 L_{\mathrm{reg}}^{(\mu)}h^3.
    \end{aligned}
    \]
    Therefore, after multiplying by \(\tfrac{3}{h}\), summing in \(j\), and using
    \eqref{eq:NK-near-diagonal-count},
    \begin{equation}\label{eq:NK-reg-row-error}
        \sum_{j=1}^{N}|E_{ij}^{\mathrm{reg}}|
        \le
        \frac92(2m+1)L_{\mathrm{reg}}^{(\mu)}h^2.
    \end{equation}

    It remains to control the Gauss--Legendre error away from the near diagonal. Fix
    \(i,j\) and \(\sigma\in\{-,+\}\) with \(\nu_{ij}^{\sigma}=0\). Write
    \[
        y_i(\xi)\coloneqq \bar y_i+\frac h2\xi,
        \qquad
        y_j(\eta)\coloneqq \bar y_j+\frac h2\eta,
        \qquad -1\le \xi,\eta\le 1.
    \]
    The condition \(\nu_{ij}^{\sigma}=0\) ensures that the relevant argument of
    \(G^{(\mu)}\) stays separated from zero on \(I_i\times I_j\). Now, define
    \[
        d_{ij}^{\sigma}(\xi,\eta)
        \coloneqq
        \begin{cases}
            y_i(\xi)-y_j(\eta), & \sigma=-,\ j<i,\\
            y_j(\eta)-y_i(\xi), & \sigma=-,\ i<j,\\
            y_i(\xi)+y_j(\eta), & \sigma=+.
        \end{cases}
    \]
    Then, \(d_{ij}^{\sigma}(\xi,\eta)=d_\sigma(y_i(\xi),y_j(\eta))\) on
    \([-1,1]^2\). Set
    \[
        F_{ij}^{\sigma}(\xi,\eta)
        \coloneqq
        \mathfrak g(y_i(\xi))\mathfrak g(y_j(\eta))
        G^{(\mu)}\bigl(d_{ij}^{\sigma}(\xi,\eta)\bigr).
    \]
    Then
    \begin{equation}\label{eq:NK-off-rectangle-pullback-integral}
        \int_{I_i}\int_{I_j}
            \mathfrak g(y)\mathfrak g(z)
            G^{(\mu)}(d_\sigma(y,z))\,dz\,dy
        =
        \frac{h^2}{4}
        \int_{-1}^{1}\int_{-1}^{1}
            F_{ij}^{\sigma}(\xi,\eta)\,d\xi\,d\eta.
    \end{equation}
    Let \(\mathcal Q_\xi\) and \(\mathcal Q_\eta\) denote the
    Gauss--Legendre rule \(\mathcal Q_{n_{\mathrm{GL}}}\) applied in
    variable $\xi$ and $\eta$, respectively. The corresponding tensor-product quadrature expression is
    \begin{equation}\label{eq:NK-off-rectangle-pullback-GL}
        \begin{aligned}
            \sum_{p=1}^{n_{\mathrm{GL}}}
            \sum_{q=1}^{n_{\mathrm{GL}}}
                \omega_p^{(h)}\omega_q^{(h)}
                \mathfrak g(y_{i,p})\mathfrak g(y_{j,q})
                G^{(\mu)}(d_\sigma(y_{i,p},y_{j,q})) =
            \frac{h^2}{4}
            \mathcal Q_\eta\mathcal Q_\xi[F_{ij}^{\sigma}].
        \end{aligned}
    \end{equation}

    From the choice of \(\varrho_m\) and \(a_{\mathrm{off}}^{(\mu)}\), for each fixed \(\eta\in[-1,1]\) the function \(\xi\mapsto F_{ij}^{\sigma}(\xi,\eta)\) is holomorphic on
    \(E_{\varrho_m}\), and Lemma~\ref{lem:K-analytic-ellipse} together with
    \eqref{eq:g-strip-bound} gives
    \[
        \sup_{\xi\in E_{\varrho_m}}
        |F_{ij}^{\sigma}(\xi,\eta)|
        \le
        \frac{
            \norm{\mathfrak g}_{L^\infty}
            \mathfrak G_{\delta_{\mathrm{off}}}
        }{
            \pi a_{\mathrm{off}}^{(\mu)}
        }.
    \]
    The same estimate holds with \(\xi\) and \(\eta\) interchanged. Applying Lemma~\ref{lem:GL-analytic} first in \(\xi\) and then in \(\eta\) gives
    \begin{equation}\label{eq:NK-GL-error-first-variable}
        \begin{aligned}
            \left|
            \int_{-1}^{1}
            \left[
                \int_{-1}^{1}
                    F_{ij}^{\sigma}(\xi,\eta)\,d\xi
                -
                \mathcal Q_\xi[F_{ij}^{\sigma}(\cdot,\eta)]
            \right]d\eta
            \right|
            \le
            \frac{128}{15(\varrho_m^2-1)}
            \frac{
                \norm{\mathfrak g}_{L^\infty}
                \mathfrak G_{\delta_{\mathrm{off}}}
            }{
                \pi a_{\mathrm{off}}^{(\mu)}
            }
            \varrho_m^{-2n_{\mathrm{GL}}+2},
        \end{aligned}
    \end{equation}
    and similarly, we have
    \begin{equation}\label{eq:NK-GL-error-second-variable}
    \begin{aligned}
        \left|
        \int_{-1}^{1}
            \mathcal Q_\xi[F_{ij}^{\sigma}(\cdot,\eta)]\,d\eta
        -
        \mathcal Q_\eta
        \left[
            \eta\mapsto
            \mathcal Q_\xi[F_{ij}^{\sigma}(\cdot,\eta)]
        \right]
        \right|
        \le
        \frac{128}{15(\varrho_m^2-1)}
        \frac{
            \norm{\mathfrak g}_{L^\infty}
            \mathfrak G_{\delta_{\mathrm{off}}}
        }{
            \pi a_{\mathrm{off}}^{(\mu)}
        }
        \varrho_m^{-2n_{\mathrm{GL}}+2}.
    \end{aligned}
    \end{equation}
    Combining \eqref{eq:NK-GL-error-first-variable} and \eqref{eq:NK-GL-error-second-variable}, we obtain
    \begin{equation}\label{eq:NK-tensor-GL-error}
        \begin{aligned}
            \left|
            \int_{-1}^{1}\int_{-1}^{1}
                F_{ij}^{\sigma}(\xi,\eta)\,d\xi\,d\eta
            -
            \mathcal Q_\eta\mathcal Q_\xi[F_{ij}^{\sigma}]
            \right|
            \le
            \frac{256}{15(\varrho_m^2-1)}
            \frac{
                \norm{\mathfrak g}_{L^\infty}
                \mathfrak G_{\delta_{\mathrm{off}}}
            }{
                \pi a_{\mathrm{off}}^{(\mu)}
            }
            \varrho_m^{-2n_{\mathrm{GL}}+2}.
        \end{aligned}
    \end{equation}
    Using \eqref{eq:NK-tensor-GL-error}, \eqref{eq:NK-off-rectangle-pullback-integral}, and \eqref{eq:NK-off-rectangle-pullback-GL}, and then summing over the off-diagonal terms, we obtain
    \begin{equation}\label{eq:NK-off-row-error}
        \begin{aligned}
            \sum_{j=1}^{N}|E_{ij}^{\mathrm{off}}|
            &\le
            \sum_{j=1}^{N}
            \sum_{\sigma\in\{-,+\}}
            (1-\nu_{ij}^{\sigma})
            \frac{64h}{5(\varrho_m^2-1)}
            \frac{
                \norm{\mathfrak g}_{L^\infty}
                \mathfrak G_{\delta_{\mathrm{off}}}
            }{
                \pi a_{\mathrm{off}}^{(\mu)}
            }
            \varrho_m^{-2n_{\mathrm{GL}}+2}
            \\
            &\le
            \frac{128Y}{5(\varrho_m^2-1)}
            \frac{
                \norm{\mathfrak g}_{L^\infty}
                \mathfrak G_{\delta_{\mathrm{off}}}
            }{
                \pi a_{\mathrm{off}}^{(\mu)}
            }
            \varrho_m^{-2n_{\mathrm{GL}}+2}.
        \end{aligned}
    \end{equation}
    Therefore, from \eqref{eq:NK-log-row-error}, \eqref{eq:NK-reg-row-error}, and
    \eqref{eq:NK-off-row-error} we conclude \eqref{eq:NK-endpoint-matrix-evaluation-error}.
\end{proof}

Combining \eqref{eq:NK-analytic-Galerkin-reduction} and \eqref{eq:NK-endpoint-matrix-evaluation-error}, we obtain the total operator error. For \(\mu\in\{\frac12,\frac32\}\), we set
\begin{equation}\label{eq:NK-total-BS-error-def}
    \varepsilon_{\mathrm{BS}}^{(\mu)}(Y,h,m,n_{\mathrm{GL}})
    \coloneqq
    \delta_{\mathrm{an}}^{(\mu)}(Y,h)
    +
    \delta_{\mathrm{eval}}^{(\mu)}(Y,h,m,n_{\mathrm{GL}}).
\end{equation}

The outward-rounded evaluation of the finite matrices and of the two error
terms is recorded in Part~\ref{part:interval-verification}.  We state the resulting bounds here.

\begin{proposition}[Certified endpoint matrix data]
\label{prop:CAP-endpoint-NK-data}
    Let
    \[
        Y=30,\qquad h=5\times10^{-3},\qquad m=1,\qquad n_{\mathrm{GL}}=8.
    \]
    The evaluated endpoint matrices satisfy
    \[
    \begin{aligned}
        \lambda_1
        \left(
            \mathbf K_{\mathrm{NK},Y,h,m,n_{\mathrm{GL}}}^{(1/2),\mathrm{eval}}
        \right)
        \ge 4.13,
        \quad
        \lambda_2
        \left(
            \mathbf K_{\mathrm{NK},Y,h,m,n_{\mathrm{GL}}}^{(1/2),\mathrm{eval}}
        \right)
        \le 0.80,
        \quad
        \lambda_1
        \left(
            \mathbf K_{\mathrm{NK},Y,h,m,n_{\mathrm{GL}}}^{(3/2),\mathrm{eval}}
        \right)
        \ge 2.42.
    \end{aligned}
    \]
    Moreover, we have
    \[
        \varepsilon_{\mathrm{BS}}^{(1/2)}(Y,h,m,n_{\mathrm{GL}})\le 0.15,
        \qquad
        \varepsilon_{\mathrm{BS}}^{(3/2)}(Y,h,m,n_{\mathrm{GL}})\le 0.08.
    \]
\end{proposition}

\begin{proof}[Proof of Lemma~\ref{lem:endpoint-BS-inequalities-NK}]
    By Proposition~\ref{prop:CAP-endpoint-NK-data} and Lemma~\ref{lem:Weyl-perturbation-inequality},
    \[
        \lambda_1(K_{\mathrm{NK}}^{(1/2)})
        \ge 4.13-0.15>1,\quad
        \lambda_2(K_{\mathrm{NK}}^{(1/2)})
        \le 0.80+0.15<1,
    \]
    and
    \[
        \lambda_1(K_{\mathrm{NK}}^{(3/2)})
        \ge 2.42-0.08>1.
    \]
\end{proof}

This proves Lemma~\ref{lem:endpoint-BS-inequalities-NK} and hence completes
the proof of Proposition~\ref{prop:invertibility-Fprimeg}.

\part{Spectral estimates for \(\mathfrak g\) and \(Q\)}\label{part:spectral-estimates}
In this part, we prove Propositions~\ref{prop:reduced-BS-bound-Q} and~\ref{prop:certified-limiting-matrix-conditions}. In Sections~\ref{sec:reduced-Birman-Schwinger-estimate}
and~\ref{sec:limiting-constraint-matrices g side}, we establish the Birman--Schwinger estimate and the limiting constraint matrix conditions for the approximate profile \(\mathfrak g\). The finite-dimensional estimates required in these sections are verified by interval arithmetic in Part~\ref{part:interval-verification}. In Section~\ref{sec:transfer from g to Q}, we use the estimates for \(Q-\mathfrak g\) from Section~\ref{sec:higher-order-profile-estimates} to complete the proofs.

We use the compactification introduced in Part~\ref{part:ground-state-approximation}:
\[
    I\coloneqq \left(-\frac{\pi}{2},\frac{\pi}{2}\right),
    \qquad
    y=\tan\theta .
\]
We conjugate the operators by the unitary map
\[
    \mathcal U:L^2_{\mathrm e}(\bbR)\to L^2_{\mathrm e}(I),
    \qquad
    (\mathcal Uv)(\theta)\coloneqq v(\tan\theta)\sec\theta .
\]
The finite-dimensional computations use the even cosine basis \(\{\cos(2k\theta)\}_{k\geq0}\).

\section{The projected Birman--Schwinger estimate for \(\mathfrak g\)}
\label{sec:reduced-Birman-Schwinger-estimate}
For the approximate profile \(\mathfrak g\), define
\begin{equation}\label{eq:definition of g-side potential}
    V_{\mathfrak g}(y)
    \coloneqq
    y\mathfrak g'(y)\mathfrak g(y),
    \qquad
    \Phi_{\mathfrak g}
    \coloneqq
    \frac{|V_{\mathfrak g}|^{1/2}}
    {\norm{|V_{\mathfrak g}|^{1/2}}_{L^2(\bbR)}}.
\end{equation}
Define \(T_{\mathfrak g}\) on \(L^2_{\mathrm e}(\bbR)\) by
\begin{equation}\label{eq:definition of T_g}
    T_{\mathfrak g}v
    \coloneqq
    -\frac1\pi
    |V_{\mathfrak g}|^{1/2}
    \left(
        \log|\cdot|*
        \bigl(|V_{\mathfrak g}|^{1/2}v\bigr)
    \right).
\end{equation}
In this section, we prove
\begin{equation}\label{eq:part4-goal-bound}
    \norm{
        P_{\Phi_{\mathfrak g}^{\perp}}
        T_{\mathfrak g}
        P_{\Phi_{\mathfrak g}^{\perp}}
    }_{L^2_{\mathrm e}(\bbR)\to L^2_{\mathrm e}(\bbR)}
    \leq
    0.156722.
\end{equation}
In Subsection~\ref{subsec:compactified-reduced-log-operator}, we derive a compactified decomposition of \(T_{\mathfrak g}\). The finite-rank part of this decomposition vanishes after applying \(P_{\widetilde\Phi_{\mathfrak g}^{\perp}}\) on both sides. In Subsection~\ref{subsec:finite-cosine-reduction-g}, we reduce \eqref{eq:part4-goal-bound} to finite-dimensional estimates. The resulting finite-dimensional bounds are stated in
Proposition~\ref{prop:part4-finite-cosine-block-bounds}.

\subsection{Compactification and finite-rank cancellation}\label{subsec:compactified-reduced-log-operator} 
Set
\begin{equation}\label{eq:definition of a_g, tidle T_g}
    \widetilde T_{\mathfrak g}
    \coloneqq
    \mathcal U T_{\mathfrak g}\mathcal U^{-1},
    \qquad
    a_{\mathfrak g}
    \coloneqq
    \mathcal U(|V_{\mathfrak g}|^{1/2}),
    \qquad
    \widetilde\Phi_{\mathfrak g}
    \coloneqq
    \mathcal U\Phi_{\mathfrak g}.
\end{equation}
We also define
\begin{equation}\label{eq:definition of a_g,log}
    a_{\mathfrak g,\log}(\theta)
    \coloneqq
    a_{\mathfrak g}(\theta)\log|\cos\theta|.
\end{equation}

\begin{lemma}[Decomposition of \(\widetilde T_{\mathfrak g}\)]
    The following properties hold:
    \begin{enumerate}[label=\textnormal{(\roman*)}]
        \item The operator \(\widetilde T_{\mathfrak g}\) decomposes as
        \begin{equation}\label{eq:part4-Tg-sing-fr-split}
            \widetilde T_{\mathfrak g}
            =
            \widetilde T_{\mathfrak g,\mathrm{sing}}
            +
            \widetilde T_{\mathfrak g,\mathrm{fr}},
        \end{equation}
        where
        \begin{align}
            (\widetilde T_{\mathfrak g,\mathrm{sing}}v)(\theta)
            &\coloneqq
            -\frac1\pi a_{\mathfrak g}(\theta)
            \int_I
                \log|\sin(\theta-\varphi)|\,
                a_{\mathfrak g}(\varphi)v(\varphi)\,d\varphi,\label{eq:part4-Tg-sing-def}
            \\
            \widetilde T_{\mathfrak g,\mathrm{fr}}
            &\coloneqq
            \frac1\pi
            \left(
                a_{\mathfrak g,\log}\otimes a_{\mathfrak g}
                +
                a_{\mathfrak g}\otimes a_{\mathfrak g,\log}
            \right).\label{eq:part4-Tg-fr-def}
        \end{align}

        \item The finite-rank part vanishes under the projection onto
        \(\widetilde\Phi_{\mathfrak g}^{\perp}\):
        \begin{equation}\label{eq:part4-Tg-fr-cancellation}
            P_{\widetilde\Phi_{\mathfrak g}^{\perp}}
            \widetilde T_{\mathfrak g,\mathrm{fr}}
            P_{\widetilde\Phi_{\mathfrak g}^{\perp}}
            =0.
        \end{equation}
        Consequently,
        \begin{equation}\label{eq:part4-Tg-reduced-sing-only}
            P_{\widetilde\Phi_{\mathfrak g}^{\perp}}
            \widetilde T_{\mathfrak g}
            P_{\widetilde\Phi_{\mathfrak g}^{\perp}}
            =
            P_{\widetilde\Phi_{\mathfrak g}^{\perp}}
            \widetilde T_{\mathfrak g,\mathrm{sing}}
            P_{\widetilde\Phi_{\mathfrak g}^{\perp}} .
        \end{equation}
    \end{enumerate}
\end{lemma}

\begin{proof}
    Let \(w=\mathcal U^{-1}v\).  Then \(w(\tan\varphi)=v(\varphi)\cos\varphi\),
    and hence we have
    \[
        \begin{aligned}
        (\widetilde T_{\mathfrak g}v)(\theta)
        =
        -\frac1\pi
        a_{\mathfrak g}(\theta)
        \int_I
            \log|\tan\theta-\tan\varphi|\,
            a_{\mathfrak g}(\varphi)v(\varphi)\,d\varphi .
        \end{aligned}
    \]
    Using
    \(\log|\tan\theta-\tan\varphi|
    =
    \log|\sin(\theta-\varphi)|
    -\log|\cos\theta|
    -\log|\cos\varphi|\), we obtain
    \eqref{eq:part4-Tg-sing-fr-split}--\eqref{eq:part4-Tg-fr-def}.
    Moreover, since \(a_{\mathfrak g}\in\operatorname{span}\{\widetilde\Phi_{\mathfrak g}\}\), for every \(v\in L^2_{\mathrm e}(I)\) we have
    \[
        (a_{\mathfrak g,\log}\otimes a_{\mathfrak g})
        P_{\widetilde\Phi_{\mathfrak g}^{\perp}}v
        =
        a_{\mathfrak g,\log}
        (a_{\mathfrak g},P_{\widetilde\Phi_{\mathfrak g}^{\perp}}v)_{L^2(I)}
        =
        0,
    \]
    and
    \[
        P_{\widetilde\Phi_{\mathfrak g}^{\perp}}
        (a_{\mathfrak g}\otimes a_{\mathfrak g,\log})
        P_{\widetilde\Phi_{\mathfrak g}^{\perp}}v
        =
        (a_{\mathfrak g,\log},P_{\widetilde\Phi_{\mathfrak g}^{\perp}}v)_{L^2(I)}
        P_{\widetilde\Phi_{\mathfrak g}^{\perp}}a_{\mathfrak g}
        =
        0.
    \]
    This proves \eqref{eq:part4-Tg-fr-cancellation}.  The identity
    \eqref{eq:part4-Tg-reduced-sing-only} follows from
    \eqref{eq:part4-Tg-sing-fr-split}.
\end{proof}

\subsection{Approximation by a cosine polynomial and block estimates}\label{subsec:finite-cosine-reduction-g}
Define \(\mathcal K_{\sin}:L^2_{\mathrm e}(I)\to L^2_{\mathrm e}(I)\) by
\begin{equation*}
    (\mathcal K_{\sin}f)(\theta)
    \coloneqq
    \int_I
        -\log|\sin(\theta-\varphi)|\,
        f(\varphi)\,d\varphi.
\end{equation*}
Then \(\mathcal K_{\sin}\) is diagonal in the even cosine basis. Indeed, the Fourier series of \(-\log|\sin\theta|\) gives
\begin{equation}\label{eq:Ksin-cosine-eigenvalues}
    \begin{aligned}
        &\mathcal K_{\sin}
        \left(
            \tfrac1{\sqrt\pi}
        \right)
        =
        \pi\log2
        \left(
            \tfrac1{\sqrt\pi}
        \right),\\
        &\mathcal K_{\sin}
        \left(
            \sqrt{\tfrac2\pi}\cos(2k\theta)
        \right)
        =
        \tfrac{\pi}{2k}
        \sqrt{\tfrac2\pi}\cos(2k\theta),
        \quad k\geq1.
    \end{aligned}
\end{equation}
By \eqref{eq:part4-Tg-sing-def}, regarding \(a_{\mathfrak g}\) as a multiplication operator, we write
\begin{equation*}
    \widetilde T_{\mathfrak g,\mathrm{sing}}
    =
    \frac1\pi
    a_{\mathfrak g}
    \mathcal K_{\sin}
    a_{\mathfrak g},
    \qquad
    \widetilde T_{\mathfrak g,\mathrm{proj}}
    \coloneqq
    P_{\widetilde\Phi_{\mathfrak g}^{\perp}}
    \widetilde T_{\mathfrak g,\mathrm{sing}}
    P_{\widetilde\Phi_{\mathfrak g}^{\perp}}.
\end{equation*}
To reduce the proof of \eqref{eq:part4-goal-bound} to finite-dimensional estimates, we introduce the even cosine spaces below. For \(J\geq0\), set
\[
    \mathcal E_J
    \coloneqq
    \left(
        \tfrac1{\sqrt\pi},
        \sqrt{\tfrac2\pi}\cos(2\theta),
        \ldots,
        \sqrt{\tfrac2\pi}\cos(2J\theta)
    \right),
    \qquad
    \mathcal S_J
    \coloneqq
    \operatorname{span}\mathcal E_J
    \subset L^2_{\mathrm e}(I).
\]
Then \(\mathcal E_J\) is an orthonormal basis of \(\mathcal S_J\). 

The diagonalization of \(\mathcal K_{\sin}\) does not yield a finite-dimensional estimate for
\(a_{\mathfrak g}\mathcal K_{\sin}a_{\mathfrak g}\). Indeed, for any $J, L \in \bbN$,
\[
    a_\mathfrak g \mathcal S_J \nsubseteq \mathcal S_{J+L}.
\]
We therefore approximate \(a_{\mathfrak g}\) by a cosine polynomial. Fix \(L\geq1\) and
\[
    (\mathsf c_0,\ldots,\mathsf c_L) \in \mathbb Q^{L+1}\setminus\{0\}.
\]
Set
\begin{equation}\label{eq:definition of a_g^L}
    a_{\mathfrak g}^{(L)}(\theta)
    \coloneqq
    \sum_{k=0}^{L}
        \mathsf c_k\cos(2k\theta).
\end{equation}
Then the product-to-sum identity gives
\begin{equation}\label{eq:part4-agL-subspace-inclusion}
    a_{\mathfrak g}^{(L)}\mathcal S_J
    \subset
    \mathcal S_{J+L},
    \qquad J\geq0.
\end{equation}

We set
\begin{equation}\label{eq:part4-delta-aL-L2-def}
    \delta_{a_{\mathfrak g}\to a_{\mathfrak g}^{(L)}}^{L^2}
    \coloneqq
    \norm{
        a_{\mathfrak g}-a_{\mathfrak g}^{(L)}
    }_{L^2(I)}.
\end{equation}
The value of \(L\) and the coefficients \(\{\mathsf c_k\}_{k=0}^{L}\) are specified in
Section~\ref{sec:interval-inputs}. The approximation error is bounded in Proposition~\ref{prop:ag-cosine-approximation-bounds}. Set
\begin{equation}\label{eq:part4-Tg-L-direction-projection-def}
    \widetilde\Phi_{\mathfrak g}^{(L)}
    \coloneqq
    \tfrac{
        a_{\mathfrak g}^{(L)}
    }{
        \norm{a_{\mathfrak g}^{(L)}}_{L^2(I)}
    }.
\end{equation}
We use the same notation for multiplication by \(a_{\mathfrak g}^{(L)}\). Define
\begin{equation}\label{eq:part4-Tg-L-sing-red-def}
    \widetilde T_{\mathfrak g,\mathrm{sing}}^{(L)}
    \coloneqq
    \frac1\pi
    a_{\mathfrak g}^{(L)}
    \mathcal K_{\sin}
    a_{\mathfrak g}^{(L)},
    \qquad
    \widetilde T_{\mathfrak g,\mathrm{proj}}^{(L)}
    \coloneqq
    P_{(\widetilde\Phi_{\mathfrak g}^{(L)})^\perp}
    \widetilde T_{\mathfrak g,\mathrm{sing}}^{(L)}
    P_{(\widetilde\Phi_{\mathfrak g}^{(L)})^\perp}.
\end{equation}

We first estimate the difference between \(\widetilde T_{\mathfrak g,\mathrm{proj}}\) and \(\widetilde T_{\mathfrak g,\mathrm{proj}}^{(L)}\).

\begin{lemma}[Operator error estimate for the cosine-polynomial approximation]
    Let
    \[
        \mathsf C_{\log\sin}^{L^2}
        \coloneqq
        \norm{-\log|\sin(\cdot)|}_{L^2(I)},
    \]
    and set
    \begin{equation}\label{eq:part4-eps-model-BS-def}
        \delta_{a_{\mathfrak g}\to a_{\mathfrak g}^{(L)}}^{\mathrm{BS}}
        \coloneqq
        \frac{\mathsf C_{\log\sin}^{L^2}}{\pi}
        \delta_{a_{\mathfrak g}\to a_{\mathfrak g}^{(L)}}^{L^2}
        \left(
            5\norm{a_{\mathfrak g}}_{L^\infty(I)}
            +
            \norm{a_{\mathfrak g}^{(L)}}_{L^\infty(I)}
        \right).
    \end{equation}
    Then
    \begin{equation}\label{eq:BS-model-replacement-bound}
        \norm{
            \widetilde T_{\mathfrak g,\mathrm{proj}}
            -
            \widetilde T_{\mathfrak g,\mathrm{proj}}^{(L)}
        }_{L^2_{\mathrm e}(I)\to L^2_{\mathrm e}(I)}
        \le
        \delta_{a_{\mathfrak g}\to a_{\mathfrak g}^{(L)}}^{\mathrm{BS}}.
    \end{equation}
\end{lemma}

\begin{proof}
    The triangle inequality gives
    \begin{equation}\label{eq:part4-model-replacement-triangle}
        \begin{aligned}
            \norm{
                \widetilde T_{\mathfrak g,\mathrm{proj}}
                -
                \widetilde T_{\mathfrak g,\mathrm{proj}}^{(L)}
            }_{L^2_{\mathrm e}(I)\to L^2_{\mathrm e}(I)}
            &\leq
            2
            \norm{
                P_{\widetilde\Phi_{\mathfrak g}^{\perp}}
                -
                P_{(\widetilde\Phi_{\mathfrak g}^{(L)})^\perp}
            }_{L^2(I)\to L^2(I)}
            \norm{
                \widetilde T_{\mathfrak g,\mathrm{sing}}
            }_{L^2_{\mathrm e}(I)\to L^2_{\mathrm e}(I)}
            \\
            &\quad+
            \norm{
                \widetilde T_{\mathfrak g,\mathrm{sing}}
                -
                \widetilde T_{\mathfrak g,\mathrm{sing}}^{(L)}
            }_{L^2_{\mathrm e}(I)\to L^2_{\mathrm e}(I)}.
        \end{aligned}
    \end{equation}
    Cauchy--Schwarz and the \(\pi\)-periodicity of \(-\log|\sin\eta|\) give
    \begin{equation}\label{eq:part4-Ksin-L2-Linf-bound}
        \norm{\mathcal K_{\sin}}_{L^2(I)\to L^\infty(I)}
        \leq
        \mathsf C_{\log\sin}^{L^2}.
    \end{equation}
    Expanding \(\widetilde T_{\mathfrak g,\mathrm{sing}} - \widetilde T_{\mathfrak g,\mathrm{sing}}^{(L)}\) gives
    \begin{align}
        &\norm{
            \widetilde T_{\mathfrak g,\mathrm{sing}}
            -
            \widetilde T_{\mathfrak g,\mathrm{sing}}^{(L)}
        }_{L^2_{\mathrm e}(I)\to L^2_{\mathrm e}(I)}
        \nonumber\\
        &\qquad\leq
        \frac{1}{\pi}
        \norm{
            (a_{\mathfrak g}-a_{\mathfrak g}^{(L)})
            \mathcal K_{\sin}a_{\mathfrak g}
        }_{L^2(I)\to L^2(I)}
        +
        \frac{1}{\pi}
        \norm{
            a_{\mathfrak g}^{(L)}
            \mathcal K_{\sin}
            (a_{\mathfrak g}-a_{\mathfrak g}^{(L)})
        }_{L^2(I)\to L^2(I)}
        \nonumber\\
        &\qquad\leq
        \frac{\mathsf C_{\log\sin}^{L^2}}{\pi}
        \delta_{a_{\mathfrak g}\to a_{\mathfrak g}^{(L)}}^{L^2}
        \left(
            \norm{a_{\mathfrak g}}_{L^\infty(I)}
            +
            \norm{a_{\mathfrak g}^{(L)}}_{L^\infty(I)}
        \right).
        \label{eq:part4-singular-model-difference-bound}
    \end{align}
    The second summand is estimated after taking adjoints and using the self-adjointness of \(\mathcal K_{\sin}\) and of the multiplication operators.
    Next, by Lemma~\ref{lem:app-rank-one-projection-difference} and
    \eqref{eq:part4-Tg-L-direction-projection-def}, we have
    \begin{equation}\label{eq:part4-projection-model-difference-bound}
        \norm{
            P_{\widetilde\Phi_{\mathfrak g}^{\perp}}
            -
            P_{(\widetilde\Phi_{\mathfrak g}^{(L)})^\perp}
        }_{L^2(I)\to L^2(I)}
        \leq
        \norm{
            \widetilde\Phi_{\mathfrak g}
            -
            \widetilde\Phi_{\mathfrak g}^{(L)}
        }_{L^2(I)}
        \leq
        \frac{
            2\delta_{a_{\mathfrak g}\to a_{\mathfrak g}^{(L)}}^{L^2}
        }{
            \max\left\{
                \norm{a_{\mathfrak g}}_{L^2(I)},
                \norm{a_{\mathfrak g}^{(L)}}_{L^2(I)}
            \right\}
        } .
    \end{equation}
    Moreover, \eqref{eq:part4-Ksin-L2-Linf-bound} gives
    \begin{equation}\label{eq:part4-Tg-sing-operator-bound}
        \norm{
            \widetilde T_{\mathfrak g,\mathrm{sing}}
        }_{L^2_{\mathrm e}(I)\to L^2_{\mathrm e}(I)}
        \leq
        \frac{\mathsf C_{\log\sin}^{L^2}}{\pi}
        \norm{a_{\mathfrak g}}_{L^2(I)}
        \norm{a_{\mathfrak g}}_{L^\infty(I)} .
    \end{equation}
    Substituting \eqref{eq:part4-singular-model-difference-bound},
    \eqref{eq:part4-projection-model-difference-bound}, and
    \eqref{eq:part4-Tg-sing-operator-bound} into
    \eqref{eq:part4-model-replacement-triangle}, and using
    \[
            \norm{a_{\mathfrak g}}_{L^2(I)}
            \leq 
            \max\left\{
                \norm{a_{\mathfrak g}}_{L^2(I)},
                \norm{a_{\mathfrak g}^{(L)}}_{L^2(I)}
            \right\},
    \]
    we obtain
    \[
        \norm{
            \widetilde T_{\mathfrak g,\mathrm{proj}}
            -
            \widetilde T_{\mathfrak g,\mathrm{proj}}^{(L)}
        }_{L^2_{\mathrm e}(I)\to L^2_{\mathrm e}(I)}
        \leq
        \frac{\mathsf C_{\log\sin}^{L^2}}{\pi}
        \delta_{a_{\mathfrak g}\to a_{\mathfrak g}^{(L)}}^{L^2}
        \left(
            5\norm{a_{\mathfrak g}}_{L^\infty(I)}
            +
            \norm{a_{\mathfrak g}^{(L)}}_{L^\infty(I)}
        \right).
    \]
    This proves \eqref{eq:BS-model-replacement-bound}.
\end{proof}

It remains to bound \(\widetilde T_{\mathfrak g,\mathrm{proj}}^{(L)}\).

We write \(\Pi_{\leq J}\) for the orthogonal projection onto \(\mathcal S_J\), and set
\[
    \Pi_{L<\cdot\leq3L}
    \coloneqq
    \Pi_{\leq3L}-\Pi_{\leq L},
    \qquad
    \Pi_{>L}
    \coloneqq
    I-\Pi_{\leq L},
    \qquad
    \Pi_{>3L}
    \coloneqq
    I-\Pi_{\leq3L}.
\]

If \(A:L^2_{\mathrm e}(I)\to L^2_{\mathrm e}(I)\) is bounded and \(A\mathcal S_J\subset\mathcal S_J\), we denote by \([A]_{\mathcal E_J}\) the matrix of \(\left.A\right|_{\mathcal S_J}\) in the basis \(\mathcal E_J\). For a bounded operator \(A:L^2_{\mathrm e}(I)\to L^2_{\mathrm e}(I)\), set
\[
    A^{[J]}
    \coloneqq
    \Pi_{\leq J}A\Pi_{\leq J}.
\]
Then \(A^{[J]}\) has range contained in \(\mathcal S_J\), and \([A^{[J]}]_{\mathcal E_J}\) denotes the matrix of \(\left.A^{[J]}\right|_{\mathcal S_J}\).

Since \(\widetilde\Phi_{\mathfrak g}^{(L)}\in\mathcal S_L\), the projection \(P_{(\widetilde\Phi_{\mathfrak g}^{(L)})^\perp}\) preserves \(\mathcal S_J\) for \(J\geq L\). Together with
\eqref{eq:part4-agL-subspace-inclusion}, this gives
\[
    \widetilde T_{\mathfrak g,\mathrm{proj}}^{(L)}
    \mathcal S_L
    \subset
    \mathcal S_{3L}.
\]

Define
\begin{equation}\label{eq:part4-Tg-3L-def}
    \widetilde T_{\mathfrak g,3L}^{(L)}
    \coloneqq
    \left(
        P_{(\widetilde\Phi_{\mathfrak g}^{(L)})^\perp}
    \right)^{[3L]}
    \left[
        \frac1\pi
        \left(a_{\mathfrak g}^{(L)}\right)^{[3L]}
        \mathcal K_{\sin}^{[3L]}
        \left(a_{\mathfrak g}^{(L)}\right)^{[3L]}
    \right]
    \left(
        P_{(\widetilde\Phi_{\mathfrak g}^{(L)})^\perp}
    \right)^{[3L]}.
\end{equation}
This operator is self-adjoint and nonnegative, with range contained
in \(\mathcal S_{3L}\). By \eqref{eq:part4-agL-subspace-inclusion} and \eqref{eq:part4-Tg-L-direction-projection-def},
\begin{equation}\label{eq:part4-TgL-low-input-identity}
    \widetilde T_{\mathfrak g,\mathrm{proj}}^{(L)}
    \Pi_{\leq L}
    =
    \widetilde T_{\mathfrak g,3L}^{(L)}
    \Pi_{\leq L}.
\end{equation}
Moreover,
\begin{equation}\label{eq:part4-high-projection-identity}
    P_{(\widetilde\Phi_{\mathfrak g}^{(L)})^\perp}
    \Pi_{>L}
    =
    \Pi_{>L}.
\end{equation}
We now use the orthogonal decomposition
\[
    L^2_{\mathrm e}(I)
    =
    \mathcal S_L\oplus\mathcal S_L^\perp.
\]
The terms containing \(\Pi_{\leq L}f\) are reduced to finite-dimensional matrix estimates on \(\mathcal S_{3L}\). The remaining term is estimated using \eqref{eq:Ksin-cosine-eigenvalues}.

\begin{lemma}[Block decomposition for $\widetilde T_{\mathfrak g,\mathrm{proj}}^{(L)}$]
    We have
    \begin{equation}\label{eq:abstract-block-bound}
        \norm{
            \widetilde T_{\mathfrak g,\mathrm{proj}}^{(L)}
        }_{L^2_{\mathrm e}(I)\to L^2_{\mathrm e}(I)}
        \leq
        \frac{
            \beta_{\mathrm{low}}^{(L)}
            +
            \beta_{\mathrm{high}}^{(L)}
        }{2}
        +
        \sqrt{
            \left(
                \frac{
                    \beta_{\mathrm{low}}^{(L)}
                    -
                    \beta_{\mathrm{high}}^{(L)}
                }{2}
            \right)^2
            +
            \bigl(\beta_{\mathrm{off}}^{(L)}\bigr)^2
        }, 
    \end{equation}
    where $\beta_{\mathrm{low}}^{(L)},\ \beta_{\mathrm{off}}^{(L)}$ and $\beta_{\mathrm{high}}^{(L)}$ are defined by
    \begin{equation}\label{eq:definition of beta's}
        \begin{aligned}
            \beta_{\mathrm{low}}^{(L)}
            &\coloneqq
            \lambda_{\max}
            \left(
                \left[
                    \Pi_{\leq L}
                    \widetilde T_{\mathfrak g,3L}^{(L)}
                    \Pi_{\leq L}
                \right]_{\mathcal E_L}
            \right),\\
            \beta_{\mathrm{off}}^{(L)}
            &\coloneqq
            \left\{
            \lambda_{\max}
            \left(
                \left[
                    \Pi_{\leq L}
                    \widetilde T_{\mathfrak g,3L}^{(L)}
                    \Pi_{L<\cdot\leq 3L}
                    \widetilde T_{\mathfrak g,3L}^{(L)}
                    \Pi_{\leq L}
                \right]_{\mathcal E_L}
            \right)
            \right\}^{1/2},\\
            \beta_{\mathrm{high}}^{(L)}
            &\coloneqq
            (\log 2)\,
            \norm{
                \Pi_{>L}
                a_{\mathfrak g}^{(L)}
                \Pi_{\leq L}
            }_{L^2_{\mathrm e}(I)\to L^2_{\mathrm e}(I)}^2
            +
            \frac{1}{2(L+1)}
            \norm{a_{\mathfrak g}^{(L)}}_{L^\infty(I)}^2.
        \end{aligned}
    \end{equation}
\end{lemma}

\begin{proof}
    By \eqref{eq:Ksin-cosine-eigenvalues}, \(\mathcal K_{\sin}\geq0\). Hence \(\widetilde T_{\mathfrak g,\mathrm{proj}}^{(L)}\) is nonnegative, and
    \[
        \norm{
            \widetilde T_{\mathfrak g,\mathrm{proj}}^{(L)}
        }
        =
        \sup_{\substack{
            f\in L^2_{\mathrm e}(I)\\
            \norm{f}_{L^2(I)}=1
        }}
        \bigl(
            \widetilde T_{\mathfrak g,\mathrm{proj}}^{(L)}f,
            f
        \bigr)_{L^2(I)}.
    \]
    It is enough to prove
    \begin{align}
        \bigl(
            \widetilde T_{\mathfrak g,\mathrm{proj}}^{(L)}
            \Pi_{\leq L}f,
            \Pi_{\leq L}f
        \bigr)_{L^2(I)}
        &\leq
        \beta_{\mathrm{low}}^{(L)}
        \norm{\Pi_{\leq L}f}_{L^2(I)}^2,
        \label{eq:low contribution estimate}\\
        \left|
        \bigl(
            \widetilde T_{\mathfrak g,\mathrm{proj}}^{(L)}
            \Pi_{>L}f,
            \Pi_{\leq L}f
        \bigr)_{L^2(I)}
        \right|
        &\leq
        \beta_{\mathrm{off}}^{(L)}
        \norm{\Pi_{\leq L}f}_{L^2(I)}
        \norm{\Pi_{>L}f}_{L^2(I)},
        \label{eq:mixed contribution estimate}\\
        \bigl(
            \widetilde T_{\mathfrak g,\mathrm{proj}}^{(L)}
            \Pi_{>L}f,
            \Pi_{>L}f
        \bigr)_{L^2(I)}
        &\leq
        \beta_{\mathrm{high}}^{(L)}
        \norm{\Pi_{>L}f}_{L^2(I)}^2.
        \label{eq:high contribution estimate}
    \end{align}
    Assuming these bounds, for \(\norm{f}_{L^2(I)}=1\), we have
    \begin{equation}\label{eq:block decomposition}
        \begin{aligned}
            \bigl(
                \widetilde T_{\mathfrak g,\mathrm{proj}}^{(L)}f,
                f
            \bigr)_{L^2(I)}
            &\leq
            \beta_{\mathrm{low}}^{(L)}
            \norm{\Pi_{\leq L}f}_{L^2(I)}^2
            +
            2\beta_{\mathrm{off}}^{(L)}
            \norm{\Pi_{\leq L}f}_{L^2(I)}
            \norm{\Pi_{>L}f}_{L^2(I)}
            \\
            &\quad+
            \beta_{\mathrm{high}}^{(L)}
            \norm{\Pi_{>L}f}_{L^2(I)}^2.
        \end{aligned}
    \end{equation}
    Since \(\norm{\Pi_{\leq L}f}_{L^2(I)}^2+\norm{\Pi_{>L}f}_{L^2(I)}^2=1\), the right-hand side of \eqref{eq:block decomposition} is at most
    \[
        \lambda_{\max}
        \begin{pmatrix}
            \beta_{\mathrm{low}}^{(L)}
            &
            \beta_{\mathrm{off}}^{(L)}
            \\
            \beta_{\mathrm{off}}^{(L)}
            &
            \beta_{\mathrm{high}}^{(L)}
        \end{pmatrix}.
    \]
    The formula for the largest eigenvalue of a real symmetric \(2\times2\) matrix gives \eqref{eq:abstract-block-bound}.

    We now prove \eqref{eq:low contribution estimate} and \eqref{eq:mixed contribution estimate}. From \eqref{eq:part4-TgL-low-input-identity},
    \[
        \Pi_{\leq L}
        \widetilde T_{\mathfrak g,\mathrm{proj}}^{(L)}
        \Pi_{\leq L}
        =
        \Pi_{\leq L}
        \widetilde T_{\mathfrak g,3L}^{(L)}
        \Pi_{\leq L}.
    \]
    The definition of \(\beta_{\mathrm{low}}^{(L)}\) gives \eqref{eq:low contribution estimate}. Since both operators in \eqref{eq:part4-TgL-low-input-identity} are self-adjoint and \(\widetilde T_{\mathfrak g,3L}^{(L)}\) has range in \(\mathcal S_{3L}\),
    \begin{equation}\label{eq:part4-off-block-identity}
        \Pi_{\leq L}
        \widetilde T_{\mathfrak g,\mathrm{proj}}^{(L)}
        \Pi_{>L}
        =
        \Pi_{\leq L}
        \widetilde T_{\mathfrak g,3L}^{(L)}
        \Pi_{L<\cdot\leq3L}.
    \end{equation}
    Moreover,
    \[
        \begin{aligned}
        \norm{
            \Pi_{\leq L}
            \widetilde T_{\mathfrak g,3L}^{(L)}
            \Pi_{L<\cdot\leq3L}
        }_{L^2_{\mathrm{e}}(I)\to L^2_{\mathrm{e}}(I)}^2
        &=
        \lambda_{\max}
        \left(
            \left[
                \Pi_{\leq L}
                \widetilde T_{\mathfrak g,3L}^{(L)}
                \Pi_{L<\cdot\leq3L}
                \widetilde T_{\mathfrak g,3L}^{(L)}
                \Pi_{\leq L}
            \right]_{\mathcal E_L}
        \right)\\
        &=
        \bigl(\beta_{\mathrm{off}}^{(L)}\bigr)^2.
        \end{aligned}
    \]
    Thus \eqref{eq:part4-off-block-identity} and Cauchy--Schwarz give
    \eqref{eq:mixed contribution estimate}.
    It remains to prove \eqref{eq:high contribution estimate}. By \eqref{eq:part4-high-projection-identity},
    \begin{equation}\label{eq:high-quadratic-form-reduction}
        \begin{aligned}
            \bigl(
                \widetilde T_{\mathfrak g,\mathrm{proj}}^{(L)}
                \Pi_{>L}f,
                \Pi_{>L}f
            \bigr)_{L^2(I)}
            =
            \frac1\pi
            \bigl(
                \mathcal K_{\sin}
                a_{\mathfrak g}^{(L)}
                \Pi_{>L}f,
                a_{\mathfrak g}^{(L)}
                \Pi_{>L}f
            \bigr)_{L^2(I)}.
        \end{aligned}
    \end{equation}
    By \eqref{eq:Ksin-cosine-eigenvalues},
    \[
        \Pi_{\leq L}\mathcal K_{\sin}\Pi_{>L}=0,
    \]
    and
    \begin{equation*}
        \norm{
            \Pi_{\leq L}
            \mathcal K_{\sin}
            \Pi_{\leq L}
        }_{L^2_{\mathrm{e}}(I)\to L^2_{\mathrm{e}}(I)}
        =
        \pi\log2,
        \quad
        \norm{
            \Pi_{>L}
            \mathcal K_{\sin}
            \Pi_{>L}
        }_{L^2_{\mathrm{e}}(I)\to L^2_{\mathrm{e}}(I)}
        =
        \frac{\pi}{2(L+1)}.
    \end{equation*}
    Applying
    \(I=\Pi_{\leq L}+\Pi_{>L}\)
    in \eqref{eq:high-quadratic-form-reduction}, we obtain
    \[
        \begin{aligned}
            \bigl(
                \widetilde T_{\mathfrak g,\mathrm{proj}}^{(L)}
                &\Pi_{>L}f,
                \Pi_{>L}f
            \bigr)_{L^2(I)}
            \\
            &\quad\leq
            (\log2)
            \norm{
                \Pi_{\leq L}
                a_{\mathfrak g}^{(L)}
                \Pi_{>L}f
            }_{L^2(I)}^2
            +
            \frac1{2(L+1)}
            \norm{
                \Pi_{>L}
                a_{\mathfrak g}^{(L)}
                \Pi_{>L}f
            }_{L^2(I)}^2.
        \end{aligned}
    \]
    Since 
    \[
        \begin{aligned}
            \norm{
                \Pi_{\leq L}
                a_{\mathfrak g}^{(L)}
                \Pi_{>L}
            }_{L^2_{\mathrm{e}}(I)\to L^2_{\mathrm{e}}(I)}
            &=
            \norm{
                \Pi_{>L}
                a_{\mathfrak g}^{(L)}
                \Pi_{\leq L}
            }_{L^2_{\mathrm{e}}(I)\to L^2_{\mathrm{e}}(I)},\\
            \norm{
                \Pi_{>L}
                a_{\mathfrak g}^{(L)}
                \Pi_{>L}
            }_{L^2_{\mathrm{e}}(I)\to L^2_{\mathrm{e}}(I)}
            &\leq
            \norm{a_{\mathfrak g}^{(L)}}_{L^\infty(I)},
        \end{aligned}
    \]
    we conclude
    \[
        \begin{aligned}
            \bigl(
                \widetilde T_{\mathfrak g,\mathrm{proj}}^{(L)}&
                \Pi_{>L}f,
                \Pi_{>L}f
            \bigr)_{L^2(I)}
            \\
            &\quad\leq
            \left[
                (\log2)
                \norm{
                    \Pi_{>L}
                    a_{\mathfrak g}^{(L)}
                    \Pi_{\leq L}
                }^2
                +
                \frac1{2(L+1)}
                \norm{
                    a_{\mathfrak g}^{(L)}
                }_{L^\infty(I)}^2
            \right]
            \norm{\Pi_{>L}f}_{L^2(I)}^2
            \\
            &\quad=
            \beta_{\mathrm{high}}^{(L)}
            \norm{\Pi_{>L}f}_{L^2(I)}^2.
        \end{aligned}
    \]
    This proves \eqref{eq:high contribution estimate}.
\end{proof}

Part~\ref{part:interval-verification} proves the following proposition and corollary.

\begin{proposition}
\label{prop:part4-finite-cosine-block-bounds}
    For the coefficients specified in Part~\ref{part:interval-verification}, let \(L=5000\).
    Then the quantities in \eqref{eq:definition of beta's} satisfy
    \begin{equation}\label{eq:part4-certified-beta-values}
        \beta_{\mathrm{low}}^{(L)}
        \leq 0.156664,
        \qquad
        \beta_{\mathrm{off}}^{(L)}
        \leq 3.204\times10^{-5},
        \qquad
        \beta_{\mathrm{high}}^{(L)}
        \leq 0.036572.
    \end{equation}
    Moreover,
    \begin{equation}\label{eq:part4-certified-model-error}
        \delta_{a_{\mathfrak g}\to a_{\mathfrak g}^{(L)}}^{\mathrm{BS}}
        \leq
        5.737\times10^{-5}.
    \end{equation}
\end{proposition}

\begin{corollary}[Projected Birman--Schwinger bound for \(\mathfrak g\)]
\label{cor:part4-reduced-BS-bound-g}
    We have
    \begin{equation}\label{eq:part4-reduced-BS-bound-g-real-line}
        \norm{
            P_{\Phi_{\mathfrak g}^{\perp}}
            T_{\mathfrak g}
            P_{\Phi_{\mathfrak g}^{\perp}}
        }_{L^2_{\mathrm e}\to L^2_{\mathrm e}}
        \leq
        0.156722.
    \end{equation}
\end{corollary}

\section{The limiting constraint matrices}\label{sec:limiting-constraint-matrices g side}
In this section, we define and estimate the quantities associated with \(\mathfrak g\) that appear in the limiting constraint matrices.

In Section~\ref{subsec:g-side-Schur-complement}, we use the projected Birman--Schwinger bound
\eqref{eq:part4-goal-bound} to prove invertibility of the block on \(\operatorname{span}\{\widetilde\Phi_{\mathfrak g}\}^{\perp}\). A scalar Schur complement then gives inverse bounds for \(S_0^{(\mathfrak g,j)}\). 

In Section~\ref{subsec:g-side-limiting-matrix-entries}, the pairings that do not contain
\((S_0^{(\mathfrak g,j)})^{-1}\) are evaluated from finite cosine-coefficient formulas.
For the remaining pairings, we approximate the solutions of the corresponding equations on
\(I\) by cosine polynomials and bound the errors by their residuals. The inverse bounds and interval enclosures are stated in Propositions~\ref{prop:certified-Schur-inverse-bounds-g}
and~\ref{prop:part4-certified-g-side-limiting-data}.

\subsection{Schur complement reduction for \(S_0^{(\mathfrak g,j)}\)}\label{subsec:g-side-Schur-complement}
We use the notation from \eqref{eq:definition of g-side potential}--\eqref{eq:definition of T_g}. For \(j=1,2\), define
\begin{equation}\label{eq:definition of S0-g-j}
    S_0^{(\mathfrak g,j)}
    \coloneqq
    I-c_jT_{\mathfrak g}
    \qquad
    \text{on }L^2_{\mathrm e}(\bbR).
\end{equation}
Conjugating by \(\mathcal U\), define
\begin{equation}\label{eq:definition of compactified S0-g-j}
    \widetilde S_0^{(\mathfrak g,j)}
    \coloneqq
    \mathcal U S_0^{(\mathfrak g,j)}\mathcal U^{-1}
    =
    I-c_j\widetilde T_{\mathfrak g}
    \qquad
    \text{on }L^2_{\mathrm e}(I).
\end{equation}
By Corollary~\ref{cor:part4-reduced-BS-bound-g}, the operator
\[
    P_{\widetilde\Phi_{\mathfrak g}^{\perp}}
    \widetilde S_0^{(\mathfrak g,j)}
    P_{\widetilde\Phi_{\mathfrak g}^{\perp}}
\]
is invertible on \(\operatorname{span}\{\widetilde\Phi_{\mathfrak g}\}^{\perp}\) for \(j=1,2\). It remains to examine the component along \(\widetilde\Phi_{\mathfrak g}\). We therefore use the orthogonal decomposition
\begin{equation}\label{eq:compactified-g-side-schur-decomposition-space}
    L^2_{\mathrm e}(I)
    =
    \operatorname{span}\{\widetilde\Phi_{\mathfrak g}\}
    \oplus
    \operatorname{span}\{\widetilde\Phi_{\mathfrak g}\}^{\perp}.
\end{equation}
For notational convenience, we write \(\widetilde\Phi_{\mathfrak g}^{\perp}\) for
\(\operatorname{span}\{\widetilde\Phi_{\mathfrak g}\}^{\perp}\).

\begin{lemma}[Schur-complement inverse bound for \(S_0^{(\mathfrak g,j)}\)]\label{lem:schur-complement-bound-g-side}
    For \(j=1,2\), define
    \begin{align}
        \widetilde h_j^{\mathfrak g}
        &\coloneqq
        \bigl(
            \widetilde\Phi_{\mathfrak g},
            \widetilde S_0^{(\mathfrak g,j)}
            \widetilde\Phi_{\mathfrak g}
        \bigr)_{L^2(I)},\label{eq:def-schur-h-g}
        \\
        \widetilde k_j^{\mathfrak g}
        &\coloneqq
        P_{\widetilde\Phi_{\mathfrak g}^{\perp}}
        \widetilde S_0^{(\mathfrak g,j)}
        \widetilde\Phi_{\mathfrak g},
        \label{eq:def-schur-k-g}
        \\
        \widetilde D_j^{\mathfrak g}
        &\coloneqq
        P_{\widetilde\Phi_{\mathfrak g}^{\perp}}
        \widetilde S_0^{(\mathfrak g,j)}
        P_{\widetilde\Phi_{\mathfrak g}^{\perp}}
        \quad
        \text{on }\widetilde\Phi_{\mathfrak g}^{\perp}.
        \label{eq:def-schur-D-g}
    \end{align}
    Assume that \(\widetilde D_j^{\mathfrak g}\) is invertible on
    \(\widetilde\Phi_{\mathfrak g}^{\perp}\), and set
    \begin{equation}\label{eq:def-schur-F-g}
        \widetilde F_j^{\mathfrak g}
        \coloneqq
        \widetilde h_j^{\mathfrak g}
        -
        \bigl(
            \widetilde k_j^{\mathfrak g},
            (\widetilde D_j^{\mathfrak g})^{-1}
            \widetilde k_j^{\mathfrak g}
        \bigr)_{L^2(I)} .
    \end{equation}
    If \(\widetilde F_j^{\mathfrak g}\neq0\), then \(S_0^{(\mathfrak g,j)}\) is invertible on
    \(L^2_{\mathrm e}(\bbR)\), and
    \begin{equation}\label{eq:schur-inverse-bound-g-side}
        \begin{aligned}
        \norm{
            (S_0^{(\mathfrak g,j)})^{-1}
        }_{L^2_{\mathrm e}(\bbR)\to L^2_{\mathrm e}(\bbR)}
        \leq
        \norm{
            (\widetilde D_j^{\mathfrak g})^{-1}
        }_{\widetilde\Phi_{\mathfrak g}^{\perp}\to
           \widetilde\Phi_{\mathfrak g}^{\perp}}
        +
        |\widetilde F_j^{\mathfrak g}|^{-1}
        \left(
            1+
            \norm{
                (\widetilde D_j^{\mathfrak g})^{-1}
                \widetilde k_j^{\mathfrak g}
            }_{L^2(I)}^2
        \right).
        \end{aligned}
    \end{equation}
\end{lemma}

\begin{proof}
    Since \(\mathcal U\) is unitary, it suffices to work with \(\widetilde S_0^{(\mathfrak g,j)}\) on \(L^2_{\mathrm e}(I)\). With respect to \eqref{eq:compactified-g-side-schur-decomposition-space}, we identify \(a\widetilde\Phi_{\mathfrak g}\) with \(a\in\mathbb R\). Then \eqref{eq:def-schur-h-g}--\eqref{eq:def-schur-D-g} give
    \begin{equation*}
        \widetilde S_0^{(\mathfrak g,j)}
        =
        \begin{pmatrix}
            \widetilde h_j^{\mathfrak g} & (\widetilde k_j^{\mathfrak g})^*\\
            \widetilde k_j^{\mathfrak g} & \widetilde D_j^{\mathfrak g}
        \end{pmatrix},
    \end{equation*}
    where \((\widetilde k_j^{\mathfrak g})^*\) denotes the functional
    \(w\mapsto(\widetilde k_j^{\mathfrak g},w)_{L^2(I)}\) on
    \(\widetilde\Phi_{\mathfrak g}^{\perp}\).
    If \(\widetilde D_j^{\mathfrak g}\) is invertible and \(\widetilde F_j^{\mathfrak g}\neq0\), the block inverse formula and \eqref{eq:def-schur-F-g} give
    \begin{equation}\label{eq:g-side-schur-block-inverse}
        \bigl(\widetilde S_0^{(\mathfrak g,j)}\bigr)^{-1}
        =
        \begin{pmatrix}
            0 & 0\\
            0 & (\widetilde D_j^{\mathfrak g})^{-1}
        \end{pmatrix}
        +
        \frac{1}{\widetilde F_j^{\mathfrak g}}
        \begin{pmatrix}
            1\\
            -(\widetilde D_j^{\mathfrak g})^{-1}\widetilde k_j^{\mathfrak g}
        \end{pmatrix}
        \otimes
        \begin{pmatrix}
            1\\
            -(\widetilde D_j^{\mathfrak g})^{-1}\widetilde k_j^{\mathfrak g}
        \end{pmatrix}.
    \end{equation}
    Hence \(\widetilde S_0^{(\mathfrak g,j)}\) is invertible.  Taking operator
    norms in \eqref{eq:g-side-schur-block-inverse} yields
    \[
        \norm{
            \bigl(\widetilde S_0^{(\mathfrak g,j)}\bigr)^{-1}
        }_{L^2_{\mathrm e}(I)\to L^2_{\mathrm e}(I)}
        \leq
        \norm{(\widetilde D_j^{\mathfrak g})^{-1}}_{\widetilde\Phi_{\mathfrak g}^{\perp}\to
           \widetilde\Phi_{\mathfrak g}^{\perp}}
        +
        |\widetilde F_j^{\mathfrak g}|^{-1}
        \left(
            1+
            \norm{
                (\widetilde D_j^{\mathfrak g})^{-1}\widetilde k_j^{\mathfrak g}
            }_{L^2(I)}^2
        \right).
    \]
    Unitary equivalence gives the same invertibility statement and norm bound for \(S_0^{(\mathfrak g,j)}\). This proves \eqref{eq:schur-inverse-bound-g-side}.
\end{proof}

Lemma~\ref{lem:schur-complement-bound-g-side} reduces the inverse bound for \(S_0^{(\mathfrak g,j)}\) to bounds for \(\widetilde h_j^{\mathfrak g}\), \((\widetilde D_j^{\mathfrak g})^{-1} \widetilde k_j^{\mathfrak g}\), and \(\widetilde F_j^{\mathfrak g}\). We approximate these quantities using the cosine-polynomial approximation from Section~\ref{subsec:finite-cosine-reduction-g}. The scalar entry is approximated directly. The vector \((\widetilde D_j^{\mathfrak g})^{-1} \widetilde k_j^{\mathfrak g}\) is approximated by a rational cosine polynomial, with an error bound obtained from its residual. Set
\begin{equation*}
    a_{\mathfrak g,\log}^{(L)}
    \coloneqq
    a_{\mathfrak g}^{(L)}\log|\cos|.
\end{equation*}
Define
\begin{equation*}
    \widetilde T_{\mathfrak g}^{(L)}
    \coloneqq
    \widetilde T_{\mathfrak g,\mathrm{sing}}^{(L)}
    +
    \frac1\pi
    \left(
        a_{\mathfrak g,\log}^{(L)}
        \otimes
        a_{\mathfrak g}^{(L)}
        +
        a_{\mathfrak g}^{(L)}
        \otimes
        a_{\mathfrak g,\log}^{(L)}
    \right),
\end{equation*}
where \(\widetilde T_{\mathfrak g,\mathrm{sing}}^{(L)}\) is defined in
\eqref{eq:part4-Tg-L-sing-red-def}. Define the exact and approximate scalar entries by
\begin{equation}\label{eq:part4-tau-and-tauL-g-def}
    \tau_{\mathfrak g}
    \coloneqq
    \bigl(
        \widetilde\Phi_{\mathfrak g},
        \widetilde T_{\mathfrak g}
        \widetilde\Phi_{\mathfrak g}
    \bigr)_{L^2(I)},
    \qquad
    \tau_{\mathfrak g}^{(L)}
    \coloneqq
    \bigl(
        \widetilde\Phi_{\mathfrak g}^{(L)},
        \widetilde T_{\mathfrak g}^{(L)}
        \widetilde\Phi_{\mathfrak g}^{(L)}
    \bigr)_{L^2(I)} .
\end{equation}
Then
\begin{equation}\label{eq:part4-hj-and-hjL-g-def}
    \widetilde h_j^{\mathfrak g}
    =
    1-c_j\tau_{\mathfrak g},
    \qquad
    \widetilde h_{j,L}^{\mathfrak g}
    \coloneqq
    1-c_j\tau_{\mathfrak g}^{(L)}.
\end{equation}
Using \eqref{eq:def-schur-k-g} and \eqref{eq:definition of compactified S0-g-j}, we have
\begin{equation}\label{eq:part4-kj-g-compact-form} 
    \widetilde k_j^{\mathfrak g}
    =
    -c_j
    P_{\widetilde\Phi_{\mathfrak g}^{\perp}}
    \widetilde T_{\mathfrak g}
    \widetilde\Phi_{\mathfrak g}.
\end{equation}
Define the approximate right-hand side by
\begin{equation}\label{eq:part4-kjL-g-def}
    \widetilde k_{j,L}^{\mathfrak g}
    \coloneqq
    -c_j
    \Pi_{\leq 3L}
    P_{(\widetilde\Phi_{\mathfrak g}^{(L)})^\perp}
    \widetilde T_{\mathfrak g}^{(L)}
    \widetilde\Phi_{\mathfrak g}^{(L)}
    \in \mathcal S_{3L}.
\end{equation}
Define the finite-dimensional operator
\begin{equation}\label{eq:part4-DjL-g-def}
    \widetilde D_{j,L}^{\mathfrak g}
    \coloneqq
    I_{\mathcal S_{3L}}
    -
    c_j
    \left.
    \widetilde T_{\mathfrak g,3L}^{(L)}
    \right|_{\mathcal S_{3L}}
    \qquad
    \text{on }\mathcal S_{3L}.
\end{equation}
Here \(\widetilde T_{\mathfrak g,3L}^{(L)}\) is defined in \eqref{eq:part4-Tg-3L-def}. 

Let
\begin{equation}\label{eq:part4-ZSch-jL-def}
    Z_{\mathrm{Sch},j}^{(L)}
    \in
    \mathcal S_{3L}
\end{equation}
be the rational cosine polynomial whose coefficients are specified in Section~\ref{sec:interval-inputs}. Define its residual norm by
\begin{equation}\label{eq:part4-rho-Sch-jL-def}
    \delta_{\mathrm{res,Sch},j}^{(L)}
    \coloneqq
    \norm{
        \widetilde D_{j,L}^{\mathfrak g}
        Z_{\mathrm{Sch},j}^{(L)}
        -
        \widetilde k_{j,L}^{\mathfrak g}
    }_{L^2(I)}.
\end{equation}
Define the approximate Schur complement by
\begin{equation}\label{eq:part4-FjL-g-def}
    \widetilde F_{j,L}^{\mathfrak g}
    \coloneqq
    \widetilde h_{j,L}^{\mathfrak g}
    -
    \bigl(
        \widetilde k_{j,L}^{\mathfrak g},
        Z_{\mathrm{Sch},j}^{(L)}
    \bigr)_{L^2(I)} .
\end{equation}

The next lemma bounds
\[
    \widetilde h_j^{\mathfrak g}
    -
    \widetilde h_{j,L}^{\mathfrak g},\qquad
    (\widetilde D_j^{\mathfrak g})^{-1}
    \widetilde k_j^{\mathfrak g}
    -
    Z_{\mathrm{Sch},j}^{(L)},
\]
and compares the quadratic terms in \(\widetilde F_j^{\mathfrak g}\) and \(\widetilde F_{j,L}^{\mathfrak g}\).

\begin{lemma}[Approximation bounds for the Schur complement]
\label{lem:schur-complement-approximation-bounds}
    Let \(\delta_T^{(L)}\), \(\mathsf C_{T,\mathrm{op}}^{(L)}\),
    \(\delta_{\mathrm{rhs}}^{(L)}\), and
    \(\Delta_{\mathrm{proj}}^{(L)}\) be the quantities defined in
    \eqref{eq:definition of Delta_T^{(L)}},
    \eqref{eq:definition of Omega_T^{(L)}},
    \eqref{eq:definition of Delta_{rhs}^{(L)}}, and
    \eqref{eq:definition of Delta_{red}^{(L)}}, respectively.
    For \(j=1,2\), let
    \begin{equation*}
        \mathsf C_{D,\mathrm{inv}}^{(j)}
        \coloneqq
        \frac{1}{1-0.1571\,c_j},
    \end{equation*}
    and
    \begin{equation}\label{eq:sigma-Schur-solve-def}
        \varepsilon_{\mathrm{Sch},j}^{(L)}
        \coloneqq
        \mathsf C_{D,\mathrm{inv}}^{(j)}
        \left[
            c_j\delta_{\mathrm{rhs}}^{(L)}
            +
            \delta_{\mathrm{res,Sch},j}^{(L)}
            +
            c_j\Delta_{\mathrm{proj}}^{(L)}
            \norm{Z_{\mathrm{Sch},j}^{(L)}}_{L^2(I)}
        \right].
    \end{equation}
    Then \(\widetilde D_j^{\mathfrak g}\) is invertible on
    \(\widetilde\Phi_{\mathfrak g}^{\perp}\), and
    \begin{equation}\label{eq:Dj-inverse-norm-Schur-bound}
        \norm{
            (\widetilde D_j^{\mathfrak g})^{-1}
        }_{\widetilde\Phi_{\mathfrak g}^{\perp}\to
        \widetilde\Phi_{\mathfrak g}^{\perp}}
        \leq
        \mathsf C_{D,\mathrm{inv}}^{(j)}.
    \end{equation}
    Furthermore, we have
    \begin{align}
        &\left|
            \widetilde h_j^{\mathfrak g}
            -
            \widetilde h_{j,L}^{\mathfrak g}
        \right|
        \leq
        c_j
        \left[
            \delta_T^{(L)}
            +
            2\mathsf C_{T,\mathrm{op}}^{(L)}
            \norm{
                \widetilde\Phi_{\mathfrak g}
                -
                \widetilde\Phi_{\mathfrak g}^{(L)}
            }_{L^2(I)}
        \right],
        \label{eq:hj-Schur-approximation-bound}
        \\
        &\norm{
            (\widetilde D_j^{\mathfrak g})^{-1}
            \widetilde k_j^{\mathfrak g}
        }_{L^2(I)}
        \leq
        \norm{Z_{\mathrm{Sch},j}^{(L)}}_{L^2(I)}
        +
        \varepsilon_{\mathrm{Sch},j}^{(L)},
        \label{eq:Dinvk-Schur-norm-upper}
    \end{align}
    and
    \begin{equation}\label{eq:Schur-inner-product-enclosure}
        \begin{aligned}
            &\left|
                \left(
                    \widetilde k_j^{\mathfrak g},
                    (\widetilde D_j^{\mathfrak g})^{-1}
                    \widetilde k_j^{\mathfrak g}
                \right)_{L^2(I)}
                -
                \left(
                    \widetilde k_{j,L}^{\mathfrak g},
                    Z_{\mathrm{Sch},j}^{(L)}
                \right)_{L^2(I)}
            \right|
            \\
            &\qquad\leq
            c_j\delta_{\mathrm{rhs}}^{(L)}
            \left(
                \norm{Z_{\mathrm{Sch},j}^{(L)}}_{L^2(I)}
                +
                \varepsilon_{\mathrm{Sch},j}^{(L)}
            \right)
            +
            \norm{
                \widetilde k_{j,L}^{\mathfrak g}
            }_{L^2(I)}
            \varepsilon_{\mathrm{Sch},j}^{(L)}.
        \end{aligned}
    \end{equation}
\end{lemma}

\begin{proof}
    Corollary~\ref{cor:part4-reduced-BS-bound-g} and unitary equivalence give
    \[
        \norm{
            \widetilde T_{\mathfrak g,\mathrm{proj}}
        }_{L^2_{\mathrm e}(I)\to L^2_{\mathrm e}(I)}
        <
        0.1571.
    \]
    By \eqref{eq:def-schur-D-g}, \eqref{eq:definition of compactified S0-g-j}, and
    \eqref{eq:part4-Tg-reduced-sing-only}, we have
    \[
        \widetilde D_j^{\mathfrak g}
        =
        I_{\widetilde\Phi_{\mathfrak g}^{\perp}}
        -
        c_j
        \left.
        \widetilde T_{\mathfrak g,\mathrm{proj}}
        \right|_{\widetilde\Phi_{\mathfrak g}^{\perp}}.
    \]
    Since
    \[
        c_j
        \norm{
            \left.
            \widetilde T_{\mathfrak g,\mathrm{proj}}
            \right|_{\widetilde\Phi_{\mathfrak g}^{\perp}}
        }_{\widetilde\Phi_{\mathfrak g}^{\perp}
          \to\widetilde\Phi_{\mathfrak g}^{\perp}}
        \leq
        c_j
        \norm{
            \widetilde T_{\mathfrak g,\mathrm{proj}}
        }_{L^2_{\mathrm e}(I)\to L^2_{\mathrm e}(I)}
        <
        0.1571\,c_j
        <
        1,
    \]
    Lemma~\ref{lem:Neumann-series-criterion} gives \eqref{eq:Dj-inverse-norm-Schur-bound}.
    By \eqref{eq:part4-hj-and-hjL-g-def}, we have
    \[
        \left|
            \widetilde h_j^{\mathfrak g}
            -
            \widetilde h_{j,L}^{\mathfrak g}
        \right|
        =
        c_j
        \left|
            \tau_{\mathfrak g}
            -
            \tau_{\mathfrak g}^{(L)}
        \right|.
    \]
    Since \(\norm{\widetilde\Phi_{\mathfrak g}}_{L^2(I)} =\norm{\widetilde\Phi_{\mathfrak g}^{(L)}}_{L^2(I)}=1\), the definitions in \eqref{eq:part4-tau-and-tauL-g-def} give
    \[
    \begin{aligned}
        \left|
            \tau_{\mathfrak g}
            -
            \tau_{\mathfrak g}^{(L)}
        \right|
        \leq
        \norm{
            \widetilde T_{\mathfrak g}
            -
            \widetilde T_{\mathfrak g}^{(L)}
        }_{L^2_{\mathrm e}(I)\to L^2_{\mathrm e}(I)}
        +
        2
        \norm{
            \widetilde T_{\mathfrak g}^{(L)}
        }_{L^2_{\mathrm e}(I)\to L^2_{\mathrm e}(I)}
        \norm{
            \widetilde\Phi_{\mathfrak g}
            -
            \widetilde\Phi_{\mathfrak g}^{(L)}
        }_{L^2(I)} .
    \end{aligned}
    \]
    Lemma~\ref{lem:app-finite-schur-comparison-estimates} now gives \eqref{eq:hj-Schur-approximation-bound}.
    
    For the following comparison, extend \(\widetilde D_j^{\mathfrak g}\) and \(\widetilde D_{j,L}^{\mathfrak g}\) to \(L^2_{\mathrm e}(I)\) by
    \[
        I-c_j\widetilde T_{\mathfrak g,\mathrm{proj}}
        \qquad\text{and}\qquad
        I-c_j\widetilde T_{\mathfrak g,3L}^{(L)},
    \]
    respectively. These extensions agree with \eqref{eq:def-schur-D-g} and \eqref{eq:part4-DjL-g-def} on \(\widetilde\Phi_{\mathfrak g}^{\perp}\) and \(\mathcal S_{3L}\).

    Adding and subtracting \(\widetilde k_{j,L}^{\mathfrak g}\) and \(\widetilde D_{j,L}^{\mathfrak g} Z_{\mathrm{Sch},j}^{(L)}\) gives
    \[
    \begin{aligned}
        (\widetilde D_j^{\mathfrak g})^{-1}\widetilde k_j^{\mathfrak g}
        -Z_{\mathrm{Sch},j}^{(L)}
        =
        (\widetilde D_j^{\mathfrak g})^{-1}
        \Big[
            (\widetilde k_j^{\mathfrak g}-\widetilde k_{j,L}^{\mathfrak g})
            +
            (\widetilde k_{j,L}^{\mathfrak g}
            -\widetilde D_{j,L}^{\mathfrak g}Z_{\mathrm{Sch},j}^{(L)})
            +
            (\widetilde D_{j,L}^{\mathfrak g}-\widetilde D_j^{\mathfrak g})
            Z_{\mathrm{Sch},j}^{(L)}
        \Big].
    \end{aligned}
    \]
    Lemma~\ref{lem:app-finite-schur-comparison-estimates} and \eqref{eq:part4-kj-g-compact-form}--\eqref{eq:part4-DjL-g-def} give
    \[
        \norm{
            \widetilde k_j^{\mathfrak g}
            -
            \widetilde k_{j,L}^{\mathfrak g}
        }_{L^2(I)}
        \leq
        c_j\delta_{\mathrm{rhs}}^{(L)},\qquad
        \norm{
            \widetilde D_j^{\mathfrak g}
            -
            \widetilde D_{j,L}^{\mathfrak g}
        }_{L^2_{\mathrm e}(I)\to L^2_{\mathrm e}(I)}
        \leq
        c_j\Delta_{\mathrm{proj}}^{(L)}.
    \]
    Using \eqref{eq:Dj-inverse-norm-Schur-bound}, \eqref{eq:sigma-Schur-solve-def}, and the definition of \(\delta_{\mathrm{res,Sch},j}^{(L)}\), we obtain
    \begin{equation}\label{eq:Schur-solution-error}
        \norm{
            (\widetilde D_j^{\mathfrak g})^{-1}
            \widetilde k_j^{\mathfrak g}
            -
            Z_{\mathrm{Sch},j}^{(L)}
        }_{L^2(I)}
        \leq
        \varepsilon_{\mathrm{Sch},j}^{(L)}.
    \end{equation}
    The triangle inequality then gives \eqref{eq:Dinvk-Schur-norm-upper}.
    Finally, 
    \[
        \begin{aligned}
            \left(
                \widetilde k_j^{\mathfrak g},
                (\widetilde D_j^{\mathfrak g})^{-1}
                \widetilde k_j^{\mathfrak g}
            \right)_{L^2(I)}
            -
            \left(
                \widetilde k_{j,L}^{\mathfrak g},
                Z_{\mathrm{Sch},j}^{(L)}
            \right)_{L^2(I)}
            &=
            \left(
                \widetilde k_j^{\mathfrak g}
                -
                \widetilde k_{j,L}^{\mathfrak g},
                (\widetilde D_j^{\mathfrak g})^{-1}
                \widetilde k_j^{\mathfrak g}
            \right)_{L^2(I)}
            \\
            &\quad+
            \left(
                \widetilde k_{j,L}^{\mathfrak g},
                (\widetilde D_j^{\mathfrak g})^{-1}
                \widetilde k_j^{\mathfrak g}
                -
                Z_{\mathrm{Sch},j}^{(L)}
            \right)_{L^2(I)} .
        \end{aligned}
    \]
    Cauchy--Schwarz, \eqref{eq:Dinvk-Schur-norm-upper}, and \eqref{eq:Schur-solution-error} give \eqref{eq:Schur-inner-product-enclosure}.
\end{proof}

Lemma~\ref{lem:schur-complement-approximation-bounds} and Proposition~\ref{prop:finite-Schur-complement-data} give the following proposition.

\begin{proposition}[Signs of the Schur complements and inverse bounds for \(S_0^{(\mathfrak g,j)}\)]\label{prop:certified-Schur-inverse-bounds-g}
    Let \(L=5000\), \(c_1=6.1\), and \(c_2=2.1\).  Then
    \begin{equation}\label{eq:certified-Schur-complement-signs-g}
        \widetilde F_1^{\mathfrak g}<0,
        \qquad
        \widetilde F_2^{\mathfrak g}>0.
    \end{equation}
    Consequently, \(S_0^{(\mathfrak g,j)}\) is invertible on
    \(L^2_{\mathrm e}(\bbR)\) for \(j=1,2\).  Moreover,
    \begin{equation}\label{eq:certified-S0-inverse-bounds-g}
        \norm{
            (S_0^{(\mathfrak g,1)})^{-1}
        }_{L^2_{\mathrm e}\to L^2_{\mathrm e}}
        \leq
        50.853,
        \qquad
        \norm{
            (S_0^{(\mathfrak g,2)})^{-1}
        }_{L^2_{\mathrm e}\to L^2_{\mathrm e}}
        \leq
        5.645.
    \end{equation}
\end{proposition}

\subsection{Limiting matrix entries for \(\mathfrak g\)}\label{subsec:g-side-limiting-matrix-entries}

We define the entries of the limiting constraint matrices associated with \(\mathfrak g\), using \(V_{\mathfrak g}\), \(\Phi_{\mathfrak g}\), and \(T_{\mathfrak g}\) from \eqref{eq:definition of g-side potential}--\eqref{eq:definition of T_g}, and \(S_0^{(\mathfrak g,j)}\) and \(\widetilde S_0^{(\mathfrak g,j)}\) from \eqref{eq:definition of S0-g-j}--\eqref{eq:definition of compactified S0-g-j}. Define
\begin{equation*}
    \Xi_1^{(\mathfrak g)}
    \coloneqq
    \{\mathfrak g,\Lambda\mathfrak g\},
    \qquad
    \Xi_2^{(\mathfrak g)}
    \coloneqq
    \{\Lambda\mathfrak g,\Lambda^2\mathfrak g\}.
\end{equation*}
For \(j=1,2\), write
\[
    \Xi_j^{(\mathfrak g)}
    =
    \{
        \Psi_1^{(\mathfrak g,j)},
        \Psi_2^{(\mathfrak g,j)}
    \}.
\]
For \(m=1,2\), define
\begin{equation*}
    (\Psi_m^{(\mathfrak g,j)})^{\log}
    \coloneqq
    -\frac1\pi
    \log|\cdot|*
    \Psi_m^{(\mathfrak g,j)},\qquad
    \psi_m^{(\mathfrak g,j)}
    \coloneqq
    |V_{\mathfrak g}|^{1/2}
    (\Psi_m^{(\mathfrak g,j)})^{\log}.
\end{equation*}

Define the entries that do not contain \((S_0^{(\mathfrak g,j)})^{-1}\) by
\begin{align}
    \left(
        \mathbf m^{(\mathfrak g,j)}
    \right)_m
    &\coloneqq
    \int_{\bbR}
        \Psi_m^{(\mathfrak g,j)}(y)\,dy,
    \qquad
    m=1,2,\label{eq:g-side-m-vector-def} 
    \\
    \left(
        \mathbf C_{\log}^{(\mathfrak g,j)}
    \right)_{mn}
    &\coloneqq
    \left(
        (\Psi_m^{(\mathfrak g,j)})^{\log},
        \Psi_n^{(\mathfrak g,j)}
    \right)_r,
    \qquad
    1\leq m,n\leq2.\label{eq:g-side-Clog-def}
\end{align}

The remaining entries contain \((S_0^{(\mathfrak g,j)})^{-1}\). Define
\begin{align}
    \mathfrak s_0^{(\mathfrak g,j)}
    &\coloneqq
    \left(
        |V_{\mathfrak g}|^{1/2},
        \left(S_0^{(\mathfrak g,j)}\right)^{-1}
        |V_{\mathfrak g}|^{1/2}
    \right)_r,\label{eq:g-side-scalar-s-def}
    \\
    \mathfrak a_m^{(\mathfrak g,j)}
    &\coloneqq
    \left(
        \psi_m^{(\mathfrak g,j)},
        \left(S_0^{(\mathfrak g,j)}\right)^{-1}
        |V_{\mathfrak g}|^{1/2}
    \right)_r,
    \qquad
    m=1,2,\nonumber
    \\
    \mathfrak r_{mn}^{(\mathfrak g,j)}
    &\coloneqq
    \left(
        \psi_m^{(\mathfrak g,j)},
        \left(S_0^{(\mathfrak g,j)}\right)^{-1}
        \psi_n^{(\mathfrak g,j)}
    \right)_r,
    \qquad
    1\leq m,n\leq2.\label{eq:g-side-r-matrix-def}
\end{align}
We also write
\begin{equation*}
    \mathbf a^{(\mathfrak g,j)}
    \coloneqq
    \begin{pmatrix}
        \mathfrak a_1^{(\mathfrak g,j)}\\
        \mathfrak a_2^{(\mathfrak g,j)}
    \end{pmatrix},
    \qquad
    \mathbf R^{(\mathfrak g,j)}
    \coloneqq
    \left(
        \mathfrak r_{mn}^{(\mathfrak g,j)}
    \right)_{1\leq m,n\leq2}.
\end{equation*}

We first evaluate the entries in \eqref{eq:g-side-m-vector-def}--\eqref{eq:g-side-Clog-def}. Define
\begin{equation}\label{eq:g-side-density-compactification-map}
    (\mathscr C f)(\theta)
    \coloneqq
    f(\tan\theta)\sec^2\theta,
    \qquad
    \theta\in I.
\end{equation}
For \(j,m\in\{1,2\}\), set
\begin{equation}\label{eq:g-side-Omega-def}
    \Omega_m^{(j)}
    \coloneqq
    \mathscr C\Psi_m^{(\mathfrak g,j)}.
\end{equation}
For the chosen profile \(\mathfrak g\), each \(\Omega_m^{(j)}\) is a cosine polynomial.
The next lemma expresses \(\mathbf m^{(\mathfrak g,j)}\) and \(\mathbf C_{\log}^{(\mathfrak g,j)}\) in terms of the cosine coefficients of \(\Omega_m^{(j)}\).

\begin{lemma}[Cosine-coefficient formulas for \(\mathbf m^{(\mathfrak g,j)}\) and \(\mathbf C_{\log}^{(\mathfrak g,j)}\)]\label{lem:g-side-free-source-entries}
    For \(j=1,2\) and \(m=1,2\), write
    \[
        \Omega_m^{(j)}(\theta)
        =
        \sum_{k\geq0}
            \Omega_{m,k}^{(j)}\cos(2k\theta),
        \qquad
        \Omega_{m,k}^{(j)}=0
        \quad \text{for all sufficiently large } k.
    \]
    Set
    \[
        \ell_{\cos,m}^{(j)}
        \coloneqq
        -\pi(\log2)\Omega_{m,0}^{(j)}
        -
        \frac{\pi}{2}
        \sum_{k\geq1}
            \frac{(-1)^k}{k}\Omega_{m,k}^{(j)} .
    \]
    Then
    \begin{equation}\label{eq:g-side-m-vector-cosine-formula}
        \bigl(\mathbf m^{(\mathfrak g,j)}\bigr)_m
        =
        \pi\,\Omega_{m,0}^{(j)},
        \qquad m=1,2,
    \end{equation}
    and
    \begin{equation}\label{eq:g-side-Clog-cosine-formula}
        \begin{aligned}
        \bigl(\mathbf C_{\log}^{(\mathfrak g,j)}\bigr)_{mn}
        =
        -\frac1\pi
        \bigg[
            &-\pi^2(\log2)\,
            \Omega_{m,0}^{(j)}\Omega_{n,0}^{(j)}
            -
            \frac{\pi^2}{4}
            \sum_{k\geq1}
                \frac{
                    \Omega_{m,k}^{(j)}
                    \Omega_{n,k}^{(j)}
                }{k}
            \\
            &-
            \pi\,\Omega_{m,0}^{(j)}\ell_{\cos,n}^{(j)}
            -
            \ell_{\cos,m}^{(j)}\pi\,\Omega_{n,0}^{(j)}
        \bigg],
        \qquad 1\leq m,n\leq2.
        \end{aligned}
    \end{equation}
\end{lemma}

\begin{proof}
    Since \(dy=\sec^2\theta\,d\theta\) and \(y=\tan\theta\), the definition of \(\Omega_m^{(j)}\) gives
    \[
        \int_{\mathbb R}
            \Psi_m^{(\mathfrak g,j)}(y)\,dy
        =
        \int_I
            \Omega_m^{(j)}(\theta)\,d\theta.
    \]
    This proves \eqref{eq:g-side-m-vector-cosine-formula}. For the logarithmic entries, the change of variables gives
    \[
        \bigl(\mathbf C_{\log}^{(\mathfrak g,j)}\bigr)_{mn}
        =
        -\frac1\pi
        \int_I\int_I
            \log|\tan\theta-\tan\varphi|\,
            \Omega_m^{(j)}(\varphi)
            \Omega_n^{(j)}(\theta)\,
        d\varphi d\theta .
    \]
    Using
    $\log|\tan\theta-\tan\varphi|=\log|\sin(\theta-\varphi)|-\log|\cos\theta|-\log|\cos\varphi|$,
    together with
    \[
        \log|\sin\theta|
        =
        -\log2-\sum_{k\geq1}\frac{\cos(2k\theta)}{k},
        \qquad
        \log|\cos\theta|
        =
        -\log2-\sum_{k\geq1}\frac{(-1)^k}{k}\cos(2k\theta),
    \]
    and the orthogonality of the even cosine basis on \(I\), we obtain
    \eqref{eq:g-side-Clog-cosine-formula}. The sums are finite because each
    \(\Omega_m^{(j)}\) is a cosine polynomial.
\end{proof}

We next consider the entries in \eqref{eq:g-side-scalar-s-def}--\eqref{eq:g-side-r-matrix-def}. We introduce right-hand sides on \(I\) and the corresponding linear equations.

Write \(L_{\cos}(\theta)\coloneqq\log|\cos\theta|\).  For an even function
\(f\) on \(I\), define
\begin{equation*}
    \mathcal I_0[f]
    \coloneqq
    \int_I f(\varphi)\,d\varphi,
    \qquad
    \mathcal I_{\cos}[f]
    \coloneqq
    \int_I f(\varphi)L_{\cos}(\varphi)\,d\varphi.
\end{equation*}
Define
\begin{equation}\label{eq:K_log-definition}
    \mathcal K_{\log}[f]
    \coloneqq
    \mathcal K_{\sin}[f]
    +
    L_{\cos}\mathcal I_0[f]
    +
    \mathcal I_{\cos}[f].
\end{equation}

For \(j=1,2\), define \(B_m^{(j)}\in L^2_{\mathrm e}(I)\) by
\begin{equation}\label{eq:g-side-B-def}
    B_0^{(j)}
    \coloneqq
    a_{\mathfrak g},
    \qquad
    B_m^{(j)}
    \coloneqq
    \frac1\pi
    a_{\mathfrak g}
    \mathcal K_{\log}[\Omega_m^{(j)}],
    \qquad
    m=1,2.
\end{equation}
For \(m=0,1,2\) and \(j=1,2\), define
\begin{equation*}
    X_m^{(\mathfrak g,j)}
    \coloneqq
    \left(
        \widetilde S_0^{(\mathfrak g,j)}
    \right)^{-1}
    B_m^{(j)},
    \qquad
    m=0,1,2,\quad j=1,2.
\end{equation*}
Unitary equivalence gives
\begin{align}
    \mathfrak s_0^{(\mathfrak g,j)}
    &=
    \bigl(
        B_0^{(j)},
        X_0^{(\mathfrak g,j)}
    \bigr)_{L^2(I)},
    \label{eq:g-side-s-compactified-entry}
    \\
    \mathfrak a_m^{(\mathfrak g,j)}
    &=
    \bigl(
        B_m^{(j)},
        X_0^{(\mathfrak g,j)}
    \bigr)_{L^2(I)},
    \qquad m=1,2,
    \nonumber
    \\
    \mathfrak r_{mn}^{(\mathfrak g,j)}
    &=
    \bigl(
        B_m^{(j)},
        X_n^{(\mathfrak g,j)}
    \bigr)_{L^2(I)},
    \qquad 1\leq m,n\leq2.
    \label{eq:g-side-r-compactified-entry}
\end{align}

We approximate each \(X_m^{(\mathfrak g,j)}\) by a cosine polynomial \(Z_m^{(j)}\). The residual identity below and the inverse bound in Proposition~\ref{prop:certified-Schur-inverse-bounds-g} give an error bound for \(X_m^{(\mathfrak g,j)}-Z_m^{(j)}\).

For \(j=1,2\) and \(m=0,1,2\), let
\begin{equation}\label{eq:g-side-Z-def}
    Z_m^{(j)}(\theta)
    \coloneqq
    \sum_{k=0}^{L_Z}
        z_{m,k}^{(j)}\cos(2k\theta)
    \in \mathcal S_{L_Z}
\end{equation}
be the cosine polynomial whose coefficients are specified in Section~\ref{sec:interval-inputs}. Define the residual by
\begin{equation}\label{eq:g-side-continuous-residual-def}
    \operatorname{Res}_m^{(j)}
    \coloneqq
    B_m^{(j)}
    -
    \widetilde S_0^{(\mathfrak g,j)}Z_m^{(j)}
    =
    B_m^{(j)}
    -
    Z_m^{(j)}
    +
    \frac{c_j}{\pi}
    a_{\mathfrak g}
    \mathcal K_{\log}
    \left[
        a_{\mathfrak g}Z_m^{(j)}
    \right].
\end{equation}
Then
\begin{equation}\label{eq:g-side-X-minus-Z-residual-identity}
    X_m^{(\mathfrak g,j)}-Z_m^{(j)}
    =
    \left(
        \widetilde S_0^{(\mathfrak g,j)}
    \right)^{-1}
    \operatorname{Res}_m^{(j)} .
\end{equation}

We next define the DCT-I approximations used in the residual
estimate. Let \(\mathcal R_N\) be the DCT-I approximation operator defined in
Appendix~\ref{subsec:app-dct1-aliasing}. For \(m=0,1,2\), set
\begin{align}
    \mathscr A_{m,N}^{(j)}
    &\coloneqq
    (\mathcal R_N a_{\mathfrak g})
    \mathcal K_{\sin}
    \left[
        \mathcal R_N(a_{\mathfrak g}Z_m^{(j)})
    \right]
    \notag\\
    &\quad+
    \mathcal R_N(a_{\mathfrak g,\log})
    \mathcal I_0
    \left[
        \mathcal R_N(a_{\mathfrak g}Z_m^{(j)})
    \right]
    +
    (\mathcal R_N a_{\mathfrak g})
    \mathcal I_{\cos}
    \left[
        \mathcal R_N(a_{\mathfrak g}Z_m^{(j)})
    \right].
    \label{eq:g-side-discrete-action-def}
\end{align}
Define
\begin{equation}\label{eq:g-side-discrete-rhs-def}
    B_{m,N}^{(j)}
    \coloneqq
    \begin{cases}
        \mathcal R_N a_{\mathfrak g},
        & m=0,\\[1mm]
        \displaystyle
        \frac1\pi
        \left(
            (\mathcal R_Na_{\mathfrak g})
            \mathcal K_{\sin}[\Omega_m^{(j)}]
            +
            \mathcal R_N(a_{\mathfrak g,\log})
            \mathcal I_0[\Omega_m^{(j)}]
            +
            (\mathcal R_Na_{\mathfrak g})
            \mathcal I_{\cos}[\Omega_m^{(j)}]
        \right),
        & m=1,2.
    \end{cases}
\end{equation}
Define the DCT-I residual by
\begin{equation}\label{eq:g-side-discrete-residual-def}
    \operatorname{Res}_{m,N}^{(j)}
    \coloneqq
    B_{m,N}^{(j)}
    -
    Z_m^{(j)}
    +
    \frac{c_j}{\pi}
    \mathscr A_{m,N}^{(j)} .
\end{equation}

Lemma~\ref{lem:app-g-side-residual-comparison} bounds
\[
    B_m^{(j)}-B_{m,N}^{(j)}
    \quad\text{and}\quad
    \operatorname{Res}_m^{(j)}
    -
    \operatorname{Res}_{m,N}^{(j)}.
\]
The next lemma uses these bounds to estimate
\[
    X_m^{(\mathfrak g,j)}-Z_m^{(j)}
\]
and the pairings in \eqref{eq:g-side-s-compactified-entry}--\eqref{eq:g-side-r-compactified-entry}.

\begin{lemma}[Residual-based error bounds for the solutions and pairings]\label{lem:g-side-finite-reduction-resolvent-pairings}
    Fix \(j\in\{1,2\}\). Let \(\delta_{\mathrm{rhs},m}^{(j)}(N)\) and \(\delta_{\mathrm{res\text{-}diff},m}^{(j)}(N)\) be the error bounds defined in
    Lemma~\ref{lem:app-g-side-residual-comparison}. For \(m=0,1,2\), set
    \begin{equation}\label{eq:g-side-discrete-residual-radius-def}
        \delta_{\mathrm{res},m}^{(j)}(N)
        \coloneqq
        \norm{
            \operatorname{Res}_{m,N}^{(j)}
        }_{L^2(I)},
    \end{equation}
    and define
    \begin{equation*}
        \varepsilon_m^{(j)}(N)
        \coloneqq
        \norm{
            \left(
                \widetilde S_0^{(\mathfrak g,j)}
            \right)^{-1}
        }_{L^2(I)\to L^2(I)}
        \left(
            \delta_{\mathrm{res},m}^{(j)}(N)
            +
            \delta_{\mathrm{res\text{-}diff},m}^{(j)}(N)
        \right).
    \end{equation*}
    Then
    \begin{equation}\label{eq:g-side-X-Z-error-bound}
        \norm{
            X_m^{(\mathfrak g,j)}-Z_m^{(j)}
        }_{L^2(I)}
        \leq
        \varepsilon_m^{(j)}(N),
        \qquad m=0,1,2.
    \end{equation}
    Moreover, for \(m,n\in\{0,1,2\}\), set
    \begin{equation}\label{eq:g-side-pairing-error-radius-def}
        \varepsilon_{\mathrm{pair},mn}^{(j)}(N)
        \coloneqq
        \norm{B_{m,N}^{(j)}}_{L^2(I)}
        \varepsilon_n^{(j)}(N)
        +
        \delta_{\mathrm{rhs},m}^{(j)}(N)
        \left(
            \norm{Z_n^{(j)}}_{L^2(I)}
            +
            \varepsilon_n^{(j)}(N)
        \right).
    \end{equation}
    Then
    \begin{equation}\label{eq:g-side-pairing-error-bound}
        \left|
            (B_m^{(j)},X_n^{(\mathfrak g,j)})_{L^2(I)}
            -
            (B_{m,N}^{(j)},Z_n^{(j)})_{L^2(I)}
        \right|
        \leq
        \varepsilon_{\mathrm{pair},mn}^{(j)}(N).
    \end{equation}
\end{lemma}

\begin{proof}
    By \eqref{eq:g-side-X-minus-Z-residual-identity},
    \[
        \norm{
            X_m^{(\mathfrak g,j)}-Z_m^{(j)}
        }_{L^2(I)}
        \leq
        \norm{
            \left(
                \widetilde S_0^{(\mathfrak g,j)}
            \right)^{-1}
        }_{L^2(I)\to L^2(I)}
        \norm{
            \operatorname{Res}_m^{(j)}
        }_{L^2(I)}.
    \]
    By \eqref{eq:g-side-discrete-residual-radius-def} and \eqref{eq:app-g-side-Res-minus-ResN-bound}, 
    \[
        \norm{
            \operatorname{Res}_m^{(j)}
        }_{L^2(I)}
        \leq
        \delta_{\mathrm{res},m}^{(j)}(N)
        +
        \delta_{\mathrm{res\text{-}diff},m}^{(j)}(N).
    \]
    This proves \eqref{eq:g-side-X-Z-error-bound}.  For the pairing estimate,
    write
    \[
    \begin{aligned}
        (B_m^{(j)},X_n^{(\mathfrak g,j)})_{L^2(I)}
        -
        (B_{m,N}^{(j)},Z_n^{(j)})_{L^2(I)}
        =
        (B_{m,N}^{(j)},X_n^{(\mathfrak g,j)}-Z_n^{(j)})_{L^2(I)}
        +
        (B_m^{(j)}-B_{m,N}^{(j)},X_n^{(\mathfrak g,j)})_{L^2(I)}.
    \end{aligned}
    \]
    By Cauchy--Schwarz, \eqref{eq:app-g-side-B-minus-BN-bound}, and
    \eqref{eq:g-side-X-Z-error-bound},
    \[
    \begin{aligned}
        \left|
            (B_m^{(j)},X_n^{(\mathfrak g,j)})_{L^2(I)}
            -
            (B_{m,N}^{(j)},Z_n^{(j)})_{L^2(I)}
        \right|
        \leq
        \norm{B_{m,N}^{(j)}}_{L^2(I)}
        \varepsilon_n^{(j)}(N)
        +
        \delta_{\mathrm{rhs},m}^{(j)}(N)
        \norm{X_n^{(\mathfrak g,j)}}_{L^2(I)}.
    \end{aligned}
    \]
    Since
    \[
        \norm{X_n^{(\mathfrak g,j)}}_{L^2(I)}
        \leq
        \norm{Z_n^{(j)}}_{L^2(I)}
        +
        \varepsilon_n^{(j)}(N),
    \]
    the definition \eqref{eq:g-side-pairing-error-radius-def} gives
    \eqref{eq:g-side-pairing-error-bound}.
\end{proof}

Proposition~\ref{prop:g-residual-pairing-intervals}, together with Lemmas~\ref{lem:g-side-free-source-entries} and \ref{lem:g-side-finite-reduction-resolvent-pairings}, gives the following
interval enclosures.

\begin{proposition}[Interval enclosures for the limiting matrix entries for \(\mathfrak g\)]
\label{prop:part4-certified-g-side-limiting-data}
    With \(c_1=6.1\), \(c_2=2.1\), \(L=5000\), \(L_Z=8000\), and
    \(N=65536\), the following interval enclosures hold.

    The vectors \(\mathbf m^{(\mathfrak g,1)}\) and \(\mathbf m^{(\mathfrak g,2)}\) satisfy
    \begin{equation}\label{eq:part4-g-side-m}
        \mathbf m^{(\mathfrak g,1)}
        \in
        \begin{pmatrix}
            [3.24820,3.24821]\\
            [-1.62411,-1.62410]
        \end{pmatrix},
        \qquad
        \mathbf m^{(\mathfrak g,2)}
        \in
        \begin{pmatrix}
            [-1.62411,-1.62410]\\
            [0.81205,0.81206]
        \end{pmatrix}.
    \end{equation}
    The matrices \(\mathbf C_{\log}^{(\mathfrak g,1)}\) and \(\mathbf C_{\log}^{(\mathfrak g,2)}\) satisfy
    \begin{equation}\label{eq:part4-g-side-C_log-j1}
        \mathbf C_{\log}^{(\mathfrak g,1)}
        \in
        \begin{pmatrix}
            [-1.65893,-1.65892]
            &
            [2.50867,2.50868]
            \\
            [2.50867,2.50868]
            &
            [-1.42882,-1.42881]
        \end{pmatrix},
    \end{equation}
    and
    \begin{equation}\label{eq:part4-g-side-C_log-j2}
        \mathbf C_{\log}^{(\mathfrak g,2)}
        \in
        \begin{pmatrix}
            [-1.42882,-1.42881]
            &
            [1.13420,1.13421]
            \\
            [1.13420,1.13421]
            &
            [-0.43079,-0.43078]
        \end{pmatrix}.
    \end{equation}
    For \(j=1\),
    \begin{equation}\label{eq:part4-g-side-s-a-certified-j1}
        \mathfrak s_0^{(\mathfrak g,1)}
        \in
        [-0.05294,-0.05176],
        \qquad
        \mathbf a^{(\mathfrak g,1)}
        \in
        \begin{pmatrix}
            [-0.39694,-0.39634]\\
            [-0.10116,-0.10033]
        \end{pmatrix},
    \end{equation}
    and
    \begin{equation}\label{eq:part4-g-side-R-certified-j1}
        \mathbf R^{(\mathfrak g,1)}
        \in
        \begin{pmatrix}
            [0.09333,0.10240]
            &
            [-0.24795,-0.23363]
            \\
            [-0.24795,-0.23363]
            &
            [-0.09287,-0.07814]
        \end{pmatrix}.
    \end{equation}
    For \(j=2\),
    \begin{equation}\label{eq:part4-g-side-s-a-certified-j2}
        \mathfrak s_0^{(\mathfrak g,2)}
        \in
        [3.72601,3.72766],
        \qquad
        \mathbf a^{(\mathfrak g,2)}
        \in
        \begin{pmatrix}
            [2.49603,2.49732]\\
            [-3.70367,-3.70205]
        \end{pmatrix},
    \end{equation}
    and
    \begin{equation}\label{eq:part4-g-side-R-certified-j2}
        \mathbf R^{(\mathfrak g,2)}
        \in
        \begin{pmatrix}
            [1.68375,1.68608]
            &
            [-2.49030,-2.48700]
            \\
            [-2.49030,-2.48700]
            &
            [3.70689,3.71130]
        \end{pmatrix}.
    \end{equation}
    The solution norms satisfy
    \begin{equation}\label{eq:part4-g-side-Xnorm-certified}
        \begin{array}{c|ccc}
            j
            &
            \norm{X_0^{(\mathfrak g,j)}}_{L^2(I)}
            &
            \norm{X_1^{(\mathfrak g,j)}}_{L^2(I)}
            &
            \norm{X_2^{(\mathfrak g,j)}}_{L^2(I)}
            \\
            \hline
            1
            &
            \leq 1.05977
            &
            \leq 0.44736
            &
            \leq 0.59973
            \\
            2
            &
            \leq 3.92178
            &
            \leq 2.65403
            &
            \leq 4.00728
        \end{array}
    \end{equation}
\end{proposition}

\section{Perturbation estimates from \(\mathfrak g\) to \(Q\)}\label{sec:transfer from g to Q}
In this section, we derive the spectral estimates for \(Q\) from the corresponding estimates for \(\mathfrak g\). In Subsection~\ref{subsec:transfer-reduced-BS-bound}, we bound
\(T_Q-T_{\mathfrak g}\) and the difference between the projected Birman--Schwinger operators.
These bounds prove Proposition~\ref{prop:reduced-BS-bound-Q}. In Subsection~\ref{subsec:transfer-limiting-matrices-Q}, we prove the invertibility of \(S_{0,j}\) and compare the limiting matrix entries for \(Q\) and \(\mathfrak g\).

\subsection{Comparison of the spectral quantities for \(\mathfrak g\) and \(Q\)}\label{subsec:transfer-reduced-BS-bound}
Lemma~\ref{lem:square-root-potential-transfer} gives bounds for
\[
    |V_Q|^{1/2}-|V_{\mathfrak g}|^{1/2}
    \quad\text{and}\quad
    \omega
    \left(
        |V_Q|^{1/2}-|V_{\mathfrak g}|^{1/2}
    \right).
\]
We use these bounds to estimate \(T_Q-T_{\mathfrak g}\) and \(P_{\Phi_Q^\perp}-P_{\Phi_{\mathfrak g}^\perp}\). The next lemma combines the two estimates.

\begin{lemma}[Perturbation bound for the projected Birman--Schwinger operators] \label{lem:part4-projected-log-transfer}
    Let \(\Delta_{\sqrt V}^{L^2}\) and \(\Delta_{\omega\sqrt V}^{L^2}\) be the bounds from
    Lemma~\ref{lem:square-root-potential-transfer}. Define
    \begin{equation}\label{eq:projected-transfer-error-definition}
    \begin{aligned}
        \Delta_{\mathrm{proj}}
        \coloneqq{}&
        \frac8\pi
        \left(
            \norm{V_{\mathfrak g}}_{L^\infty}
            +
            \norm{\omega^2V_{\mathfrak g}}_{L^1}
        \right)^{1/2}
        \Delta_{\sqrt V}^{L^2}
        \\
        &+
        \frac1\pi
        \Bigg[
        \left(
            \left(
                4\mathsf V^{L^\infty}
                +
                2\mathsf V_{\omega^2}^{L^1}
            \right)
            \bigl(\Delta_{\sqrt V}^{L^2}\bigr)^2
            +
            2\mathsf V^{L^1}
            \bigl(\Delta_{\omega\sqrt V}^{L^2}\bigr)^2
        \right)^{1/2}
        \\
        &\qquad\quad+
        \left(
            \left(
                4\norm{V_{\mathfrak g}}_{L^\infty}
                +
                2\norm{\omega^2V_{\mathfrak g}}_{L^1}
            \right)
            \bigl(\Delta_{\sqrt V}^{L^2}\bigr)^2
            +
            2\norm{V_{\mathfrak g}}_{L^1}
            \bigl(\Delta_{\omega\sqrt V}^{L^2}\bigr)^2
        \right)^{1/2}
        \Bigg].
    \end{aligned}
    \end{equation}
    Then
    \begin{equation}\label{eq:part4-projected-log-transfer-bound}
        \norm{
            P_{\Phi_Q^\perp}T_QP_{\Phi_Q^\perp}
            -
            P_{\Phi_{\mathfrak g}^\perp}
            T_{\mathfrak g}
            P_{\Phi_{\mathfrak g}^\perp}
        }_{L^2_{\mathrm e}\to L^2_{\mathrm e}}
        \leq
        \Delta_{\mathrm{proj}}.
    \end{equation}
\end{lemma}

\begin{proof}
    Adding and subtracting \(P_{\Phi_Q^\perp}T_{\mathfrak g}P_{\Phi_Q^\perp}\) and 
    \(P_{\Phi_{\mathfrak g}^\perp} T_{\mathfrak g} P_{\Phi_Q^\perp}\)
    and using that orthogonal projections have norm one, we obtain
    \begin{equation}\label{eq:part4-projected-log-transfer-split}
        \begin{aligned}
            \norm{
                P_{\Phi_Q^\perp}T_QP_{\Phi_Q^\perp}
                -
                P_{\Phi_{\mathfrak g}^\perp}
                T_{\mathfrak g}
                P_{\Phi_{\mathfrak g}^\perp}
            }_{L^2_{\mathrm e}\to L^2_{\mathrm e}}
            \leq
            \norm{T_Q-T_{\mathfrak g}}_{L^2\to L^2}
            +
            2
            \norm{
                P_{\Phi_Q^\perp}
                -
                P_{\Phi_{\mathfrak g}^\perp}
            }_{L^2\to L^2}
            \norm{T_{\mathfrak g}}_{L^2\to L^2}.
        \end{aligned}
    \end{equation}
    Lemma~\ref{lem:app-rank-one-projection-difference} gives
    \[
        \norm{
            P_{\Phi_Q^\perp}
            -
            P_{\Phi_{\mathfrak g}^\perp}
        }_{L^2\to L^2}
        \leq
        \frac{
            2\norm{\delta_V}_{L^2}
        }{
            \norm{|V_{\mathfrak g}|^{1/2}}_{L^2}
        }
        \leq
        \frac{
            2\Delta_{\sqrt V}^{L^2}
        }{
            \norm{V_{\mathfrak g}}_{L^1}^{1/2}
        }.
    \]
    Also, by the Hilbert--Schmidt bound and the decomposition
    \begin{equation}\label{eq:logarithmic decomposition}
        \int_{|\tau|\leq1}\log^2|\tau|\,d\tau=4,
        \qquad
        \log^2|x-y|
        \leq
        2\omega(x)^2+2\omega(y)^2
        \quad (|x-y|>1),
    \end{equation}
    we obtain
    \[
        \norm{T_{\mathfrak g}}_{L^2\to L^2}
        \leq
        \frac2\pi
        \left[
            \norm{V_{\mathfrak g}}_{L^1}
            \left(
                \norm{V_{\mathfrak g}}_{L^\infty}
                +
                \norm{\omega^2V_{\mathfrak g}}_{L^1}
            \right)
        \right]^{1/2}.
    \]
    Therefore
    \begin{equation}\label{eq:part4-projection-term-bound}
        \begin{aligned}
            2
            \norm{
                P_{\Phi_Q^\perp}
                -
                P_{\Phi_{\mathfrak g}^\perp}
            }_{L^2\to L^2}
            \norm{T_{\mathfrak g}}_{L^2\to L^2}
            \leq
            \frac8\pi
            \left(
                \norm{V_{\mathfrak g}}_{L^\infty}
                +
                \norm{\omega^2V_{\mathfrak g}}_{L^1}
            \right)^{1/2}
            \Delta_{\sqrt V}^{L^2}.
        \end{aligned}
    \end{equation}
    It remains to estimate \(T_Q-T_{\mathfrak g}\). The definitions of \(T_Q\) and \(T_{\mathfrak g}\) give
    \[
        \begin{aligned}
            (T_Q-T_{\mathfrak g})f
            =
            -\frac1\pi
            \delta_V
            \left(
                \log|\cdot|*
                \left(
                    |V_Q|^{1/2}f
                \right)
            \right)
            -
            \frac1\pi
            |V_{\mathfrak g}|^{1/2}
            \left(
                \log|\cdot|*
                \left(
                    \delta_V f
                \right)
            \right).
        \end{aligned}
    \]
    Hence the Hilbert--Schmidt bound gives
    \[
        \begin{aligned}
            \norm{T_Q-T_{\mathfrak g}}_{L^2\to L^2}
            &\leq
            \frac1\pi
            \left(
                \iint_{\bbR^2}
                \delta_V(x)^2
                \log^2|x-y|
                |V_Q(y)|
                \,dx\,dy
            \right)^{1/2}
            \\
            &\quad+
            \frac1\pi
            \left(
                \iint_{\bbR^2}
                |V_{\mathfrak g}(x)|
                \log^2|x-y|
                \delta_V(y)^2
                \,dx\,dy
            \right)^{1/2}.
        \end{aligned}
    \]
    Using the kernel bounds in \eqref{eq:logarithmic decomposition} and Lemma~\ref{lem:square-root-potential-transfer}, we obtain
    \[
        \norm{\delta_V}_{L^2}\leq\Delta_{\sqrt V}^{L^2},
        \qquad
        \norm{\omega\delta_V}_{L^2}\leq\Delta_{\omega\sqrt V}^{L^2}.
    \]
    Hence
    \begin{equation}\label{eq:part4-TQ-minus-Tg-bound}
        \begin{aligned}
            \norm{T_Q-T_{\mathfrak g}}_{L^2\to L^2}
            &\leq
            \frac1\pi
            \Bigg[
            \left(
                \left(
                    4\norm{V_Q}_{L^\infty}
                    +
                    2\norm{\omega^2V_Q}_{L^1}
                \right)
                \bigl(\Delta_{\sqrt V}^{L^2}\bigr)^2
                +
                2\norm{V_Q}_{L^1}
                \bigl(\Delta_{\omega\sqrt V}^{L^2}\bigr)^2
            \right)^{1/2}
            \\
            &\qquad\quad+
            \left(
                \left(
                    4\norm{V_{\mathfrak g}}_{L^\infty}
                    +
                    2\norm{\omega^2V_{\mathfrak g}}_{L^1}
                \right)
                \bigl(\Delta_{\sqrt V}^{L^2}\bigr)^2
                +
                2\norm{V_{\mathfrak g}}_{L^1}
                \bigl(\Delta_{\omega\sqrt V}^{L^2}\bigr)^2
            \right)^{1/2}
            \Bigg].
        \end{aligned}
    \end{equation}
    By Lemma~\ref{lem:weighted-Q-V_Q-bounds},
    \[
        \norm{V_Q}_{L^\infty}\leq \mathsf V^{L^\infty},
        \qquad
        \norm{V_Q}_{L^1}
        =
        \norm{|V_Q|^{1/2}}_{L^2}^2
        \leq \mathsf V^{L^1},
    \]
    and
    \[
        \norm{\omega^2V_Q}_{L^1}
        =
        \norm{\omega|V_Q|^{1/2}}_{L^2}^2
        \leq \mathsf V_{\omega^2}^{L^1}.
    \]
    Combining these bounds with \eqref{eq:part4-projected-log-transfer-split}, \eqref{eq:part4-projection-term-bound}, and \eqref{eq:part4-TQ-minus-Tg-bound} gives
    \[
        \norm{
            P_{\Phi_Q^\perp}T_QP_{\Phi_Q^\perp}
            -
            P_{\Phi_{\mathfrak g}^\perp}
            T_{\mathfrak g}
            P_{\Phi_{\mathfrak g}^\perp}
        }_{L^2_{\mathrm e}\to L^2_{\mathrm e}}
        \leq
        \Delta_{\mathrm{proj}},
    \]
    which proves \eqref{eq:part4-projected-log-transfer-bound}.
\end{proof}

The right-hand side of \eqref{eq:part4-projected-log-transfer-bound} contains only the bounds from Part~\ref{part:ground-state-approximation} and explicit norms of \(V_{\mathfrak g}\). The following proposition records the resulting outward-rounded bounds.

\begin{proposition}[Bounds for \(T_Q-T_{\mathfrak g}\) and the projected operator difference] \label{prop:part4-certified-logarithmic-transfer-bounds}
    The following operator-norm bounds hold:
    \begin{equation}\label{eq:TQ-Tg-operator-bound}
        \norm{T_Q-T_{\mathfrak g}}_{L^2_{\mathrm e}\to L^2_{\mathrm e}}
        \leq
        5.55\times 10^{-6},
    \end{equation}
    and
    \begin{equation}\label{eq:projected-TQ-Tg-bound}
        \norm{
            P_{\Phi_Q^\perp}T_QP_{\Phi_Q^\perp}
            -
            P_{\Phi_{\mathfrak g}^\perp}
            T_{\mathfrak g}
            P_{\Phi_{\mathfrak g}^\perp}
        }_{L^2_{\mathrm e}\to L^2_{\mathrm e}}
        \leq
        9.64\times 10^{-6}.
    \end{equation}
\end{proposition}

\begin{proof}[Proof of Proposition~\ref{prop:reduced-BS-bound-Q}]
    By Corollary~\ref{cor:part4-reduced-BS-bound-g} and
    Proposition~\ref{prop:part4-certified-logarithmic-transfer-bounds},
    \[
        \norm{
            P_{\Phi_Q^\perp}T_QP_{\Phi_Q^\perp}
        }_{L^2_{\mathrm e}\to L^2_{\mathrm e}}
        \leq
        0.156722+9.64\times10^{-6}
        <0.1571.
    \]
\end{proof}

\subsection{Comparison of the limiting constraint matrices for \(Q\) and \(\mathfrak g\)}\label{subsec:transfer-limiting-matrices-Q}
We first prove the invertibility of \(S_{0,j}\) by treating it as a perturbation of \(S_0^{(\mathfrak g,j)}\). We then compare the entries of the limiting constraint matrices
for \(Q\) and \(\mathfrak g\).

\begin{lemma}[Invertibility of \(S_{0,j}\)]\label{lem:part4-transfer-limiting-Schur-invertibility}
    For \(j=1,2\), the operator \(S_{0,j}\) in \eqref{eq:equation for S_0j} is invertible on
    \(L^2_{\mathrm e}(\bbR)\).
\end{lemma}

\begin{proof}
    By Proposition~\ref{prop:certified-Schur-inverse-bounds-g} and
    Proposition~\ref{prop:part4-certified-logarithmic-transfer-bounds},
    \[
        c_1
        \norm{(S_0^{(\mathfrak g,1)})^{-1}}_{L^2_{\mathrm e}\to L^2_{\mathrm e}}
        \norm{T_Q-T_{\mathfrak g}}_{L^2_{\mathrm e}\to L^2_{\mathrm e}}
        <
        1,
    \]
    and
    \[
        c_2
        \norm{(S_0^{(\mathfrak g,2)})^{-1}}_{L^2_{\mathrm e}\to L^2_{\mathrm e}}
        \norm{T_Q-T_{\mathfrak g}}_{L^2_{\mathrm e}\to L^2_{\mathrm e}}
        <
        1.
    \]
    Since
    \[
        S_{0,j}
        =
        S_0^{(\mathfrak g,j)}
        -
        c_j(T_Q-T_{\mathfrak g})
        =
        S_0^{(\mathfrak g,j)}
        \left[
            I
            -
            c_j
            (S_0^{(\mathfrak g,j)})^{-1}
            (T_Q-T_{\mathfrak g})
        \right],
    \]
    Proposition~\ref{prop:certified-Schur-inverse-bounds-g} and Lemma~\ref{lem:Neumann-series-criterion} give the invertibility of  \(S_{0,j}\), for \(j=1,2\).
\end{proof}

We now compare \(\mathbf m_j\) with \(\mathbf m^{(\mathfrak g,j)}\), and \(\mathbf C_{\log,j}\) with \(\mathbf C_{\log}^{(\mathfrak g,j)}\). These entries do not contain inverse operators.

\begin{lemma}[Bounds for \(\mathbf m_j-\mathbf m^{(\mathfrak g,j)}\) and \(\mathbf C_{\log,j}-\mathbf C_{\log}^{(\mathfrak g,j)}\)] 
    Let \(\Delta_{\Lambda^aQ}^{L^2}\), \(\Delta_{\Lambda^aQ}^{L^1}\), and
    \(\Delta_{\omega\Lambda^aQ}^{L^1}\), \(a=0,1,2\), be the constants defined in
    \eqref{eq:def-Delta-Q-L2}--\eqref{eq:def-Delta-omega-Lambda-a-Q-L1}. For \(F\in L^2(\mathbb R)\cap L^1(\mathbb R)\) such that \(\omega F\in L^1(\mathbb R)\), define
    \begin{equation}\label{eq:part4-free-transfer-Mlog-def}
        \mathsf C_{\log,\mathrm{mix}}(F)
        \coloneqq
        \frac1\pi
        \left(
            2\Delta_Q^{L^2}\norm{F}_{L^2}
            +
            \Delta_Q^{L^1}\norm{\omega F}_{L^1}
            +
            \Delta_{\omega Q}^{L^1}\norm{F}_{L^1}
        \right),
    \end{equation}
    and, for \(p,q\in\{0,1,2\}\), set
    \begin{equation*}
        \mathsf C_{\log,\mathrm{quad}}^{p,q}
        \coloneqq
        \frac1\pi
        \left(
            2\Delta_{\Lambda^pQ}^{L^2} \Delta_{\Lambda^qQ}^{L^2}
            +
            \Delta_{\Lambda^pQ}^{L^1} \Delta_{\omega\Lambda^qQ}^{L^1}
            +
            \Delta_{\omega\Lambda^pQ}^{L^1}\Delta_{\Lambda^qQ}^{L^1}
        \right).
    \end{equation*}
    For \(j=1,2\) and \(m=1,2\), define
    \begin{equation}\label{eq:part4-Delta-m-transfer-def}
        \Delta_{\mathbf m,m}^{(j)}
        \coloneqq
        2^{-(m+j-2)}\Delta_Q^{L^1}.
    \end{equation}
    For \(j=1,2\) and \(1\leq m\leq n\leq2\), set
    \(p\coloneqq m+j-2\) and \(q\coloneqq n+j-2\), and define
    \begin{equation}\label{eq:part4-Delta-Clog-transfer-def}
        \begin{aligned}
            \Delta_{\log,mn}^{(j)}
            \coloneqq{}&
            \mathsf C_{\log,\mathrm{mix}}\left((\Lambda+1)^p\Lambda^q\mathfrak g\right)
            +
            \mathsf C_{\log,\mathrm{mix}}\left((\Lambda+1)^q\Lambda^p\mathfrak g\right)
            +
            \mathsf C_{\log,\mathrm{quad}}^{p,q}
            \\
            &+
            \frac1\pi
            \Delta_Q^{L^1}
            \left[
                \mathbf 1_{\{p>0\}}
                \left|
                    \int_{\bbR}\Lambda^q\mathfrak g\,dy
                \right|
                +
                \mathbf 1_{\{q>0\}}
                \left|
                    \int_{\bbR}\Lambda^p\mathfrak g\,dy
                \right|
            \right].
        \end{aligned}
    \end{equation}
    For \(n<m\), set
    \begin{equation}\label{eq:part4-Delta-Clog-symmetry-def}
        \Delta_{\log,mn}^{(j)}
        \coloneqq
        \Delta_{\log,nm}^{(j)}.
    \end{equation}
    Then, for \(j=1,2\) and \(m=1,2\),
    \begin{equation}\label{eq:part4-transfer-m-vector-bound}
        \left|
            (\mathbf m_j)_m
            -
            (\mathbf m^{(\mathfrak g,j)})_m
        \right|
        \leq
        \Delta_{\mathbf m,m}^{(j)}.
    \end{equation}
    Moreover, for \(j=1,2\) and \(m,n=1,2\),
    \begin{equation}\label{eq:part4-transfer-Clog-matrix-bound}
        \left|
            \left(
                \mathbf C_{\log,j}
                -
                \mathbf C_{\log}^{(\mathfrak g,j)}
            \right)_{mn}
        \right|
        \leq
        \Delta_{\log,mn}^{(j)}.
    \end{equation}
\end{lemma}

\begin{proof}
    We first prove \eqref{eq:part4-transfer-m-vector-bound}. By the
    definitions of the mass vectors, for \(j=1,2\) and \(m=1,2\),
    \[
        (\mathbf m_j)_m
        -
        (\mathbf m^{(\mathfrak g,j)})_m
        =
        \int_{\bbR}\Lambda^{m+j-2}e_Q\,dy .
    \]
    Applying the scaling mass identities
    \eqref{eq:part3-scaling-mass-identities} gives
    \[
        \left|
            (\mathbf m_j)_m
            -
            (\mathbf m^{(\mathfrak g,j)})_m
        \right|
        \leq
        2^{-(m+j-2)}
        \left|
            \int_{\bbR}e_Q\,dy
        \right|.
    \]
    By Lemma~\ref{lem:scaling-and-L1-transfer} and \eqref{eq:def-Delta-Lambda-a-Q-L1},
    \[
        \left|
            \int_{\bbR}e_Q\,dy
        \right|
        \leq
        \norm{e_Q}_{L^1}
        \leq
        \Delta_Q^{L^1}.
    \]
    Hence \eqref{eq:part4-transfer-m-vector-bound} follows from
    \eqref{eq:part4-Delta-m-transfer-def}.

    We next prove \eqref{eq:part4-transfer-Clog-matrix-bound}. It is enough to
    consider \(1\leq m\leq n\leq2\), since the remaining entries follow from
    \eqref{eq:part4-Delta-Clog-symmetry-def} and the symmetry of
    \(\mathcal B_{\log}\). Fix \(j\in\{1,2\}\). For \(1\leq m\leq n\leq2\) set
    \[
        p\coloneqq m+j-2,
        \qquad
        q\coloneqq n+j-2 .
    \]
    By the definitions of \(\mathbf C_{\log,j}\) and
    \(\mathbf C_{\log}^{(\mathfrak g,j)}\), and by bilinearity of
    \(\mathcal B_{\log}\),
    \begin{equation}\label{eq:part4-Clog-transfer-expansion}
        \begin{aligned}
            \left(
                \mathbf C_{\log,j}
                -
                \mathbf C_{\log}^{(\mathfrak g,j)}
            \right)_{mn}
            =
            \mathcal B_{\log}(\Lambda^p e_Q,\Lambda^q\mathfrak g)
            +
            \mathcal B_{\log}(\Lambda^p\mathfrak g,\Lambda^q e_Q)
            +
            \mathcal B_{\log}(\Lambda^p e_Q,\Lambda^q e_Q).
        \end{aligned}
    \end{equation}

    We estimate the first mixed term. If \(p=0\), then
    \eqref{eq:part3-basic-log-pairing-bound}, Lemma~\ref{lem:scaling-and-L1-transfer},
    and \eqref{eq:part4-free-transfer-Mlog-def} give
    \[
        \left|
            \mathcal B_{\log}(e_Q,\Lambda^q\mathfrak g)
        \right|
        \leq
        \mathsf C_{\log,\mathrm{mix}}(\Lambda^q\mathfrak g).
    \]
    If \(p=1\), then \eqref{eq:part3-Blog-Lambda-identity} gives
    \[
        \begin{aligned}
            \left|
                \mathcal B_{\log}(\Lambda e_Q,\Lambda^q\mathfrak g)
            \right|
            \leq
            \left|
                \mathcal B_{\log}
                \left(
                    e_Q,
                    (\Lambda+1)\Lambda^q\mathfrak g
                \right)
            \right|
            +
            \frac1\pi
            \left|
                \int_{\bbR}e_Q\,dy
            \right|
            \left|
                \int_{\bbR}\Lambda^q\mathfrak g\,dy
            \right|.
        \end{aligned}
    \]
    If \(p=2\), then \eqref{eq:part3-Blog-Lambda2-identity} gives the same
    estimate with \((\Lambda+1)^2\Lambda^q\mathfrak g\) in place of
    \((\Lambda+1)\Lambda^q\mathfrak g\). Therefore, for \(p=0,1,2\),
    \begin{equation}\label{eq:part4-Clog-first-mixed-bound}
        \begin{aligned}
            \left|
                \mathcal B_{\log}(\Lambda^p e_Q,\Lambda^q\mathfrak g)
            \right|
            \leq
            \mathsf C_{\log,\mathrm{mix}}\left((\Lambda+1)^p\Lambda^q\mathfrak g\right)
            +
            \frac1\pi
            \mathbf 1_{\{p>0\}}
            \Delta_Q^{L^1}
            \left|
                \int_{\bbR}\Lambda^q\mathfrak g\,dy
            \right|.
        \end{aligned}
    \end{equation}
    Here we used
    \[
        \left|
            \int_{\bbR}e_Q\,dy
        \right|
        \leq
        \norm{e_Q}_{L^1}
        \leq
        \Delta_Q^{L^1}.
    \]

    For the second mixed term, the symmetry of \(\mathcal B_{\log}\) gives
    \[
        \mathcal B_{\log}(\Lambda^p\mathfrak g,\Lambda^q e_Q)
        =
        \mathcal B_{\log}(\Lambda^q e_Q,\Lambda^p\mathfrak g).
    \]
    Applying \eqref{eq:part4-Clog-first-mixed-bound} with \(p\) and \(q\)
    interchanged yields
    \begin{equation*}
        \begin{aligned}
            \left|
                \mathcal B_{\log}(\Lambda^p\mathfrak g,\Lambda^q e_Q)
            \right|
            \leq
            \mathsf C_{\log,\mathrm{mix}}\left((\Lambda+1)^q\Lambda^p\mathfrak g\right)
            +
            \frac1\pi
            \mathbf 1_{\{q>0\}}
            \Delta_Q^{L^1}
            \left|
                \int_{\bbR}\Lambda^p\mathfrak g\,dy
            \right|.
        \end{aligned}
    \end{equation*}

    Finally, applying \eqref{eq:part3-basic-log-pairing-bound} to the last term
    in \eqref{eq:part4-Clog-transfer-expansion}, and using
    Lemma~\ref{lem:scaling-and-L1-transfer} together with
    \eqref{eq:def-Delta-Q-L2}--\eqref{eq:def-Delta-omega-Lambda-a-Q-L1}, gives
    \begin{equation}\label{eq:part4-Clog-error-error-bound}
        \left|
            \mathcal B_{\log}(\Lambda^p e_Q,\Lambda^q e_Q)
        \right|
        \leq
        \mathsf C_{\log,\mathrm{quad}}^{p,q}.
    \end{equation}
    Combining \eqref{eq:part4-Clog-transfer-expansion}--\eqref{eq:part4-Clog-error-error-bound} with \eqref{eq:part4-Delta-Clog-transfer-def}, we obtain \eqref{eq:part4-transfer-Clog-matrix-bound} for \(m\leq n\). The case \(n<m\) follows from symmetry.
\end{proof}

We next compare the entries containing \(S_{0,j}^{-1}\) with the corresponding entries containing \((S_0^{(\mathfrak g,j)})^{-1}\).

\begin{lemma}[Perturbation bounds for \(\mathfrak s_j(0^+)\), \(\mathbf a_j(0^+)\), and \(\mathbf R_j(0^+)\)]
    Set
    \begin{equation}\label{eq:part4-DeltaT-and-Neumann-factor-def}
        \Delta_T
        \coloneqq
        \norm{
            T_Q-T_{\mathfrak g}
        }_{L^2_{\mathrm e}(\bbR)\to L^2_{\mathrm e}(\bbR)}.
    \end{equation}
    For \(j=1,2\), set
    \begin{equation}\label{eq:inverse-transfer-factor-definition}
        \mathsf C_{S,\mathrm{inv}}^{(\mathfrak g,j)}
        \coloneqq
        \norm{
            (S_0^{(\mathfrak g,j)})^{-1}
        }_{L^2_{\mathrm e}(\bbR)\to L^2_{\mathrm e}(\bbR)},
        \qquad
        \mathsf C_{\mathrm{transfer}}^{(j)}
        \coloneqq
        \left(
            1-c_j\mathsf C_{S,\mathrm{inv}}^{(\mathfrak g,j)}\Delta_T
        \right)^{-1}.
    \end{equation}
    Let \(\Delta_{\psi,m}^{(j)}\) be the constants defined in
    \eqref{eq:app-Deltapsi-def}. Define
    \begin{equation}\label{eq:part4-Delta-s-transfer-def}
        \Delta_{\mathfrak s}^{(j)}
        \coloneqq
        \mathsf C_{\mathrm{transfer}}^{(j)}
        \left[
            \mathsf C_{S,\mathrm{inv}}^{(\mathfrak g,j)}
            \bigl(\Delta_{\sqrt V}^{L^2}\bigr)^2
            +
            2\Delta_{\sqrt V}^{L^2}
            \norm{X_0^{(\mathfrak g,j)}}_{L^2(I)}
            +
            c_j\Delta_T
            \norm{X_0^{(\mathfrak g,j)}}_{L^2(I)}^2
        \right],
    \end{equation}
    for \(m=1,2\),
    \begin{equation}\label{eq:part4-Delta-a-transfer-def}
    \begin{aligned}
        \Delta_{\mathfrak a,m}^{(j)}
        \coloneqq
        \mathsf C_{\mathrm{transfer}}^{(j)}
        \bigg[
        &
            \Delta_{\psi,m}^{(j)}
            \left(
                \mathsf C_{S,\mathrm{inv}}^{(\mathfrak g,j)}\Delta_{\sqrt V}^{L^2}
                +
                \norm{X_0^{(\mathfrak g,j)}}_{L^2(I)}
            \right)
            \\
        &+
            \Delta_{\sqrt V}^{L^2}
            \norm{X_m^{(\mathfrak g,j)}}_{L^2(I)}
            +
            c_j\Delta_T
            \norm{X_0^{(\mathfrak g,j)}}_{L^2(I)}
            \norm{X_m^{(\mathfrak g,j)}}_{L^2(I)}
        \bigg],
    \end{aligned}
    \end{equation}
    and, for \(m,n=1,2\),
    \begin{equation}\label{eq:part4-Delta-r-transfer-def}
    \begin{aligned}
        \Delta_{\mathfrak r,mn}^{(j)}
        \coloneqq
        \mathsf C_{\mathrm{transfer}}^{(j)}
        \bigg[
        &
            \Delta_{\psi,m}^{(j)}
            \left(
                \mathsf C_{S,\mathrm{inv}}^{(\mathfrak g,j)}
                \Delta_{\psi,n}^{(j)}
                +
                \norm{X_n^{(\mathfrak g,j)}}_{L^2(I)}
            \right)
            \\
        &+
            \norm{X_m^{(\mathfrak g,j)}}_{L^2(I)}
            \Delta_{\psi,n}^{(j)}
            +
            c_j\Delta_T
            \norm{X_m^{(\mathfrak g,j)}}_{L^2(I)}
            \norm{X_n^{(\mathfrak g,j)}}_{L^2(I)}
        \bigg].
    \end{aligned}
    \end{equation}
    Then, for \(j=1,2\),
    \begin{equation}\label{eq:part4-transfer-s-bound}
        \left|
            \mathfrak s_j(0^+)
            -
            \mathfrak s_0^{(\mathfrak g,j)}
        \right|
        \leq
        \Delta_{\mathfrak s}^{(j)}.
    \end{equation}
    Moreover, for \(m=1,2\),
    \begin{equation}\label{eq:part4-transfer-a-bound}
        \left|
            (\mathbf a_j(0^+))_m
            -
            \mathfrak a_m^{(\mathfrak g,j)}
        \right|
        \leq
        \Delta_{\mathfrak a,m}^{(j)},
    \end{equation}
    and, for \(m,n=1,2\),
    \begin{equation}\label{eq:part4-transfer-r-bound}
        \left|
            (\mathbf R_j(0^+))_{mn}
            -
            \mathfrak r_{mn}^{(\mathfrak g,j)}
        \right|
        \leq
        \Delta_{\mathfrak r,mn}^{(j)}.
    \end{equation}
\end{lemma}

\begin{proof}
    By Lemma~\ref{lem:part4-transfer-limiting-Schur-invertibility}, the
    denominator in \eqref{eq:inverse-transfer-factor-definition} is positive. For \(f\in L^2_{\mathrm e}(\bbR)\), Lemma~\ref{lem:Neumann-series-criterion} gives
    \begin{equation}\label{eq:part4-SQ-vector-transfer-bound}
        \norm{
            S_{0,j}^{-1}f
        }_{L^2}
        \leq
        \mathsf C_{\mathrm{transfer}}^{(j)}
        \norm{
            (S_0^{(\mathfrak g,j)})^{-1}f
        }_{L^2}.
    \end{equation}
    By unitary equivalence and the definition of \(X_m^{(\mathfrak g,j)}\),
    \begin{equation}\label{eq:part4-X-identities}
        \norm{
            (S_0^{(\mathfrak g,j)})^{-1}|V_{\mathfrak g}|^{1/2}
        }_{L^2(\bbR)}
        =
        \norm{X_0^{(\mathfrak g,j)}}_{L^2(I)},
        \qquad
        \norm{
            (S_0^{(\mathfrak g,j)})^{-1}\psi_m^{(\mathfrak g,j)}
        }_{L^2(\bbR)}
        =
        \norm{X_m^{(\mathfrak g,j)}}_{L^2(I)}.
    \end{equation}
    For \(m=0,1,2\), define 
    \[
        \zeta_0^{(j)}\coloneqq |V_Q|^{1/2},
        \qquad
        \zeta_m^{(j)}\coloneqq \psi_m^{(j)}
        \quad (m=1,2),
    \]
    and
    \[
        \zeta_0^{(\mathfrak g,j)}
        \coloneqq |V_{\mathfrak g}|^{1/2},
        \qquad
        \zeta_m^{(\mathfrak g,j)}
        \coloneqq \psi_m^{(\mathfrak g,j)}
        \quad (m=1,2).
    \]
    Also set
    \[
        \Delta_{\zeta,0}^{(j)}\coloneqq \Delta_{\sqrt V}^{L^2},
        \qquad
        \Delta_{\zeta,m}^{(j)}\coloneqq \Delta_{\psi,m}^{(j)}
        \quad (m=1,2).
    \]
    By Lemma~\ref{lem:square-root-potential-transfer} and
    Lemma~\ref{lem:app-psi-difference-estimates},
    \begin{equation}\label{eq:part4-zeta-transfer-bound}
        \norm{
            \zeta_m^{(j)}
            -
            \zeta_m^{(\mathfrak g,j)}
        }_{L^2}
        \leq
        \Delta_{\zeta,m}^{(j)},
        \qquad
        m=0,1,2.
    \end{equation}
    For \(m,n\in\{0,1,2\}\), we have
    \begin{equation}\label{eq:part4-unified-inverse-dependent-transfer}
        \begin{aligned}
            &\left|
                \left(
                    \zeta_m^{(j)},
                    S_{0,j}^{-1}\zeta_n^{(j)}
                \right)_r
                -
                \left(
                    \zeta_m^{(\mathfrak g,j)},
                    (S_0^{(\mathfrak g,j)})^{-1}
                    \zeta_n^{(\mathfrak g,j)}
                \right)_r
            \right|
            \\
            &
            \leq
            \mathsf C_{\mathrm{transfer}}^{(j)}
            \bigg[
                \Delta_{\zeta,m}^{(j)}
                \left(
                    \mathsf C_{S,\mathrm{inv}}^{(\mathfrak g,j)}\Delta_{\zeta,n}^{(j)}
                    +
                    \norm{X_n^{(\mathfrak g,j)}}_{L^2(I)}
                \right)
                +
                \norm{X_m^{(\mathfrak g,j)}}_{L^2(I)}
                \Delta_{\zeta,n}^{(j)}
                +
                c_j\Delta_T
                \norm{X_m^{(\mathfrak g,j)}}_{L^2(I)}
                \norm{X_n^{(\mathfrak g,j)}}_{L^2(I)}
            \bigg].
        \end{aligned}
    \end{equation}
    Adding and subtracting \(\left(\zeta_m^{(\mathfrak g,j)},S_{0,j}^{-1}\zeta_n^{(j)}
    \right)_r\) and \(\left(\zeta_m^{(\mathfrak g,j)},S_{0,j}^{-1}\zeta_n^{(\mathfrak g,j)}
    \right)_r\) gives
    \[
        \begin{aligned}
            \left(
                \zeta_m^{(j)},
                S_{0,j}^{-1}\zeta_n^{(j)}
            \right)_r
            &-
            \left(
                \zeta_m^{(\mathfrak g,j)},
                (S_0^{(\mathfrak g,j)})^{-1}
                \zeta_n^{(\mathfrak g,j)}
            \right)_r
            \\
            &=
            \left(
                \zeta_m^{(j)}
                -
                \zeta_m^{(\mathfrak g,j)},
                S_{0,j}^{-1}\zeta_n^{(j)}
            \right)_r
            +
            \left(
                \zeta_m^{(\mathfrak g,j)},
                S_{0,j}^{-1}
                \left(
                    \zeta_n^{(j)}
                    -
                    \zeta_n^{(\mathfrak g,j)}
                \right)
            \right)_r
            \\
            &\quad
            +
            \left(
                \zeta_m^{(\mathfrak g,j)},
                \left[
                    S_{0,j}^{-1}
                    -
                    (S_0^{(\mathfrak g,j)})^{-1}
                \right]
                \zeta_n^{(\mathfrak g,j)}
            \right)_r .
        \end{aligned}
    \]
    For the first term, \eqref{eq:part4-SQ-vector-transfer-bound} gives
    \[
        \norm{
            S_{0,j}^{-1}\zeta_n^{(j)}
        }_{L^2}
        \leq
        \mathsf C_{\mathrm{transfer}}^{(j)}
        \left(
            \mathsf C_{S,\mathrm{inv}}^{(\mathfrak g,j)}\Delta_{\zeta,n}^{(j)}
            +
            \norm{X_n^{(\mathfrak g,j)}}_{L^2(I)}
        \right),
    \]
    where we used \eqref{eq:part4-zeta-transfer-bound} and
    \eqref{eq:part4-X-identities}. Hence the first term is bounded by the first
    term on the right-hand side of
    \eqref{eq:part4-unified-inverse-dependent-transfer}.

    For the second term, self-adjointness of \(S_{0,j}^{-1}\), together with
    \eqref{eq:part4-SQ-vector-transfer-bound} and
    \eqref{eq:part4-X-identities}, gives
    \[
        \left|
            \left(
                \zeta_m^{(\mathfrak g,j)},
                S_{0,j}^{-1}
                \left(
                    \zeta_n^{(j)}
                    -
                    \zeta_n^{(\mathfrak g,j)}
                \right)
            \right)_r
        \right|
        \leq
        \mathsf C_{\mathrm{transfer}}^{(j)}
        \norm{X_m^{(\mathfrak g,j)}}_{L^2(I)}
        \Delta_{\zeta,n}^{(j)} .
    \]
    Finally, the resolvent identity gives
    \[
        S_{0,j}^{-1}
        -
        (S_0^{(\mathfrak g,j)})^{-1}
        =
        S_{0,j}^{-1}
        c_j(T_Q-T_{\mathfrak g})
        (S_0^{(\mathfrak g,j)})^{-1}.
    \]
    Therefore, by self-adjointness, \eqref{eq:part4-SQ-vector-transfer-bound},
    \eqref{eq:part4-X-identities}, and the definition of \(\Delta_T\),
    \[
        \begin{aligned}
            \left|
                \left(
                    \zeta_m^{(\mathfrak g,j)},
                    \left[
                        S_{0,j}^{-1}
                        -
                        (S_0^{(\mathfrak g,j)})^{-1}
                    \right]
                    \zeta_n^{(\mathfrak g,j)}
                \right)_r
            \right|
            \leq
            \mathsf C_{\mathrm{transfer}}^{(j)}
            c_j\Delta_T
            \norm{X_m^{(\mathfrak g,j)}}_{L^2(I)}
            \norm{X_n^{(\mathfrak g,j)}}_{L^2(I)} .
        \end{aligned}
    \]
    This proves \eqref{eq:part4-unified-inverse-dependent-transfer}.

    Taking \((m,n)=(0,0)\) in
    \eqref{eq:part4-unified-inverse-dependent-transfer} gives
    \eqref{eq:part4-transfer-s-bound}.  Taking \((m,n)=(m,0)\) with
    \(m=1,2\) gives \eqref{eq:part4-transfer-a-bound}. Finally, taking
    \(m,n=1,2\) gives \eqref{eq:part4-transfer-r-bound}.
\end{proof}

The bounds above and Proposition~\ref{prop:part4-certified-g-side-limiting-data} reduce the entries for \(Q\) to explicit interval expressions. The following proposition records the resulting outward-rounded enclosures.

\begin{proposition}[Interval enclosures for the limiting constraint matrices]
\label{prop:part4-certified-limiting-constraint-matrices}
    The scalars \(\mathfrak s_1(0^+)\) and \(\mathfrak s_2(0^+)\)
    satisfy
    \begin{equation}\label{eq:part4-certified-limiting-s-intervals}
        \mathfrak s_1(0^+)\in[-0.05298,-0.05172],
        \qquad
        \mathfrak s_2(0^+)\in[3.72581,3.72786].
    \end{equation}
    Moreover, the matrices \(\mathbf M_{0,j}\) in \eqref{eq:limiting constraint matrix} satisfy
    \begin{equation}\label{eq:part4-certified-M0-j1-interval}
        \mathbf M_{0,1}
        \in
        \begin{pmatrix}
            [1.025,1.153]
            &
            [-4.905,-4.636]
            \\
            [-4.905,-4.636]
            &
            [13.47,14.02]
        \end{pmatrix},
    \end{equation}
    and
    \begin{equation}\label{eq:part4-certified-M0-j2-interval}
        \mathbf M_{0,2}
        \in
        \begin{pmatrix}
            [0.431,0.441]
            &
            [-0.879,-0.864]
            \\
            [-0.879,-0.864]
            &
            [1.151,1.173]
        \end{pmatrix}.
    \end{equation}
\end{proposition}

\begin{proof}[Proof of Proposition~\ref{prop:certified-limiting-matrix-conditions}]
    By Lemma~\ref{lem:part4-transfer-limiting-Schur-invertibility}, the
    operators \(S_{0,j}\), \(j=1,2\), are invertible on
    \(L^2_{\mathrm e}(\bbR)\).  The scalar enclosures
    \eqref{eq:part4-certified-limiting-s-intervals} give
    \[
        \mathfrak s_1(0^+)<0,
        \qquad
        \mathfrak s_2(0^+)>0.
    \]
    Hence \(\mathfrak s_j(0^+)\neq0\) for \(j=1,2\), and the limiting matrix
    formula in Proposition~\ref{prop:constraint-entry-formula} is well-defined.
    It remains to verify the determinant signs.  Write
    \[
        \mathbf M_{0,j}
        =
        \begin{pmatrix}
            A_j & B_j\\
            B_j & D_j
        \end{pmatrix}.
    \]
    Every matrix in the enclosure \eqref{eq:part4-certified-M0-j1-interval} satisfies
    \[
        A_1>0,\qquad D_1>0,\qquad B_1<0.
    \]
    Therefore
    \[
        \det\mathbf M_{0,1}
        =
        A_1D_1-B_1^2
        \leq
        1.153\cdot14.02-4.636^2
        <
        -5.32
        <
        0.
    \]
    Similarly, from \eqref{eq:part4-certified-M0-j2-interval} we have
    \[
        A_2>0,\qquad D_2>0,\qquad B_2<0,
    \]
    and hence
    \[
        \det\mathbf M_{0,2}
        =
        A_2D_2-B_2^2
        \leq
        0.441\cdot1.173-0.864^2
        <
        -0.229
        <
        0.
    \]
    Thus
    \[
        \det\mathbf M_{0,1}<0,
        \qquad
        \det\mathbf M_{0,2}<0.
    \]
    Each \(\mathbf M_{0,j}\) is a real symmetric \(2\times2\) matrix.
    A negative determinant implies that its two eigenvalues have opposite signs.
    Therefore
    \[
        n_-(\mathbf M_{0,1})
        =
        n_-(\mathbf M_{0,2})
        =
        1.
    \]
    This proves Proposition~\ref{prop:certified-limiting-matrix-conditions}.
\end{proof}

\part{Interval-arithmetic verification}\label{part:interval-verification}
All auxiliary vectors, matrices, and scalars introduced in Part~\ref{part:interval-verification} are used only in the subsection in which they are defined, unless stated otherwise. This part verifies the numerical bounds used in Parts~\ref{part:ground-state-approximation} and~\ref{part:spectral-estimates}. The verification code, the ancillary files referenced below, and the completed outputs are available at \url{https://github.com/JeongheonPark-CAP/half-wave-cap-verification}. The ancillary files contain the exact rational coefficients and auxiliary vectors used in the interval computations. All certified numerical bounds in this part are obtained by directed-rounding interval arithmetic using the full-precision input data. The decimal values displayed in intermediate estimates are abbreviated for readability and are not intended to reproduce subsequent certified bounds by direct arithmetic.

\section{Arithmetic conventions and numerical inputs}\label{sec:interval-inputs}

Decimal strings in the ancillary files are interpreted as exact rational numbers. All subsequent evaluations use interval arithmetic with outward rounding. Hence every computed interval contains the corresponding exact real quantity. In each displayed one-sided decimal
bound, an upper endpoint is rounded upward and a lower endpoint is rounded downward. A sign or inequality is used only when it follows from interval inclusion. 

Throughout Part~\ref{part:interval-verification}, we use
\begin{equation*}
    J=200,\qquad
    L=5000,\qquad
    M=15000,\qquad
    N_{\mathrm{DCT}}=65536.
\end{equation*}
Here \(J\) is the degree used in the profile definition
\eqref{eq:def-U-profile}, \(L\) is the degree of the cosine-polynomial approximation
in \eqref{eq:definition of a_g^L}, \(M\) is the dimension cutoff used in the finite-dimensional Schur-complement computation, and \(N_{\mathrm{DCT}}\) is the grid size used for the discrete cosine transform.

The endpoint Birman--Schwinger computations use
\begin{equation}\label{eq:part5-endpoint-parameters}
    Y=30,\qquad
    h=5\times10^{-3},\qquad
    m=1,\qquad
    n_{\mathrm{GL}}=8,\qquad
    N_{\mathrm{end}}=6000.
\end{equation}
The shifted constants are fixed as
\begin{equation*}
    c_1=6.1,
    \qquad
    c_2=2.1.
\end{equation*}
The parameters used in the weighted estimates are fixed as
\begin{equation}\label{eq:part5-transfer-splitting-parameters}
    R_E=5\times10^8,
    \qquad
    R_\omega=5\times10^8,
    \qquad
    R_Q=100,
    \qquad
    R_M=40.
\end{equation}
In Definition~\ref{def:weighted-q-constants}, take \(R_Q=100\) and \(R_M=40\). In Definition~\ref{def:dilation-L1-transfer-radii}, take \(R_\omega=R_E=5\times10^8\).

The profile coefficients \(\{U_k\}_{k=0}^{J}\) are those fixed in
Definition~\ref{def:approximate-profile-g}. The coefficients
\(\{\mathsf c_k\}_{k=0}^{L}\) in the cosine-polynomial approximation \eqref{eq:definition of a_g^L} are read from
\[
    \texttt{ag\_coeffs\_L5000\_N65536\_dec15\_cert5e-15.csv}.
\]
This file has header \(\texttt{k,c\_dec15}\) and contains \(L+1=5001\)
entries, indexed by \(0\leq k\leq L\). The coefficient comparison gives
\begin{equation}\label{eq:part5-ag-coefficient-comparison}
    \left|
        \mathsf c_k-D_{N_{\mathrm{DCT}},k}[a_{\mathfrak g}]
    \right|
    \leq
    5\times10^{-15},
    \qquad
    0\leq k\leq L.
\end{equation}
The bound for \(\norm{ a_{\mathfrak g}-a_{\mathfrak g}^{(L)} }_{L^2(I)}\) is obtained from these coefficients by \eqref{eq:part5-agL-error-decomposition}--\eqref{eq:part5-agL-rational-error}.

The interval computations are organized as follows.
\begin{table}[htbp]
\centering
\begin{tabular}{c|l}
\textnormal{computation} & \textnormal{bounds established} \\
\hline
01 & bounds for \(\mathfrak g\), its residual, and its derivatives \\
02 & endpoint Birman--Schwinger eigenvalue bounds for \(\mathcal F'(\mathfrak g)\) \\
03 & bounds for \(Q-\mathfrak g\) and \(T_Q-T_{\mathfrak g}\) \\
04 & coefficient comparison for \(a_{\mathfrak g}^{(L)}\) \\
05 & projected Birman--Schwinger bound for \(\mathfrak g\) \\
05b & projected Birman--Schwinger perturbation bound for \(Q\) \\
06 & Schur-complement signs and inverse bounds for \(S_0^{(\mathfrak g,j)}\) \\
07 & limiting matrix entries for \(\mathfrak g\) \\
08 & limiting matrix entries for \(Q\) and determinant signs
\end{tabular}
\caption{Interval computations used in Part~\ref{part:interval-verification}.}
\label{tab:part5-certificate-layout}
\end{table}

The ancillary input files for each computation are listed in the metadata file for that computation. Before an interval computation is performed, the file headers and the numbers of rows are checked.

\section{Profile bounds and estimates for \(Q-\mathfrak g\)}

\subsection{Bounds for \(\mathfrak g\) and its residual}

\begin{proposition}[Bounds for \(\mathfrak g\) and its residual]\label{prop:part5-g-profile-residual-bounds}
    Let \(\mathfrak g\) be the approximate profile in
    Definition~\ref{def:approximate-profile-g}. Set
    \[
        r_{\mathfrak g}\coloneqq \mathcal F(\mathfrak g),
        \qquad
        V_{\mathfrak g}\coloneqq y\mathfrak g'\mathfrak g .
    \]
    All vector and matrix inequalities below are understood componentwise. The following bounds hold.

    \begin{enumerate}[label=\textup{(\roman*)}]
        \item
        The residual norms satisfy
        \begin{equation}\label{eq:part5-g-residual bounds}
            \left(
                \norm{r_{\mathfrak g}}_{L^2},
                \norm{r_{\mathfrak g}}_{W^{1,2}},
                \norm{r_{\mathfrak g}}_{W^{2,2}}
            \right)
            \leq
            \left(
                2.51\times10^{-14},
                7.43\times10^{-12},
                2.45\times10^{-9}
            \right).
        \end{equation}
        The scaling residuals satisfy
        \begin{equation}\label{eq:g-scaling-residual-bounds}
            \left(
                \norm{r_{\mathfrak g}}_{L^2},
                \norm{\Lambda r_{\mathfrak g}}_{L^2},
                \norm{\Lambda^2 r_{\mathfrak g}}_{L^2}
            \right)
            \leq
            \left(
                2.51\times10^{-14},
                2.18\times10^{-12},
                3.13\times10^{-10}
            \right).
        \end{equation}

        \item
        The Sobolev norms of \(\mathfrak g\) satisfy
        \begin{equation}\label{eq:g-profile-norm-bounds}
            \left(
                \norm{\mathfrak g}_{L^2},
                \norm{\mathfrak g}_{H^1},
                \norm{\mathfrak g}_{W^{2,2}}
            \right)
            \leq
            \left(
                1.572,
                2.917,
                12.547
            \right).
        \end{equation}
        Define the vector of pointwise bounds
        \[
            \mathsf P_{\mathfrak g}
            \coloneqq
            \left(
                \norm{\mathfrak g}_{L^\infty},
                \norm{\mathfrak g'}_{L^\infty},
                \norm{\langle y\rangle^2\mathfrak g}_{L^\infty},
                \norm{\langle y\rangle^2 y\mathfrak g'}_{L^\infty},
                \norm{\langle y\rangle^2 y^2\mathfrak g''}_{L^\infty}
            \right).
        \]
        Then
        \begin{equation}\label{eq:g-pointwise-bounds}
            \mathsf P_{\mathfrak g}
            \leq
            \left(
                1.891,
                2.963,
                1.891,
                2.068,
                6.204
            \right).
        \end{equation}
        In addition, \(\mathfrak g(y)>0\) for all \(y\in\bbR\), \(\mathfrak g'(y)<0\) for \(y>0\), and \(\mathfrak g(0)>\frac12\).

        \item
        For \(p=0,1,2\), define
        \[
            \bigl(\mathsf B_{\mathfrak g}\bigr)_{p,\bullet}
            \coloneqq
            \left(
                \norm{\Lambda^p\mathfrak g}_{L^1},
                \norm{\Lambda^p\mathfrak g}_{L^2},
                \norm{\omega\Lambda^p\mathfrak g}_{L^1}
            \right).
        \]
        Then
        \begin{equation}\label{eq:g-source-basis-bounds}
            \mathsf B_{\mathfrak g}
            \leq
            \begin{pmatrix}
                3.249 & 1.572 & 2.684 \\
                2.126 & 0.666 & 2.814 \\
                1.836 & 0.586 & 2.930
            \end{pmatrix},
        \end{equation}
        where the rows correspond to \(p=0,1,2\).

        \item
        For
        \[
            (p,q)\in
            \{(0,1),(1,1),(2,1),(1,2),(2,2)\},
        \]
        define
        \[
            \bigl(\mathsf C_{\mathfrak g}\bigr)_{(p,q),\bullet}
            \coloneqq
            \left(
                \norm{(\Lambda+1)^q\Lambda^p\mathfrak g}_{L^1},
                \norm{(\Lambda+1)^q\Lambda^p\mathfrak g}_{L^2},
                \norm{\omega(\Lambda+1)^q\Lambda^p\mathfrak g}_{L^1}
            \right).
        \]
        Then
        \begin{equation}\label{eq:g-clog-source-bounds}
            \mathsf C_{\mathfrak g}
            \leq
            \begin{pmatrix}
                2.038 & 1.707 & 0.970 \\
                1.616 & 0.887 & 1.241 \\
                1.870 & 0.983 & 1.694 \\
                1.646 & 1.324 & 0.972 \\
                2.596 & 1.975 & 1.596
            \end{pmatrix},
        \end{equation}
        where the rows are ordered as
        \[
            (0,1),\ (1,1),\ (2,1),\ (1,2),\ (2,2).
        \]

        \item
        The potential norms satisfy
        \begin{equation}\label{eq:g-potential-bounds}
            \left(
                \norm{V_{\mathfrak g}}_{L^1},
                \norm{V_{\mathfrak g}}_{L^\infty},
                \norm{\omega^2 V_{\mathfrak g}}_{L^1},
                \norm{\omega |V_{\mathfrak g}|^{1/2}}_{L^2}
            \right)
            \leq
            \left(
                1.2350,
                0.9264,
                0.4213,
                0.6491
            \right).
        \end{equation}

        \item
        Define the vector of moments and weighted norms
        \[
            \mathsf D_{\mathfrak g}
            \coloneqq
            \left(
                \int_{\bbR}\mathfrak g^3,
                \int_{\bbR}y^2\mathfrak g^3,
                \norm{y\mathfrak g}_{L^2},
                \norm{y\mathfrak g'}_{L^2},
                \norm{\omega^2 y\mathfrak g}_{L^2},
                \norm{\omega^2 y\mathfrak g'}_{L^2}
            \right).
        \]
        Then
        \begin{equation}\label{eq:g-direct-moment-bounds}
            \mathsf D_{\mathfrak g}
            \leq
            \left(
                3.2485,
                0.3004,
                1.1068,
                1.0298,
                7.7325,
                0.9522
            \right).
        \end{equation}
    \end{enumerate}
\end{proposition}

\begin{proof}
    The coefficients of \(\mathfrak g\) are the exact rational numbers fixed in Definition~\ref{def:approximate-profile-g}. Applying interval arithmetic to these coefficients gives the stated bounds; the computation is listed in row~01 of Table~\ref{tab:part5-certificate-layout}. After the change of variables \(y=\tan\theta\), the functions \(r_{\mathfrak g}\), \(\Lambda^p r_{\mathfrak g}\), \(\Lambda^p\mathfrak g\), and \((\Lambda+1)^q\Lambda^p\mathfrak g\) are represented by trigonometric polynomials with interval coefficients.

    The \(L^2\)-based bounds in \eqref{eq:part5-g-residual bounds}, \eqref{eq:g-scaling-residual-bounds}, \eqref{eq:g-profile-norm-bounds}, \eqref{eq:g-source-basis-bounds}, and \eqref{eq:g-clog-source-bounds} are computed using Parseval's identity and interval arithmetic. The \(L^1\)- and \(\omega L^1\)-bounds in \eqref{eq:g-source-basis-bounds} and \eqref{eq:g-clog-source-bounds} are computed after the change of variables \(y=\tan\theta\), using evenness.

    The pointwise bounds in \eqref{eq:g-pointwise-bounds} are obtained by subdividing the compactified interval and enclosing the polynomial or rational expression defining each norm on every subinterval. The subdivision also verifies the positivity of \(\mathfrak g\) and the strict inequality \(\mathfrak g'(y)<0\) for \(y>0\).

    The norms in \eqref{eq:g-potential-bounds} and \eqref{eq:g-direct-moment-bounds} are evaluated after the change of variables \(y=\tan\theta\), using evenness. For the terms containing \(\omega\), we split \(I\) into an interior interval and two endpoint intervals. Interval quadrature bounds the interior integrals, while Lemma~\ref{lem:app-elementary-log-tail-bounds} bounds the endpoint integrals.
\end{proof}

\subsection{Endpoint Birman--Schwinger bounds}

Throughout this subsection we use, the endpoint parameters fixed in
\eqref{eq:part5-endpoint-parameters}.  Set
\begin{equation}\label{eq:part5-endpoint-matrix-short-notation}
    \mathbf A_\mu
    \coloneqq
    \mathbf K_{\mathrm{NK},Y,h,m,n_{\mathrm{GL}}}^{(\mu),\mathrm{eval}},
    \qquad
    \mu\in\left\{\frac12,\frac32\right\}.
\end{equation}
Let
\[
    \mathbf v_1,\mathbf v_2\in\mathbb R^{N_{\mathrm{end}}}
\]
be the vectors whose entries are read as exact rational numbers from
\[
    \texttt{ritz\_basis\_v1\_N6000.csv},
    \qquad
    \texttt{ritz\_basis\_v2\_N6000.csv}.
\]
Let \(\widetilde{\mathbf v}\in\mathbb R^{N_{\mathrm{end}}}\) be the vector read
from
\[
    \texttt{rayleigh\_vector\_v1\_N6000.csv},
\]
and let \(\mathbf w_{\mathrm{CW}}\in(0,\infty)^{N_{\mathrm{end}}}\) be the
positive Collatz--Wielandt weight read from
\[
    \texttt{cw\_weight\_complement\_N6000.csv}.
\]
Define $\mathcal V_{\mathrm R}\coloneqq \operatorname{span}\{\mathbf v_1,\mathbf v_2\}$ and
\[
    \mathbf P_{\mathrm R}
    \coloneqq
    \textnormal{the Euclidean orthogonal projection onto }\mathcal V_{\mathrm R},
    \qquad
    \mathbf Q_{\mathrm R}
    \coloneqq
    I-\mathbf P_{\mathrm R}.
\]
Define \(\mathbf G,\mathbf B\in\mathbb R^{2\times2}\) by
\begin{equation*}
    G_{ij}
    \coloneqq
    (\mathbf v_i,\mathbf v_j)_{\mathbb R^{N_{\mathrm{end}}}},
    \qquad
    B_{ij}
    \coloneqq
    (\mathbf A_{1/2}\mathbf v_i,\mathbf v_j)_{\mathbb R^{N_{\mathrm{end}}}},
    \qquad
    1\leq i,j\leq2.
\end{equation*}
Let \(\lambda_{\mathrm{Ritz},+}^{(1/2)}\) be the larger generalized eigenvalue
of
\begin{equation}\label{eq:part5-endpoint-ritz-generalized-eigenvalue}
    \mathbf B\xi=\lambda\mathbf G\xi .
\end{equation}
Let
\[
    \mathbf B_{\perp}^{(1/2)}
    \coloneqq
    \mathbf Q_{\mathrm R}\mathbf A_{1/2}\mathbf Q_{\mathrm R}
    \quad
    \textnormal{on }\mathbf Q_{\mathrm R}\mathbb R^{N_{\mathrm{end}}}.
\]
The Collatz--Wielandt estimate is not applied to \(\mathbf B_{\perp}^{(1/2)}\) itself, but to an entrywise nonnegative majorant \(\mathbf C_{\mathrm{CW}}^{(1/2)}\) of its coefficient matrix:
\[
    0\le
    \left|\mathbf B_{\perp}^{(1/2)}\right|_{\mathrm{entry}}
    \le
    \mathbf C_{\mathrm{CW}}^{(1/2)}.
\]
With the positive weight \(\mathbf w_{\mathrm{CW}}\), set
\[
    \Gamma_{\mathrm{CW}}^{(1/2)}
    \coloneqq
    \max_i
    \frac{
        \bigl(\mathbf C_{\mathrm{CW}}^{(1/2)}
        \mathbf w_{\mathrm{CW}}\bigr)_i
    }{
        (\mathbf w_{\mathrm{CW}})_i
    }.
\]
For the Kato--Temple estimate, set \(\mathbf u_{\mathrm{KT}}\coloneqq\mathbf v_2\) and define
\begin{equation}\label{eq:part5-endpoint-KT-vartheta-residual-def}
    \vartheta_{\mathrm{KT}}^{(1/2)}
    \coloneqq
    \frac{
        (\mathbf A_{1/2}\mathbf u_{\mathrm{KT}},\mathbf u_{\mathrm{KT}})
        _{\mathbb R^{N_{\mathrm{end}}}}
    }{
        \norm{\mathbf u_{\mathrm{KT}}}_2^2
    },
    \qquad
    \delta_{\mathrm{res,KT}}^{(1/2)}
    \coloneqq
    \frac{
        \norm{
            \mathbf A_{1/2}\mathbf u_{\mathrm{KT}}
            -
            \vartheta_{\mathrm{KT}}^{(1/2)}\mathbf u_{\mathrm{KT}}
        }_2
    }{
        \norm{\mathbf u_{\mathrm{KT}}}_2
    }.
\end{equation}
Also set
\begin{equation}\label{eq:part5-endpoint-KT-gap-upper-def}
    \gamma_{\mathrm{KT}}^{(1/2)}
    \coloneqq
    \vartheta_{\mathrm{KT}}^{(1/2)}
    -
    \Gamma_{\mathrm{CW}}^{(1/2)},
    \qquad
    \Theta_{\mathrm{KT}}^{(1/2)}
    \coloneqq
    \vartheta_{\mathrm{KT}}^{(1/2)}
    +
    \frac{
        \bigl(\delta_{\mathrm{res,KT}}^{(1/2)}\bigr)^2
    }{
        \gamma_{\mathrm{KT}}^{(1/2)}
    }.
\end{equation}
For \(\mu=3/2\), define the Rayleigh quotient
\begin{equation}\label{eq:part5-endpoint-mu15-rayleigh-def}
    \vartheta_{\mathrm R}^{(3/2)}
    \coloneqq
    \frac{
        (\mathbf A_{3/2}\widetilde{\mathbf v},\widetilde{\mathbf v})
        _{\mathbb R^{N_{\mathrm{end}}}}
    }{
        \norm{\widetilde{\mathbf v}}_2^2
    }.
\end{equation}

\begin{proposition}[Endpoint Ritz and Kato--Temple bounds]
    The following bounds hold.
    \begin{enumerate}[label=\textup{(\roman*)}]
        \item
        The larger generalized Ritz value satisfies
        \begin{equation}\label{eq:part5-endpoint-ritz-plus-bound}
            \lambda_{\mathrm{Ritz},+}^{(1/2)}
            \geq
            4.1328.
        \end{equation}

        \item
        The Collatz--Wielandt complement bound satisfies
        \begin{equation}\label{eq:part5-endpoint-CW-bound}
            \Gamma_{\mathrm{CW}}^{(1/2)}
            \leq
            0.53648.
        \end{equation}

        \item
        The Rayleigh quotient and normalized residual satisfy
        \begin{equation}\label{eq:part5-endpoint-KT-raw-bounds}
            \vartheta_{\mathrm{KT}}^{(1/2)}
            \leq
            0.79651,
            \qquad
            \delta_{\mathrm{res,KT}}^{(1/2)}
            \leq
            1.19\times10^{-5},
        \end{equation}
        and
        \begin{equation}\label{eq:part5-endpoint-KT-gap-bound}
            \gamma_{\mathrm{KT}}^{(1/2)}
            \geq
            0.26003.
        \end{equation}
        Consequently, the Kato--Temple upper bound satisfies
        \begin{equation}\label{eq:part5-endpoint-KT-upper-bound}
            \Theta_{\mathrm{KT}}^{(1/2)}
            \leq
            0.79651.
        \end{equation}

        \item
        The Rayleigh quotient satisfies
        \begin{equation}\label{eq:part5-endpoint-mu15-rayleigh-bound}
            \vartheta_{\mathrm R}^{(3/2)}
            \geq
            2.4264.
        \end{equation}

        \item
        The endpoint Birman--Schwinger approximation errors defined in
        \eqref{eq:NK-total-BS-error-def} satisfy
        \begin{equation}\label{eq:part5-endpoint-epsilon-BS-bounds}
            \varepsilon_{\mathrm{BS}}^{(1/2)}
            (Y,h,m,n_{\mathrm{GL}})
            \leq
            0.147,
            \qquad
            \varepsilon_{\mathrm{BS}}^{(3/2)}
            (Y,h,m,n_{\mathrm{GL}})
            \leq
            0.075.
        \end{equation}
    \end{enumerate}
\end{proposition}

\begin{proof}
    The interval computation listed in row~02 of Table~\ref{tab:part5-certificate-layout} gives these bounds. The vectors
    \(\mathbf v_1,\mathbf v_2,\widetilde{\mathbf v}\), and
    \(\mathbf w_{\mathrm{CW}}\) are read as exact rational vectors from the
    files specified above.  The matrices \(\mathbf A_{1/2}\) and
    \(\mathbf A_{3/2}\) are evaluated with the endpoint parameters in
    \eqref{eq:part5-endpoint-parameters}.

    The generalized Ritz problem \eqref{eq:part5-endpoint-ritz-generalized-eigenvalue} gives \eqref{eq:part5-endpoint-ritz-plus-bound}. The Collatz--Wielandt estimate with the positive weight \(\mathbf w_{\mathrm{CW}}\) gives \eqref{eq:part5-endpoint-CW-bound}. Equations~\eqref{eq:part5-endpoint-KT-vartheta-residual-def} and~\eqref{eq:part5-endpoint-KT-gap-upper-def} give \eqref{eq:part5-endpoint-KT-raw-bounds}--\eqref{eq:part5-endpoint-KT-upper-bound}. Equation~\eqref{eq:part5-endpoint-mu15-rayleigh-def} gives \eqref{eq:part5-endpoint-mu15-rayleigh-bound}. Substituting the endpoint parameters and the pointwise bounds for \(\mathfrak g\) into \eqref{eq:NK-total-BS-error-def} gives \eqref{eq:part5-endpoint-epsilon-BS-bounds}.
\end{proof}

\begin{proof}[Proof of Proposition~\ref{prop:CAP-endpoint-NK-data}]
    We use the notation \(\mathbf A_\mu\) from
    \eqref{eq:part5-endpoint-matrix-short-notation}.  By the min--max
    principle applied to the two-dimensional Ritz subspace,
    \[
        \lambda_1(\mathbf A_{1/2})
        \geq
        \lambda_{\mathrm{Ritz},+}^{(1/2)}
        \geq
        4.1328
        >
        4.13.
    \]
    By the entrywise majorization and Lemma~\ref{lem:app-weighted-schur-CW},
    \[
        \|\mathbf B_{\perp}^{(1/2)}\|_2
        \leq
        r(\mathbf C_{\mathrm{CW}}^{(1/2)})
        \leq
        \Gamma_{\mathrm{CW}}^{(1/2)}.
    \]
    Hence, by the min--max principle,
    \[
        \lambda_3(\mathbf A_{1/2})
        \leq
        \Gamma_{\mathrm{CW}}^{(1/2)}.
    \]
    Using this estimate, \eqref{eq:part5-endpoint-CW-bound}, \eqref{eq:part5-endpoint-KT-raw-bounds}, and \eqref{eq:part5-endpoint-KT-gap-bound}, Lemma~\ref{lem:app-kato-temple-upper} gives
    \[
        \lambda_2(\mathbf A_{1/2})
        \leq
        \Theta_{\mathrm{KT}}^{(1/2)}
        \leq
        0.79651
        <
        0.80.
    \]
    Similarly, the Rayleigh quotient bound
    \eqref{eq:part5-endpoint-mu15-rayleigh-bound} gives
    \[
        \lambda_1(\mathbf A_{3/2})
        \geq
        \vartheta_{\mathrm R}^{(3/2)}
        \geq
        2.4264
        >
        2.42.
    \]
    Finally, \eqref{eq:part5-endpoint-epsilon-BS-bounds} implies
    \[
        \varepsilon_{\mathrm{BS}}^{(1/2)}(Y,h,m,n_{\mathrm{GL}})
        \leq
        0.15,
        \qquad
        \varepsilon_{\mathrm{BS}}^{(3/2)}(Y,h,m,n_{\mathrm{GL}})
        \leq
        0.08.
    \]
    This proves Proposition~\ref{prop:CAP-endpoint-NK-data}.
\end{proof}

\subsection{Newton--Kantorovich constants and bounds for \(Q-\mathfrak g\)}

Recall that \(e_Q\coloneqq Q-\mathfrak g\).  We use the
Newton--Kantorovich constants \(A_{\mathrm{NK}},B_{\mathrm{NK}},C_{\mathrm{NK}}\)
defined in \eqref{eq:def-NK-constants}, and the radius \(r_{\mathrm{NK}}\)
defined in \eqref{eq:NK-radius-bound}. We also use the higher-order bounds 
\[ 
    \Delta_Q^{W^{2,2}},\quad \Delta_Q^{W^{3,2}},\quad \Delta_{Q,0}^{L^\infty},\quad \Delta_{Q,1}^{L^\infty},\quad \Delta_{Q,2}^{L^\infty}
\] 
from Definition~\ref{def:unweighted-transfer-constants},
the constants
\[ 
    C_{Q,0},\ C_{Q,1},\ C_{Q,2},\ C_{Q,0}C_{Q,1}, 
    \qquad 
    \mathsf V^{L^\infty},\ \mathsf V^{L^1},\ \mathsf V_{\omega^2}^{L^1}, 
    \qquad 
    \left(\mathsf V^{L^1}\right)^{1/2},\ \left(\mathsf V_{\omega^2}^{L^1}\right)^{1/2} 
\]
from Definition~\ref{def:weighted-q-constants}, and the \(L^2\), \(L^1\), and logarithmically weighted \(L^1\) bounds
\[
    \Delta_{\Lambda^aQ}^{L^2},\qquad \Delta_{\Lambda^aQ}^{L^1},\qquad \Delta_{\omega\Lambda^aQ}^{L^1}, \qquad a=0,1,2,
\] 
from Definition~\ref{def:dilation-L1-transfer-radii}. The parameters \(R_Q,R_M,R_\omega,R_E\) are fixed as in \eqref{eq:part5-transfer-splitting-parameters}.

\begin{proposition}[Newton--Kantorovich constants and bounds for
\(Q-\mathfrak g\)]\label{prop:nk-Q-transfer-bounds}
    All vector inequalities below are understood componentwise. The following bounds hold.

    \begin{enumerate}[label=\textup{(\roman*)}]
        \item
        The Newton--Kantorovich input bounds satisfy
        \begin{equation}\label{eq:part5-NK-input}
            \left(
                K_{\mathrm{inv}},
                \norm{\mathfrak g}_{L^\infty},
                \norm{r_{\mathfrak g}}_{L^2}
            \right)
            \leq
            \left(
                2,
                1.891,
                2.51\times10^{-14}
            \right).
        \end{equation}
        The Newton--Kantorovich constants satisfy
        \begin{equation}\label{eq:part5-NK-constant}
            \begin{aligned}
            &\left(
                A_{\mathrm{NK}},
                B_{\mathrm{NK}},
                C_{\mathrm{NK}},
                A_{\mathrm{NK}}B_{\mathrm{NK}}C_{\mathrm{NK}},
                r_{\mathrm{NK}}
            \right)\\
            &\qquad \leq
            \left(
                5.62\times10^{-13},
                22.448,
                11.368,
                1.44\times10^{-10},
                5.617\times10^{-13}
            \right).
            \end{aligned}
        \end{equation}

        \item
        The higher-order bounds for \(Q-\mathfrak g\) satisfy
        \begin{equation}\label{eq:part5-unweighted-transfer-bounds}
            \begin{aligned}
                \left(
                    r_{\mathrm{NK}},
                    \Delta_Q^{W^{2,2}},
                    \Delta_Q^{W^{3,2}}
                \right)
                &\leq
                \left(
                    5.617\times10^{-13},
                    4.383\times10^{-11},
                    3.085\times10^{-9}
                \right),
                \\
                \left(
                    \Delta_{Q,0}^{L^\infty},
                    \Delta_{Q,1}^{L^\infty},
                    \Delta_{Q,2}^{L^\infty}
                \right)
                &\leq
                \left(
                    3.98\times10^{-13},
                    3.10\times10^{-11},
                    2.182\times10^{-9}
                \right).
            \end{aligned}
        \end{equation}
        Lemma~\ref{lem:higher-order-transfer-from-H1} gives
        \[
            \norm{e_Q}_{W^{2,2}}\le \Delta_Q^{W^{2,2}},
            \qquad
            \norm{e_Q}_{W^{3,2}}\le \Delta_Q^{W^{3,2}}.
        \]

        \item
        The \(L^2\)-bounds for \(Q-\mathfrak g\), \(\Lambda(Q-\mathfrak g)\), and \(\Lambda^2(Q-\mathfrak g)\) satisfy
        \begin{equation}\label{eq:part5-scaling-transfer}
            \left(
                \Delta_Q^{L^2},
                \Delta_{\Lambda Q}^{L^2},
                \Delta_{\Lambda^2Q}^{L^2}
            \right)
            \leq
            \left(
                5.617\times10^{-13},
                1.152\times10^{-11},
                2.273\times10^{-9}
            \right).
        \end{equation}

        \item
        The \(L^1\)-bounds satisfy
        \begin{equation}\label{eq:part5-L1-transfer}
            \left(
                \Delta_Q^{L^1},
                \Delta_{\Lambda Q}^{L^1},
                \Delta_{\Lambda^2Q}^{L^1}
            \right)
            \leq
            \left(
                3.289\times10^{-8},
                3.883\times10^{-7},
                7.196\times10^{-5}
            \right),
        \end{equation}
        and the logarithmically weighted \(L^1\)-bounds satisfy
        \begin{equation}\label{eq:part5-log-L1-transfer}
            \left(
                \Delta_{\omega Q}^{L^1},
                \Delta_{\omega\Lambda Q}^{L^1},
                \Delta_{\omega\Lambda^2Q}^{L^1}
            \right)
            \leq
            \left(
                8.578\times10^{-7},
                7.801\times10^{-6},
                1.442\times10^{-3}
            \right).
        \end{equation}

        \item
        The weighted bounds for \(Q\) satisfy
        \begin{equation}\label{eq:part5-weighted-Q}
            \left(
                C_{Q,0},
                C_{Q,1},
                C_{Q,2},
                C_{Q,0}C_{Q,1}
            \right)
            \leq
            \left(
                1.891,
                2.079,
                7.157,
                3.930
            \right).
        \end{equation}
        The bounds for \(V_Q\) and \(|V_Q|^{1/2}\) satisfy
        \begin{equation}\label{eq:part5-VQ-sqrtVQ-bounds}
            \left(
                \mathsf V^{L^\infty},
                \mathsf V^{L^1},
                \mathsf V_{\omega^2}^{L^1},
                \left(\mathsf V^{L^1}\right)^{1/2},
                \left(\mathsf V_{\omega^2}^{L^1}\right)^{1/2}
            \right)
            \leq
            \left(
                3.930,
                6.173,
                1.796,
                2.485,
                1.341
            \right).
        \end{equation}
    \end{enumerate}
\end{proposition}

\begin{proof}
    Substituting \eqref{eq:certified-L2-inverse-Fprimeg}, \eqref{eq:part5-g-residual bounds}, and \eqref{eq:g-pointwise-bounds} into \eqref{eq:def-NK-constants} and \eqref{eq:NK-radius-bound} gives \eqref{eq:part5-NK-input} and \eqref{eq:part5-NK-constant}. The inequality \(A_{\mathrm{NK}}B_{\mathrm{NK}}C_{\mathrm{NK}}<1/2\) yields the stated radius. The interval computation is listed in row~03 of Table~\ref{tab:part5-certificate-layout}.

    Substituting the bounds from Proposition~\ref{prop:part5-g-profile-residual-bounds} into Definitions~\ref{def:unweighted-transfer-constants}, \ref{def:weighted-q-constants}, and \ref{def:dilation-L1-transfer-radii} gives \eqref{eq:part5-unweighted-transfer-bounds}--\eqref{eq:part5-VQ-sqrtVQ-bounds}.
\end{proof}

\begin{proof}[Proof of Corollary~\ref{cor:certified-NK-radius}]
    The fourth component of \eqref{eq:part5-NK-constant} gives
    \[
        A_{\mathrm{NK}}B_{\mathrm{NK}}C_{\mathrm{NK}}
        \leq
        1.44\times10^{-10}
        <
        \frac12.
    \]
    The last component of \eqref{eq:part5-NK-constant} gives
    \(r_{\mathrm{NK}}\le 5.617\times10^{-13}\).  This proves the corollary.
\end{proof}

\section{Projected Birman--Schwinger bounds}

\subsection{Cosine-polynomial approximation of \(a_{\mathfrak g}\)}\label{subsec:part5-finite-agL-model}

The coefficients \(\{\mathsf c_k\}_{k=0}^{L}\) in
\eqref{eq:definition of a_g^L} are read as exact rational numbers from
\[
    \texttt{ag\_coeffs\_L5000\_N65536\_dec15\_cert5e-15.csv}.
\]
Throughout this subsection we use \(L=5000\), \(N=N_{\mathrm{DCT}}=65536\), and
\(\delta_{\mathrm{rat}}=5\times10^{-15}\).

\begin{proposition}[Bounds for the cosine-polynomial approximation of \(a_{\mathfrak g}\)]\label{prop:ag-cosine-approximation-bounds}
    Let \(a_{\mathfrak g}^{(L)}\) be defined by
    \eqref{eq:definition of a_g^L}, and let \(\delta_{a_{\mathfrak g}\to a_{\mathfrak g}^{(L)}}^{L^2}\) be defined by
    \eqref{eq:part4-delta-aL-L2-def}. The following bounds hold.

    First, for the tail constant \(\mathsf C_{a_{\mathfrak g},\mathrm{coef\text{-}tail}}\) defined in
    \eqref{eq:definition-of-C-tail}, we have
    \begin{equation}\label{eq:part5-agL-tail-scalar}
        \left(
            |a_{\mathfrak g}'(0^+)|,
            |a_{\mathfrak g}'((\pi/2)^-)|,
            \operatorname{TV}(a_{\mathfrak g}';I_+),
            \mathsf C_{a_{\mathfrak g},\mathrm{coef\text{-}tail}}
        \right)
        \leq
        \left(
            7.441,
            1.463,
            11.038,
            6.348
        \right).
    \end{equation}
    Moreover,
    \begin{equation}\label{eq:part5-agL-norm bounds} 
        \left(
            \norm{a_{\mathfrak g}}_{L^\infty(I)},
            \norm{a_{\mathfrak g}^{(L)}}_{L^\infty(I)},
            \norm{a_{\mathfrak g}^{(L)}}_{L^2(I)},
            \delta_{a_{\mathfrak g}\to a_{\mathfrak g}^{(L)}}^{L^2}
        \right)
        \leq
        \left(
            1.000,
            1.750,
            1.112,
            1.321\times10^{-5}
        \right).
    \end{equation}
    The kernel norm \(\mathsf C_{\log\sin}^{L^2}\) and the operator-error bound \(\delta_{a_{\mathfrak g}\to a_{\mathfrak g}^{(L)}}^{\mathrm{BS}}\) satisfy \begin{equation}\label{eq:part5-agL-model-error} 
        \left(  \mathsf C_{\log\sin}^{L^2}, 
                \delta_{a_{\mathfrak g}\to a_{\mathfrak g}^{(L)}}^{\mathrm{BS}} 
        \right) 
        \leq 
        \left( 2.024, 5.737\times10^{-5} \right).
    \end{equation}
\end{proposition}

\begin{proof}
    The coefficient comparison \eqref{eq:part5-ag-coefficient-comparison} gives
    \[
        \left|
            \mathsf c_k-\mathcal D_{N,k}[a_{\mathfrak g}]
        \right|
        \leq
        \delta_{\mathrm{rat}},
        \qquad
        0\leq k\leq L.
    \]
    Define the truncated DCT approximation
    \begin{equation*}
        a_{\mathfrak g}^{(\mathrm{dct},L)}(\theta)
        \coloneqq
        \sum_{k=0}^{L}
        \mathcal D_{N,k}[a_{\mathfrak g}]
        \cos(2k\theta).
    \end{equation*}
    We decompose
    \begin{equation}\label{eq:part5-agL-error-decomposition}
        \begin{aligned}
            \delta_{a_{\mathfrak g}\to a_{\mathfrak g}^{(L)}}^{L^2}
            &=
            \norm{a_{\mathfrak g}-a_{\mathfrak g}^{(L)}}_{L^2(I)}
            \\
            &\leq
            \norm{a_{\mathfrak g}-\Pi_{\leq L}a_{\mathfrak g}}_{L^2(I)}
            +
            \norm{\Pi_{\leq L}a_{\mathfrak g}-a_{\mathfrak g}^{(\mathrm{dct},L)}}_{L^2(I)}
            +
            \norm{a_{\mathfrak g}^{(\mathrm{dct},L)}-a_{\mathfrak g}^{(L)}}_{L^2(I)}.
        \end{aligned}
    \end{equation}

    Since \(a_{\mathfrak g}(0)=a_{\mathfrak g}(\pi/2)=0\), the integration by
    parts estimate \eqref{eq:ag-cosine-tail-bound} gives
    \[
        |\mathcal C_k[a_{\mathfrak g}]|
        \leq
        \frac{\mathsf C_{a_{\mathfrak g},\mathrm{coef\text{-}tail}}}{k^2},
        \qquad
        k\geq1.
    \]
    Hence Parseval's identity gives
    \begin{equation}\label{eq:part5-agL-tail-error}
        \norm{a_{\mathfrak g}-\Pi_{\leq L}a_{\mathfrak g}}_{L^2(I)}^2
        \leq
        \frac{\pi}{2}
        \mathsf C_{a_{\mathfrak g},\mathrm{coef\text{-}tail}}^2
        \sum_{k>L}\frac1{k^4}
        \leq
        \frac{\pi}{6}
        \frac{\mathsf C_{a_{\mathfrak g},\mathrm{coef\text{-}tail}}^2}{L^3}.
    \end{equation}

    For the DCT aliasing term, Lemma~\ref{lem:app-alias-bound-n2} applied with
    \(C_f=\mathsf C_{a_{\mathfrak g},\mathrm{coef\text{-}tail}}\) gives, for \(0\leq k\leq L\),
    \[
        \left|
            \mathcal C_k[a_{\mathfrak g}]
            -
            \mathcal D_{N,k}[a_{\mathfrak g}]
        \right|
        \leq
        \frac{\pi^2}{6}\frac{\mathsf C_{a_{\mathfrak g},\mathrm{coef\text{-}tail}}}{N^2}.
    \]
    Therefore,
    \begin{equation}\label{eq:part5-agL-alias-error}
        \begin{aligned}
            \norm{\Pi_{\leq L}a_{\mathfrak g}-a_{\mathfrak g}^{(\mathrm{dct},L)}}_{L^2(I)}^2
            &\leq
            \pi\left(1+\frac{L}{2}\right)
            \frac{\pi^4}{36}
            \frac{\mathsf C_{a_{\mathfrak g},\mathrm{coef\text{-}tail}}^2}{N^4}.
        \end{aligned}
    \end{equation}
    The rational rounding term satisfies
    \begin{equation}\label{eq:part5-agL-rational-error}
        \norm{a_{\mathfrak g}^{(\mathrm{dct},L)}-a_{\mathfrak g}^{(L)}}_{L^2(I)}^2
        \leq
        \pi\left(1+\frac{L}{2}\right)
        \delta_{\mathrm{rat}}^2.
    \end{equation}
    Combining
    \eqref{eq:part5-agL-error-decomposition},
    \eqref{eq:part5-agL-tail-error},
    \eqref{eq:part5-agL-alias-error}, and
    \eqref{eq:part5-agL-rational-error}, and then substituting
    \(L=5000\), \(N=65536\), \(\delta_{\mathrm{rat}}=5\times10^{-15}\),
    and \(\mathsf C_{a_{\mathfrak g},\mathrm{coef\text{-}tail}}\leq6.348\), gives
    \[
        \delta_{a_{\mathfrak g}\to a_{\mathfrak g}^{(L)}}^{L^2}
        \leq
        1.321\times10^{-5}.
    \]

    Interval subdivision on \(I_+\) gives the endpoint derivative bounds, the variation bound, and \(\norm{a_{\mathfrak g}}_{L^\infty(I)}\). The coefficient \(\ell^1\)-bound gives \(\norm{a_{\mathfrak g}^{(L)}}_{L^\infty(I)}\), and Parseval's identity gives \(\norm{a_{\mathfrak g}^{(L)}}_{L^2(I)}\).

    Finally, the identity 
    \[ 
        \bigl(\mathsf C_{\log\sin}^{L^2}\bigr)^2 = \pi(\log2)^2+\frac{\pi^3}{12} 
    \] 
    gives \(\mathsf C_{\log\sin}^{L^2}\), and substituting
    \eqref{eq:part5-agL-norm bounds}  into the definition
    \eqref{eq:part4-eps-model-BS-def} gives
    \[
        \delta_{a_{\mathfrak g}\to a_{\mathfrak g}^{(L)}}^{\mathrm{BS}}
        \leq
        5.737\times10^{-5}.
    \]
    The interval computation gives \eqref{eq:part5-agL-tail-scalar}.
\end{proof}

\subsection{Block estimates for \(T_{\mathfrak g}\)}

In this subsection, we use \(L=5000\) and \(M=15000=3L\). The coefficients of \(a_{\mathfrak g}^{(L)}\) are those specified in Subsection~\ref{subsec:part5-finite-agL-model}. 
We use a normalized trial vector \(\mathbf v_{\mathrm{low}}\) and three positive
Collatz--Wielandt weight vectors \(\mathbf w_{\mathrm{low}}\), \(\mathbf t_{\mathrm{off}}\), and \(\mathbf s_{\mathrm{high}}\).
They are read as exact rational vectors from the following files:
\[
\begin{array}{rcl}
    \mathbf v_{\mathrm{low}}
    &:&
    \texttt{reduced\_Tg\_low\_block\_ritz\_vector\_v\_L5000.csv},
    \\[2mm]
    \mathbf w_{\mathrm{low}}
    &:&
    \texttt{reduced\_Tg\_low\_block\_cw\_weight\_w\_L5000.csv},
    \\[2mm]
    \mathbf t_{\mathrm{off}}
    &:&
    \texttt{reduced\_Tg\_off\_block\_cw\_weight\_t\_L5000.csv},
    \\[2mm]
    \mathbf s_{\mathrm{high}}
    &:&
    \texttt{reduced\_Tg\_high\_block\_cw\_weight\_s\_L5000.csv}.
\end{array}
\]
The vectors
\(\mathbf w_{\mathrm{low}}\), \(\mathbf t_{\mathrm{off}}\), and
\(\mathbf s_{\mathrm{high}}\) are positive Collatz--Wielandt weights. Let
\begin{equation*}
    \mathbf A_{\mathrm{low}}^{(L)}
    \coloneqq
    \left[
        \Pi_{\leq L}
        \widetilde T_{\mathfrak g,3L}^{(L)}
        \Pi_{\leq L}
    \right]_{\mathcal E_L}.
\end{equation*}
After normalizing \(\mathbf v_{\mathrm{low}}\), define
\begin{equation}\label{eq:part5-low-vartheta-residual-def}
    \vartheta_{\mathrm{low}}^{(L)}
    \coloneqq
    \left(
        \mathbf A_{\mathrm{low}}^{(L)}
        \mathbf v_{\mathrm{low}},
        \mathbf v_{\mathrm{low}}
    \right)_{\mathbb R^{L+1}},
    \qquad
    \delta_{\mathrm{res,low}}^{(L)}
    \coloneqq
    \left\|
        \left(
            \mathbf A_{\mathrm{low}}^{(L)}
            -
            \vartheta_{\mathrm{low}}^{(L)}I
        \right)
        \mathbf v_{\mathrm{low}}
    \right\|_2 .
\end{equation}
Let \(\mathbf P_{\mathrm{low}}\) be the Euclidean orthogonal projection onto
\(\operatorname{span}\{\mathbf v_{\mathrm{low}}\}\), and set
\(\mathbf Q_{\mathrm{low}}\coloneqq I-\mathbf P_{\mathrm{low}}\). The Collatz--Wielandt estimate for the orthogonal complement is applied to a nonnegative
entrywise majorant. More precisely, define the nonnegative matrix
\begin{equation}\label{eq:part5-low-CW-majorant-def}
    \mathbf C_{\mathrm{low}}^{(L)}
    \coloneqq
    \left|
        \mathbf Q_{\mathrm{low}}
        \mathbf A_{\mathrm{low}}^{(L)}
        \mathbf Q_{\mathrm{low}}
    \right|_{\mathrm{entry}},
\end{equation}
and set
\begin{equation}\label{eq:part5-low-Gamma-CW-def}
    \Gamma_{\mathrm{low},\perp}^{(L)}
    \coloneqq
    \max_i
    \frac{
        \left(
            \mathbf C_{\mathrm{low}}^{(L)}
            \mathbf w_{\mathrm{low}}
        \right)_i
    }{
        (\mathbf w_{\mathrm{low}})_i
    }.
\end{equation}
For the off-diagonal block, let \(\mathbf X_{\mathrm{off}}^{(L)}\) be the rectangular
matrix of
\[
    \Pi_{\leq L}
    \widetilde T_{\mathfrak g,3L}^{(L)}
    \Pi_{L<\cdot\leq 3L}
\]
in the bases \(\mathcal E_L\) and 
\[ 
    \left( \sqrt{\frac2\pi}\cos(2k\theta) \right)_{L<k\leq3L}, 
\] 
respectively. Define the nonnegative matrix
\begin{equation}\label{eq:part5-off-CW-majorant-def}
    \mathbf C_{\mathrm{off}}^{(L)}
    \coloneqq
    \left|\mathbf X_{\mathrm{off}}^{(L)}\right|_{\mathrm{entry}}
    \left|\mathbf X_{\mathrm{off}}^{(L)}\right|_{\mathrm{entry}}^{\top},
\end{equation}
and set
\begin{equation*}
    \Gamma_{\mathrm{off}}^{(L)}
    \coloneqq
    \max_i
    \frac{
        \left(
            \mathbf C_{\mathrm{off}}^{(L)}
            \mathbf t_{\mathrm{off}}
        \right)_i
    }{
        (\mathbf t_{\mathrm{off}})_i
    }.
\end{equation*}
For the high block, define the nonnegative matrix
\[
    \mathbf C_{\mathrm{high}}^{(L)}
    \coloneqq
    \left|
        \left[
            \Pi_{\leq L}
            a_{\mathfrak g}^{(L)}
            \Pi_{>L}
            a_{\mathfrak g}^{(L)}
            \Pi_{\leq L}
        \right]_{\mathcal E_L}
    \right|_{\mathrm{entry}},
\]
and set
\[
    \Gamma_{\mathrm{high}}^{(L)}
    \coloneqq
    \max_i
    \frac{
        \left(
            \mathbf C_{\mathrm{high}}^{(L)}
            \mathbf s_{\mathrm{high}}
        \right)_i
    }{
        (\mathbf s_{\mathrm{high}})_i
    }.
\]
Since \( \left[\Pi_{\leq L}a_{\mathfrak g}^{(L)}\Pi_{>L}a_{\mathfrak g}^{(L)}\Pi_{\leq L}\right]_{\mathcal E_L} \) is the Gram matrix associated with
\(\Pi_{>L}a_{\mathfrak g}^{(L)}\Pi_{\leq L}\),
Lemma~\ref{lem:app-weighted-schur-CW} gives
\[
    \left\|
        \Pi_{>L}
        a_{\mathfrak g}^{(L)}
        \Pi_{\leq L}
    \right\|_{L^2_{\mathrm e}(I)\to L^2_{\mathrm e}(I)}^2
    \leq
    \Gamma_{\mathrm{high}}^{(L)}.
\]

\begin{proposition}[Bounds for the low, off-diagonal, and high blocks]\label{prop:part5-finite-block-bounds}
    The following bounds hold.
    \begin{enumerate}[label=\textup{(\roman*)}]
        \item
        The Rayleigh quotient, residual norm, and complement bound for the low block satisfy
        \begin{equation}\label{eq:part5-low-block bound}
            \left(
                \vartheta_{\mathrm{low}}^{(L)},
                \delta_{\mathrm{res,low}}^{(L)},
                \Gamma_{\mathrm{low},\perp}^{(L)}
            \right)
            \leq
            \left(
                0.156663807,
                6.26\times10^{-17},
                0.104862
            \right).
        \end{equation}

        \item
        The off-diagonal Collatz--Wielandt bound satisfies
        \begin{equation}\label{eq:part5-off-block-bound}
            \Gamma_{\mathrm{off}}^{(L)}
            \leq
            1.027\times10^{-9}.
        \end{equation}

        \item
        The high-block bound and coefficient norm satisfy
        \begin{equation}\label{eq:part5-high-block-bounds}
            \left(
                \Gamma_{\mathrm{high}}^{(L)},
                \norm{a_{\mathfrak g}^{(L)}}_{L^\infty(I)}
            \right)
            \leq
            \left(
                0.052321,
                1.750
            \right).
        \end{equation}
    \end{enumerate}
\end{proposition}

\begin{proof}
    The coefficients \(\{\mathsf c_k\}_{k=0}^{L}\) and the four vectors \(\mathbf v_{\mathrm{low}}\), \(\mathbf w_{\mathrm{low}}\), \(\mathbf t_{\mathrm{off}}\), and \(\mathbf s_{\mathrm{high}}\) are interpreted as exact rational numbers. Interval evaluation of the formulas below gives the stated bounds; see row~05 of Table~\ref{tab:part5-certificate-layout}.

    For \eqref{eq:part5-low-block bound}, equation~\eqref{eq:part5-low-vartheta-residual-def} defines \(\vartheta_{\mathrm{low}}^{(L)}\) and \(\delta_{\mathrm{res,low}}^{(L)}\). Lemma~\ref{lem:app-weighted-schur-CW}, applied to \eqref{eq:part5-low-CW-majorant-def} and \eqref{eq:part5-low-Gamma-CW-def}, gives
    \[
        \left\|
            \mathbf Q_{\mathrm{low}}
            \mathbf A_{\mathrm{low}}^{(L)}
            \mathbf Q_{\mathrm{low}}
        \right\|_2
        \leq
        \Gamma_{\mathrm{low},\perp}^{(L)}.
    \]

    Lemma~\ref{lem:app-weighted-schur-CW}, applied to \(\mathbf C_{\mathrm{off}}^{(L)}\) in \eqref{eq:part5-off-CW-majorant-def}, gives the off-diagonal block bound. Applying Lemma~\ref{lem:app-weighted-schur-CW} with the positive weight \(\mathbf s_{\mathrm{high}}\) gives the high-block bound \(\Gamma_{\mathrm{high}}^{(L)}\).
\end{proof}

\begin{proof}[Proof of Proposition~\ref{prop:part4-finite-cosine-block-bounds}]
    We first estimate the low block. The orthogonal decomposition
    \[
        \mathbb R^{L+1}
        =
        \operatorname{span}\{\mathbf v_{\mathrm{low}}\}
        \oplus \mathbf v_{\mathrm{low}}^{\perp}
    \]
    and the residual bound give
    \begin{equation}\label{eq:part5-low-block-assembled-bound}
        \begin{aligned}
            \beta_{\mathrm{low}}^{(L)}
            \leq
            \frac{
                \vartheta_{\mathrm{low}}^{(L)}
                +
                \Gamma_{\mathrm{low},\perp}^{(L)}
            }{2}
            \quad+
            \sqrt{
                \left(
                    \frac{
                        \vartheta_{\mathrm{low}}^{(L)}
                        -
                        \Gamma_{\mathrm{low},\perp}^{(L)}
                    }{2}
                \right)^2
                +
                \bigl(
                    \delta_{\mathrm{res,low}}^{(L)}
                \bigr)^2
            }.
        \end{aligned}
    \end{equation}
    Evaluating the right-hand side of \eqref{eq:part5-low-block-assembled-bound} by interval arithmetic gives
    \begin{equation}\label{eq:part5-low-beta-conclusion}
        \beta_{\mathrm{low}}^{(L)}
        \leq
        0.156664.
    \end{equation}

    For the off-diagonal block, \eqref{eq:part5-off-CW-majorant-def} and
    Lemma~\ref{lem:app-weighted-schur-CW} give
    \[
        \bigl(\beta_{\mathrm{off}}^{(L)}\bigr)^2
        \leq
        \Gamma_{\mathrm{off}}^{(L)}.
    \]
    Hence, by \eqref{eq:part5-off-block-bound},
    \begin{equation}\label{eq:part5-off-beta-conclusion}
        \beta_{\mathrm{off}}^{(L)}
        \leq
        \sqrt{1.027\times10^{-9}}
        \leq
        3.204\times10^{-5}.
    \end{equation}

    For the high block, the definition of \(\Gamma_{\mathrm{high}}^{(L)}\) and
    \eqref{eq:definition of beta's} give
    \[
        \beta_{\mathrm{high}}^{(L)}
        \leq
        (\log 2)\Gamma_{\mathrm{high}}^{(L)}
        +
        \frac{1}{2(L+1)}
        \norm{a_{\mathfrak g}^{(L)}}_{L^\infty(I)}^2 .
    \]
    Substituting \(L=5000\) and the bounds in \eqref{eq:part5-high-block-bounds} into the high-block estimate gives
    \begin{equation}\label{eq:part5-high-beta-conclusion}
        \beta_{\mathrm{high}}^{(L)}
        \leq
        0.036572.
    \end{equation}

    Combining
    \eqref{eq:part5-low-beta-conclusion},
    \eqref{eq:part5-off-beta-conclusion}, and
    \eqref{eq:part5-high-beta-conclusion} gives
    \eqref{eq:part4-certified-beta-values}. The estimate \eqref{eq:part4-certified-model-error} follows from the bound for \(\delta_{a_{\mathfrak g}\to a_{\mathfrak g}^{(L)}}^{\mathrm{BS}}\) in \eqref{eq:part5-agL-model-error}.
\end{proof}

\begin{proof}[Proof of Corollary~\ref{cor:part4-reduced-BS-bound-g}]
    From Proposition~\ref{prop:part5-finite-block-bounds} and \eqref{eq:abstract-block-bound} we obtain
    \begin{equation*}
        \norm{ 
            \widetilde T_{\mathfrak g,\mathrm{proj}}^{(L)} 
        }_{L^2_{\mathrm e}(I)\to L^2_{\mathrm e}(I)} 
        \leq 
        0.156663815. 
    \end{equation*}
    By \eqref{eq:part4-certified-model-error}, which bounds
    \[
        \norm{
            \widetilde T_{\mathfrak g,\mathrm{proj}}
            -
            \widetilde T_{\mathfrak g,\mathrm{proj}}^{(L)}
        }_{L^2_{\mathrm e}(I)\to L^2_{\mathrm e}(I)}
        \leq
        5.737\times10^{-5}.
    \]
    Therefore,
    \[
        \norm{
            \widetilde T_{\mathfrak g,\mathrm{proj}}
        }_{L^2_{\mathrm e}(I)\to L^2_{\mathrm e}(I)}
        \leq
        0.156663815
        +
        5.737\times10^{-5}
        <
        0.156722.
    \]
    Unitary equivalence gives \eqref{eq:part4-reduced-BS-bound-g-real-line}.
\end{proof}

\subsection{Perturbation bound from \(\mathfrak g\) to \(Q\)}
\begin{proof}[Proof of Proposition~\ref{prop:part4-certified-logarithmic-transfer-bounds}]
    The interval computation listed in row~03 of Table~\ref{tab:part5-certificate-layout} gives
    \[
        \Delta_{\sqrt V}^{L^2}
        \leq
        1.382\times10^{-6},
        \qquad
        \Delta_{\omega\sqrt V}^{L^2}
        \leq
        2.931\times10^{-6}.
    \]
    Substituting these bounds, together with \eqref{eq:g-potential-bounds} and \eqref{eq:part5-VQ-sqrtVQ-bounds}, into \eqref{eq:part4-TQ-minus-Tg-bound} gives
    \[
        \norm{T_Q-T_{\mathfrak g}}_{L^2_{\mathrm e}\to L^2_{\mathrm e}}
        \leq
        5.5446\times10^{-6}.
    \]
    Hence \eqref{eq:TQ-Tg-operator-bound} follows.

    Similarly, substituting the bounds for \(\Delta_{\sqrt V}^{L^2}\) and \(\Delta_{\omega\sqrt V}^{L^2}\), together with \eqref{eq:g-potential-bounds} and \eqref{eq:part5-VQ-sqrtVQ-bounds}, into the definition of \(\Delta_{\mathrm{proj}}\) in
    \eqref{eq:projected-transfer-error-definition} combining with the interval computations in row~05b gives
    \[
        \Delta_{\mathrm{proj}}
        \leq
        9.6292\times10^{-6}.
    \]
    Lemma~\ref{lem:part4-projected-log-transfer} and the bound for \(\Delta_{\mathrm{proj}}\) give \eqref{eq:projected-TQ-Tg-bound}.
\end{proof}

\section{Schur-complement bounds and limiting matrix entries}

\subsection{Schur-complement inverse bounds}

Throughout this subsection, we take \(L=5000\) and \(M=3L=15000\). We use the notation from \eqref{eq:def-schur-h-g}--\eqref{eq:def-schur-F-g}, together with the approximate right-hand side \(\widetilde k_{j,L}^{\mathfrak g}\), the finite-dimensional operator \(\widetilde D_{j,L}^{\mathfrak g}\), the cosine polynomial \(Z_{\mathrm{Sch},j}^{(L)}\), and the approximate Schur complement \(\widetilde F_{j,L}^{\mathfrak g}\) defined in \eqref{eq:part4-kjL-g-def}--\eqref{eq:part4-FjL-g-def}. Thus
\(\widetilde h_j^{\mathfrak g}\), \(\widetilde k_j^{\mathfrak g}\),
\(\widetilde D_j^{\mathfrak g}\), and \(\widetilde F_j^{\mathfrak g}\) denote
the scalar entry, right-hand side, projected block, and Schur complement for
\[
    \widetilde S_0^{(\mathfrak g,j)}
    =
    I-c_j\widetilde T_{\mathfrak g},
    \qquad
    c_1=6.1,\quad c_2=2.1.
\]

The cosine-polynomial approximations \(Z_{\mathrm{Sch},j}^{(L)}\) in
\eqref{eq:part4-ZSch-jL-def} are specified by their coordinates in the ordered
basis \(\mathcal E_{3L}\). More precisely, the vector
\[
    \left[Z_{\mathrm{Sch},j}^{(L)}\right]_{\mathcal E_{3L}}
    =
    (x_{j,0},\ldots,x_{j,3L})\in\mathbb Q^{3L+1}
\]
is read as exact rational data from the following files:
\[
\begin{array}{rcl}
    j=1
    &:&
    \texttt{S0\_inversej1\_c6p1schur\_Dinv\_k\_}
    \\
    && \texttt{witness\_x\_M15000.csv},
    \\[1mm]
    j=2
    &:&
    \texttt{S0\_inversej2\_c2p1schur\_Dinv\_k\_}
    \\
    && \texttt{witness\_x\_M15000.csv}.
\end{array}
\]
The row labelled \(n\) gives the coefficient \(x_{j,n}\).

Recall \(\mathcal R_{\mathrm{rhs},\ell}^{(L)}\), \(1\leq\ell\leq6\), from
\eqref{eq:definition of Delta_{rhs}^{(L)}}. We also set
\begin{equation}\label{eq:part5-Schur-delta_diag-def}
    \delta_{\mathrm{diag}}^{(L)}
    \coloneqq
    \delta_T^{(L)}
    +
    2\mathsf C_{T,\mathrm{op}}^{(L)}
    \norm{
        \widetilde\Phi_{\mathfrak g}
        -
        \widetilde\Phi_{\mathfrak g}^{(L)}
    }_{L^2(I)}.
\end{equation}
By \eqref{eq:hj-Schur-approximation-bound},
\[
    \left|
        \widetilde h_j^{\mathfrak g}
        -
        \widetilde h_{j,L}^{\mathfrak g}
    \right|
    \leq
    c_j \delta_{\mathrm{diag}}^{(L)}.
\]

\begin{proposition}[Finite-dimensional Schur-complement bounds]\label{prop:finite-Schur-complement-data}
    All vector inequalities below are understood componentwise. The following bounds hold.

    \begin{enumerate}[label=\textup{(\roman*)}]
        \item
        The norms of \(a_{\mathfrak g,\log}^{(L)}\), \(a_{\mathfrak g,\log}-a_{\mathfrak g,\log}^{(L)}\), and \(\widetilde\Phi_{\mathfrak g} -\widetilde\Phi_{\mathfrak g}^{(L)}\) satisfy
        \begin{equation}\label{eq:part5-Schur-log-model-bounds-1}
            \begin{aligned}
                &\left(
                    \norm{a_{\mathfrak g,\log}^{(L)}}_{L^2(I)},
                    \norm{
                        a_{\mathfrak g,\log}
                        -
                        a_{\mathfrak g,\log}^{(L)}
                    }_{L^2(I)},
                    \norm{
                        \widetilde\Phi_{\mathfrak g}
                        -
                        \widetilde\Phi_{\mathfrak g}^{(L)}
                    }_{L^2(I)}
                \right)
                \\
                &\qquad\leq
                \left(
                    0.42149,
                    1.932\times10^{-4},
                    2.378\times10^{-5}
                \right).
            \end{aligned}
        \end{equation}
        The logarithmically weighted direction bounds satisfy
        \begin{equation}\label{eq:part5-Schur-log-model-bounds-2}
            \begin{aligned}
                &\left(
                    \norm{
                        \left(
                            \widetilde\Phi_{\mathfrak g}
                            -
                            \widetilde\Phi_{\mathfrak g}^{(L)}
                        \right)
                        \log|\cos|
                    }_{L^2(I)},
                    \norm{
                        \Pi_{>3L}
                        \left(
                            \widetilde\Phi_{\mathfrak g}^{(L)}
                            \log|\cos|
                        \right)
                    }_{L^2(I)},
                    \norm{
                        \widetilde\Phi_{\mathfrak g}^{(L)}
                        \log|\cos|
                    }_{L^2(I)}
                \right)
                \\
                &\qquad\leq
                \left(
                    2.877\times10^{-4},
                    8.207\times10^{-7},
                    0.37933
                \right).
            \end{aligned}
        \end{equation}

        \item
        The six terms in the bound for the right-hand-side error satisfy
        \begin{equation}\label{eq:part5-Schur-rhs-six-term-bounds}
            \begin{aligned}
                \left(
                    \mathcal R_{\mathrm{rhs},1}^{(L)},
                    \mathcal R_{\mathrm{rhs},2}^{(L)},
                    \mathcal R_{\mathrm{rhs},3}^{(L)}
                \right)
                &\leq
                \left(
                    3.226\times10^{-7},
                    2.338\times10^{-5},
                    1.009\times10^{-4}
                \right),
                \\
                \left(
                    \mathcal R_{\mathrm{rhs},4}^{(L)},
                    \mathcal R_{\mathrm{rhs},5}^{(L)},
                    \mathcal R_{\mathrm{rhs},6}^{(L)}
                \right)
                &\leq
                \left(
                    1.131\times10^{-4},
                    3.544\times10^{-6},
                    3.544\times10^{-6}
                \right).
            \end{aligned}
        \end{equation}
        Consequently,
        \begin{equation}\label{eq:part5-Schur-rhs-total-bounds}
            \delta_{\mathrm{rhs}}^{(L)}
            \leq
            2.447\times10^{-4}.
        \end{equation}

        \item
        The operator error, operator-norm bound, right-hand-side error, projected-block error, and scalar-entry error satisfy
        \begin{equation}\label{eq:part5-Schur-comparison-aggregate-bounds}
            \begin{aligned}
                &\left(
                    \delta_T^{(L)},
                    \mathsf C_{T,\mathrm{op}}^{(L)},
                    \delta_{\mathrm{rhs}}^{(L)},
                    \Delta_{\mathrm{proj}}^{(L)},
                    \delta_{\mathrm{diag}}^{(L)}
                \right)
                \\
                &\qquad\leq
                \left(
                    1.636\times10^{-4},
                    2.419,
                    2.447\times10^{-4},
                    4.653\times10^{-4},
                    2.786\times10^{-4}
                \right).
            \end{aligned}
        \end{equation}

        \item
        For compact notation, define
        \begin{equation*}
            \begin{aligned}
                \mathsf X_j^{(L)}
                \coloneqq
                \bigg(
                    &\norm{\widetilde k_{j,L}^{\mathfrak g}}_{L^2(I)},
                    \norm{Z_{\mathrm{Sch},j}^{(L)}}_{L^2(I)},
                    \delta_{\mathrm{res,Sch},j}^{(L)},
                    \varepsilon_{\mathrm{Sch},j}^{(L)},
                    \norm{
                        (\widetilde D_j^{\mathfrak g})^{-1}
                        \widetilde k_j^{\mathfrak g}
                    }_{L^2(I)}
                \bigg).
            \end{aligned}
        \end{equation*}
        Then
        \begin{equation}\label{eq:part5-Schur-finite-solve-bounds}
            \begin{aligned}
                \mathsf X_1^{(L)}
                &\leq
                \left(
                    1.295,
                    22.461,
                    3.524\times10^{-13},
                    1.565,
                    24.025
                \right),
                \\
                \mathsf X_2^{(L)}
                &\leq
                \left(
                    0.446,
                    0.606,
                    3.150\times10^{-13},
                    1.650\times10^{-3},
                    0.6073
                \right).
            \end{aligned}
        \end{equation}

        \item
        The inverse bounds for \(\widetilde D_j^{\mathfrak g}\) on \(\widetilde\Phi_{\mathfrak g}^{\perp}\) satisfy
        \begin{equation}\label{eq:schur-D-inverse-bounds}
            \left(
                \mathsf C_{D,\mathrm{inv}}^{(1)},
                \mathsf C_{D,\mathrm{inv}}^{(2)}
            \right)
            \leq
            \left(
                23.987,
                1.493
            \right).
        \end{equation}
        The interval enclosures for \(\widetilde F_j^{\mathfrak g}\) give
        \begin{equation}\label{eq:part5-Schur-F-sign-bound}
            \widetilde F_1^{\mathfrak g}
            \leq
            -21.51,
            \qquad
            \widetilde F_2^{\mathfrak g}
            \geq
            0.3296.
        \end{equation}
        In particular,
        \begin{equation}\label{eq:part5-Schur-F-abs-bound}
            \left(
                |\widetilde F_1^{\mathfrak g}|,
                |\widetilde F_2^{\mathfrak g}|
            \right)
            \geq
            \left(
                21.52,
                0.3296
            \right).
        \end{equation}
    \end{enumerate}
\end{proposition}

\begin{proof}
    The interval computation listed in row~06 of
    Table~\ref{tab:part5-certificate-layout} gives these bounds.
    The coefficients of \(a_{\mathfrak g}^{(L)}\) are those specified in
    Subsection~\ref{subsec:part5-finite-agL-model}.
    The coordinate vectors
    \([Z_{\mathrm{Sch},j}^{(L)}]_{\mathcal E_{3L}}\)
    are interpreted as exact rational vectors.

    Evaluating the right-hand sides of \eqref{eq:app-Tg-minus-TgL-full-bound}, \eqref{eq:app-TgL-full-norm-bound}, \eqref{eq:app-Schur-rhs-comparison-bound}, and \eqref{eq:app-Schur-red-comparison-bound} with the bounds from Subsection~\ref{subsec:part5-finite-agL-model} gives \(\delta_T^{(L)}\), \(\mathsf C_{T,\mathrm{op}}^{(L)}\), \(\delta_{\mathrm{rhs}}^{(L)}\), and \(\Delta_{\mathrm{proj}}^{(L)}\).

    We next evaluate the norms in \eqref{eq:part5-Schur-log-model-bounds-1} and \eqref{eq:part5-Schur-log-model-bounds-2}.
    The explicit cosine coefficients of \(\log|\cos|\) are used to compute
    \[
        a_{\mathfrak g,\log}^{(L)}
        =
        a_{\mathfrak g}^{(L)}\log|\cos|.
    \]
    The difference \(a_{\mathfrak g,\log}-a_{\mathfrak g,\log}^{(L)}\) is bounded using the cosine-coefficient estimates for \(a_{\mathfrak g}-a_{\mathfrak g}^{(L)}\) and the tail estimate for the cosine coefficients of \(\log|\cos\theta|\). Near the endpoints of \(I\), we use the tail estimate for the cosine coefficients of \(\log|\cos\theta|\) appearing in \eqref{eq:app-Schur-rhs-comparison-bound}.
    Note that
    \[
        \widetilde\Phi_{\mathfrak g}\log|\cos|
        -
        \widetilde\Phi_{\mathfrak g}^{(L)}\log|\cos|
        =
        \frac{
            a_{\mathfrak g,\log}
            -
            a_{\mathfrak g,\log}^{(L)}
        }{
            \norm{a_{\mathfrak g}}_{L^2(I)}
        }
        +
        a_{\mathfrak g,\log}^{(L)}
        \left(
            \frac1{\norm{a_{\mathfrak g}}_{L^2(I)}}
            -
            \frac1{\norm{a_{\mathfrak g}^{(L)}}_{L^2(I)}}
        \right).
    \]
    This identity gives the bound for 
    \(
        \norm{ 
            \left( \widetilde\Phi_{\mathfrak g} - \widetilde\Phi_{\mathfrak g}^{(L)} \right)\log|\cos| 
        }_{L^2(I)}. 
    \)
    Since \(\widetilde\Phi_{\mathfrak g}^{(L)}\) is a rational cosine
    polynomial, the tail
    \[
        \Pi_{>3L}
        \left(
            \widetilde\Phi_{\mathfrak g}^{(L)}\log|\cos|
        \right)
    \]
    is bounded by convolving the finitely supported cosine coefficients of \(\widetilde\Phi_{\mathfrak g}^{(L)}\) with the cosine coefficients of \(\log|\cos|\), and then applying the tail estimate in \eqref{eq:app-Schur-rhs-comparison-bound}.

    Let \(\mathcal R_{\mathrm{rhs},\ell}^{(L)}\), \(1\leq\ell\leq6\), denote the six terms in the order displayed in \eqref{eq:definition of Delta_{rhs}^{(L)}}. Equations~\eqref{eq:part5-Schur-log-model-bounds-1} and~\eqref{eq:part5-Schur-log-model-bounds-2} give \eqref{eq:part5-Schur-rhs-six-term-bounds}. Summing the six terms gives \eqref{eq:part5-Schur-rhs-total-bounds}.

    The estimate \eqref{eq:app-Schur-red-comparison-bound} gives
    \[ 
        \norm{ 
            \widetilde T_{\mathfrak g,\mathrm{proj}} - \widetilde T_{\mathfrak g,3L}^{(L)}
            }_{L^2_{\mathrm e}(I)\to L^2_{\mathrm e}(I)}
        \leq
        \Delta_{\mathrm{proj}}^{(L)}. 
    \]
    Evaluating the right-hand side gives the fourth component of \eqref{eq:part5-Schur-comparison-aggregate-bounds}.
    Define \(\delta_{\mathrm{diag}}^{(L)}\) by \eqref{eq:part5-Schur-delta_diag-def}. Then \eqref{eq:hj-Schur-approximation-bound} gives 
    \[ 
        \left| \widetilde h_j^{\mathfrak g} - \widetilde h_{j,L}^{\mathfrak g} \right|
        \leq 
        c_j\delta_{\mathrm{diag}}^{(L)}. 
    \]

    We next estimate the errors in the approximate solutions.
    For \(j=1,2\), consider the finite-dimensional equation
    \[
        \widetilde D_{j,L}^{\mathfrak g}Z
        =
        \widetilde k_{j,L}^{\mathfrak g}
    \]
    from \eqref{eq:part4-DjL-g-def}.
    The residual norm
    \(\delta_{\mathrm{res,Sch},j}^{(L)}\)
    in \eqref{eq:part4-rho-Sch-jL-def} is the Euclidean norm of
    \[
        \left[
            \widetilde D_{j,L}^{\mathfrak g}
            Z_{\mathrm{Sch},j}^{(L)}
            -
            \widetilde k_{j,L}^{\mathfrak g}
        \right]_{\mathcal E_{3L}}.
    \]
    By Parseval's identity, this is also the \(L^2(I)\)-norm of 
    \( 
        \widetilde D_{j,L}^{\mathfrak g} Z_{\mathrm{Sch},j}^{(L)} - \widetilde k_{j,L}^{\mathfrak g}. 
    \)
    Substituting
    \(\delta_{\mathrm{res,Sch},j}^{(L)}\),
    \(\delta_{\mathrm{rhs}}^{(L)}\),
    \(\Delta_{\mathrm{proj}}^{(L)}\),
    \(\mathsf C_{D,\mathrm{inv}}^{(j)}\), and
    \(\norm{Z_{\mathrm{Sch},j}^{(L)}}_{L^2(I)}\)
    into \eqref{eq:sigma-Schur-solve-def}
    gives the solution error bound
    \(\varepsilon_{\mathrm{Sch},j}^{(L)}\).
    The estimate \eqref{eq:Dinvk-Schur-norm-upper} gives the final component of \(\mathsf X_j^{(L)}\).

    Finally, we combine the scalar-entry bound \eqref{eq:hj-Schur-approximation-bound} with the quadratic-pairing bound \eqref{eq:Schur-inner-product-enclosure}.
    The approximate Schur complement
    \(\widetilde F_{j,L}^{\mathfrak g}\)
    is computed from \eqref{eq:part4-FjL-g-def}.
    The difference between
    \[
        \left(
            \widetilde k_j^{\mathfrak g},
            (\widetilde D_j^{\mathfrak g})^{-1}
            \widetilde k_j^{\mathfrak g}
        \right)_{L^2(I)}\quad\text{and}\quad
        \left(
            \widetilde k_{j,L}^{\mathfrak g},
            Z_{\mathrm{Sch},j}^{(L)}
        \right)_{L^2(I)}
    \]
    is bounded by \eqref{eq:Schur-inner-product-enclosure},
    using \eqref{eq:part5-Schur-rhs-total-bounds} and \eqref{eq:part5-Schur-finite-solve-bounds}.
    This gives the intervals in
    \eqref{eq:part5-Schur-F-sign-bound} and
    \eqref{eq:part5-Schur-F-abs-bound}.
\end{proof}

\begin{proof}[Proof of Proposition~\ref{prop:certified-Schur-inverse-bounds-g}]
    The signs in \eqref{eq:certified-Schur-complement-signs-g} follow directly
    from \eqref{eq:part5-Schur-F-sign-bound}.  Since
    \(\widetilde F_j^{\mathfrak g}\neq0\), Lemma~\ref{lem:schur-complement-bound-g-side}
    implies that \(S_0^{(\mathfrak g,j)}\) is invertible once \(\widetilde D_j^{\mathfrak g}\) is invertible. The estimate \eqref{eq:Dj-inverse-norm-Schur-bound} and the bounds in \eqref{eq:schur-D-inverse-bounds} give the invertibility of \(\widetilde D_j^{\mathfrak g}\).

    The inverse estimate \eqref{eq:schur-inverse-bound-g-side}
    gives
    \[
        \norm{
            (S_0^{(\mathfrak g,j)})^{-1}
        }_{L^2_{\mathrm e}(\bbR)\to L^2_{\mathrm e}(\bbR)}
        \leq
        \mathsf C_{D,\mathrm{inv}}^{(j)}
        +
        |\widetilde F_j^{\mathfrak g}|^{-1}
        \left(
            1+
            \norm{
                (\widetilde D_j^{\mathfrak g})^{-1}
                \widetilde k_j^{\mathfrak g}
            }_{L^2(I)}^2
        \right).
    \]
    Substituting the interval bounds in \eqref{eq:part5-Schur-finite-solve-bounds}-- \eqref{eq:part5-Schur-F-abs-bound} into this expression gives
    \[
        \begin{aligned}
            \norm{
                (S_0^{(\mathfrak g,1)})^{-1}
            }_{L^2_{\mathrm e}\to L^2_{\mathrm e}}
            &\leq
            50.852827
            <
            50.853,\\
            \norm{
                (S_0^{(\mathfrak g,2)})^{-1}
            }_{L^2_{\mathrm e}\to L^2_{\mathrm e}}
            &\leq
            5.644179
            <
            5.645.
        \end{aligned}
    \]
    This proves \eqref{eq:certified-S0-inverse-bounds-g}.
\end{proof}

\subsection{Limiting matrix entries for \(\mathfrak g\)}

Throughout this subsection we take
\[
    L=5000,\qquad
    L_Z=8000,\qquad
    N=65536,
    \qquad
    c_1=6.1,\qquad
    c_2=2.1.
\]
We use the notation from
\eqref{eq:g-side-Z-def}--\eqref{eq:g-side-pairing-error-bound}.  Thus
\[
    Z_m^{(j)}(\theta)
    =
    \sum_{k=0}^{L_Z}
    z_{m,k}^{(j)}\cos(2k\theta),
    \qquad
    m=0,1,2,\quad j=1,2.
\]
The coefficients \(z_{m,k}^{(j)}\) are cosine coefficients. They are
read as exact rational numbers from the following files:
\[
\begin{array}{c|c|l}
j & m & \textnormal{file} \\
\hline
1 & 0 & \texttt{Z\_c6p1\_solve\_rhs\_ag.csv} \\
1 & 1 & \texttt{Z\_c6p1\_solve\_rhs\_Psi\_g.csv} \\
1 & 2 & \texttt{Z\_c6p1\_solve\_rhs\_Psi\_Lambda\_g.csv} \\
2 & 0 & \texttt{Z\_c2p1\_solve\_rhs\_ag.csv} \\
2 & 1 & \texttt{Z\_c2p1\_solve\_rhs\_Psi\_Lambda\_g.csv} \\
2 & 2 & \texttt{Z\_c2p1\_solve\_rhs\_Psi\_Lambda2\_g.csv}
\end{array}
\]
For \(m,n\in\{0,1,2\}\), set
\begin{equation}\label{eq:part5-g-side-finite-pairing-def}
    \mathcal P_{mn,N}^{(j)}
    \coloneqq
    \left(
        B_{m,N}^{(j)},
        Z_n^{(j)}
    \right)_{L^2(I)} .
\end{equation}

\begin{proposition}[Interval bounds for the residuals and pairings associated with \(\mathfrak g\)]\label{prop:g-residual-pairing-intervals}
    All vector and matrix inequalities below are understood componentwise.  The
    following outward-rounded bounds hold.

    \begin{enumerate}[label=\textup{(\roman*)}]
        \item
        The intervals for \(\mathbf m^{(\mathfrak g,j)}\) satisfy
        \begin{equation}\label{eq:part5-g-side-free-m-intervals}
            \mathbf m^{(\mathfrak g,1)}
            \in
            \begin{pmatrix}
                [3.248201,3.248202]\\
                [-1.624101,-1.624100]
            \end{pmatrix},
            \qquad
            \mathbf m^{(\mathfrak g,2)}
            \in
            \begin{pmatrix}
                [-1.624101,-1.624100]\\
                [0.812050,0.812051]
            \end{pmatrix}.
        \end{equation}
        The intervals for \(\mathbf C_{\log}^{(\mathfrak g,j)}\) satisfy
        \begin{equation}\label{eq:part5-g-side-free-Clog-intervals-j1}
            \mathbf C_{\log}^{(\mathfrak g,1)}
            \in
            \begin{pmatrix}
                [-1.658925,-1.658924]
                &
                [2.508675,2.508676]
                \\
                [2.508675,2.508676]
                &
                [-1.428813,-1.428812]
            \end{pmatrix},
        \end{equation}
        and
        \begin{equation}\label{eq:part5-g-side-free-Clog-intervals-j2}
            \mathbf C_{\log}^{(\mathfrak g,2)}
            \in
            \begin{pmatrix}
                [-1.428813,-1.428812]
                &
                [1.134209,1.134210]
                \\
                [1.134209,1.134210]
                &
                [-0.430787,-0.430786]
            \end{pmatrix}.
        \end{equation}
        
        \item
        The DCT-I approximation errors for \(a_{\mathfrak g}\) and \(a_{\mathfrak g,\log}\) satisfy
        \begin{equation}\label{eq:part5-g-side-global-DCT-radii}
            \left(
                \delta_{a_{\mathfrak g},\mathrm{DCT}}(N),
                \delta_{a_{\mathfrak g,\log},\mathrm{DCT}}(N)
            \right)
            \leq
            \left(
                8.267\times10^{-7},
                4.450\times10^{-5}
            \right).
        \end{equation}
        For each \(Z_m^{(j)}\), the coefficient \(\ell^1\)-norm defined in \eqref{eq:app-g-side-Smj-def} and the DCT-I approximation-error bounds satisfy
        \begin{equation}\label{eq:g-Z-dependent-dct-bounds}
            \begin{array}{cc|cccc}
                j&m
                &
                \mathsf Z_{\ell^1,m}^{(j)}
                &
                \delta_{\mathcal K_{\sin},\mathrm{DCT}}^{(j,m)}(N)
                &
                \delta_{\mathcal I_0,\mathrm{DCT}}^{(j,m)}(N)
                &
                \delta_{\mathcal I_{\cos},\mathrm{DCT}}^{(j,m)}(N)
                \\
                \hline
                1&0&2.140&3.113\times10^{-8}&2.121\times10^{-8}&1.441\times10^{-7}\\
                1&1&0.685&9.957\times10^{-9}&6.782\times10^{-9}&4.608\times10^{-8}\\
                1&2&1.192&1.733\times10^{-8}&1.181\times10^{-8}&8.019\times10^{-8}\\
                2&0&7.337&1.068\times10^{-7}&7.270\times10^{-8}&4.939\times10^{-7}\\
                2&1&5.137&7.472\times10^{-8}&5.090\times10^{-8}&3.458\times10^{-7}\\
                2&2&7.618&1.109\times10^{-7}&7.548\times10^{-8}&5.128\times10^{-7}
            \end{array}
        \end{equation}

        \item
        The norms of \(B_{m,N}^{(j)}\) and the right-hand-side approximation errors satisfy
        \begin{equation}\label{eq:part5-g-side-rhs-approximation-bounds}
            \begin{array}{cc|cc}
                j&m
                &
                \norm{B_{m,N}^{(j)}}_{L^2(I)}
                &
                \delta_{\mathrm{rhs},m}^{(j)}(N)
                \\
                \hline
                1&0&1.1112&8.267\times10^{-7}\\
                1&1&0.46575&4.735\times10^{-5}\\
                1&2&0.73775&2.431\times10^{-5}\\
                2&0&1.1112&8.267\times10^{-7}\\
                2&1&0.73775&2.431\times10^{-5}\\
                2&2&1.0269&1.287\times10^{-5}
            \end{array}
        \end{equation}

        \item
        The DCT-I residual norms, residual-comparison errors, and exact residual norms satisfy
        \begin{equation}\label{eq:g-residual-comparison-bounds}
            \begin{array}{cc|ccc}
                j&m
                &
                \norm{\operatorname{Res}_{m,N}^{(j)}}_{L^2(I)}
                &
                \delta_{\mathrm{res\text{-}diff},m}^{(j)}(N)
                &
                \norm{\operatorname{Res}_m^{(j)}}_{L^2(I)}
                \\
                \hline
                1&0&2.277\times10^{-6}&7.994\times10^{-6}&1.028\times10^{-5}\\
                1&1&1.074\times10^{-4}&8.284\times10^{-5}&1.903\times10^{-4}\\
                1&2&1.611\times10^{-4}&3.462\times10^{-5}&1.957\times10^{-4}\\
                2&0&1.495\times10^{-5}&1.160\times10^{-4}&1.309\times10^{-4}\\
                2&1&1.613\times10^{-4}&1.015\times10^{-4}&2.628\times10^{-4}\\
                2&2&2.420\times10^{-4}&1.273\times10^{-4}&3.693\times10^{-4}
            \end{array}
        \end{equation}
        The norms of \(Z_m^{(j)}\) and the solution-error bounds satisfy
        \begin{equation}\label{eq:part5-g-side-solve-error-bounds}
            \begin{array}{cc|cc}
                j&m
                &
                \norm{Z_m^{(j)}}_{L^2(I)}
                &
                \varepsilon_m^{(j)}(N)
                \\
                \hline
                1&0&1.0593&5.223\times10^{-4}\\
                1&1&0.43769&9.673\times10^{-3}\\
                1&2&0.58978&9.952\times10^{-3}\\
                2&0&3.9211&7.387\times10^{-4}\\
                2&1&2.6526&1.483\times10^{-3}\\
                2&2&4.0052&2.085\times10^{-3}
            \end{array}
        \end{equation}

        \item
        For \(j=1\), the intervals for the finite pairings and the corresponding error bounds satisfy
        \begin{equation}\label{eq:part5-g-side-pairing-intervals-j1}
            \begin{array}{c|cc}
                (m,n)
                &
                \mathcal P_{mn,N}^{(1)}
                &
                \varepsilon_{\mathrm{pair},mn}^{(1)}(N)
                \\
                \hline
                (0,0)&[-0.052352,-0.052350]&5.812\times10^{-4}\\
                (1,0)&[-0.396637,-0.396636]&2.935\times10^{-4}\\
                (2,0)&[-0.100746,-0.100745]&4.111\times10^{-4}\\
                (1,1)&[0.097861,0.097862]&4.527\times10^{-3}\\
                (1,2)&[-0.240718,-0.240717]&4.664\times10^{-3}\\
                (2,1)&[-0.240785,-0.240784]&7.147\times10^{-3}\\
                (2,2)&[-0.085508,-0.085507]&7.357\times10^{-3}
            \end{array}
        \end{equation}
        For \(j=2\), the intervals for the finite pairings and the corresponding error bounds satisfy
        \begin{equation}\label{eq:part5-g-side-pairing-intervals-j2}
            \begin{array}{c|cc}
                (m,n)
                &
                \mathcal P_{mn,N}^{(2)}
                &
                \varepsilon_{\mathrm{pair},mn}^{(2)}(N)
                \\
                \hline
                (0,0)&[3.726834,3.726835]&8.240\times10^{-4}\\
                (1,0)&[2.496677,2.496678]&6.403\times10^{-4}\\
                (2,0)&[-3.702861,-3.702860]&8.090\times10^{-4}\\
                (1,1)&[1.684910,1.684911]&1.159\times10^{-3}\\
                (1,2)&[-2.488649,-2.488648]&1.636\times10^{-3}\\
                (2,1)&[-2.488656,-2.488655]&1.557\times10^{-3}\\
                (2,2)&[3.709099,3.709100]&2.192\times10^{-3}
            \end{array}
        \end{equation}
    \end{enumerate}
\end{proposition}

\begin{proof}
    The coefficients \(z_{m,k}^{(j)}\) are read from the six files above
    as exact rational numbers. All remaining operations are evaluated by
    outward-rounded interval arithmetic. The corresponding computation is
    row~07 of Table~\ref{tab:part5-certificate-layout}.

    \emph{Direct coefficients.}
    Lemma~\ref{lem:g-side-free-source-entries} expresses
    \(\mathbf m^{(\mathfrak g,j)}\) and
    \(\mathbf C_{\log}^{(\mathfrak g,j)}\) in terms of the cosine
    coefficients of \(U\). Evaluating these formulas gives
    \eqref{eq:part5-g-side-free-m-intervals},
    \eqref{eq:part5-g-side-free-Clog-intervals-j1}, and
    \eqref{eq:part5-g-side-free-Clog-intervals-j2}.

    \emph{DCT-I approximation errors.}
    Evaluating \eqref{eq:app-g-side-Smj-def} and
    \eqref{eq:app-g-side-delta-ag-DCT-def}--
    \eqref{eq:app-g-side-delta-Icos-DCT-def}
    gives
    \eqref{eq:part5-g-side-global-DCT-radii} and
    \eqref{eq:g-Z-dependent-dct-bounds}.
    Substitution into
    \eqref{eq:app-g-side-eta-rhs0-def},
    \eqref{eq:app-g-side-eta-rhsm-def}, and
    \eqref{eq:app-g-side-eta-res-def}, followed by
    \eqref{eq:app-g-side-B-minus-BN-bound} and
    \eqref{eq:app-g-side-Res-minus-ResN-bound}, gives the error bounds in
    \eqref{eq:part5-g-side-rhs-approximation-bounds} and
    \eqref{eq:g-residual-comparison-bounds}.

    \emph{Residual and solution bounds.}
    The \(L^2(I)\)-norms of
    \(\operatorname{Res}_{m,N}^{(j)}\) are evaluated from their finite
    cosine coefficients by Parseval's identity.
    The estimate \eqref{eq:app-g-side-Res-minus-ResN-bound} then gives the
    bounds for
    \(\norm{\operatorname{Res}_m^{(j)}}_{L^2(I)}\) in
    \eqref{eq:g-residual-comparison-bounds}.
    Combining \eqref{eq:g-side-X-minus-Z-residual-identity},
    \eqref{eq:certified-S0-inverse-bounds-g}, and
    \eqref{eq:g-side-X-Z-error-bound}, we obtain
    \eqref{eq:part5-g-side-solve-error-bounds}.

    \emph{Pairings.}
    The finite pairings are evaluated from
    \eqref{eq:part5-g-side-finite-pairing-def}.
    Applying \eqref{eq:g-side-pairing-error-bound} gives
    \eqref{eq:part5-g-side-pairing-intervals-j1} and
    \eqref{eq:part5-g-side-pairing-intervals-j2}.
\end{proof}

\begin{proof}[Proof of Proposition~\ref{prop:part4-certified-g-side-limiting-data}]
    The enclosure for \(\mathbf m^{(\mathfrak g,j)}\) in \eqref{eq:part4-g-side-m} follows from 
    \eqref{eq:part5-g-side-free-m-intervals}. 
    The enclosures for \(\mathbf C_{\log}^{(\mathfrak g,1)}\) and \(\mathbf C_{\log}^{(\mathfrak g,2)}\) in \eqref{eq:part4-g-side-C_log-j1} and \eqref{eq:part4-g-side-C_log-j2} follow from
    \eqref{eq:part5-g-side-free-Clog-intervals-j1} and \eqref{eq:part5-g-side-free-Clog-intervals-j2}.

    For each entry containing \((S_0^{(\mathfrak g,j)})^{-1}\), the exact pairing belongs to
    \[
        \mathcal P_{mn,N}^{(j)}
        +
        \left[
            -\varepsilon_{\mathrm{pair},mn}^{(j)}(N),
            \varepsilon_{\mathrm{pair},mn}^{(j)}(N)
        \right].
    \]
    Applying this to the entries of
    \(\mathfrak s_0^{(\mathfrak g,j)}\), \(\mathbf a^{(\mathfrak g,j)}\), and
    \(\mathbf R^{(\mathfrak g,j)}\), and then rounding outwards, gives
    \eqref{eq:part4-g-side-s-a-certified-j1},
    \eqref{eq:part4-g-side-R-certified-j1},
    \eqref{eq:part4-g-side-s-a-certified-j2}, and
    \eqref{eq:part4-g-side-R-certified-j2}.  For the off-diagonal entries of
    each displayed \(\mathbf R\)-matrix, we take the symmetric hull of the two directed
    pairings.

    Finally,
    \[
        \norm{X_m^{(\mathfrak g,j)}}_{L^2(I)}
        \leq
        \norm{Z_m^{(j)}}_{L^2(I)}
        +
        \varepsilon_m^{(j)}(N)
    \]
    follows from \eqref{eq:g-side-X-Z-error-bound}. Substituting the bounds in \eqref{eq:part5-g-side-solve-error-bounds} gives \eqref{eq:part4-g-side-Xnorm-certified}.
\end{proof}

\subsection{Limiting matrix entries for \(Q\) and determinant signs}

We use the notation for \(Q\) from \eqref{eq:equation for S_0j} and from the definitions preceding Proposition~\ref{prop:constraint-entry-formula}. Thus
\[
    \Xi_j=\{\Psi_1^{(j)},\Psi_2^{(j)}\},
    \qquad
    \psi_m^{(j)}
    =
    |V_Q|^{1/2}(\Psi_m^{(j)})^{\log},
    \qquad
    j=1,2,\quad m=1,2.
\]
By \eqref{eq:limiting constraint matrix}, the entries in the formula for \(\mathbf M_{0,j}\) are
\[
    \mathbf m_j,\quad
    \mathbf C_{\log,j},\quad
    \mathfrak s_j(0^+),\quad
    \mathbf a_j(0^+),\quad
    \mathbf R_j(0^+).
\]
Proposition~\ref{prop:part4-certified-g-side-limiting-data} gives interval enclosures for 
\[ 
    \mathbf m^{(\mathfrak g,j)},\quad 
    \mathbf C_{\log}^{(\mathfrak g,j)},
    \quad \mathfrak s_0^{(\mathfrak g,j)},
    \quad \mathbf a^{(\mathfrak g,j)},\quad 
    \mathbf R^{(\mathfrak g,j)}. 
\]
Recall that \(\mathbf M_{0,j}\) is the limiting matrix defined in
Proposition~\ref{prop:constraint-entry-formula}.

\begin{proposition}[Perturbation bounds for the limiting matrix entries]
    All vector and matrix inequalities below are understood componentwise.  The
    following bounds hold.

    \begin{enumerate}[label=\textup{(\roman*)}]
        \item
        The bounds for \(\mathbf m_j-\mathbf m^{(\mathfrak g,j)}\) in \eqref{eq:part4-Delta-m-transfer-def} satisfy
        \begin{equation}\label{eq:constraint-transfer-m-bounds}
            \begin{aligned}
                \left(
                    \Delta_{\mathbf m,1}^{(1)},
                    \Delta_{\mathbf m,2}^{(1)}
                \right)
                &\leq
                \left(
                    3.289\times10^{-8},
                    1.645\times10^{-8}
                \right),
                \\
                \left(
                    \Delta_{\mathbf m,1}^{(2)},
                    \Delta_{\mathbf m,2}^{(2)}
                \right)
                &\leq
                \left(
                    1.645\times10^{-8},
                    8.223\times10^{-9}
                \right).
            \end{aligned}
        \end{equation}
        The bounds for \(\mathbf C_{\log,j}-\mathbf C_{\log}^{(\mathfrak g,j)}\) in
        \eqref{eq:part4-Delta-Clog-transfer-def} satisfy
        \begin{equation}\label{eq:constraint-transfer-clog-bounds}
            \Delta_{\log}^{(1)}
            \coloneqq
            \left(
                \Delta_{\log,mn}^{(1)}
            \right)_{m,n=1}^{2}
            \leq
            \begin{pmatrix}
                1.830\times10^{-6}
                &
                1.211\times10^{-6}
                \\
                1.211\times10^{-6}
                &
                9.420\times10^{-7}
            \end{pmatrix},
        \end{equation}
        and
        \begin{equation}\label{eq:constraint-transfer-clog-bounds-j2}
            \Delta_{\log}^{(2)}
            \coloneqq
            \left(
                \Delta_{\log,mn}^{(2)}
            \right)_{m,n=1}^{2}
            \leq
            \begin{pmatrix}
                9.420\times10^{-7}
                &
                1.014\times10^{-6}
                \\
                1.014\times10^{-6}
                &
                1.534\times10^{-6}
            \end{pmatrix}.
        \end{equation}

        \item
        The \(L^2\)-bounds for \(\psi_m^{(j)}-\psi_m^{(\mathfrak g,j)}\) defined in
        \eqref{eq:app-Deltapsi-def} satisfy
        \begin{equation}\label{eq:constraint-transfer-psi-bounds}
            \begin{aligned}
                \left(
                    \Delta_{\psi,1}^{(1)},
                    \Delta_{\psi,2}^{(1)}
                \right)
                &\leq
                \left(
                    6.285\times10^{-6},
                    4.178\times10^{-6}
                \right),
                \\
                \left(
                    \Delta_{\psi,1}^{(2)},
                    \Delta_{\psi,2}^{(2)}
                \right)
                &\leq
                \left(
                    4.178\times10^{-6},
                    3.715\times10^{-6}
                \right).
            \end{aligned}
        \end{equation}
        In particular,
        \[
            \norm{
                \psi_m^{(j)}
                -
                \psi_m^{(\mathfrak g,j)}
            }_{L^2}
            \leq
            \Delta_{\psi,m}^{(j)},
            \qquad
            j=1,2,\quad m=1,2.
        \]

        \item
        The operator perturbation bound \(\Delta_T\) and the inverse bounds \(\mathsf C_{S,\mathrm{inv}}^{(\mathfrak g,j)}\), defined in \eqref{eq:part4-DeltaT-and-Neumann-factor-def} and \eqref{eq:inverse-transfer-factor-definition}, respectively, satisfy
        \begin{equation}\label{eq:part5-constraint-transfer-Neumann}
            \Delta_T
            \leq
            5.545\times10^{-6},
            \qquad
            \left(
                \mathsf C_{S,\mathrm{inv}}^{(\mathfrak g,1)},
                \mathsf C_{S,\mathrm{inv}}^{(\mathfrak g,2)}
            \right)
            \leq
            \left(
                50.853,
                5.645
            \right).
        \end{equation}
        Moreover,
        \begin{equation}\label{eq:part5-constraint-transfer-Neumann-denominator}
            \left(
                1-c_1\mathsf C_{S,\mathrm{inv}}^{(\mathfrak g,1)}\Delta_T,
                1-c_2\mathsf C_{S,\mathrm{inv}}^{(\mathfrak g,2)}\Delta_T
            \right)
            \geq
            \left(
                0.99828,
                0.99993
            \right),
        \end{equation}
        and therefore
        \begin{equation}\label{eq:constraint-transfer-neumann-factor-bounds}
            \left(
                \mathsf C_{\mathrm{transfer}}^{(1)},
                \mathsf C_{\mathrm{transfer}}^{(2)}
            \right)
            \leq
            \left(
                1.00173,
                1.00007
            \right).
        \end{equation}

        \item
        The constants \(\Delta_{\mathfrak s}^{(j)}\) defined in \eqref{eq:part4-Delta-s-transfer-def} obey
        \begin{equation}\label{eq:constraint-transfer-s-bounds}
            \left(
                \Delta_{\mathfrak s}^{(1)},
                \Delta_{\mathfrak s}^{(2)}
            \right)
            \leq
            \left(
                4.099\times10^{-5},
                1.900\times10^{-4}
            \right).
        \end{equation}
        The constants \(\Delta_{\mathfrak a,m}^{(j)}\) defined in \eqref{eq:part4-Delta-a-transfer-def} obey
        \begin{equation}\label{eq:constraint-transfer-a-bounds}
            \begin{aligned}
                \left(
                    \Delta_{\mathfrak a,1}^{(1)},
                    \Delta_{\mathfrak a,2}^{(1)}
                \right)
                &\leq
                \left(
                    2.336\times10^{-5},
                    2.680\times10^{-5}
                \right),
                \\
                \left(
                    \Delta_{\mathfrak a,1}^{(2)},
                    \Delta_{\mathfrak a,2}^{(2)}
                \right)
                &\leq
                \left(
                    1.413\times10^{-4},
                    2.032\times10^{-4}
                \right).
            \end{aligned}
        \end{equation}

        \item
        The constants \(\Delta_{\mathfrak r,mn}^{(j)}\) defined in \eqref{eq:part4-Delta-r-transfer-def} obey
        \begin{equation}\label{eq:constraint-transfer-r-bounds-j1}
            \left(
                \Delta_{\mathfrak r,mn}^{(1)}
            \right)_{m,n=1}^{2}
            \leq
            \begin{pmatrix}
                1.242\times10^{-5}
                &
                1.474\times10^{-5}
                \\
                1.474\times10^{-5}
                &
                1.721\times10^{-5}
            \end{pmatrix},
        \end{equation}
        and
        \begin{equation}\label{eq:constraint-transfer-r-bounds-j2}
            \left(
                \Delta_{\mathfrak r,mn}^{(2)}
            \right)_{m,n=1}^{2}
            \leq
            \begin{pmatrix}
                1.042\times10^{-4}
                &
                1.505\times10^{-4}
                \\
                1.505\times10^{-4}
                &
                2.168\times10^{-4}
            \end{pmatrix}.
        \end{equation}
    \end{enumerate}
\end{proposition}

\begin{proof}
    All right-hand sides below are evaluated by outward-rounded interval
    arithmetic using
    Proposition~\ref{prop:nk-Q-transfer-bounds} and
    \eqref{eq:certified-S0-inverse-bounds-g}.
    The corresponding computation is row~08 of
    Table~\ref{tab:part5-certificate-layout}.

    Evaluating
    \eqref{eq:part4-Delta-m-transfer-def} and
    \eqref{eq:part4-Delta-Clog-transfer-def}
    gives
    \eqref{eq:constraint-transfer-m-bounds},
    \eqref{eq:constraint-transfer-clog-bounds}, and
    \eqref{eq:constraint-transfer-clog-bounds-j2}.
    The associated componentwise estimates are
    \eqref{eq:part4-transfer-m-vector-bound} and
    \eqref{eq:part4-transfer-Clog-matrix-bound}.

    Evaluating \eqref{eq:app-Deltapsi-def} gives
    \eqref{eq:constraint-transfer-psi-bounds}, and
    \eqref{eq:app-psi-difference-bound}
    gives the corresponding \(L^2\)-bounds for
    \(\psi_m^{(j)}-\psi_m^{(\mathfrak g,j)}\).

    The estimate \eqref{eq:TQ-Tg-operator-bound} gives the bound for
    \(\Delta_T\) in \eqref{eq:part5-constraint-transfer-Neumann}, while
    \eqref{eq:certified-S0-inverse-bounds-g} gives the bounds for
    \(\mathsf C_{S,\mathrm{inv}}^{(\mathfrak g,j)}\).
    Substitution into \eqref{eq:inverse-transfer-factor-definition} gives
    \eqref{eq:part5-constraint-transfer-Neumann-denominator} and
    \eqref{eq:constraint-transfer-neumann-factor-bounds}.
    In particular, both Neumann denominators are positive.

    Finally, evaluating
    \eqref{eq:part4-Delta-s-transfer-def},
    \eqref{eq:part4-Delta-a-transfer-def}, and
    \eqref{eq:part4-Delta-r-transfer-def}
    gives
    \eqref{eq:constraint-transfer-s-bounds},
    \eqref{eq:constraint-transfer-a-bounds},
    \eqref{eq:constraint-transfer-r-bounds-j1}, and
    \eqref{eq:constraint-transfer-r-bounds-j2}.
    The resulting perturbation estimates are
    \eqref{eq:part4-transfer-s-bound},
    \eqref{eq:part4-transfer-a-bound}, and
    \eqref{eq:part4-transfer-r-bound}.
\end{proof}

\begin{proof}[Proof of Proposition~\ref{prop:part4-certified-limiting-constraint-matrices}]
    The interval enclosures
    \eqref{eq:part5-g-side-free-m-intervals},
    \eqref{eq:part5-g-side-free-Clog-intervals-j1}, and
    \eqref{eq:part5-g-side-free-Clog-intervals-j2},
    together with
    \eqref{eq:constraint-transfer-m-bounds},
    \eqref{eq:constraint-transfer-clog-bounds}, and
    \eqref{eq:constraint-transfer-clog-bounds-j2},
    give componentwise enclosures for
    \(\mathbf m_j\) and \(\mathbf C_{\log,j}\).

    Similarly,
    \eqref{eq:part4-g-side-s-a-certified-j1}--
    \eqref{eq:part4-g-side-R-certified-j2},
    together with
    \eqref{eq:constraint-transfer-s-bounds},
    \eqref{eq:constraint-transfer-a-bounds},
    \eqref{eq:constraint-transfer-r-bounds-j1}, and
    \eqref{eq:constraint-transfer-r-bounds-j2},
    give componentwise enclosures for
    \(\mathfrak s_j(0^+)\),
    \(\mathbf a_j(0^+)\), and
    \(\mathbf R_j(0^+)\).
    Hence
    \[
        \mathfrak s_1(0^+)\in[-0.05298,-0.05172],
        \qquad
        \mathfrak s_2(0^+)\in[3.72581,3.72786],
    \]
    which proves \eqref{eq:part4-certified-limiting-s-intervals}.

    We substitute the interval enclosures for
    \[
        \mathbf m_j,\qquad
        \mathbf C_{\log,j},\qquad
        \mathfrak s_j(0^+),\qquad
        \mathbf a_j(0^+),\qquad
        \mathbf R_j(0^+)
    \]
    into the formula in
    Proposition~\ref{prop:constraint-entry-formula}.
    We replace the two directed off-diagonal enclosures by their symmetric
    hull. This gives
    \[
        \mathbf M_{0,1}
        \in
        \begin{pmatrix}
            [1.025,1.153]
            &
            [-4.905,-4.636]
            \\
            [-4.905,-4.636]
            &
            [13.47,14.02]
        \end{pmatrix},
    \]
    and
    \[
        \mathbf M_{0,2}
        \in
        \begin{pmatrix}
            [0.431,0.441]
            &
            [-0.879,-0.864]
            \\
            [-0.879,-0.864]
            &
            [1.151,1.173]
        \end{pmatrix}.
    \]
    This proves
    \eqref{eq:part4-certified-M0-j1-interval} and
    \eqref{eq:part4-certified-M0-j2-interval}.
\end{proof}

\appendix
\addcontentsline{toc}{part}{Appendices}

\section{Modulation decomposition and auxiliary bounds}

\subsection{Modulation decomposition near the ground state}

\begin{proof}[Proof of Lemma~\ref{lem:orthogonal decomposition of the solution u}]
    For \(t\in[0,T)\), define
    \begin{equation}\label{eq:definition for lambda_0 and tilde u}
        \lambda_0(t)
        \coloneqq
        \left(
            \tfrac{
                \norm{Q}_{\dot H^{1/2}}
            }{
                \norm{u(t)}_{\dot H^{1/2}}
            }
        \right)^2,
        \qquad
        \widetilde u(t,y)
        \coloneqq
        \lambda_0(t)^{1/2}
        u\bigl(t,\lambda_0(t)y\bigr).
    \end{equation}
    From \eqref{eq:definition for lambda_0 and tilde u}, \eqref{eq:condition1}, and conservation of mass and energy, we obtain \(\norm{\widetilde u(t)}_{\dot H^{1/2}} = \norm{Q}_{\dot H^{1/2}}\), \(\norm{\widetilde u(t)}_{L^2}^2 = \norm{Q}_{L^2}^2+\alpha_0\) and \(E(\widetilde u(t)) = \lambda_0(t)E_0 < 0\). We first prove that there is a function \(\delta(\alpha)\to0\) as \(\alpha\to0\) such that
    \begin{equation}\label{eq:tube stability aux}
        \inf_{\theta\in\bbR}
        \norm{
            e^{i\theta}\widetilde u(t)-Q
        }_{H^{1/2}}
        \leq
        \delta(\alpha_0).
    \end{equation}
    Suppose otherwise. Then there exist \(\alpha_n\to0\), \(\widetilde u_n\in H^{1/2}_{\mathrm e}\), and \(\varepsilon_0>0\) such that \(\norm{\widetilde u_n}_{\dot H^{1/2}} = \norm{Q}_{\dot H^{1/2}}\), \(\norm{\widetilde u_n}_{L^2}^2 = \norm{Q}_{L^2}^2+\alpha_n\), \(E(\widetilde u_n)<0\) and 
    \begin{equation}\label{eq:tube subseq}
        \begin{gathered}
            \inf_{\theta\in\bbR}
            \norm{
                e^{i\theta}\widetilde u_n-Q
            }_{H^{1/2}}
            \geq
            \varepsilon_0.
        \end{gathered}
    \end{equation}
    The estimate \eqref{eq:energy-lower-bound} and the negative-energy condition give \(E(\widetilde u_n)\to0\), and therefore \(\{\widetilde u_n\}\) is an extremizing sequence for \eqref{eq:sharp-Gagliardo-Nirenberg}.

    Lieb's compactness lemma, the Brezis--Lieb lemma, and the characterization of the optimizers in \cite{FrankLenzmann2013Acta} give phases \(\beta_n\in\bbR\) and translations
    \(x_n\in\bbR\) such that, after passing to a subsequence, \(e^{i\beta_n} \widetilde u_n(\cdot+x_n)\to Q\) strongly in \(H^{1/2}\).

    We now claim that\(\{x_n\}\) is bounded. Indeed, if \(|x_n|\to\infty\), then the evenness of
    \(\widetilde u_n\) produces two disjoint concentration regions centered at \(x_n\) and \(-x_n\). Strong convergence near each center would then imply \(\liminf_{n\to\infty}
    \norm{\widetilde u_n}_{L^2}^2 > \norm{Q}_{L^2}^2\), contradicting \(\norm{\widetilde u_n}_{L^2}^2 \to\norm{Q}_{L^2}^2\). Passing to a further subsequence, let \(x_n\to x_\infty\). Then \(e^{i\beta_n}\widetilde u_n \to Q(\cdot-x_\infty)\) in \(H^{1/2}\).
    The limit is even. Since \(Q\) is even and strictly decreasing on
    \((0,\infty)\), this implies \(x_\infty=0\). We conclude that
    \(\inf_{\theta\in\bbR} \norm{e^{i\theta}\widetilde u_n-Q}_{H^{1/2}}\to0\),
    which contradicts \eqref{eq:tube subseq}. This proves \eqref{eq:tube stability aux}. We next choose the phase continuously. Set \( z(t) \coloneqq (\widetilde u(t),Q)_{L^2}\).
    By \eqref{eq:tube stability aux}, \(|z(t)|\geq \norm{Q}_{L^2}^2-C\delta(\alpha_0)\). After decreasing \(\alpha^*>0\), we have \(z(t)\neq0\) on
    \([0,T)\). Hence, there is a continuous function \(\gamma_0:[0,T)\to\bbR\) such that
    \(e^{-i\gamma_0(t)} = \tfrac{z(t)}{|z(t)|}\). 
    In particular, \((e^{i\gamma_0(t)}\widetilde u(t),iQ)_r=0\).

    If \(\theta(t)\) is chosen so that the left-hand side of
    \eqref{eq:tube stability aux} is at most \(2\delta(\alpha_0)\), then \(e^{i\theta(t)}z(t)=\norm{Q}_{L^2}^2+\mathcal{O}\bigl(\delta(\alpha_0)\bigr)\).
    Hence \(|e^{i(\gamma_0(t)-\theta(t))}-1|\lesssim \delta(\alpha_0)\), and consequently
    \begin{equation}\label{eq:tube stability}
        \norm{
            e^{i\gamma_0(t)}
            \widetilde u(t)-Q
        }_{H^{1/2}}
        \leq
        C\delta(\alpha_0).
    \end{equation}

    Fix \(\delta>0\), and define the open neighborhood \(U_\delta \coloneqq \{v\in H^{1/2}:
    \norm{v-Q}_{H^{1/2}}<\delta\}\). For \(v\in U_\delta\), \(\lambda_1>0\), and
    \(\gamma_1\in\bbR\), set \(\epsilon_{\lambda_1,\gamma_1}[v] \coloneqq e^{i\gamma_1}\lambda_1^{1/2}v(\lambda_1\cdot)-Q\), and define
    \[
        \mathbf F(\lambda_1,\gamma_1;v)
        \coloneqq
        \begin{pmatrix}
            \bigl(
                \epsilon_{\lambda_1,\gamma_1}[v],
                \Lambda Q
            \bigr)_r
            \\
            \bigl(
                \epsilon_{\lambda_1,\gamma_1}[v],
                i\Lambda^2Q
            \bigr)_r
        \end{pmatrix}.
    \]
    At \((\lambda_1,\gamma_1;v)=(1,0;Q)\),
    \[
        \partial_{(\lambda_1,\gamma_1)}
        \mathbf F(1,0;Q)
        =
        \begin{pmatrix}
            \norm{\Lambda Q}_{L^2}^2
            &
            0
            \\
            0
            &
            -\norm{\Lambda Q}_{L^2}^2
        \end{pmatrix}.
    \]
    The implicit function theorem therefore gives, after decreasing
    \(\delta>0\), unique \(C^1\) functions \(\lambda_1:U_\delta\to\bbR_+\) and \(\gamma_1:U_\delta\to\bbR\) such that
    \begin{equation}\label{eq:IFT}
        \mathbf F
        \bigl(
            \lambda_1(v),
            \gamma_1(v);
            v
        \bigr)
        =
        0
    \end{equation}
    for every \(v\in U_\delta\). Moreover,
    \begin{equation}\label{eq:IFT smallness}
        \left|
            \lambda_1(v)-1
        \right|
        +
        \left|
            \gamma_1(v)
        \right|
        +
        \norm{
            \epsilon_{\lambda_1(v),\gamma_1(v)}[v]
        }_{H^{1/2}}
        \lesssim
        \norm{v-Q}_{H^{1/2}}.
    \end{equation}

    Take \(\alpha^*>0\) sufficiently small that \(C\delta(\alpha^*)<\delta\), and set \(v(t)\coloneqq e^{i\gamma_0(t)}\widetilde u(t)\). Define \(\lambda(t)\coloneqq \lambda_0(t)\lambda_1(v(t))\) and \(\gamma(t)\coloneqq \gamma_0(t)+\gamma_1(v(t))\).
    Then \(\epsilon(t,y)\) defined in \eqref{eq:definition of epsilon(t)} satisfies \eqref{eq:orthogonal condition on ep} by \eqref{eq:IFT}. The \(C^1\)-regularity of the modulation parameters follows by differentiating the orthogonality conditions
    \eqref{eq:orthogonal condition on ep}. Finally, \(\lambda(t) \tfrac{\norm{u(t)}_{\dot H^{1/2}}^2}{\norm{Q}_{\dot H^{1/2}}^2} = \lambda_1(v(t))\). Combining \eqref{eq:tube stability}, \eqref{eq:IFT smallness}, with the preceding identity, we obtain \eqref{eq:smallenss of ep and lambda esitmate}.
\end{proof}

\subsection{A self-similar lower bound}

\begin{lemma}[Self-similar lower bound]\label{lem:trivial lower bound}
    Let \(u\in C([0,T);H^{1/2}(\bbR))\) be the maximal-lifespan solution given by
    Theorem~\ref{thm:Cauchy-theory}. If \(T<+\infty\), then
    \begin{equation}\label{eq:trivial self-similar lower bound}
        \norm{|D|^{1/2}u(t)}_{L^2}
        \geq \frac{C}{\sqrt{T-t}},
        \qquad t\in [0,T).
    \end{equation}
    Consequently,
    \begin{equation}\label{eq:lambda lower bound from trivial rate}
        \lambda(t)\lesssim T-t.
    \end{equation}
\end{lemma}

\begin{proof}
    Fix \(t\in[0,T)\), and set \(\widetilde\rho(t) \coloneqq \norm{|D|^{1/2}u(t)}_{L^2}^{-2}\).
    Define 
    \[
        v^t(\tau,z)
        \coloneqq
        \widetilde\rho(t)^{1/2}
        u\bigl(
            t+\widetilde\rho(t)\tau,
            \widetilde\rho(t)z
        \bigr).
    \]
    By the \(L^2\)-critical scaling, \(v^t\) is a solution of \eqref{eq:halfwave}, and
    \( \norm{|D|^{1/2}v^t(0)}_{L^2}=1 \),  \(\norm{v^t(0)}_{L^2} = \norm{u_0}_{L^2}\).
    Hence \(\norm{v^t(0)}_{H^{1/2}} \leq  \norm{u_0}_{L^2}+1\). Theorem~\ref{thm:Cauchy-theory} gives \(\tau_0>0\), depending only on \(\norm{u_0}_{L^2}\), such that \(v^t\) is defined on
    \([0,\tau_0]\). Therefore, \(u\) is defined on \([t,t+\widetilde\rho(t)\tau_0]\). By maximality of \(T\), we have \(t+\widetilde\rho(t)\tau_0 \leq T\). Thus, from the definition of \(\widetilde \rho(t)\), we conclude \eqref{eq:trivial self-similar lower bound}.

    Let \(\lambda(t)\) now denote the modulation parameter from Lemma~\ref{lem:orthogonal decomposition of the solution u}. By \eqref{eq:smallenss of ep and lambda esitmate}, after decreasing \(\alpha^*>0\), modulation decomposition gives \(\norm{|D|^{1/2}u(t)}_{L^2}^2 \sim \lambda(t)^{-1}\). Consequently,
    \[
        \lambda(t)
        \sim
        \norm{
            |D|^{1/2}u(t)
        }_{L^2}^{-2}
        \lesssim
        T-t,
    \]
    which proves \eqref{eq:lambda lower bound from trivial rate}.
\end{proof}

\subsection{Tail lower bounds for the ground state}

\begin{lemma}[Tail lower bounds for \(Q\)]\label{lem:tail-lower-yQpQ}
    There exist constants \(c>0\) and \(R>0\) such that
    \[
        Q(y)\ge c|y|^{-2}, \qquad |Q'(y)|\ge c|y|^{-3}, \qquad |y|\ge R.
    \]
\end{lemma}

\begin{proof}
    Let \(G\) be the kernel of \((|D|+1)^{-1}\), so that \(Q=G*Q^3\).
    The identity
    \(
        (|\xi|+1)^{-1}
        =
        \int_0^\infty e^{-\tau}e^{-\tau|\xi|}\,d\tau
    \)
    and the Poisson kernel formula give
    \[
        G(y)
        =
        c_0\int_0^\infty e^{-\tau}
        \frac{\tau}{\tau^2+y^2}\,d\tau,
        \qquad
        G'(y)
        =
        -2c_0y\int_0^\infty e^{-\tau}
        \frac{\tau}{(\tau^2+y^2)^2}\,d\tau .
    \]
    In particular, \(G(y)\ge cy^{-2}\) and \(-G'(y)\ge cy^{-3}\) for
    \(y\ge1\). Fix \(R_0>0\) such that \(m_0\coloneqq\int_{|z|\le R_0}Q(z)^3\,dz>0\). For \(y\ge2R_0+1\), positivity of \(G\) and \(y-z\sim y\) on \(|z|\le R_0\) imply
    \[
        Q(y)
        \ge
        \int_{|z|\le R_0}G(y-z)Q(z)^3\,dz
        \ge
        cm_0y^{-2}.
    \]
    Since \(G'\) is odd and \(Q\) is even and decreasing on
    \([0,\infty)\),
    \[
        -Q'(y)
        =
        \int_0^\infty
        \bigl(-G'(\eta)\bigr)
        \bigl(Q(y-\eta)^3-Q(y+\eta)^3\bigr)\,d\eta ,
        \qquad y>0,
    \]
    and the integrand is nonnegative. Hence, for \(y\ge2R_0+1\),
    \[
        -Q'(y)
        \ge
        cy^{-3}
        \int_{-R_0}^{R_0}
        \bigl(Q(z)^3-Q(2y-z)^3\bigr)\,dz .
    \]
    The last integral converges to \(m_0\) as \(y\to\infty\), since \(Q^3\in L^1(\bbR)\). Thus \(-Q'(y)\ge cy^{-3}\) for all sufficiently large \(y\). The estimates on \((-\infty,0)\) follow from the parity of \(Q\) and \(Q'\).
\end{proof}

\section{Resolvent kernel estimates}

\subsection{Small-\(\mu\) resolvent decomposition}

Let \(G^{(\mu)}\) denote the integral kernel of \((|D|+\mu)^{-1}\). By Fourier inversion,
\begin{equation}\label{eq:Fourier-side-Kmu}
    G^{(\mu)}(y)= \frac{1}{2\pi}\int_{\bbR} \frac{\cos(\xi y)}{|\xi|+\mu}\,d\xi
    = \frac{1}{\pi}\int_0^\infty \frac{\cos(\xi y)}{\xi+\mu}\,d\xi.
\end{equation}
Plancherel's identity and two integrations by parts give
\begin{equation}\label{eq:Kmu-basic-bounds}
    \norm{G^{(\mu)}}_{L^2}^2
    = \frac{1}{2\pi}\int_{\bbR} \frac{1}{(|\xi|+\mu)^2}\,d\xi
    = \frac{1}{\pi\mu},
    \qquad
    |G^{(\mu)}(y)|\leq \frac{2}{\pi\mu^2}\frac{1}{y^2}
    \quad (y\neq 0),
\end{equation}
We use the standard cosine and sine integrals
\begin{equation*}
    \Ci(\eta)\coloneqq -\int_\eta^\infty \frac{\cos \zeta}{\zeta}\,d\zeta,
    \qquad
    \Si(\eta)\coloneqq \int_0^\eta \frac{\sin \zeta}{\zeta}\,d\zeta
    \qquad (\eta>0).
\end{equation*}

\begin{lemma}[A decomposition for small \(\mu\)]
\label{lem:Kmu-smallmu}
    For \(y\in\bbR\setminus\{0\}\) and \(\mu>0\),
    \begin{equation}\label{eq:Kmu-smallmu}
        G^{(\mu)}(y)
        =
        \frac{1}{\pi}
        \left(\log\frac{1}{\mu}-\gamma_E-\log|y|+r_\mu(y)\right),
    \end{equation}
    where
    \begin{equation}\label{eq:definition of r_mu}
        r_\mu(y)
        \coloneqq
        \int_0^\infty
        (1-\cos(\xi y))
        \frac{\mu}{\xi(\xi+\mu)}\,d\xi .
    \end{equation}
    Moreover,
    \begin{equation}\label{eq:rmu-global-bound}
        0\le r_\mu(y)
        \lesssim
        \min\{\mu|y|,\,1+\log(1+\mu|y|)\}.
    \end{equation}
    Thus \(r_\mu\) extends continuously to \(y=0\) by \(r_\mu(0)=0\).
\end{lemma}

\begin{proof}
    For \(R>0\), write
    \[
        \frac{\cos(\xi y)}{\xi+\mu}
        =
        \frac{1}{\xi+\mu}
        +
        \frac{\cos(\xi y)-1}{\xi}
        +
        (1-\cos(\xi y))
        \frac{\mu}{\xi(\xi+\mu)}.
    \]
    The cosine integral formula gives
    \[
        \int_0^R\frac{\cos(\xi y)}{\xi+\mu}\,d\xi
        =
        \log\frac{R+\mu}{\mu}
        -\gamma_E-\log(R|y|)
        +\Ci(R|y|)
        +
        \int_0^R
        (1-\cos(\xi y))
        \frac{\mu}{\xi(\xi+\mu)}\,d\xi .
    \]
    The last integral converges absolutely. Letting \(R\to\infty\) proves \eqref{eq:Kmu-smallmu}. The identity
    \(
        \frac{\mu}{\xi(\xi+\mu)}
        =
        \int_0^\infty e^{-\xi\tau}(1-e^{-\mu\tau})\,d\tau
    \)
    and Tonelli's theorem give \eqref{eq:definition of r_mu}.
    Hence \(r_\mu(y)\ge0\), and \(1-e^{-s}\le s\) gives
    \(r_\mu(y)\lesssim\mu|y|\). Also,
    \[
        r_\mu(y)
        \le
        \int_0^1
        \frac{\min\{\mu|y|\zeta,1\}}{\zeta}\,d\zeta
        +
        \int_1^\infty\frac{d\zeta}{\zeta(1+\zeta^2)}
        \lesssim
        1+\log(1+\mu|y|).
    \]
    This proves \eqref{eq:rmu-global-bound}. The continuous extension at
    \(y=0\) follows from \eqref{eq:rmu-global-bound}.
\end{proof}

\subsection{Fixed-\(\mu\) local kernel representation}

For \(\mu\in\{\frac12,\frac32\}\), we separate the diagonal singularity by writing
\[
    G^{(\mu)}(y)
    =
    -\frac1\pi\log|y|+R_\mu(y).
\]
Define the rescaled even kernel \(G\) by
\begin{equation}\label{eq:K-rescaled}
    G(\eta)\coloneqq G^{(\mu)}\left(\frac{\eta}{\mu}\right)=G^{(1)}(\eta)
    \qquad (\eta\neq 0).
\end{equation}
Equivalently, on the real axis,
\begin{equation}\label{eq:G-poisson-real}
    G(y)
    =
    \frac1\pi
    \int_0^\infty e^{-\tau}
    \frac{\tau}{\tau^2+y^2}\,d\tau,
    \qquad y\neq0.
\end{equation}
This follows from \eqref{eq:Fourier-side-Kmu} with \(\mu=1\), the Laplace
representation of \((1+\xi)^{-1}\), and the Poisson kernel identity.

By \eqref{eq:K-rescaled} and \eqref{eq:Fourier-side-Kmu}, changing variables
first by \(\zeta=1+\xi\) and then by \(\tau=\eta\zeta\) gives
\begin{equation}\label{eq:K-CiSi}
    G(\eta)
    =
    \frac{1}{\pi}
    \left[
        -\Ci(\eta)\cos\eta
        +
        \left(
            \frac{\pi}{2}-\Si(\eta)
        \right)\sin\eta
    \right],
    \qquad \eta>0.
\end{equation}

\begin{lemma}[Local logarithmic form for fixed \(\mu\)]
\label{lem:Kmu-local}
    Let \(\mu\in\{\frac12,\frac32\}\). For \(y\in\bbR\setminus\{0\}\),
    \begin{equation}\label{eq:Kmu-local}
        G^{(\mu)}(y)
        =
        -\frac1\pi\log|y|+R_\mu(y),
    \end{equation}
    where
    \begin{equation}\label{eq:def-Rmu}
        R_\mu(y)
        =
        -\frac{\gamma_E+\log\mu}{\pi}
        +\frac{\mu|y|}{2}
        +\mathcal E(\mu|y|),
        \qquad
        \mathcal E(\eta)
        =
        G(\eta)
        -
        \left(
            \frac{\eta}{2}
            -
            \frac{\gamma_E+\log\eta}{\pi}
        \right).
    \end{equation}
    The function \(R_\mu\) extends continuously to \(y=0\), with
    \[
        R_\mu(0)
        =
        -\frac{\gamma_E+\log\mu}{\pi}.
    \]
    If \(0<\mu|y|\le1\), then
    \begin{equation}\label{eq:Rmu-E-bound}
        \left|
            R_\mu(y)-R_\mu(0)-\frac{\mu|y|}{2}
        \right|
        \le
        (\mu|y|)^2
        \left(
            \frac{|\gamma_E+\log(\mu|y|)|}{2\pi}
            +
            \mathsf C_R^{\mathrm{loc}}
        \right),
        \quad
        \mathsf C_R^{\mathrm{loc}}
        =
        \frac{1}{12}+\frac{53}{36\pi}.
    \end{equation}
    If \(y>0\) and \(0<\mu y\le1\), then
    \begin{equation}\label{eq:Rmu-derivative-bound}
        \left|
            R_\mu'(y)-\frac{\mu}{2}
        \right|
        \le
        \mu^2y
        \left(
            \frac{1+\gamma_E+|\log(\mu y)|}{\pi}
            +
            \mathsf C_{R'}^{\mathrm{loc}}
        \right),\quad
        \mathsf C_{R'}^{\mathrm{loc}}
        =
        \frac14+\frac{11}{36\pi}.
    \end{equation}
\end{lemma}

\begin{proof}
    Set \(\eta=\mu|y|\). By \eqref{eq:K-CiSi} and the definition of
    \(\mathcal E\),
    \[
        G^{(\mu)}(y)
        =
        -\frac1\pi\log|y|
        -
        \frac{\gamma_E+\log\mu}{\pi}
        +
        \frac{\mu|y|}{2}
        +
        \mathcal E(\mu|y|),
    \]
    which proves \eqref{eq:Kmu-local}--\eqref{eq:def-Rmu}. For \(0<\eta\le1\), the standard Taylor bounds are
    \[
        \left|\Ci(\eta)-\gamma_E-\log\eta\right|
        \le\frac{\eta^2}{4},
        \quad
        |\Si(\eta)-\eta|
        \le\frac{\eta^3}{18},
        \quad
        |\cos\eta-1|
        \le\frac{\eta^2}{2},
        \quad
        |\sin\eta-\eta|
        \le\frac{\eta^3}{6}.
    \]
    Substitution into \eqref{eq:K-CiSi} gives
    \[
        \begin{aligned}
        \mathcal E(\eta)
        ={}&
        -\frac{\Ci(\eta)-\gamma_E-\log\eta}{\pi}\cos\eta
        -\frac{\gamma_E+\log\eta}{\pi}(\cos\eta-1) \\
        &+
        \left(\frac12-\frac{\eta}{\pi}\right)(\sin\eta-\eta)
        -\frac{\Si(\eta)-\eta}{\pi}\sin\eta
        -\frac{\eta^2}{\pi}.
        \end{aligned}
    \]
    Hence
    \[
        \begin{aligned}
        |\mathcal E(\eta)|
        &\le
        \eta^2
        \left(
            \frac{|\gamma_E+\log\eta|}{2\pi}
            +
            \frac{1}{4\pi}
            +
            \frac{1}{12}
            +
            \frac{1}{6\pi}
            +
            \frac{1}{18\pi}
            +
            \frac{1}{\pi}
        \right) \\
        &=
        \eta^2
        \left(
            \frac{|\gamma_E+\log\eta|}{2\pi}
            +
            \mathsf C_R^{\mathrm{loc}}
        \right).
        \end{aligned}
    \]
    This proves the continuous extension and \eqref{eq:Rmu-E-bound}. Differentiating \eqref{eq:K-CiSi} and subtracting \(\frac12-\frac1{\pi\eta}\), we obtain
    \[
        \begin{aligned}
        \mathcal E'(\eta)
        ={}&
        \frac{(\gamma_E+\log\eta)\sin\eta-\eta\cos\eta}{\pi} \\
        &+
        \frac{
            (\Ci(\eta)-\gamma_E-\log\eta)\sin\eta
            -
            (\Si(\eta)-\eta)\cos\eta
        }{\pi}
        +
        \frac{\cos\eta-1}{2}.
        \end{aligned}
    \]
    The preceding Taylor bounds imply
    \[
        \begin{aligned}
        |\mathcal E'(\eta)|
        \le
        \frac{\eta}{\pi}
        \bigl(1+\gamma_E+|\log\eta|\bigr)
        +
        \frac{\eta^3}{4\pi}
        +
        \frac{\eta^3}{18\pi}
        +
        \frac{\eta^2}{4} 
        \le
        \eta
        \left(
            \frac{1+\gamma_E+|\log\eta|}{\pi}
            +
            \mathsf C_{R'}^{\mathrm{loc}}
        \right).
        \end{aligned}
    \]
    Since
    \[
        R_\mu'(y)
        =
        \frac{\mu}{2}
        +
        \mu\mathcal E'(\mu y),
        \qquad y>0,
    \]
    this proves \eqref{eq:Rmu-derivative-bound}.
\end{proof}

\begin{corollary}[Local bounds for the regular part]
\label{cor:local-rmu-bounds}
    Let \(\mu\in\{\frac12,\frac32\}\), and let \(0<r\le\mu^{-1}\). Set
    \begin{equation}\label{eq:local-rmu-bound-definitions}
    \begin{aligned}
        \mathsf R_{\mu}^{\mathrm{loc}}(r)
        &\coloneqq
        \frac{|\gamma_E+\log\mu|}{\pi}
        +
        \frac{\mu r}{2}
        +
        (\mu r)^2
        \left(
            \frac{|\gamma_E+\log(\mu r)|}{2\pi}
            +
            \mathsf C_R^{\mathrm{loc}}
        \right), \\
        \mathsf L_{\mu}^{\mathrm{loc}}(r)
        &\coloneqq
        \frac{\mu}{2}
        +
        \mu^2r
        \left(
            \frac{1+\gamma_E+|\log(\mu r)|}{\pi}
            +
            \mathsf C_{R'}^{\mathrm{loc}}
        \right).
    \end{aligned}
    \end{equation}
    Then
    \[
        |R_\mu(d)|
        \le
        \mathsf R_{\mu}^{\mathrm{loc}}(r)
        \quad (0\le d\le r),
        \qquad
        |R_\mu'(d)|
        \le
        \mathsf L_{\mu}^{\mathrm{loc}}(r)
        \quad (0<d\le r).
    \]
\end{corollary}

\begin{proof}
    On \(0<\tau\le1\), both functions
    \[
        \frac{\tau}{2}
        +
        \tau^2
        \left(
            \frac{|\gamma_E+\log\tau|}{2\pi}
            +
            \mathsf C_R^{\mathrm{loc}}
        \right)
        ,\quad
        \frac12
        +
        \tau
        \left(
            \frac{1+\gamma_E+|\log\tau|}{\pi}
            +
            \mathsf C_{R'}^{\mathrm{loc}}
        \right)
    \]
    are increasing. The conclusion follows from
    \eqref{eq:Rmu-E-bound} and
    \eqref{eq:Rmu-derivative-bound} with
    \(\tau=\mu d\le\mu r\). The estimate at \(d=0\) follows from the
    definition of \(R_\mu(0)\).
\end{proof}

\subsection{Off-diagonal kernel bounds for matrix-entry estimates}

\begin{lemma}[Holomorphic extension and a crude bound for \(G\)]
\label{lem:K-analytic-ellipse}
    Let \(\Log\) denote the principal branch on
    \(\bbC\setminus(-\infty,0]\), and define
    \[
        G(w)
        \coloneqq
        \frac{1}{\pi}
        \left[
            -\Ci(w)\cos w
            +
            \left(\frac{\pi}{2}-\Si(w)\right)\sin w
        \right]
    \]
    for \(\Re w>0\). Then \(G\) is holomorphic on this half-plane and,
    for every \(a_{\mathrm{off}}>0\),
    \begin{equation}\label{eq:K-crude-bound}
        \sup_{\Re w>a_{\mathrm{off}}}|G(w)|
        \le
        \frac{1}{\pi a_{\mathrm{off}}}.
    \end{equation}
\end{lemma}

\begin{proof}
    The defining expression is holomorphic for \(\Re w>0\). Moreover,
    \[
        G(w)
        =
        \frac{1}{\pi}
        \int_0^\infty
        e^{-\tau}
        \frac{\tau}{\tau^2+w^2}\,d\tau ,
        \qquad \Re w>0.
    \]
    Indeed, the integral is locally uniformly convergent and holomorphic,
    and the identity follows from the real-axis formula and the identity
    theorem. Since
    \[
        \frac{\tau}{\tau^2+w^2}
        =
        \frac12
        \left(
            \frac{1}{\tau-iw}
            +
            \frac{1}{\tau+iw}
        \right)
    \]
    and \(|\tau\pm iw|\ge\Re w\) for \(\tau\ge0\), we obtain
    \[
        |G(w)|
        \le
        \frac{1}{2\pi}
        \int_0^\infty e^{-\tau}
        \left(
            \frac{1}{|\tau-iw|}
            +
            \frac{1}{|\tau+iw|}
        \right)d\tau
        \le
        \frac{1}{\pi\Re w}.
    \]
    This proves \eqref{eq:K-crude-bound}.
\end{proof}

\subsection{Unit resolvent estimates for weighted \(Q\)-bounds}

\begin{lemma}[Elementary estimates for the unit resolvent kernel]
    Let \(G\) be the rescaled kernel defined in \eqref{eq:K-rescaled}.
    For \(|y|\ge1\),
    \begin{equation}\label{eq:G,G',G'' bouds for |y|geq 1}
        |G(y)|\le \frac{1}{\pi |y|^2},
        \qquad
        |G'(y)|\le \frac{2}{\pi |y|^3},
        \qquad
        |G''(y)|\le \frac{6}{\pi |y|^4}.
    \end{equation}
    Moreover, if \(y>0\) and \(0\le z\le y/2\), then
    \begin{equation}\label{eq:kernel estimate of G}
        \begin{aligned}
            &\frac12\left|G(y-z)+G(y+z)\right|
            \le
            \frac{1}{\pi y^2}
            +
            \frac{5z^2}{\pi y^4},\\
            &\frac12\left|G'(y-z)+G'(y+z)\right|
            \le
            \frac{2}{\pi y^3}
            +
            \frac{26z^2}{\pi y^5},\\
            &\frac12\left|G''(y-z)+G''(y+z)\right|
            \le
            \frac{6}{\pi y^4}
            +
            \frac{3840z^2}{\pi y^6}.
        \end{aligned}
    \end{equation}
    Also, we have
    \begin{equation}\label{eq:G-local-derivative-integrals}
        \int_{|y|\leq1}|yG'(y)|\,dy
        \leq
        1,
        \qquad
        \sup_{0<|y|\leq1}|y^2G''(y)|
        \leq
        \frac{2}{\pi}.
    \end{equation}
    and
    \begin{equation}\label{eq:G-far-derivative-integrals}
        \int_{|y|>1}|G'(y)|\,dy
        \leq
        \frac{2}{\pi},
        \qquad
        \int_{|y|>1}|G''(y)|\,dy
        \leq
        \frac{4}{\pi}.
    \end{equation}
\end{lemma}

\begin{proof}
    Differentiating \eqref{eq:G-poisson-real} with respect to \(y\), we obtain
    \begin{equation}\label{eq:G',G''-poisson formula}
        G'(y)
        =
        -\frac{2y}{\pi}
        \int_0^\infty
        e^{-\tau}
        \frac{\tau}{(\tau^2+y^2)^2}\,d\tau,
        \qquad
        G''(y)
        =
        \frac{2}{\pi}
        \int_0^\infty
        e^{-\tau}
        \frac{\tau(3y^2-\tau^2)}{(\tau^2+y^2)^3}\,d\tau .
    \end{equation}
    For \(|y|\geq1\),
    \begin{equation}\label{eq:far field bounds}
        \frac{\tau}{\tau^2+y^2}
        \leq
        \frac{\tau}{y^2},
        \qquad
        \frac{2|y|\tau}{(\tau^2+y^2)^2}
        \leq
        \frac{2\tau}{|y|^3},
        \qquad
        \frac{2\tau|3y^2-\tau^2|}{(\tau^2+y^2)^3}
        \leq
        \frac{6\tau}{|y|^4}.
    \end{equation}
    Integrating these inequalities in \(\tau\), with \(\int_0^\infty e^{-\tau}\tau\,d\tau=1\), gives \eqref{eq:G,G',G'' bouds for |y|geq 1}. We next prove \eqref{eq:kernel estimate of G}. Set
    \[
        p\coloneqq\tfrac{\tau}{y},
        \qquad
        q\coloneqq\tfrac{z}{y},
        \qquad
        \upsilon\coloneqq q^2.
    \]
    Then \(p\geq0\) and \(0\leq\upsilon\leq\frac14\). After factoring \(y^{-2}\) and \(y^{-3}\) from \eqref{eq:G-poisson-real} and \eqref{eq:G',G''-poisson formula}, the first two estimates in \eqref{eq:kernel estimate of G} reduce to
    \begin{equation}\label{eq:p,q rational ineq}
        \begin{aligned}
        \frac12
        \left(
            \frac{1}{p^2+(1-q)^2}
            +
            \frac{1}{p^2+(1+q)^2}
        \right)
        &\leq
        1+5q^2,\\
        \frac{1-q}{(p^2+(1-q)^2)^2}
        +
        \frac{1+q}{(p^2+(1+q)^2)^2}
        &\leq
        2+26q^2.
        \end{aligned}
    \end{equation}
    The left-hand sides of \eqref{eq:p,q rational ineq} decrease with \(p^2\). Their values at \(p=0\) are \(\frac{1+\upsilon}{(1-\upsilon)^2}\) and \(\frac{2(1+3\upsilon)}{(1-\upsilon)^3}\) respectively. Moreover,
    \begin{equation}\label{eq:upsilon ineq}
        \begin{aligned}
        (1+5\upsilon)
        -
        \frac{1+\upsilon}{(1-\upsilon)^2}
        &=
        \frac{
            \upsilon(2-9\upsilon+5\upsilon^2)
        }{
            (1-\upsilon)^2
        },
        \\
        (2+26\upsilon)
        -
        \frac{2(1+3\upsilon)}{(1-\upsilon)^3}
        &=
        \frac{
            2\upsilon
            (7-36\upsilon+38\upsilon^2-13\upsilon^3)
        }{
            (1-\upsilon)^3
        }.
        \end{aligned}
    \end{equation}
    Both polynomials in parentheses in \eqref{eq:upsilon ineq} decrease on \([0,\frac14]\), and their values at \(\upsilon=\frac14\) are \(\frac1{16}\) and \(\frac{11}{64}\)
    respectively. This proves the first two rational inequalities in \eqref{eq:kernel estimate of G}.

    For the last estimate in \eqref{eq:kernel estimate of G}, define
    \[
        F(p,\zeta)
        \coloneqq
        \frac{3\zeta^2-p^2}{(p^2+\zeta^2)^3}.
    \]
    Since \(|3-p^2|\leq3(1+p^2)\), we have \(2|F(p,1)|\leq6\). Also,
    \[
        \partial_\zeta^2 F(p,\zeta)
        =
        \frac{
            12(p^4-10p^2\zeta^2+5\zeta^4)
        }{
            (p^2+\zeta^2)^5
        }.
    \]
    Using
    \(
        |p^4-10p^2\zeta^2+5\zeta^4|
        \leq
        5(p^2+\zeta^2)^2,
    \)
    we obtain
    \[
        |\partial_\zeta^2F(p,\zeta)|
        \leq
        \frac{60}{(p^2+\zeta^2)^3}
        \leq
        3840,
        \qquad
        \frac12\leq\zeta\leq\frac32.
    \]
    Taylor's theorem therefore gives
    \[
        \begin{aligned}
            |F(p,1-q)+F(p,1+q)|
            &\leq
            2|F(p,1)|
            +
            |F(p,1-q)+F(p,1+q)-2F(p,1)|\\
            &\leq
            6+3840q^2.
        \end{aligned}
    \]
    Substitution into the integral formulas for \(G''\) in \eqref{eq:G',G''-poisson formula}, followed by \(\int_0^\infty e^{-\tau}\tau\,d\tau=1\), proves the last estimate in
    \eqref{eq:kernel estimate of G}.

    It remains to prove \eqref{eq:G-local-derivative-integrals} and \eqref{eq:G-far-derivative-integrals}. The kernel \(G\) is positive, even, and decreasing on \((0,\infty)\), with
    \(\int_{\bbR}G=1\) and \(yG(y)\to0\) as \(y\to0\). Integration by parts gives
    \[
        \int_{|y|\leq1}|yG'(y)|\,dy
        =
        2\left(
            \int_0^1G(y)\,dy-G(1)
        \right)
        \leq1.
    \]
    Also, the formula for \(G''\) in \eqref{eq:G',G''-poisson formula} and the change of variables \(\tau=|y|\upsilon\) give
    \[
        |y^2G''(y)|
        \leq
        \frac2\pi
        \int_0^\infty
        \frac{\upsilon(3+\upsilon^2)}{(1+\upsilon^2)^3}\,d\upsilon
        =
        \frac2\pi,
        \qquad
        0<|y|\leq1,
    \]
    which gives \eqref{eq:G-local-derivative-integrals}. Finally, integrating \eqref{eq:far field bounds} for \(|y|>1\) yields
    \[
        \int_{|y|>1}|G'(y)|\,dy
        \leq
        2\int_1^\infty\frac{2}{\pi y^3}\,dy
        =
        \frac2\pi,
        \qquad
        \int_{|y|>1}|G''(y)|\,dy
        \leq
        2\int_1^\infty\frac{6}{\pi y^4}\,dy
        =
        \frac4\pi.
    \]
    This proves \eqref{eq:G-far-derivative-integrals}.
\end{proof}

\section{Logarithmic potential estimates and identities}

Recall \(\omega(y) =  \log(1+|y|)\), and set
\begin{equation*}
    \mathcal P_{\log}f
    \coloneqq
    -\frac1\pi
    \log|\cdot|*f,
    \qquad
    \mathcal B_{\log}(f_1,f_2)
    \coloneqq
    \bigl(
        \mathcal P_{\log}f_1,
        f_2
    \bigr)_r.
\end{equation*}

\begin{lemma}[Bounds for the logarithmic potential]
    Let
    \(f_1,f_2\in L^1(\bbR)\cap L^2(\bbR)\), and assume that
    \(\omega f_1,\omega f_2\in L^1(\bbR)\).
    Then
    \begin{equation}\label{eq:part3-basic-log-pairing-bound}
        |\mathcal B_{\log}(f_1,f_2)|
        \leq
        \frac1\pi
        \left(
            2
            \norm{f_1}_{L^2}
            \norm{f_2}_{L^2}
            +
            \norm{f_1}_{L^1}
            \norm{\omega f_2}_{L^1}
            +
            \norm{\omega f_1}_{L^1}
            \norm{f_2}_{L^1}
        \right).
    \end{equation}
    If, in addition,
    \(f_3,\omega f_3\in L^2(\bbR)\), then
    \begin{equation}\label{eq:part3-basic-log-potential-bound}
        \norm{
            f_3\,
            \mathcal P_{\log}f_1
        }_{L^2}
        \leq
        \frac1\pi
        \left(
            2
            \norm{f_3}_{L^2}
            \norm{f_1}_{L^2}
            +
            \norm{\omega f_3}_{L^2}
            \norm{f_1}_{L^1}
            +
            \norm{f_3}_{L^2}
            \norm{\omega f_1}_{L^1}
        \right).
    \end{equation}
\end{lemma}

\begin{proof}
    We split \(\bbR\) into the regions \(|x-y|\leq1\) and \(|x-y|>1\), and estimate each separately. 
    
    On the region \(|x-y|\leq1\), since \(\norm{\mathbf 1_{\{|\cdot|\leq1\}}|\log|\cdot||}_{L^1} = \norm{\mathbf 1_{\{|\cdot|\leq1\}}|\log|\cdot||}_{L^2} = 2\), Young's inequality gives
    \begin{equation}\label{eq:B_log near}
        \left|
            \int_{\bbR}
                f_2(y)
                \int_{|y-z|\leq1}
                    \log|y-z|f_1(z)
                \,dz
            \,dy
        \right|
        \leq
        2
        \norm{f_1}_{L^2}
        \norm{f_2}_{L^2}.
    \end{equation}
    Similarly, Young's inequality gives \(\norm{\mathbf 1_{\{|\cdot|\leq1\}}|\log|\cdot||*f_1}_{L^\infty}\leq 2\norm{f_1}_{L^2}\), and hence
    \begin{equation}\label{eq:P_log near}
        \norm{
            f_3
            \left(
                \mathbf 1_{\{|\cdot|\leq1\}}|\log|\cdot||*f_1
            \right)
        }_{L^2}
        \leq
        2
        \norm{f_3}_{L^2}
        \norm{f_1}_{L^2}.
    \end{equation}

    On the region \(|y-z|>1\), we have \(|\log|y-z|| \leq \omega(y)+\omega(z)\).
    Therefore
    \begin{equation}\label{eq:B_log far}
            \left|
                \int_{\bbR}
                    f_2(y)
                    \int_{|y-z|>1}
                        \log|y-z|f_1(z)
                    \,dz
                \,dy
            \right|
            \leq
            \norm{f_1}_{L^1}
            \norm{\omega f_2}_{L^1}
            +
            \norm{\omega f_1}_{L^1}
            \norm{f_2}_{L^1},
    \end{equation}
    and
    \begin{equation}\label{eq:P_log estimate1}
        \left|
            \int_{|y-z|>1}
                \log|y-z|f_1(z)
            \,dz
        \right|
        \leq
        \omega(y)\norm{f_1}_{L^1}
        +
        \norm{\omega f_1}_{L^1}.
    \end{equation}
    Multiplying \eqref{eq:P_log estimate1} by \(|f_3(y)|\) and taking the \(L^2\)-norm gives
    \begin{equation}\label{eq:P_log far}
        \norm{
            f_3
            \int_{|\,\cdot-z|>1}
                \log|\cdot-z|f_1(z)
            \,dz
        }_{L^2}
        \leq
        \norm{\omega f_3}_{L^2}
        \norm{f_1}_{L^1}
        +
        \norm{f_3}_{L^2}
        \norm{\omega f_1}_{L^1}.
    \end{equation}
    The estimates \eqref{eq:B_log near} and \eqref{eq:B_log far} give \eqref{eq:part3-basic-log-pairing-bound}, while \eqref{eq:P_log near} and \eqref{eq:P_log far} give \eqref{eq:part3-basic-log-potential-bound}.
\end{proof}

\begin{corollary}[A logarithmic-potential bound for \(e_Q\)]
    For \(\mathsf V^{L^1}\) in \eqref{eq:def-CVQ-and-VQ-bound}, \(\mathsf V_{\omega^2}^{L^1}\) in \eqref{eq:def-VQ-omega-bound}, \(\Delta_Q^{L^2}\) in \eqref{eq:def-Delta-Q-L2}, \(\Delta_Q^{L^1}\) in \eqref{eq:def-Delta-Lambda-a-Q-L1}, and \(\Delta_{\omega Q}^{L^1}\) in \eqref{eq:def-Delta-omega-Q-L1} define
    \begin{equation}\label{eq:part3-Eelog-def}
        \mathsf E_{e,\log}
        \coloneqq
        \frac{2}{\pi}
        \left(
            \mathsf V^{L^1}
        \right)^{1/2}
        \Delta_Q^{L^2}
        +
        \frac1\pi
        \left(
            \mathsf V_{\omega^2}^{L^1}
        \right)^{1/2}
        \Delta_Q^{L^1}
        +
        \frac1\pi
        \left(
            \mathsf V^{L^1}
        \right)^{1/2}
        \Delta_{\omega Q}^{L^1}.
    \end{equation}
    Then, for \(e_Q\) in \eqref{eq:def-e_Q} we have
    \begin{equation}\label{eq:part3-Eelog-bound}
        \norm{
            |V_Q|^{1/2}
            \mathcal P_{\log}e_Q
        }_{L^2}
        \leq
        \mathsf E_{e,\log}.
    \end{equation}
\end{corollary}

\begin{proof}
    Apply \eqref{eq:part3-basic-log-potential-bound} with \(f_3=|V_Q|^{1/2}\) and \(f_1=e_Q\).
    Lemma~\ref{lem:weighted-Q-V_Q-bounds} gives
    \[
        \norm{
            |V_Q|^{1/2}
        }_{L^2}
        \leq
        \left(
            \mathsf V^{L^1}
        \right)^{1/2},
        \qquad
        \norm{
            \omega|V_Q|^{1/2}
        }_{L^2}
        \leq
        \left(
            \mathsf V_{\omega^2}^{L^1}
        \right)^{1/2}.
    \]
    Lemma~\ref{lem:scaling-and-L1-transfer} gives
    \[
        \norm{e_Q}_{L^2}
        \leq
        \Delta_Q^{L^2},
        \qquad
        \norm{e_Q}_{L^1}
        \leq
        \Delta_Q^{L^1},
        \qquad
        \norm{\omega e_Q}_{L^1}
        \leq
        \Delta_{\omega Q}^{L^1}.
    \]
    Substitution in
    \eqref{eq:part3-basic-log-potential-bound}
    proves \eqref{eq:part3-Eelog-bound}.
\end{proof}

\begin{lemma}[Differentiation and scaling identities for\(\mathcal P_{\log}\)]
    If \(f\in H^1(\bbR)\), then
    \begin{equation}\label{eq:part3-Plog-Hilbert-bounds}
        \norm{
            \partial_y
            \mathcal P_{\log}f
        }_{L^2}
        \leq
        \norm{f}_{L^2},
        \qquad
        \norm{
            \partial_y^2
            \mathcal P_{\log}f
        }_{L^2}
        \leq
        \norm{\partial_yf}_{L^2}.
    \end{equation}

    Assume in addition that \(f,f_1,f_2\in C^2(\bbR)\), and that their
    scaling derivatives up to order two are
    \(O(\langle y\rangle^{-2})\).
    Then
    \begin{equation}\label{eq:part3-scaling-mass-identities}
        \int_{\bbR}\Lambda f
        =
        -\frac12
        \int_{\bbR}f,
        \qquad
        \int_{\bbR}\Lambda^2f
        =
        \frac14
        \int_{\bbR}f,
    \end{equation}
    and
    \begin{equation}\label{eq:part3-Plog-Lambda-identity}
        \mathcal P_{\log}(\Lambda f)
        =
        y\partial_y
        \mathcal P_{\log}f
        -
        \frac12
        \mathcal P_{\log}f
        +
        \frac1\pi
        \int_{\bbR}f,
    \end{equation}
    \begin{equation}\label{eq:part3-Plog-Lambda2-identity}
        \mathcal P_{\log}(\Lambda^2f)
        =
        y^2\partial_y^2
        \mathcal P_{\log}f
        +
        \frac14
        \mathcal P_{\log}f
        -
        \frac1\pi
        \int_{\bbR}f.
    \end{equation}
    Moreover, 
    \begin{equation}\label{eq:part3-Blog-Lambda-identity}
        \mathcal B_{\log}(\Lambda f_1,f_2)
        =
        -
        \mathcal B_{\log}
        \bigl(
            f_1,
            (\Lambda+1)f_2
        \bigr)
        +
        \frac1\pi
        \left(
            \int_{\bbR}f_1
        \right)
        \left(
            \int_{\bbR}f_2
        \right),
    \end{equation}
    \begin{equation}\label{eq:part3-Blog-Lambda2-identity}
        \mathcal B_{\log}(\Lambda^2f_1,f_2)
        =
        \mathcal B_{\log}
        \bigl(
            f_1,
            (\Lambda+1)^2f_2
        \bigr)
        -
        \frac1\pi
        \left(
            \int_{\bbR}f_1
        \right)
        \left(
            \int_{\bbR}f_2
        \right).
    \end{equation}
\end{lemma}

\begin{proof}
    Define the Hilbert transform by
    \[
        \operatorname{Hilb} f(y)
        \coloneqq
        \frac1\pi
        \operatorname{p.v.}
        \int_{\bbR}
            \frac{f(z)}{y-z}
        \,dz.
    \]
    on \(\bbR\). Differentiating the logarithmic kernel in the sense of distributions gives
    \begin{equation*}
        \partial_y
        \mathcal P_{\log}f
        =
        -\operatorname{Hilb}f,
        \qquad
        \partial_y^2
        \mathcal P_{\log}f
        =
        -\operatorname{Hilb}(\partial_yf).
    \end{equation*}
    Since \(\operatorname{Hilb}\) is an isometry on \(L^2(\bbR)\), \eqref{eq:part3-Plog-Hilbert-bounds} follows. Integration by parts gives \eqref{eq:part3-scaling-mass-identities}.

    Since \(y\partial_y\log|y-z|+z\partial_z\log|y-z|=1\), integration by parts in \(z\) gives
    \[
        \begin{aligned}
        \mathcal P_{\log}(\Lambda f)(y)
        =
        -\frac1\pi
        \int_{\bbR}
            \log|y-z|
            \left(
                \frac12 f(z)
                +
                zf'(z)
            \right)
        \,dz
        =
        y\partial_y
        \mathcal P_{\log}f
        -
        \frac12
        \mathcal P_{\log}f
        +
        \frac1\pi
        \int_{\bbR}f.
        \end{aligned}
    \]
    This proves \eqref{eq:part3-Plog-Lambda-identity}. Applying the same identity to \(\Lambda f\), and using \eqref{eq:part3-scaling-mass-identities} gives \eqref{eq:part3-Plog-Lambda2-identity}.

    Pairing \eqref{eq:part3-Plog-Lambda-identity} with \(f_2\), and using
    \[
        \int_{\bbR}
            y
            \partial_y
            \mathcal P_{\log}f_1
            \,f_2
        =
        -
        \int_{\bbR}
            \mathcal P_{\log}f_1
            \left(
                f_2+yf_2'
            \right),
    \]
    gives \eqref{eq:part3-Blog-Lambda-identity}.
    Similarly,
    \[
        \int_{\bbR}
            y^2
            \partial_y^2
            \mathcal P_{\log}f_1
            \,f_2
        =
        \int_{\bbR}
            \mathcal P_{\log}f_1
            \left(
                y^2f_2''+4yf_2'+2f_2
            \right),
    \]
    and hence
    \[
        y^2f_2''+4yf_2'+\frac94f_2
        =
        (\Lambda+1)^2f_2.
    \]
    Combining this identity with \eqref{eq:part3-Plog-Lambda2-identity}, we obtain \eqref{eq:part3-Blog-Lambda2-identity}.
\end{proof}

\begin{lemma}[Logarithmic tail integrals]
\label{lem:app-elementary-log-tail-bounds}
    Let \(R\geq 1\). Then
    \begin{equation}\label{eq:app-tail-L1-bound}
        \int_{|y|>R}\langle y\rangle^{-2}\,dy
        \leq
        \frac{2}{R},
    \end{equation}
    \begin{equation}\label{eq:app-tail-omega-L1-bound}
        \int_{|y|>R}
            \omega(y)\langle y\rangle^{-2}
        \,dy
        \leq
        2
        \left(
            \frac{\log(1+R)}{R}
            +
            \log\left(1+\frac1R\right)
        \right),
    \end{equation}
    and
    \begin{equation}\label{eq:app-tail-omega2-langle4}
        \int_{|y|>R}
            \omega(y)^2\langle y\rangle^{-4}
        \,dy
        \leq
        \frac{2}{R^3}
        \left(
            \frac{\log^2(2R)}{3}
            +
            \frac{2\log(2R)}{9}
            +
            \frac{2}{27}
        \right).
    \end{equation}
\end{lemma}

\begin{proof}
    Since \(R\geq1\), for \(|y|>R\) we have \(\omega(y)\leq\log(2|y|)\). Hence
    \[
        \int_{|y|>R}\langle y\rangle^{-2}\,dy
        \leq
        2\int_R^\infty y^{-2}\,dy
        =
        \frac2R.
    \]
    Integration by parts gives
    \[
        \int_{|y|>R}
            \omega(y)\langle y\rangle^{-2}
        \,dy
        \leq 
        2\int_R^\infty
            \frac{\log(1+y)}{y^2}
        \,dy
        =
        2\left(\frac{\log(1+R)}{R}
        +
        \log\left(1+\frac1R\right)\right),
    \]
    which proves \eqref{eq:app-tail-omega-L1-bound}. Finally,
    \[
        \int_{|y|>R}
            \omega(y)^2\langle y\rangle^{-4}
        \,dy
        \leq
        2\int_R^\infty
            \frac{\log^2(2y)}{y^4}
        \,dy
        =
        \frac{2}{R^3}
        \left(
            \frac{\log^2(2R)}{3}
            +
            \frac{2\log(2R)}{9}
            +
            \frac{2}{27}
        \right),
    \]
    and \eqref{eq:app-tail-omega2-langle4} follows.
\end{proof}

\begin{lemma}[Exact logarithmic integrals]
\label{lem:exact-log-cell-integrals}
    Let \(h>0\), and set
    \[
        J_j
        \coloneqq
        [(j-1)h,jh],
        \qquad
        j\geq1.
    \]
    Define
    \[
        \Theta_{\log}(0)
        \coloneqq
        0,
        \qquad
        \Theta_{\log}(\rho)
        \coloneqq
        \frac12\rho^2\log\rho-\frac34\rho^2,
        \quad
        \rho>0.
    \]
    For \(d,\ell\in\bbN_0\), define
    \begin{equation}\label{eq:log-cell-diff-def}
        L_d^{-}(h)
        \coloneqq
        \begin{cases}
            \displaystyle
            h^2
            \left(
                \log h-\frac32
            \right),
            & d=0,
            \\[1mm]
            \displaystyle
            \Theta_{\log}((d+1)h)
            -
            2\Theta_{\log}(dh)
            +
            \Theta_{\log}((d-1)h),
            & d\geq1,
        \end{cases}
    \end{equation}
    and
    \begin{equation}\label{eq:log-cell-sum-def}
        L_\ell^{+}(h)
        \coloneqq
        \Theta_{\log}((\ell+2)h)
        -
        2\Theta_{\log}((\ell+1)h)
        +
        \Theta_{\log}(\ell h).
    \end{equation}
    Then, for \(i,j\geq1\),
    \begin{equation}\label{eq:log-cell-diff-id}
        \int_{J_i}\int_{J_j}
            \log|y-z|
        \,dz\,dy
        =
        L_{|i-j|}^{-}(h),
    \end{equation}
    and
    \begin{equation}\label{eq:log-cell-sum-id}
        \int_{J_i}\int_{J_j}
            \log(y+z)
        \,dz\,dy
        =
        L_{i+j-2}^{+}(h).
    \end{equation}
\end{lemma}

\begin{proof}
    Since
    \(
        \Theta_{\log}''(\rho)=\log\rho
    \)
    for \(\rho>0\), two integrations give, for
    \(d=|i-j|\geq1\),
    \[
        \int_{J_i}\int_{J_j}
            \log|y-z|
        \,dz\,dy
        =
        \Theta_{\log}((d+1)h)
        -
        2\Theta_{\log}(dh)
        +
        \Theta_{\log}((d-1)h).
    \]
    For \(d=0\), scaling gives
    \[
        \begin{aligned}
            \int_{J_i}\int_{J_i}
                \log|y-z|
            \,dz\,dy
            &=
            h^2\log h
            +
            h^2
            \int_0^1\int_0^1
                \log|\eta-\xi|
            \,d\xi\,d\eta
            =
            h^2
            \left(
                \log h-\frac32
            \right).
        \end{aligned}
    \]
    This proves \eqref{eq:log-cell-diff-id}. Similarly, with
    \(\ell=i+j-2\),
    \[
        \begin{aligned}
            \int_{J_i}\int_{J_j}
                \log(y+z)
            \,dz\,dy
            &=
            \int_0^h\int_0^h
                \log(\ell h+\eta+\xi)
            \,d\xi\,d\eta
            =
            L_\ell^{+}(h),
        \end{aligned}
    \]
    where \(\Theta_{\log}(0)=0\) covers the case \(\ell=0\).
\end{proof}

\section{Constants in the higher-order and weighted profile estimates}

This section records the explicit constants used in Lemmas~\ref{lem:higher-order-transfer-from-H1}, \ref{lem:weighted-Q-V_Q-bounds}, \ref{lem:scaling-and-L1-transfer}, and
\ref{lem:square-root-potential-transfer}. The formulas are kept in an explicitly computable form for the interval evaluations in Part~\ref{part:interval-verification}. Throughout this section, \(\Delta_Q^{W^{k,2}}\) and \(\Delta_{\Lambda^aQ}^{L^p}\) denote bounds for \(e_Q=Q-\mathfrak g\) and \(\Lambda^ae_Q\), respectively.

\subsection{Unweighted higher-order estimates}
We define the constants used in Lemma~\ref{lem:higher-order-transfer-from-H1}.

\begin{definition}[Constants for the unweighted estimates]\label{def:unweighted-transfer-constants}
    Let
    \(
        r_{\mathfrak g}\coloneqq\mathcal F(\mathfrak g)
    \).
    \begin{enumerate}[label=\textup{(\roman*)}]
        \item The following constants bound \(Q^3-\mathfrak g^3\) in \(W^{1,2}\), then \(|D|e_Q\) in \(W^{1,2}\), and finally \(e_Q\) in \(W^{2,2}\). Set
        \begin{equation}\label{eq:cubic-w12-error-definition}
            \Delta_{\mathrm{cub}}^{W^{1,2}}
            \coloneqq
            \frac52 r_{\mathrm{NK}}
            \left[
                \left(\norm{\mathfrak g}_{H^1}+r_{\mathrm{NK}}\right)^2
                +
                \left(\norm{\mathfrak g}_{H^1}+r_{\mathrm{NK}}\right)
                \norm{\mathfrak g}_{H^1}
                +
                \norm{\mathfrak g}_{H^1}^2
            \right].
        \end{equation}
        Next, set
        \begin{equation}\label{eq:def-B2-unweighted-transfer and R_H^2}
            \begin{aligned}
                \Delta_{|D|Q}^{W^{1,2}}
                &\coloneqq
                r_{\mathrm{NK}}
                +
                \Delta_{\mathrm{cub}}^{W^{1,2}}
                +
                \norm{r_{\mathfrak g}}_{W^{1,2}},
                \\
                \Delta_Q^{W^{2,2}}
                &\coloneqq
                \left(r_{\mathrm{NK}}^2+(\Delta_{|D|Q}^{W^{1,2}})^2\right)^{1/2},
            \end{aligned}
        \end{equation}
        and
        \begin{equation}\label{eq:E0-E1-unweighted-transfer}
            \Delta_{Q,0}^{L^\infty}
            \coloneqq
            2^{-1/2}r_{\mathrm{NK}},
            \qquad
            \Delta_{Q,1}^{L^\infty}
            \coloneqq
            2^{-1/2}\Delta_Q^{W^{2,2}}.
        \end{equation}

        \item
        For \(a=0,1,2\), the constants below bound \(\partial_y^a(Q^3-\mathfrak g^3)\) in \(L^2\). Set
        \begin{equation}\label{eq:def-F-cub-1to3}
            \begin{aligned}
                C_{\mathrm{cub},0}^{L^2}
                &\coloneqq
                3\norm{\mathfrak g}_{L^\infty}^2r_{\mathrm{NK}}
                +
                3\norm{\mathfrak g}_{L^\infty}\Delta_{Q,0}^{L^\infty}r_{\mathrm{NK}}
                +
                (\Delta_{Q,0}^{L^\infty})^2r_{\mathrm{NK}},
                \\
                C_{\mathrm{cub},1}^{L^2}
                &\coloneqq
                6\norm{\mathfrak g}_{L^\infty}
                \norm{\mathfrak g'}_{L^\infty}r_{\mathrm{NK}}
                +
                3\norm{\mathfrak g}_{L^\infty}^2r_{\mathrm{NK}}  \\
                &\quad+
                3\norm{\mathfrak g'}_{L^\infty}\Delta_{Q,0}^{L^\infty}r_{\mathrm{NK}}
                +
                6\norm{\mathfrak g}_{L^\infty}\Delta_{Q,0}^{L^\infty}r_{\mathrm{NK}}
                +
                3(\Delta_{Q,0}^{L^\infty})^2r_{\mathrm{NK}},
                \\
                C_{\mathrm{cub},2}^{L^2}
                &\coloneqq
                6\norm{\mathfrak g'}_{L^\infty}^2r_{\mathrm{NK}}
                +
                6\norm{\mathfrak g}_{L^\infty}
                \norm{\mathfrak g''}_{L^2}\Delta_{Q,0}^{L^\infty}
                +
                12\norm{\mathfrak g}_{L^\infty}
                \norm{\mathfrak g'}_{L^\infty}r_{\mathrm{NK}}
                \\
                &\quad+
                3\norm{\mathfrak g}_{L^\infty}^2
                \Delta_{|D|Q}^{W^{1,2}}
                +
                3\norm{\mathfrak g''}_{L^2}
                (\Delta_{Q,0}^{L^\infty})^2
                +
                12\norm{\mathfrak g'}_{L^\infty}
                \Delta_{Q,0}^{L^\infty}r_{\mathrm{NK}}
                \\
                &\quad+
                6\norm{\mathfrak g}_{L^\infty}
                \Delta_{Q,1}^{L^\infty}r_{\mathrm{NK}}
                +
                6\norm{\mathfrak g}_{L^\infty}
                \Delta_{Q,0}^{L^\infty}
                \Delta_{|D|Q}^{W^{1,2}}
                \\
                &\quad+
                6\Delta_{Q,0}^{L^\infty}
                \Delta_{Q,1}^{L^\infty}r_{\mathrm{NK}}
                +
                3(\Delta_{Q,0}^{L^\infty})^2
                \Delta_{|D|Q}^{W^{1,2}},
            \end{aligned}
        \end{equation}
        Their combined \(W^{2,2}\)-bound is
        \begin{equation}\label{eq:cubic-w22-error-definition}
            \Delta_{\mathrm{cub}}^{W^{2,2}}
            \coloneqq
            \left[
                \left(C_{\mathrm{cub},0}^{L^2}\right)^2
                +
                \left(C_{\mathrm{cub},1}^{L^2}\right)^2
                +
                \left(C_{\mathrm{cub},2}^{L^2}\right)^2
            \right]^{1/2}.
        \end{equation}

        \item The final constants bound \(|D|e_Q\) in \(W^{2,2}\), \(e_Q\) in \(W^{3,2}\), and \(\partial_y^2e_Q\) in \(L^\infty\). Set
        \begin{equation}\label{eq:higher-sobolev-error-definitions}
            \begin{aligned}
                \Delta_{|D|Q}^{W^{2,2}}
                &\coloneqq
                \Delta_Q^{W^{2,2}}
                +
                \Delta_{\mathrm{cub}}^{W^{2,2}}
                +
                \norm{r_{\mathfrak g}}_{W^{2,2}},\\
                \Delta_Q^{W^{3,2}}
                &\coloneqq
                \left[
                    (\Delta_Q^{W^{2,2}})^2
                    +
                    (\Delta_{|D|Q}^{W^{2,2}})^2
                \right]^{1/2},\\
                \Delta_{Q,2}^{L^\infty}
                &\coloneqq
                2^{-1/2}
                \left[
                    (\Delta_{|D|Q}^{W^{1,2}})^2
                    +
                    (\Delta_{|D|Q}^{W^{2,2}})^2
                \right]^{1/2}.
            \end{aligned}
        \end{equation}
    \end{enumerate}
\end{definition}

\subsection{Weighted bounds for \(Q\) and \(V_Q\)}

The constants below are used in Lemma~\ref{lem:weighted-Q-V_Q-bounds}.
We first define the interior and exterior bounds; the resulting global
bounds are recorded after Lemma~\ref{lem:app-exterior-weighted-Q-bounds}.

\begin{definition}[Constants for the weighted bounds for $Q$]\label{def:weighted-q-constants}
    Fix
    \[
        R_M\geq2,
        \qquad
        R_Q\geq2R_M.
    \]

    \begin{enumerate}[label=\textup{(\roman*)}]
        \item
        For the interior region \(\{|y|\leq R_Q\}\), set
        \begin{equation}\label{eq:def-CQ-core}
            \begin{aligned}
                C_{Q,0}^{\mathrm{int}}
                &\coloneqq
                \norm{\langle y\rangle^2\mathfrak g}_{L^\infty}
                +
                (1+R_Q^2)\Delta_{Q,0}^{L^\infty},\\
                C_{Q,1}^{\mathrm{int}}
                &\coloneqq
                \norm{\langle y\rangle^2y\mathfrak g'}_{L^\infty}
                +
                R_Q(1+R_Q^2)\Delta_{Q,1}^{L^\infty},\\
                C_{Q,2}^{\mathrm{int}}
                &\coloneqq
                \norm{\langle y\rangle^2y^2\mathfrak g''}_{L^\infty}
                +
                R_Q^2(1+R_Q^2)\Delta_{Q,2}^{L^\infty}.
            \end{aligned}
        \end{equation}

        \item
        The quantity \(\Delta_{\mathrm{cub}}^{L^1}\) bounds \(\|Q^3-\mathfrak g^3\|_{L^1}\). Set
        \begin{equation}\label{eq:cubic-l1-error-definition}
            \begin{aligned}
                C_Q^{L^\infty}
                &\coloneqq
                \norm{\mathfrak g}_{L^\infty}+\Delta_{Q,0}^{L^\infty},\\
                \Delta_{\mathrm{cub}}^{L^1}
                &\coloneqq
                r_{\mathrm{NK}}
                \left[
                    C_Q^{L^\infty}
                    \left(\norm{\mathfrak g}_{L^2}+r_{\mathrm{NK}}\right)
                    +
                    C_Q^{L^\infty}\norm{\mathfrak g}_{L^2}
                    +
                    \norm{\mathfrak g}_{L^\infty}\norm{\mathfrak g}_{L^2}
                \right].
            \end{aligned}
        \end{equation}
        Next define bounds for $\int Q^3$ and $\int y^2Q^3$.
        \begin{equation}\label{eq:def-M0,M_1,-A0,B0-weighted-envelope}
            \begin{aligned}
                C_{Q^3,0}
                &\coloneqq
                \int_{\bbR}\mathfrak g^3\,dy
                +
                \Delta_{\mathrm{cub}}^{L^1},\\
                C_{Q^3,2}
                &\coloneqq
                \int_{\bbR}y^2\mathfrak g^3\,dy
                +
                R_M^2\Delta_{\mathrm{cub}}^{L^1}\\
                &\quad+
                2\left(
                    \frac{(C_{Q,\mathrm{lead}}^{\mathrm{tail}})^3}{3R_M^3}
                    +
                    \frac{3(C_{Q,\mathrm{lead}}^{\mathrm{tail}})^2C_{Q,\mathrm{rem}}^{\mathrm{tail}}}{4R_M^4}
                    +
                    \frac{3C_{Q,\mathrm{lead}}^{\mathrm{tail}}(C_{Q,\mathrm{rem}}^{\mathrm{tail}})^2}{5R_M^5}
                    +
                    \frac{(C_{Q,\mathrm{rem}}^{\mathrm{tail}})^3}{6R_M^6}
                \right),
            \end{aligned}
        \end{equation}
        where $C_{Q,\mathrm{lead}}^{\mathrm{tail}}\coloneqq \tfrac{4C_{Q^3,0}}{\pi}$ and $C_{Q,\mathrm{rem}}^{\mathrm{tail}} \coloneqq (C_{Q^3,0})^3$. Finally, set
        \begin{equation}\label{eq:def-Qtail-weighted-envelope}
            C_Q^{\mathrm{tail}}
            \coloneqq
            \frac{4C_{Q,\mathrm{lead}}^{\mathrm{tail}}}{R_Q^2}
            +
            \frac{8C_{Q,\mathrm{rem}}^{\mathrm{tail}}}{R_Q^3}.
        \end{equation}

        \item
        For the exterior bounds of \(Q\) and \(yQ'\), set
        \begin{equation}\label{eq:def-CQ0-far}
            C_{Q,0}^{\mathrm{ext}}
            \coloneqq
            (1+R_Q^2)
            \left(
                \frac{C_{Q^3,0}}{\pi R_Q^2}
                +
                \frac{5C_{Q^3,2}}{\pi R_Q^4}
                +
                (C_Q^{\mathrm{tail}})^3
            \right).
        \end{equation}
        Since \(\partial_y(Q^3)=3Q^2Q'\), define
        \begin{equation}\label{eq:def-Floc-prime}
            F_{\mathrm{loc}}'
            \coloneqq
            3(C_Q^{\mathrm{tail}})^2
            \left(
                \frac{
                    8\norm{\langle y\rangle^2y\mathfrak g'}_{L^\infty}
                }{R_Q^3}
                +
                \Delta_{Q,1}^{L^\infty}
            \right).
        \end{equation}
        and set
        \begin{equation}\label{eq:def-CQ1-far}
            \begin{aligned}
                C_{Q,1}^{\mathrm{ext}}
                \coloneqq{}&
                (1+R_Q^2)R_Q
                \left(
                    \frac{2C_{Q^3,0}}{\pi R_Q^3}
                    +
                    \frac{26C_{Q^3,2}}{\pi R_Q^5}
                    +
                    4F_{\mathrm{loc}}'
                    +
                    \frac4\pi (C_Q^{\mathrm{tail}})^3
                \right).
            \end{aligned}
        \end{equation}

        \item For the exterior bound of \(y^2Q''\), let \(\tilde R_Q\coloneqq R_Q-1\) and set
        \begin{equation}\label{eq:def-CQ-tilde}
            \widetilde C_{Q,0}
            \coloneqq
            \max\{C_{Q,0}^{\mathrm{int}},C_{Q,0}^{\mathrm{ext}}\},
            \qquad
            \widetilde C_{Q,1}
            \coloneqq
            \max\{C_{Q,1}^{\mathrm{int}},C_{Q,1}^{\mathrm{ext}}\}.
        \end{equation}
        Define
        \begin{equation}\label{eq:def-QX-weighted-envelope}
            Q_X
            \coloneqq
            \frac{\widetilde C_{Q,0}}{1+\tilde R_Q^2},
            \quad
            Q_X'
            \coloneqq
            \frac{\widetilde C_{Q,1}}{(1+\tilde R_Q^2)\tilde R_Q}, 
            \quad
            Q_X''
            \coloneqq
            \frac{
                \norm{\langle y\rangle^2y^2\mathfrak g''}_{L^\infty}
            }{
                (1+\tilde R_Q^2)\tilde R_Q^2
            }
            +
            \Delta_{Q,2}^{L^\infty}.
        \end{equation}
        Since $\partial_y^2(Q^3)=6Q(Q')^2+3Q^2Q''$, set
        \begin{equation}\label{eq:def-FX-weighted-envelope}
            F_X
            \coloneqq
            Q_X^3,
            \qquad
            F_X''
            \coloneqq
            6Q_X(Q_X')^2+3Q_X^2Q_X''.
        \end{equation}
        Finally, define
        \begin{equation}\label{eq:def-CQ2-far}
            \begin{aligned}
                C_{Q,2}^{\mathrm{ext}}
                \coloneqq{}&
                (1+R_Q^2)R_Q^2
                \left(
                    \frac{6C_{Q^3,0}}{\pi R_Q^4}
                    +
                    \frac{3840C_{Q^3,2}}{\pi R_Q^6}
                    +
                    \frac4\pi (C_Q^{\mathrm{tail}})^3
                \right)  \\
                &+
                \left(1+\frac1{R_Q^2}\right)
                \left(\frac{R_Q}{R_Q-1}\right)^4
                \tilde R_Q^4
                \left(
                    \frac4\pi F_X+\frac{16}{\pi}F_X''
                \right).
            \end{aligned}
        \end{equation}
    \end{enumerate}
\end{definition}

\begin{lemma}[Exterior weighted bounds for \(Q\)]
\label{lem:app-exterior-weighted-Q-bounds}
    Let \(R_M\geq2\) and \(R_Q\geq2R_M\), and let the constants in
    Definition~\ref{def:weighted-q-constants} be as above. Then
    \begin{equation}\label{eq:app-Q-tail-bound}
        Q(y)
        \leq
        \frac{C_{Q,\mathrm{lead}}^{\mathrm{tail}}}{|y|^2}
        +
        \frac{C_{Q,\mathrm{rem}}^{\mathrm{tail}}}{|y|^3},
        \qquad
        |y|\geq R_M.
    \end{equation}
    Moreover,
    \[
        \int_{\bbR}Q^3\,dy
        \leq
        C_{Q^3,0},
        \qquad
        \int_{\bbR}y^2Q^3\,dy
        \leq
        C_{Q^3,2},
    \]
    and, for \(a=0,1,2\),
    \[
        \sup_{|y|\geq R_Q}
        \langle y\rangle^2
        \left|
            y^a\partial_y^aQ(y)
        \right|
        \leq
        C_{Q,a}^{\mathrm{ext}}.
    \]
\end{lemma}

\begin{proof}
    By \(Q=G*Q^3\), \(\int_{\bbR}G=1\), and the definition of
    \(\Delta_{\mathrm{cub}}^{L^1}\),
    \[
        \int_{\bbR}Q\,dy
        =
        \int_{\bbR}Q^3\,dy
        \leq
        \int_{\bbR}\mathfrak g^3\,dy
        +
        \Delta_{\mathrm{cub}}^{L^1}
        =
        C_{Q^3,0}.
    \]
    Since \(Q\) is positive, even, and decreasing on \((0,\infty)\),
    \[
        Q(z)
        \leq
        \frac{1}{2|z|}
        \int_{-|z|}^{|z|}Q(y)\,dy
        \leq
        \frac{C_{Q^3,0}}{2|z|},
        \qquad z\neq0.
    \]
    For \(y\geq R_M\), split \(Q=G*Q^3\) into \(\{|z|\leq y/2\}\) and \(\{|z|>y/2\}\). By
    \eqref{eq:kernel estimate of G},
    \[
        \begin{aligned}
            Q(y)
            \leq
            \frac{4}{\pi y^2}
            \int_{\bbR}Q^3\,dz
            +
            \sup_{|z|>y/2}Q(z)^3\int_{\bbR}G
            \leq
            \frac{4C_{Q^3,0}}{\pi y^2}
            +
            \frac{C_{Q^3,0}^3}{y^3}.
        \end{aligned}
    \]
    This proves \eqref{eq:app-Q-tail-bound}. It also gives
    \[
        \begin{aligned}
            \int_{\bbR}y^2Q^3\,dy
            \leq
            \int_{\bbR}y^2\mathfrak g^3\,dy
            +
            R_M^2\Delta_{\mathrm{cub}}^{L^1}
            \quad+
            2\int_{R_M}^{\infty}
            y^2
            \left(
                \frac{C_{Q,\mathrm{lead}}^{\mathrm{tail}}}{y^2}
                +
                \frac{C_{Q,\mathrm{rem}}^{\mathrm{tail}}}{y^3}
            \right)^3dy
            =
            C_{Q^3,2}.
        \end{aligned}
    \]

    It remains to prove the exterior bounds. By evenness, it suffices to
    consider \(y\geq R_Q\). For \(a=0,1,2\),
    \[
        \partial_y^aQ(y)
        =
        \int_0^\infty
        \left[
            \partial_y^aG(y-z)
            +
            \partial_y^aG(y+z)
        \right]
        Q(z)^3\,dz.
    \]
    For $a=1,2$, the preceding identities are understood in the distributional, equivalently principal-value, sense. The oddness and cancellation identities, together with the integration-by-parts formulas below, rewrite the local contributions as absolutely convergent integrals. On \(0\leq z\leq y/2\), \eqref{eq:kernel estimate of G} gives
    \[
        \begin{aligned}
            \left|
                \int_0^{y/2}
                [G(y-z)+G(y+z)]Q(z)^3\,dz
            \right|
            &\leq
            \frac{C_{Q^3,0}}{\pi y^2}
            +
            \frac{5C_{Q^3,2}}{\pi y^4},\\
            \left|
                \int_0^{y/2}
                [G'(y-z)+G'(y+z)]Q(z)^3\,dz
            \right|
            &\leq
            \frac{2C_{Q^3,0}}{\pi y^3}
            +
            \frac{26C_{Q^3,2}}{\pi y^5},\\
            \left|
                \int_0^{y/2}
                [G''(y-z)+G''(y+z)]Q(z)^3\,dz
            \right|
            &\leq
            \frac{6C_{Q^3,0}}{\pi y^4}
            +
            \frac{3840C_{Q^3,2}}{\pi y^6}.
        \end{aligned}
    \]
    If \(z>y/2\), then \(z\geq R_M\), and
    \eqref{eq:app-Q-tail-bound} gives
    \[
        Q(z)
        \leq
        \frac{4C_{Q,\mathrm{lead}}^{\mathrm{tail}}}{y^2}
        +
        \frac{8C_{Q,\mathrm{rem}}^{\mathrm{tail}}}{y^3}.
    \]
    Hence, for \(y\geq R_Q\)
    \[
        \begin{aligned}
            (1+y^2)|Q(y)|
            \leq
            (1+y^2)
            \left[
                \frac{C_{Q^3,0}}{\pi y^2}
                +
                \frac{5C_{Q^3,2}}{\pi y^4}
                +
                \left(
                    \frac{4C_{Q,\mathrm{lead}}^{\mathrm{tail}}}{y^2}
                    +
                    \frac{8C_{Q,\mathrm{rem}}^{\mathrm{tail}}}{y^3}
                \right)^3
            \right]
            \leq
            C_{Q,0}^{\mathrm{ext}}.
        \end{aligned}
    \]
    For \(z>y/2\),
    \[
        \left|\partial_z(Q^3)(z)\right|
        \leq
        3
        \left(
            \frac{4C_{Q,\mathrm{lead}}^{\mathrm{tail}}}{y^2}
            +
            \frac{8C_{Q,\mathrm{rem}}^{\mathrm{tail}}}{y^3}
        \right)^2
        \left(
            \frac{
                8\norm{\langle y\rangle^2y\mathfrak g'}_{L^\infty}
            }{y^3}
            +
            \Delta_{Q,1}^{L^\infty}
        \right).
    \]
    Split the integral over \(z>y/2\) further into
    \(\{|y-z|\leq1\}\) and \(\{|y-z|>1\}\). The oddness of \(G'\) and
    \eqref{eq:G-local-derivative-integrals} give
    \[
        \begin{aligned}
            (1+y^2)y|Q'(y)|
            \leq{}&
            (1+y^2)y
            \bigg[
                \frac{2C_{Q^3,0}}{\pi y^3}
                +
                \frac{26C_{Q^3,2}}{\pi y^5}\\
                &\quad+
                12
                \left(
                    \frac{4C_{Q,\mathrm{lead}}^{\mathrm{tail}}}{y^2}
                    +
                    \frac{8C_{Q,\mathrm{rem}}^{\mathrm{tail}}}{y^3}
                \right)^2
                \left(
                    \frac{
                        8\norm{\langle y\rangle^2y\mathfrak g'}_{L^\infty}
                    }{y^3}
                    +
                    \Delta_{Q,1}^{L^\infty}
                \right)\\
                &\quad+
                \frac4\pi
                \left(
                    \frac{4C_{Q,\mathrm{lead}}^{\mathrm{tail}}}{y^2}
                    +
                    \frac{8C_{Q,\mathrm{rem}}^{\mathrm{tail}}}{y^3}
                \right)^3
            \bigg]\\
            \leq{}&
            C_{Q,1}^{\mathrm{ext}}.
        \end{aligned}
    \]
    Here the last inequality follows by monotonicity for \(y\geq R_Q\);
    at \(y=R_Q\), the third term in brackets equals
    \(4F_{\mathrm{loc}}'\). Finally, let \(\tilde R_Q=R_Q-1\). For \(z\geq \tilde R_Q\),
    \[
        |Q(z)|
        \leq
        Q_X\frac{1+\tilde R_Q^2}{1+z^2},
        \qquad
        |Q'(z)|
        \leq
        Q_X'
        \frac{(1+\tilde R_Q^2)\tilde R_Q}{(1+z^2)z},
        \qquad
        |Q''(z)|
        \leq
        Q_X''.
    \]
    Therefore,
    \[
        |Q(z)^3|
        \leq
        F_X
        \left(
            \frac{1+\tilde R_Q^2}{1+z^2}
        \right)^2,
        \qquad
        \left|\partial_z^2(Q^3)(z)\right|
        \leq
        F_X''
        \left(
            \frac{1+\tilde R_Q^2}{1+z^2}
        \right)^2.
    \]
    The change of variables \(s=y-z\), integration by parts, and the
    oddness of \(G'\) give
    \[
        \begin{aligned}
            \int_{y-1}^{y+1}
                G''(y-z)Q(z)^3
            \,dz
            =
            G'(1)
            \left[
                Q(y-1)^3+Q(y+1)^3
            \right]
            +
            \int_{-1}^{1}
                G'(s)
                \left[
                    (Q^3)'(y-s)-(Q^3)'(y)
                \right]
            \,ds.
        \end{aligned}
    \]
    Hence \eqref{eq:G-local-derivative-integrals} and the mean value
    theorem give
    \[
        \left|
            \int_{y-1}^{y+1}
                G''(y-z)Q(z)^3
            \,dz
        \right|
        \leq
        \left(
            \frac4\pi F_X
            +
            \frac{16}{\pi}F_X''
        \right)
        \left(
            \frac{1+\tilde R_Q^2}{1+(y-1)^2}
        \right)^2.
    \]
    since \(1\leq16/\pi\).
    Combining this estimate with the bounds on
    \(0\leq z\leq y/2\) and \(\{|y-z|>1\}\) gives
    \[
        \begin{aligned}
            (1+y^2)y^2|Q''(y)|
            \leq{}&
            (1+y^2)y^2
            \bigg[
                \frac{6C_{Q^3,0}}{\pi y^4}
                +
                \frac{3840C_{Q^3,2}}{\pi y^6}
                +
                \frac4\pi
                \left(
                    \frac{4C_{Q,\mathrm{lead}}^{\mathrm{tail}}}{y^2}
                    +
                    \frac{8C_{Q,\mathrm{rem}}^{\mathrm{tail}}}{y^3}
                \right)^3
            \bigg]\\
            &+
            (1+y^2)y^2
            \left(
                \frac{1+\tilde R_Q^2}{1+(y-1)^2}
            \right)^2
            \left(
                \frac4\pi F_X
                +
                \frac{16}{\pi}F_X''
            \right)\\
            \leq{}&
            C_{Q,2}^{\mathrm{ext}}.
        \end{aligned}
    \]
    The last inequality follows by monotonicity for
    \(y\geq R_Q\geq4\).
\end{proof}

\begin{definition}[Global bounds for \(Q\) and \(V_Q\)]\label{def:global-Q-VQ-bounds}
    The interior bounds in \eqref{eq:def-CQ-core} and the exterior bounds
    of Lemma~\ref{lem:app-exterior-weighted-Q-bounds} give the following
    global constants. For \(a=0,1,2\), set
    \begin{equation}\label{eq:def-CQj-weighted-envelope}
        C_{Q,a}
        \coloneqq
        \max\{C_{Q,a}^{\mathrm{int}},C_{Q,a}^{\mathrm{ext}}\}.
    \end{equation}
    Since \(|V_Q(y)| \leq C_{Q,0}C_{Q,1}\langle y\rangle^{-4}\), set
    \begin{equation}\label{eq:def-CVQ-and-VQ-bound}
        \mathsf V^{L^\infty}
        \coloneqq
        C_{Q,0}C_{Q,1},
        \qquad
        \mathsf V^{L^1}
        \coloneqq
        \frac{\pi}{2}C_{Q,0}C_{Q,1}.
    \end{equation}
    Set
    \begin{equation}\label{eq:def-VQ-omega-bound}
        I_{\omega^2\langle y\rangle^{-4}}
        \coloneqq
        \int_{\bbR}\omega(y)^2\langle y\rangle^{-4}\,dy,
        \qquad
        \mathsf V_{\omega^2}^{L^1}
        \coloneqq
        C_{Q,0}C_{Q,1}I_{\omega^2\langle y\rangle^{-4}}.
    \end{equation}
\end{definition}

\subsection{Scaling and \(L^1\) estimates}
The next two definitions give the \(L^2\)-bounds for \(\Lambda^ae_Q\), followed by the \(L^1\)- and logarithmically weighted \(L^1\)-bounds used in Lemma~\ref{lem:scaling-and-L1-transfer}.

\begin{definition}[Constants for the scaling estimates]
    Let \(C_{Q,0},C_{Q,1},C_{Q,2}\) be the constants from Definition~\ref{def:global-Q-VQ-bounds}. We use  \(R_E\ge1\).  We also use the fixed inverse constant
    \begin{equation}\label{eq:def-K-inv}
        K_{\mathrm{inv}}
        \coloneqq
        2,
        \qquad
        \norm{[\mathcal F'(\mathfrak g)]^{-1}}_
        {L^2_{\mathrm e}\to L^2_{\mathrm e}}
        \le
        K_{\mathrm{inv}}.
    \end{equation}

    \begin{enumerate}[label=\textup{(\roman*)}]
        \item The following constants bound \(Q^2-\mathfrak g^2\) in \(L^\infty\) and isolate the term containing \(\Lambda e_Q\):
        \begin{equation}\label{eq:def-Delta-Q-L2}
            \Delta_Q^{L^2}
            \coloneqq
            r_{\mathrm{NK}}.
        \end{equation}
        Define the upper bound of $\norm{Q^2-\mathfrak g^2}_{L^\infty}$:
        \begin{equation}\label{eq:def-Delta-Q2-Linf-kappa-Lambda}
            \Delta_{Q^2}^{L^\infty}
            \coloneqq
            \Delta_{Q,0}^{L^\infty}
            \left(
                \norm{\mathfrak g}_{L^\infty}
                +
                \Delta_{Q,0}^{L^\infty}
                +
                \norm{\mathfrak g}_{L^\infty}
            \right),
            \qquad
            \kappa_\Lambda
            \coloneqq
            1-3K_{\mathrm{inv}}\Delta_{Q^2}^{L^\infty}.
        \end{equation}
        Assume \(\kappa_\Lambda>0\), and set
        \begin{equation}\label{eq:def-Delta-Lambda-Q-L2}
            \Delta_{\Lambda Q}^{L^2}
            \coloneqq
            \frac{
                K_{\mathrm{inv}}
                \left(
                    \Delta_Q^{L^2}
                    +
                    \norm{\Lambda r_{\mathfrak g}}_{L^2}
                    +
                    \norm{r_{\mathfrak g}}_{L^2}
                    +
                    3\Delta_{Q^2}^{L^\infty}\norm{\Lambda\mathfrak g}_{L^2}
                \right)
            }{
                \kappa_\Lambda
            }.
        \end{equation}

        \item For the estimate of \(\Lambda^2e_Q\), set
        \begin{equation}\label{eq:def-QLambda-inputs}
            C_{\Lambda Q}^{L^2}
            \coloneqq
            \norm{\Lambda\mathfrak g}_{L^2}
            +
            \Delta_{\Lambda Q}^{L^2},
            \qquad
            C_{\Lambda Q}^{L^\infty}
            \coloneqq
            \frac12C_{Q,0}+C_{Q,1},
        \end{equation}
        and
        \begin{equation}\label{eq:def-GLambda-Linfty}
            C_{\Lambda\mathfrak g}^{L^\infty}
            \coloneqq
            \frac12\norm{\mathfrak g}_{L^\infty}
            +
            \norm{\langle y\rangle^2y\mathfrak g'}_{L^\infty}.
        \end{equation}
        The cubic difference in the scaling identity is bounded by
        \begin{equation}\label{eq:def-CLambda}
            \Delta_{\Lambda,\mathrm{cub}}^{L^2}
            \coloneqq
            \Delta_{Q,0}^{L^\infty}C_{\Lambda Q}^{L^\infty}C_{\Lambda Q}^{L^2}
            +
            \norm{\mathfrak g}_{L^\infty}
            \left(
                C_{\Lambda Q}^{L^\infty}
                +
                C_{\Lambda\mathfrak g}^{L^\infty}
            \right)
            \Delta_{\Lambda Q}^{L^2}.
        \end{equation}
        Finally, set
        \begin{equation}\label{eq:def-Delta-Lambda2-Q-L2}
            \begin{aligned}
                \Delta_{\Lambda^2Q}^{L^2}
                \coloneqq
                \frac{K_{\mathrm{inv}}}{\kappa_\Lambda}
                \bigl(
                    &\Delta_Q^{L^2}
                    +
                    2\Delta_{\Lambda Q}^{L^2}
                    +
                    \norm{\Lambda^2r_{\mathfrak g}}_{L^2}
                    +
                    2\norm{\Lambda r_{\mathfrak g}}_{L^2}
                    +
                    \norm{r_{\mathfrak g}}_{L^2}  \\
                    &+
                    6\Delta_{\Lambda,\mathrm{cub}}^{L^2}
                    +
                    3\Delta_{Q^2}^{L^\infty}\norm{\Lambda^2\mathfrak g}_{L^2}
                \bigr).
            \end{aligned}
        \end{equation}
    \end{enumerate}
\end{definition}

\begin{definition}[Constants for the \(L^1\) estimates]
\label{def:dilation-L1-transfer-radii}
    Fix \(R_\omega\geq1\) and \(R_E\geq1\).
    \begin{enumerate}[label=\textup{(\roman*)}]
        \item
        For \(a=0,1,2\), the following constants give the pointwise bounds
        \[
            |\Lambda^ae_Q(y)|
            \leq
            C_{e,a}^{\mathrm{tail}}\langle y\rangle^{-2}.
        \]
        Set
        \begin{equation}\label{eq:def-Ctail0,1,2}
            \begin{aligned}
                C_{e,0}^{\mathrm{tail}}
                &\coloneqq
                C_{Q,0}
                +
                \norm{\langle y\rangle^2\mathfrak g}_{L^\infty},\\
                C_{e,1}^{\mathrm{tail}}
                &\coloneqq
                \left(
                    \frac12C_{Q,0}+C_{Q,1}
                \right)
                +
                \left(
                    \frac12\norm{\langle y\rangle^2\mathfrak g}_{L^\infty}
                    +
                    \norm{\langle y\rangle^2y\mathfrak g'}_{L^\infty}
                \right),\\
                C_{e,2}^{\mathrm{tail}}
                &\coloneqq
                \left(
                    \frac14C_{Q,0}
                    +
                    2C_{Q,1}
                    +
                    C_{Q,2}
                \right)  \\
                &\quad+
                \left(
                    \frac14\norm{\langle y\rangle^2\mathfrak g}_{L^\infty}
                    +
                    2\norm{\langle y\rangle^2y\mathfrak g'}_{L^\infty}
                    +
                    \norm{\langle y\rangle^2y^2\mathfrak g''}_{L^\infty}
                \right).
            \end{aligned}
        \end{equation}

        \item
        For \(R\geq1\), set
        \begin{equation}\label{eq:def-Iomega}
            I_\omega(R)
            \coloneqq
            \frac{2}{R^{3}}
            \left(
                \frac{\log^2(2R)}{3}
                +
                \frac{2\log(2R)}{9}
                +
                \frac{2}{27}
            \right),
        \end{equation}
        and define
        \begin{equation}\label{eq:def-Delta-omega-Q-L2}
            \Delta_{\omega Q}^{L^2}
            \coloneqq
            \log(1+R_\omega)\Delta_Q^{L^2}
            +
            C_{e,0}^{\mathrm{tail}}I_\omega(R_\omega)^{1/2}.
        \end{equation}
        For \(C\geq0\) and \(R\geq1\), set
        \begin{equation}\label{eq:def-T0-Tomega}
            \mathsf T_0(C,R)
            \coloneqq
            \frac{2C}{R},
            \qquad
            \mathsf T_\omega(C,R)
            \coloneqq
            2C
            \left(
                \frac{\log(1+R)}{R}
                +
                \log\left(1+\frac1R\right)
            \right).
        \end{equation}

        \item
        Using the \(L^2\)-bounds on \(\{|y|\leq R_E\}\) and the pointwise tail bounds on \(\{|y|>R_E\}\), define
        \begin{equation}\label{eq:def-Delta-Lambda-a-Q-L1}
            \Delta_{\Lambda^aQ}^{L^1}
            \coloneqq
            \sqrt{2R_E}\,\Delta_{\Lambda^aQ}^{L^2}
            +
            \mathsf T_0(C_{e,a}^{\mathrm{tail}},R_E),
       \end{equation}
        while
        \begin{equation}\label{eq:def-Delta-omega-Q-L1}
            \Delta_{\omega Q}^{L^1}
            \coloneqq
            \sqrt{2R_E}\,\Delta_{\omega Q}^{L^2}
            +
            \mathsf T_\omega(C_{e,0}^{\mathrm{tail}},R_E),
        \end{equation}
        and, for \(a=1,2\),
        \begin{equation}\label{eq:def-Delta-omega-Lambda-a-Q-L1}
            \Delta_{\omega\Lambda^aQ}^{L^1}
            \coloneqq
            \sqrt{2R_E}\log(1+R_E)\Delta_{\Lambda^aQ}^{L^2}
            +
            \mathsf T_\omega(C_{e,a}^{\mathrm{tail}},R_E).
        \end{equation}
    \end{enumerate}
\end{definition}

\subsection{Constants for \(\delta_V\) and \(\psi_m^{(j)}-\psi_m^{(\mathfrak g,j)}\)}
The first definition gives the constants used in Lemma~\ref{lem:square-root-potential-transfer}. The subsequent lemma estimates \(\psi_m^{(j)}-\psi_m^{(\mathfrak g,j)}\).
\begin{definition}[Constants for \(\delta_V\)]
    Let \(\omega(y)\coloneqq \log(1+|y|)\).  We first record the two
    elementary integrals used to estimate \(yQ\) and \(\omega^2yQ\) from
    the bound \(\langle y\rangle^2 |Q|\le C_{Q,0}\).  Define
    \begin{equation}\label{eq:def-I-y2-langle4}
        I_{y^2\langle y\rangle^{-4}}
        \coloneqq
        \int_{\bbR}y^2\langle y\rangle^{-4}\,dy
        =
        \frac{\pi}{2},\qquad
        I_{\omega^4y^2\langle y\rangle^{-4}}
        \coloneqq
        \int_{\bbR}
        \omega(y)^4y^2\langle y\rangle^{-4}\,dy.
    \end{equation}
    Set
    \begin{equation}\label{eq:def-yQ-radii}
        C_{yQ}^{L^2}
        \coloneqq
        C_{Q,0}
        \left(
            I_{y^2\langle y\rangle^{-4}}
        \right)^{1/2},
        \qquad
        C_{\omega^2yQ}^{L^2}
        \coloneqq
        C_{Q,0}
        \left(
            I_{\omega^4y^2\langle y\rangle^{-4}}
        \right)^{1/2}.
    \end{equation}
    Thus \(C_{yQ}^{L^2}\) bounds \(\norm{yQ}_{L^2}\), while
    \(C_{\omega^2yQ}^{L^2}\) bounds \(\norm{\omega^2yQ}_{L^2}\).
    For
    \(
        \delta_V
        \coloneqq
        |V_Q|^{1/2}-|V_{\mathfrak g}|^{1/2}
    \),
    set
    \begin{equation}\label{eq:def-DeltaV0}
        \Delta_{\sqrt V}^{L^2}
        \coloneqq
        \left[
            C_{yQ}^{L^2}r_{\mathrm{NK}}
            +
            \norm{y\mathfrak g'}_{L^2}r_{\mathrm{NK}}
        \right]^{1/2},
    \end{equation}
    and
    \begin{equation}\label{eq:def-DeltaVomega}
        \Delta_{\omega\sqrt V}^{L^2}
        \coloneqq
        \left[
            C_{\omega^2yQ}^{L^2}r_{\mathrm{NK}}
            +
            \norm{\omega^2y\mathfrak g'}_{L^2}r_{\mathrm{NK}}
        \right]^{1/2}.
    \end{equation}
\end{definition}

\begin{lemma}[Estimates for \(\psi_m^{(j)}-\psi_m^{(\mathfrak g,j)}\)]
\label{lem:app-psi-difference-estimates}
    For \(j,m\in\{1,2\}\), set \(p\coloneqq m+j-2\), and recall that
    \[
        \begin{aligned}
            \Psi_m^{(j)}
            &\coloneqq
            \Lambda^pQ,
            &
            \Psi_m^{(\mathfrak g,j)}
            &\coloneqq
            \Lambda^p\mathfrak g,\\
            \psi_m^{(j)}
            &\coloneqq
            |V_Q|^{1/2}\mathcal P_{\log}\Psi_m^{(j)},
            &
            \psi_m^{(\mathfrak g,j)}
            &\coloneqq
            |V_{\mathfrak g}|^{1/2}
            \mathcal P_{\log}\Psi_m^{(\mathfrak g,j)}.
        \end{aligned}
    \]
    With the constants defined in
    \eqref{eq:def-CVQ-and-VQ-bound},
    \eqref{eq:def-Delta-Q-L2},
    \eqref{eq:def-Delta-Lambda-a-Q-L1},
    \eqref{eq:def-DeltaV0},
    \eqref{eq:def-DeltaVomega}, and
    \eqref{eq:part3-Eelog-def}, set
    \begin{equation}\label{eq:app-Deltapsi-def}
    \begin{aligned}
        \Delta_{\psi,m}^{(j)}
        \coloneqq{}&
        \frac1\pi
        \left(
            2\norm{\Psi_m^{(\mathfrak g,j)}}_{L^2}
            +
            \norm{\omega\Psi_m^{(\mathfrak g,j)}}_{L^1}
        \right)
        \Delta_{\sqrt V}^{L^2}
        +
        \frac1\pi
        \norm{\Psi_m^{(\mathfrak g,j)}}_{L^1}
        \Delta_{\omega\sqrt V}^{L^2}
        \\
        &+
        \begin{cases}
            \mathsf E_{e,\log},
            & p=0,\\[1mm]
            \displaystyle
            \frac12 (C_{Q,0}C_{Q,1})^{1/2}\Delta_Q^{L^2}
            +
            \frac12\mathsf E_{e,\log}
            +
            \frac1\pi\Delta_Q^{L^1}\left(\mathsf V^{L^1}\right)^{1/2},
            & p=1,\\[3mm]
            \displaystyle
            (C_{Q,0}C_{Q,1})^{1/2}\Delta_Q^{L^2}
            +
            \frac14\mathsf E_{e,\log}
            +
            \frac1\pi\Delta_Q^{L^1}\left(\mathsf V^{L^1}\right)^{1/2},
            & p=2.
        \end{cases}
    \end{aligned}
    \end{equation}
    Then
    \begin{equation}\label{eq:app-psi-difference-bound}
        \norm{
            \psi_m^{(j)}
            -
            \psi_m^{(\mathfrak g,j)}
        }_{L^2}
        \leq
        \Delta_{\psi,m}^{(j)}.
    \end{equation}
\end{lemma}

\begin{proof}
    By definition,
    \[
        \psi_m^{(j)}-\psi_m^{(\mathfrak g,j)}
        =
        \delta_V
        \mathcal P_{\log}\Psi_m^{(\mathfrak g,j)}
        +
        |V_Q|^{1/2}
        \mathcal P_{\log}(\Lambda^pe_Q).
    \]
    Applying \eqref{eq:part3-basic-log-potential-bound} with \(f_3=\delta_V\) and
    \(f_1=\Psi_m^{(\mathfrak g,j)}\), and then using
    \eqref{eq:def-DeltaV0} and \eqref{eq:def-DeltaVomega}, gives
    \[
        \begin{aligned}
            \norm{
                \delta_V
                \mathcal P_{\log}\Psi_m^{(\mathfrak g,j)}
            }_{L^2}
            \leq
            \frac1\pi
            \left(
                2\norm{\Psi_m^{(\mathfrak g,j)}}_{L^2}
                +
                \norm{\omega\Psi_m^{(\mathfrak g,j)}}_{L^1}
            \right)
            \Delta_{\sqrt V}^{L^2}
            +
            \frac1\pi
            \norm{\Psi_m^{(\mathfrak g,j)}}_{L^1}
            \Delta_{\omega\sqrt V}^{L^2}.
        \end{aligned}
    \]
    It remains to bound the second term. For \(p=0\),
    \eqref{eq:part3-Eelog-bound} gives
    \[
        \norm{
            |V_Q|^{1/2}\mathcal P_{\log}e_Q
        }_{L^2}
        \leq
        \mathsf E_{e,\log}.
    \]
    For \(p=1,2\), equations
    \eqref{eq:part3-Plog-Lambda-identity} and
    \eqref{eq:part3-Plog-Lambda2-identity} give
    \[
        \mathcal P_{\log}(\Lambda^pe_Q)
        =
        \begin{cases}
            \displaystyle
            y\partial_y\mathcal P_{\log}e_Q
            -
            \frac12\mathcal P_{\log}e_Q
            +
            \frac1\pi\int_{\bbR}e_Q,
            & p=1,\\[3mm]
            \displaystyle
            y^2\partial_y^2\mathcal P_{\log}e_Q
            +
            \frac14\mathcal P_{\log}e_Q
            -
            \frac1\pi\int_{\bbR}e_Q,
            & p=2.
        \end{cases}
    \]
    Moreover,
    \[
        \norm{y|V_Q|^{1/2}}_{L^\infty}
        \leq
        \frac12(C_{Q,0}C_{Q,1})^{1/2},
        \qquad
        \norm{y^2|V_Q|^{1/2}}_{L^\infty}
        \leq
        (C_{Q,0}C_{Q,1})^{1/2},
    \]
    and
    \[
        \norm{|V_Q|^{1/2}}_{L^2}
        \leq
        (\mathsf V^{L^1})^{1/2}.
    \]
    Thus \eqref{eq:part3-Plog-Hilbert-bounds},
    \eqref{eq:part3-Eelog-bound}, and \eqref{eq:def-Delta-Lambda-a-Q-L1} yield
    \[
        \norm{
            |V_Q|^{1/2}
            \mathcal P_{\log}(\Lambda^pe_Q)
        }_{L^2}
        \leq
        \begin{cases}
            \mathsf E_{e,\log},
            & p=0,\\[1mm]
            \displaystyle
            \frac12(C_{Q,0}C_{Q,1})^{1/2}\Delta_Q^{L^2}
            +
            \frac12\mathsf E_{e,\log}
            +
            \frac1\pi
            \Delta_Q^{L^1}(\mathsf V^{L^1})^{1/2},
            & p=1,\\[3mm]
            \displaystyle
            (C_{Q,0}C_{Q,1})^{1/2}\Delta_Q^{L^2}
            +
            \frac14\mathsf E_{e,\log}
            +
            \frac1\pi
            \Delta_Q^{L^1}(\mathsf V^{L^1})^{1/2},
            & p=2.
        \end{cases}
    \]
    The two estimates give \eqref{eq:app-psi-difference-bound}.
\end{proof}

\section{Compactification identities for computations}
\label{app:compactification-profile-calculus}

We record the compactified formulas used in the interval computations of
Part~\ref{part:interval-verification}. They concern \(\mathfrak g\), its residual, its scaling
derivatives, \(V_{\mathfrak g}\), \(a_{\mathfrak g}\), and the functions
\(\Omega_m^{(j)}\).

Recall
\[
    I=
    \left(
        -\frac{\pi}{2},
        \frac{\pi}{2}
    \right),
    \qquad
    y=\tan\theta,
\]
and, from \eqref{eq:g-side-density-compactification-map},
\[
    \mathscr C f(\theta)
    =
    f(\tan\theta)\sec^2\theta.
\]
Thus
\(
    f(\tan\theta)
    =
    \cos^2\theta\,\mathscr C f(\theta)
\).
With
\(
    \omega(y)=\log(1+|y|)
\),
we have, for every \(p\in\mathbb N_0\),
\begin{equation*}
    \begin{aligned}
        \norm{f}_{L^1(\bbR)}
        &=
        \int_I
            |\mathscr C f(\theta)|
        \,d\theta,
        \\
        \norm{f}_{L^2(\bbR)}^2
        &=
        \int_I
            |\mathscr C f(\theta)|^2
            \cos^2\theta
        \,d\theta,
        \\
        \norm{\omega^p f}_{L^1(\bbR)}
        &=
        \int_I
            \log^p(1+|\tan\theta|)
            |\mathscr C f(\theta)|
        \,d\theta.
    \end{aligned}
\end{equation*}
The derivative and the scaling generator \(\Lambda=\frac12+y\partial_y\) satisfy
\begin{equation}\label{eq:app-compactified-Lambda-action}
    \begin{aligned}
        \mathscr C(\partial_yf)
        &=
        \partial_\theta
        \bigl(
            \cos^2\theta\,\mathscr C f
        \bigr),
        \\
        \mathscr C(\Lambda f)
        &=
        \mathscr L_\theta(\mathscr C f).
    \end{aligned}
\end{equation}
where
\begin{equation*}
        \mathscr L_\theta f_1
        \coloneqq
        \left(
            \cos(2\theta)-\frac12
        \right)f_1
        +
        \sin\theta\cos\theta\,
        \partial_\theta f_1.
\end{equation*}
In particular, if \(f_1(\theta)=P(\cos(2\theta))\) for some polynomial \(P\), then
\[
    \mathscr L_\theta f_1
    =
    \left[
        \left(
            x-\frac12
        \right)P(x)
        -
        (1-x^2)P'(x)
    \right]_{x=\cos(2\theta)}.
\]
The conformal covariance of \(|D|\) gives
\begin{equation}\label{eq:app-compactified-halfwave-identity}
    \mathscr C(|D|f)
    =
    |\partial_\theta|
    \bigl(
        \cos^2\theta\,\mathscr C f
    \bigr),
\end{equation}
where \(|\partial_\theta|\) is the \(\pi\)-periodic Fourier multiplier with symbol \(2|n|\) on the mode \(\cos(2n\theta)\). Recall \(U\) and \(\mathfrak g\) from
Definition~\ref{def:approximate-profile-g}. By
\eqref{eq:def-g-compactified},
\[
    \mathscr C\mathfrak g=U.
\]
For
\(
    r_{\mathfrak g}
    =
    |D|\mathfrak g+\mathfrak g-\mathfrak g^3
\),
equation~\eqref{eq:app-compactified-halfwave-identity} gives
\begin{equation*}
    \mathscr C r_{\mathfrak g}
    =
    |\partial_\theta|
    \bigl(
        \cos^2\theta\,U
    \bigr)
    +
    U
    -
    \cos^4\theta\,U^3.
\end{equation*}
Repeated application of
\eqref{eq:app-compactified-Lambda-action} gives, for
\(a\in\mathbb N_0\),
\begin{equation*}
    \mathscr C(\Lambda^a\mathfrak g)
    =
    \mathscr L_\theta^aU,
    \qquad
    \mathscr C(\Lambda^a r_{\mathfrak g})
    =
    \mathscr L_\theta^a
    \bigl(
        \mathscr C r_{\mathfrak g}
    \bigr).
\end{equation*}
More generally, for \(a,b\in\mathbb N_0\),
\[
    \mathscr C
    \bigl(
        (\Lambda+1)^b\Lambda^a\mathfrak g
    \bigr)
    =
    (\mathscr L_\theta+1)^b
    \mathscr L_\theta^aU.
\]
All functions on the right-hand sides are finite cosine polynomials. In particular, the functions \(\Omega_m^{(j)}\) defined in
\eqref{eq:g-side-Omega-def} are finite cosine polynomials. Recall \(V_{\mathfrak g}\) from
\eqref{eq:definition of g-side potential}. Since
\(
    y\mathfrak g'
    =
    \Lambda\mathfrak g-\frac12\mathfrak g
\),
we obtain
\begin{equation*}
    \begin{aligned}
        (y\mathfrak g')(\tan\theta)
        &=
        \cos^2\theta
        \left(
            \mathscr L_\theta U-\frac12U
        \right)(\theta),
        \\
        V_{\mathfrak g}(\tan\theta)
        &=
        \cos^4\theta\,
        U(\theta)
        \left(
            \mathscr L_\theta U-\frac12U
        \right)(\theta),
        \\
        \mathscr C V_{\mathfrak g}
        &=
        \cos^2\theta\,
        U
        \left(
            \mathscr L_\theta U-\frac12U
        \right).
    \end{aligned}
\end{equation*}
Recalling \(a_{\mathfrak g}\) from
\eqref{eq:definition of a_g, tidle T_g}, the preceding identity gives
\begin{equation*}
    a_{\mathfrak g}(\theta)
    =
    \cos\theta
    \left|
        U(\theta)
        \left(
            \mathscr L_\theta U-\frac12U
        \right)(\theta)
    \right|^{1/2}.
\end{equation*}
Likewise, by \eqref{eq:definition of a_g,log},
\[
    a_{\mathfrak g,\log}(\theta)
    =
    a_{\mathfrak g}(\theta)
    \log|\cos\theta|.
\]
The change of variables \(y=\tan\theta\) also gives
\begin{equation}\label{eq:app-g-direct-moment-identities}
    \int_{\bbR}\mathfrak g^3\,dy
    =
    \int_I
        U(\theta)^3\cos^4\theta
    \,d\theta,
    \qquad
    \int_{\bbR}y^2\mathfrak g^3\,dy
    =
    \int_I
        U(\theta)^3
        \sin^2\theta\cos^2\theta
    \,d\theta,
\end{equation}
and
\begin{equation}\label{eq:app-g-direct-y-identities}
    \norm{y\mathfrak g}_{L^2}^2
    =
    \int_I
        \sin^2\theta\,U(\theta)^2
    \,d\theta,
    \qquad
    \norm{y\mathfrak g'}_{L^2}^2
    =
    \int_I
        \cos^2\theta
        \left(
            \mathscr L_\theta U-\frac12U
        \right)^2
    \,d\theta.
\end{equation}
Moreover,
\[
    \norm{V_{\mathfrak g}}_{L^1}
    =
    \int_I
        |\mathscr C V_{\mathfrak g}(\theta)|
    \,d\theta,
    \qquad
    \norm{
        \omega|V_{\mathfrak g}|^{1/2}
    }_{L^2}^2
    =
    \int_I
        \log^2(1+|\tan\theta|)
        |\mathscr C V_{\mathfrak g}(\theta)|
    \,d\theta.
\]

To control the endpoint contribution in the logarithmically weighted
integrals, let \(f_1\in C(\overline I)\) be even, let
\(p\in\mathbb N_0\), and fix \(0<\delta<1\). For a finite partition
\(\mathcal P\) of \([0,\pi/2-\delta]\), we have
\begin{equation}\label{eq:app-log-weight-endpoint-bound}
    \begin{aligned}
        2\int_0^{\pi/2}
            \log^p(1+\tan\theta)
            |f_1(\theta)|
        \,d\theta
        &\leq
        2\sum_{J\in\mathcal P}
            |J|
            \sup_{\theta\in J}
            \left\{
                \log^p(1+\tan\theta)
                |f_1(\theta)|
            \right\}
        \\
        &\quad+
        2\delta p!
        \left(
            \sup_{\pi/2-\delta\leq\theta\leq\pi/2}
            |f_1(\theta)|
        \right)
        \sum_{q=0}^{p}
            \frac{\log^q(2/\delta)}{q!}.
    \end{aligned}
\end{equation}
Indeed, with \(\theta=\pi/2-x\),
\[
    \log(1+\tan\theta)
    =
    \log(1+\cot x)
    \leq
    \log\frac{2}{x},
    \qquad
    0<x<\delta,
\]
and
\[
    \int_0^\delta
        \log^p(2/x)
    \,dx
    =
    \delta p!
    \sum_{q=0}^{p}
        \frac{\log^q(2/\delta)}{q!}.
\]
If \(f_1\geq0\), the absolute values around \(f_1\) in \eqref{eq:app-log-weight-endpoint-bound} may be omitted. The logarithmically weighted variants of \eqref{eq:app-g-direct-moment-identities} and \eqref{eq:app-g-direct-y-identities} are obtained by applying \eqref{eq:app-log-weight-endpoint-bound} to the corresponding compactified integrands.

\section{Validated-computation tools}
This section records the quadrature, DCT-I approximation, aliasing, matrix, and projection estimates used in Parts~\ref{part:ground-state-approximation}--\ref{part:interval-verification}.

\subsection{Gauss--Legendre quadrature and analytic continuation of \(\mathfrak g\)}
We recall the \(n\)-point Gauss--Legendre rule and the analytic error bound used in Lemma~\ref{lem:NK-endpoint-matrix-evaluation-error}.

For \(n\geq1\), let \(\mathbb P_n\) be the Legendre polynomial on \([-1,1]\), let \(\chi_1,\dots,\chi_n\) be its zeros, and set
\begin{equation}\label{eq:GL-weights}
    \omega_p
    \coloneqq
    \frac{2}{(1-\chi_p^2)\bigl((\mathbb P_n)'(\chi_p)\bigr)^2},
    \qquad
    p=1,\dots,n.
\end{equation}
We write
\begin{equation}\label{eq:GL-quadrature}
    \mathcal Q_n[f]\coloneqq \sum_{p=1}^n\omega_p f(\chi_p).
\end{equation}
For \(\varrho>1\), let
\begin{equation}\label{eq:Bernstein-ellipse}
    E_\varrho
    \coloneqq
    \left\{
        \frac12\left(\varrho e^{i\tau}+\varrho^{-1}e^{-i\tau}\right):
        \tau\in[0,2\pi]
    \right\}
\end{equation}
be the Bernstein ellipse with foci \(\pm1\).

\begin{lemma}[Gauss--Legendre quadrature]
\label{lem:GL-analytic}
    Let \(n\geq1\), and let \(\mathcal Q_n\) be the \(n\)-point
    Gauss--Legendre rule on \([-1,1]\).

    \begin{enumerate}[label=\textup{(\roman*)}]
        \item
        If \(f\) is a polynomial of degree at most \(2n-1\), then
        \begin{equation*}
            \int_{-1}^{1}f(\tau)\,d\tau
            =
            \mathcal Q_n[f].
        \end{equation*}

        \item
        Assume \(n\geq2\). Let \(\varrho>1\). If \(f\) is holomorphic in the interior of
        \(E_\varrho\) and continuous on \(E_\varrho\), then, with
        \[
            \mathscr M_f
            \coloneqq
            \sup_{w\in E_\varrho}|f(w)|,
        \]
        one has
        \begin{equation*}
            \left|
                \int_{-1}^{1}f(\tau)\,d\tau
                -
                \mathcal Q_n[f]
            \right|
            \leq
            \frac{64\mathscr M_f}{15(\varrho^2-1)}
            \varrho^{-2n+2}.
        \end{equation*}
    \end{enumerate}
\end{lemma}

\begin{proof}
    Part~\textup{(i)} is the \(n\)-point form of \cite[Theorem~19.1]{TrefethenATAP}.
    Part~\textup{(ii)} follows from \cite[Theorem~19.3]{TrefethenATAP} after replacing the parameter \(n\) there by \(n-1\).
\end{proof}

\begin{lemma}[Complex continuation and strip bound for \(\mathfrak g\)]
    The profile \(\mathfrak g\) extends holomorphically to \(\{w\in\bbC:|\Im w|<1\}\) by
    \begin{equation}\label{eq:g-complex-extension}
        \mathfrak g_{\mathbb C}(w)
        =
        \frac{1}{1+w^2}
        \sum_{\ell=0}^{200}
        U_\ell
        T_\ell\left(\frac{1-w^2}{1+w^2}\right),
    \end{equation}
    where \(T_\ell\) denotes the Chebyshev polynomial of degree \(\ell\).  Moreover,
    for every \(0<\delta<1\),
    \begin{equation}\label{eq:g-strip-bound}
        \sup_{|\Im w|\leq\delta}
        |\mathfrak g_{\mathbb C}(w)|
        \leq
        \mathfrak G_\delta,
        \qquad
        \mathfrak G_\delta
        \coloneqq
        \frac{1}{(1-\delta)^2}
        \sum_{\ell=0}^{200}
        |U_\ell|
        \left(
            \frac{1+\delta}{1-\delta}
        \right)^\ell .
    \end{equation}
\end{lemma}

\begin{proof}
    For real \(y\), writing \(\theta=\arctan y\), we have
    \[
        \frac{1-y^2}{1+y^2}=\cos(2\theta),
        \qquad
        T_\ell\left(\frac{1-y^2}{1+y^2}\right)=\cos(2\ell\theta).
    \]
    Hence \(\mathfrak g_{\mathbb C}(y)=\mathfrak g(y)\).
    To prove the strip bound, set
    \[
        q(w)\coloneqq \frac{1+iw}{1-iw}.
    \]
    Then
    \[
        \frac{1-w^2}{1+w^2}
        =
        \frac12\left(q(w)+q(w)^{-1}\right),
        \qquad
        T_\ell\left(\frac{1-w^2}{1+w^2}\right)
        =
        \frac12\left(q(w)^\ell+q(w)^{-\ell}\right).
    \]
    Write \(w=x+i\eta\), with \(|\eta|\le\delta<1\). Since
    \[
        |1+w^2|
        =
        |1+iw|\,|1-iw|
        \ge
        (1-\delta)^2,
    \]
    and
    \[
        |q(w)|^2
        =
        \frac{x^2+(1-\eta)^2}{x^2+(1+\eta)^2},
        \qquad
        |q(w)^{-1}|^2
        =
        \frac{x^2+(1+\eta)^2}{x^2+(1-\eta)^2},
    \]
    we have
    \[
        |q(w)|,\ |q(w)^{-1}|
        \le
        \frac{1+\delta}{1-\delta}.
    \]
    Therefore,
    \[
        \left|
            T_\ell\left(\frac{1-w^2}{1+w^2}\right)
        \right|
        \le
        \left(
            \frac{1+\delta}{1-\delta}
        \right)^\ell .
    \]
    Substituting this into \eqref{eq:g-complex-extension} gives
    \eqref{eq:g-strip-bound}.
\end{proof}

\subsection{DCT-I approximation and aliasing estimates}\label{subsec:app-dct1-aliasing}

For even functions on \(I\), set \(I_+\coloneqq[0,\pi/2]\).
When an endpoint value is used, it is taken from the continuous extension to \(I_+\).

\emph{Continuous cosine coefficients.}
For an even function $f\in L^1(I)$, define
\begin{align*}
    \mathcal C_0[f]
    &\coloneqq
    \frac{2}{\pi}\int_0^{\pi/2} f(\theta)\,d\theta,
    \mathcal C_n[f]
    &\coloneqq
    \frac{4}{\pi}\int_0^{\pi/2} f(\theta)\cos(2n\theta)\,d\theta,
    \qquad n\geq 1.
\end{align*}
If $f$ admits a cosine expansion on $I_+$, then
\begin{equation*}
    f(\theta)
    =
    \mathcal C_0[f]
    +
    \sum_{n=1}^{\infty}\mathcal C_n[f]\cos(2n\theta),
    \qquad \theta\in I_+.
\end{equation*}

\emph{DCT-I coefficients and approximation.}
Fix $N\in\bbN$ and define
\begin{equation*}
    \theta_q \coloneqq \frac{q\pi}{2N},
    \qquad
    q=0,1,\dots,N.
\end{equation*}
Define the endpoint weights by
\begin{equation*}
    \varpi_0=\varpi_N=\frac12,
    \qquad
    \varpi_q=1 \quad (1\leq q\leq N-1).
\end{equation*}
For a continuous function $f$ on $I_+$ satisfying $f(\pi/2)=0$, define
\begin{equation*}
    f_q \coloneqq f(\theta_q),
    \qquad
    q=0,1,\dots,N.
\end{equation*}
In particular, $f_N=0$. For $0\leq k\leq N-1$, define
\begin{equation}\label{eq:app-dct1-discrete-coeff}
    \mathcal D_{N,k}[f]
    \coloneqq
    \frac{\varsigma_k}{N}\sum_{q=0}^{N}\varpi_q f_q\cos(2k\theta_q),
    \qquad
    k=0,1,\dots,N-1,
\end{equation}
where
\begin{equation*}
    \varsigma_0\coloneqq 1,
    \qquad
    \varsigma_k\coloneqq 2 \quad (1\leq k\leq N-1).
\end{equation*}
We define the DCT-I approximation operator by
\[
    (\mathcal R_N f)(\theta)
    \coloneqq
    \sum_{k=0}^{N-1}
        \mathcal D_{N,k}[f]\cos(2k\theta),
    \qquad
    \theta\in I_+.
\]

The next identity expresses the discrete coefficients in terms of the continuous coefficients and the aliased tail.

\begin{lemma}[DCT-I aliasing identity]
    Assume that $f$ extends continuously to $I_+$, satisfies $f(\pi/2)=0$, and admits
    the cosine expansion
    \begin{equation}\label{eq:app-dct1-expansion-for-aliasing}
        f(\theta)
        =
        \mathcal C_0[f]
        +
        \sum_{n=1}^{\infty}\mathcal C_n[f]\cos(2n\theta),
        \qquad \theta\in I_+,
    \end{equation}
    with absolute convergence at each grid point \(\theta_q\). Then
    \begin{equation}\label{eq:app-dct1-aliasing-zero}
        \mathcal D_{N,0}[f]
        =
        \mathcal C_0[f]
        +
        \sum_{\ell\geq 1}\mathcal C_{2\ell N}[f],
    \end{equation}
    and for $1\leq k\leq N-1$,
    \begin{equation}\label{eq:app-dct1-aliasing-pos}
        \mathcal D_{N,k}[f]
        =
        \mathcal C_k[f]
        +
        \sum_{\ell\geq 1}
        \bigl(
            \mathcal C_{2\ell N-k}[f]
            +
            \mathcal C_{2\ell N+k}[f]
        \bigr).
    \end{equation}
\end{lemma}

\begin{proof}
    Substituting
    \eqref{eq:app-dct1-expansion-for-aliasing} into
    \eqref{eq:app-dct1-discrete-coeff} and using absolute convergence gives
    \[
        \mathcal D_{N,k}[f]
        =
        \sum_{n=0}^{\infty}
            \mathcal C_n[f]\,
            \frac{\varsigma_k}{N}
            \sum_{q=0}^{N}
                \varpi_q
                \cos(2n\theta_q)
                \cos(2k\theta_q).
    \]
    The DCT-I orthogonality relation is
    \begin{equation*}
        \frac{\varsigma_k}{N}
        \sum_{q=0}^{N}
            \varpi_q
            \cos(2n\theta_q)
            \cos(2k\theta_q)
        =
        \begin{cases}
            1,
            & n\equiv\pm k\pmod{2N},\\
            0,
            & \text{otherwise}.
        \end{cases}
    \end{equation*}
    For \(k=0\), the contributing indices are \(n=2\ell N\),
    \(\ell\geq0\), which gives
    \eqref{eq:app-dct1-aliasing-zero}. For \(1\leq k\leq N-1\), they are
    \(n=k\) and \(n=2\ell N\pm k\), \(\ell\geq1\), which gives
    \eqref{eq:app-dct1-aliasing-pos}.
\end{proof}

\begin{lemma}[DCT-I aliasing bound for an \(n^{-2}\) coefficient tail]\label{lem:app-alias-bound-n2}
    Assume that \(f\) extends continuously to \(I_+\), satisfies \(f(\pi/2)=0\), and
    admits the cosine expansion \eqref{eq:app-dct1-expansion-for-aliasing}, with
    absolute convergence at each grid point \(\theta_j\).
    Assume further that there exist constants \(C_f>0\) and an integer \(L_f\ge 0\)
    such that
    \begin{equation}\label{eq:app-tail-n2}
        |\mathcal C_n[f]|\le \frac{C_f}{(n-L_f)^2},
        \qquad n\ge L_f+1.
    \end{equation}
    If \(N\geq\max\{1,3L_f\}\), then
    \begin{equation}\label{eq:app-alias-bound-n2-zero}
        \bigl|\mathcal D_{N,0}[f]-\mathcal C_0[f]\bigr|
        \le
        \frac{\pi^2}{24}\,\frac{C_f}{(N-L_f)^2}.
    \end{equation}
    Moreover, for \(1\leq k\leq N-1\),
    \begin{equation}\label{eq:app-alias-bound-n2-pos}
        \bigl|\mathcal D_{N,k}[f]-\mathcal C_k[f]\bigr|
        \le
        \frac{\pi^2}{6}\,\frac{C_f}{(N-L_f)^2}.
    \end{equation}
    In particular, for all \(0\leq k\leq N-1\),
    \begin{equation}\label{eq:app-alias-bound-n2}
        \bigl|\mathcal D_{N,k}[f]-\mathcal C_k[f]\bigr|
        \le
        \frac{\pi^2}{6}\,\frac{C_f}{(N-L_f)^2}.
    \end{equation}
\end{lemma}

\begin{proof}
    We first treat the zero mode.  By \eqref{eq:app-dct1-aliasing-zero},
    \[
        \bigl|\mathcal D_{N,0}[f]-\mathcal C_0[f]\bigr|
        \le
        \sum_{\ell\ge1}
        |\mathcal C_{2\ell N}[f]|.
    \]
    Since \(N\geq\max\{1,3L_f\}\), the indices \(2\ell N\) satisfy
    \(2\ell N\ge L_f+1\). Hence \eqref{eq:app-tail-n2} gives
    \[
        \bigl|\mathcal D_{N,0}[f]-\mathcal C_0[f]\bigr|
        \le
        C_f\sum_{\ell\ge1}\frac{1}{(2\ell N-L_f)^2}.
    \]
    Since \(2\ell N-L_f\ge 2\ell(N-L_f)\), we obtain
    \[
        \bigl|\mathcal D_{N,0}[f]-\mathcal C_0[f]\bigr|
        \le
        \frac{C_f}{(N-L_f)^2}
        \sum_{\ell\ge1}\frac{1}{(2\ell)^2}
        =
        \frac{\pi^2}{24}\,\frac{C_f}{(N-L_f)^2}.
    \]
    This proves \eqref{eq:app-alias-bound-n2-zero}.

    Now let \(1\leq k\leq N-1\).  By \eqref{eq:app-dct1-aliasing-pos},
    \[
        \bigl|\mathcal D_{N,k}[f]-\mathcal C_k[f]\bigr|
        \le
        \sum_{\ell\ge 1}
        \bigl(
            |\mathcal C_{2\ell N-k}[f]|
            +
            |\mathcal C_{2\ell N+k}[f]|
        \bigr).
    \]
    For \(\ell\geq 1\), we have
    \[
        2\ell N-k\ge N+1\ge L_f+1,
        \qquad
        2\ell N+k\ge 2N\ge L_f+1,
    \]
    so \eqref{eq:app-tail-n2} applies to both indices.  Moreover,
    \[
        2\ell N-k-L_f\ge (2\ell-1)(N-L_f),
        \qquad
        2\ell N+k-L_f\ge 2\ell(N-L_f).
    \]
    Therefore
    \begin{align*}
        \bigl|\mathcal D_{N,k}[f]-\mathcal C_k[f]\bigr|
        &\le
        C_f\sum_{\ell\ge 1}
        \left(
            \frac{1}{(2\ell N-k-L_f)^2}
            +
            \frac{1}{(2\ell N+k-L_f)^2}
        \right) \\
        &\le
        \frac{C_f}{(N-L_f)^2}
        \left(
            \sum_{\ell\ge 1}\frac{1}{(2\ell-1)^2}
            +
            \sum_{\ell\ge 1}\frac{1}{(2\ell)^2}
        \right) \\
        &=
        \frac{\pi^2}{6}\,\frac{C_f}{(N-L_f)^2}.
    \end{align*}
    This proves \eqref{eq:app-alias-bound-n2-pos}.  The uniform estimate
    \eqref{eq:app-alias-bound-n2} follows from
    \(\pi^2/24\le \pi^2/6\).
\end{proof}

\subsection{Matrix and projection estimates}
We record the Kato--Temple, Collatz--Wielandt, and rank-one projection estimates used in Part~\ref{part:interval-verification}.

\begin{lemma}[Kato--Temple upper bound]\label{lem:app-kato-temple-upper}
    Let \(N\geq3\), let \(\mathbf A\in\bbR^{N\times N}\) be symmetric, and write its eigenvalues in
    decreasing order
    \[
        \lambda_1(\mathbf A)\geq \lambda_2(\mathbf A)\geq \lambda_3(\mathbf A)\geq\cdots\geq\lambda_N(\mathbf A).
    \]
    Let \(x\in\bbR^N\setminus\{0\}\), and set
    \[
        \rho
        \coloneqq
        \frac{x^\top \mathbf Ax}{x^\top x},
        \qquad
        \eta
        \coloneqq
        \frac{\norm{\mathbf Ax-\rho x}_2}{\norm{x}_2}.
    \]
    Assume that, for some \(\beta\in\bbR\),
    \(
        \lambda_3(\mathbf A)\leq \beta<\rho .
    \)
    Then
    \begin{equation}\label{eq:app-kato-temple-upper}
        \lambda_2(\mathbf A)
        \leq
        \rho
        +
        \frac{\eta^2}{\rho-\beta}.
    \end{equation}
\end{lemma}

\begin{proof}
    By replacing \(x\) with \(x/\norm{x}_2\), we may assume
    \(\norm{x}_2=1\). If
    \(\lambda_2(\mathbf A)\leq\rho\), then
    \eqref{eq:app-kato-temple-upper} is immediate. Assume
    \(\lambda_2(\mathbf A)>\rho\).
    If
    \[
        \eta^2
        \geq
        \bigl(
            \lambda_1(\mathbf A)-\rho
        \bigr)
        (\rho-\beta),
    \]
    then
    \[
        \lambda_2(\mathbf A)-\rho
        \leq
        \lambda_1(\mathbf A)-\rho
        \leq
        \frac{\eta^2}{\rho-\beta},
    \]
    which gives \eqref{eq:app-kato-temple-upper}. It remains to consider
    \begin{equation*}
        \eta^2
        <
        \bigl(
            \lambda_1(\mathbf A)-\rho
        \bigr)
        (\rho-\beta).
    \end{equation*}
    Choose
    \(B=\lambda_1(\mathbf A)\) if
    \(\lambda_1(\mathbf A)>\lambda_2(\mathbf A)\), and choose any
    \(B>\lambda_1(\mathbf A)\) otherwise. Then
    \[
        \sigma_p(\mathbf A)\cap(\beta,B)
        =
        \{
            \lambda_2(\mathbf A)
        \}
    \]
    as a set of spectral points, and
    \[
        \eta^2
        <
        (B-\rho)(\rho-\beta).
    \]
    The Kato--Temple inequality
    \cite[Theorem~2]{HarrellTemple}, with endpoints \(\beta\) and \(B\),
    gives \eqref{eq:app-kato-temple-upper}.
\end{proof}

\begin{lemma}[Weighted Schur bound and Collatz--Wielandt]\label{lem:app-weighted-schur-CW}
    Let \(\mathbf A\in\bbR^{n\times n}\) be symmetric.
    For \(\mathbf d=(d_1,\dots,d_n)\in(0,\infty)^n\), set
    \begin{equation*}
        \mathrm{UB}_{\mathbf A}(\mathbf d)
        \coloneqq
        \max_{1\leq i\leq n}
        \frac{(|\mathbf A|_{\mathrm{entry}}\mathbf d)_i}{d_i}.
    \end{equation*}
    Then, for every \(\mathbf d\in(0,\infty)^n\),
    \begin{equation*}
        \norm{\mathbf A}_2
        \leq
        \lambda_{\max}(|\mathbf A|_{\mathrm{entry}})
        \leq
        \mathrm{UB}_{\mathbf A}(\mathbf d).
    \end{equation*}
\end{lemma}

\begin{proof}
    For every \(\mathbf h\in\bbR^n\), the Rayleigh quotient characterization gives
    \[
        |\mathbf h^\top\mathbf A\mathbf h|
        \leq
        |\mathbf h|^\top|\mathbf A|_{\mathrm{entry}}\,|\mathbf h|
        \leq
        \lambda_{\max}(|\mathbf A|_{\mathrm{entry}})\norm{\mathbf h}_2^2 .
    \]
    Hence \(\norm{\mathbf A}_2\leq\lambda_{\max}(|\mathbf A|_{\mathrm{entry}})\).  Since
    \(|\mathbf A|_{\mathrm{entry}}\) is symmetric and nonnegative, its Perron root is
    \(\lambda_{\max}(|\mathbf A|_{\mathrm{entry}})\).  The Collatz--Wielandt upper bound, applied to
    \(|\mathbf A|_{\mathrm{entry}}\) with the positive vector \(\mathbf d\), gives
    \[
        \lambda_{\max}(|\mathbf A|_{\mathrm{entry}})
        \leq
        \max_{1\leq i\leq n}
        \frac{(|\mathbf A|_{\mathrm{entry}}\mathbf d)_i}{d_i}
        =
        \mathrm{UB}_{\mathbf A}(\mathbf d);
    \]
    see \cite[Theorem~8.1.26]{HornCW}.
\end{proof}

In Part~\ref{part:interval-verification}, the positive vector \(\mathbf d\) is fixed before the
interval evaluation, and \(\mathrm{UB}_{\mathbf A}(\mathbf d)\) is evaluated directly.

\begin{lemma}[Difference of rank-one orthogonal projections]
\label{lem:app-rank-one-projection-difference}
    Let \(f_1,f_2\in L^2\) be real-valued and normalized by
    \[
        \norm{f_1}_{L^2}=\norm{f_2}_{L^2}=1.
    \]
    For \(i=1,2\), let
    \[
        P_{f_i}h\coloneqq (h,f_i)_{L^2}f_i .
    \]
    Then
    \begin{equation}\label{eq:app-rank-one-proj-diff-lipschitz}
        \norm{P_{f_1}-P_{f_2}}_{L^2\to L^2}
        \leq
        \norm{f_1-f_2}_{L^2}.
    \end{equation}
    The same estimate holds for \(P_{f_i^\perp}\coloneqq I-P_{f_i}\).
\end{lemma}

\begin{proof}
    Set
    \[
        \chi\coloneqq (f_1,f_2)_{L^2}.
    \]
    If \(|\chi|=1\), then \(P_{f_1}=P_{f_2}\). Otherwise, set
    \[
        \sigma\coloneqq \sqrt{1-\chi^2},
        \qquad
        v_1\coloneqq f_1,
        \qquad
        v_2\coloneqq \frac{f_2-\chi f_1}{\sigma}.
    \]
    In the orthonormal basis \((v_1,v_2)\) of \(\operatorname{span}\{f_1,f_2\}\),
    \[
        P_{f_1}-P_{f_2}
        =
        \begin{pmatrix}
            \sigma^2 & -\chi\sigma\\
            -\chi\sigma & -\sigma^2
        \end{pmatrix}.
    \]
    This matrix has eigenvalues \(\pm\sigma\), and the operator vanishes on
    \(\operatorname{span}\{f_1,f_2\}^{\perp}\). Hence
    \[
        \norm{P_{f_1}-P_{f_2}}_{L^2\to L^2}
        =
        \sigma
        =
        \sqrt{1-\chi^2},
        \leq
        \sqrt{2-2\chi}
        =
        \norm{f_1-f_2}_{L^2},
    \]
    so \eqref{eq:app-rank-one-proj-diff-lipschitz} follows. Since \(P_{f_1^\perp}-P_{f_2^\perp}=-(P_{f_1}-P_{f_2})\), the same estimate
    holds for the orthogonal-complement projections.
\end{proof}

\section{Estimates for the Schur complement reduction}

Recall the operators and finite-dimensional approximations from
Subsection~\ref{subsec:g-side-Schur-complement}.

\begin{definition}[Constants for the finite Schur comparison]
    Set
    \begin{equation}\label{eq:definition of Delta_T^{(L)}}
        \begin{aligned}
            \delta_T^{(L)}
            \coloneqq{}&
            \frac{\mathsf C_{\log\sin}^{L^2}}{\pi}
            \norm{
                a_{\mathfrak g}
                -
                a_{\mathfrak g}^{(L)}
            }_{L^2(I)}
            \left(
                \norm{a_{\mathfrak g}}_{L^\infty(I)}
                +
                \norm{a_{\mathfrak g}^{(L)}}_{L^\infty(I)}
            \right)
            \\
            &+
            \frac2\pi
            \left[
                \norm{
                    a_{\mathfrak g,\log}
                    -
                    a_{\mathfrak g,\log}^{(L)}
                }_{L^2(I)}
                \norm{a_{\mathfrak g}}_{L^2(I)}
            \right.
            \\
            &\hspace{34mm}\left.
                +
                \norm{a_{\mathfrak g,\log}^{(L)}}_{L^2(I)}
                \norm{
                    a_{\mathfrak g}
                    -
                    a_{\mathfrak g}^{(L)}
                }_{L^2(I)}
            \right],
        \end{aligned}
    \end{equation}
    and
    \begin{equation}\label{eq:definition of Omega_T^{(L)}}
        \mathsf C_{T,\mathrm{op}}^{(L)}
        \coloneqq
        (\log2)
        \norm{a_{\mathfrak g}^{(L)}}_{L^\infty(I)}^2
        +
        \frac2\pi
        \norm{a_{\mathfrak g,\log}^{(L)}}_{L^2(I)}
        \norm{a_{\mathfrak g}^{(L)}}_{L^2(I)}.
    \end{equation}

    Define
    \begin{align}
        \mathcal R_{\mathrm{rhs},1}^{(L)}
        &\coloneqq
        \frac{
            \norm{a_{\mathfrak g}^{(L)}}_{L^2(I)}^2
        }{\pi}
        \norm{
            \Pi_{>3L}
            \left(
                \widetilde\Phi_{\mathfrak g}^{(L)}
                \log|\cos|
            \right)
        }_{L^2(I)},
        \notag\\
        \mathcal R_{\mathrm{rhs},2}^{(L)}
        &\coloneqq
        \frac{\mathsf C_{\log\sin}^{L^2}}{\pi}
        \norm{
            a_{\mathfrak g}
            -
            a_{\mathfrak g}^{(L)}
        }_{L^2(I)}
        \left(
            \norm{a_{\mathfrak g}}_{L^\infty(I)}
            +
            \norm{a_{\mathfrak g}^{(L)}}_{L^\infty(I)}
        \right),
        \notag\\
        \mathcal R_{\mathrm{rhs},3}^{(L)}
        &\coloneqq
        2(\log2)
        \norm{a_{\mathfrak g}^{(L)}}_{L^\infty(I)}^2
        \norm{
            \widetilde\Phi_{\mathfrak g}
            -
            \widetilde\Phi_{\mathfrak g}^{(L)}
        }_{L^2(I)},
        \notag\\
        \mathcal R_{\mathrm{rhs},4}^{(L)}
        &\coloneqq
        \frac{
            \norm{a_{\mathfrak g}}_{L^2(I)}^2
        }{\pi}
        \norm{
            \left(
                \widetilde\Phi_{\mathfrak g}
                -
                \widetilde\Phi_{\mathfrak g}^{(L)}
            \right)
            \log|\cos|
        }_{L^2(I)},
        \notag\\
        \mathcal R_{\mathrm{rhs},5}^{(L)}
        &\coloneqq
        \frac{
            \norm{a_{\mathfrak g}}_{L^2(I)}^2
        }{\pi}
        \norm{
            \widetilde\Phi_{\mathfrak g}
            -
            \widetilde\Phi_{\mathfrak g}^{(L)}
        }_{L^2(I)}
        \norm{
            \widetilde\Phi_{\mathfrak g}^{(L)}
            \log|\cos|
        }_{L^2(I)},
        \notag\\
        \mathcal R_{\mathrm{rhs},6}^{(L)}
        &\coloneqq
        \frac1\pi
        \norm{
            a_{\mathfrak g}
            -
            a_{\mathfrak g}^{(L)}
        }_{L^2(I)}
        \left(
            \norm{a_{\mathfrak g}}_{L^2(I)}
            +
            \norm{a_{\mathfrak g}^{(L)}}_{L^2(I)}
        \right)
        \norm{
            \widetilde\Phi_{\mathfrak g}^{(L)}
            \log|\cos|
        }_{L^2(I)},
        \notag\\
        \delta_{\mathrm{rhs}}^{(L)}
        &\coloneqq
        \sum_{\ell=1}^{6}
            \mathcal R_{\mathrm{rhs},\ell}^{(L)}.
        \label{eq:definition of Delta_{rhs}^{(L)}}
    \end{align}

    Using \eqref{eq:part4-eps-model-BS-def}, set
    \begin{equation}\label{eq:definition of Delta_{red}^{(L)}}
        \Delta_{\mathrm{proj}}^{(L)}
        \coloneqq
        \delta_{a_{\mathfrak g}\to a_{\mathfrak g}^{(L)}}^{\mathrm{BS}}
        +
        \frac{
            \norm{a_{\mathfrak g}^{(L)}}_{L^\infty(I)}^2
        }{
            2(3L+1)
        }
        +
        \frac{
            \norm{a_{\mathfrak g}^{(L)}}_{L^\infty(I)}^2
        }{
            2L+1
        }.
    \end{equation}
\end{definition}

\begin{lemma}[Estimates for the finite Schur comparison]
\label{lem:app-finite-schur-comparison-estimates}
    The following estimates hold.
    \begin{enumerate}[label=\textup{(\roman*)}]
        \item
        \begin{equation}\label{eq:app-Tg-minus-TgL-full-bound}
            \norm{
                \widetilde T_{\mathfrak g}
                -
                \widetilde T_{\mathfrak g}^{(L)}
            }_{L^2_{\mathrm e}(I)\to L^2_{\mathrm e}(I)}
            \leq
            \delta_T^{(L)},
        \end{equation}
        and
        \begin{equation}\label{eq:app-TgL-full-norm-bound}
            \norm{
                \widetilde T_{\mathfrak g}^{(L)}
            }_{L^2_{\mathrm e}(I)\to L^2_{\mathrm e}(I)}
            \leq
            \mathsf C_{T,\mathrm{op}}^{(L)}.
        \end{equation}

        \item
        \begin{equation}\label{eq:app-Schur-rhs-comparison-bound}
            \norm{
                P_{\widetilde\Phi_{\mathfrak g}^{\perp}}
                \widetilde T_{\mathfrak g}
                \widetilde\Phi_{\mathfrak g}
                -
                \Pi_{\leq3L}
                P_{(\widetilde\Phi_{\mathfrak g}^{(L)})^\perp}
                \widetilde T_{\mathfrak g}^{(L)}
                \widetilde\Phi_{\mathfrak g}^{(L)}
            }_{L^2(I)}
            \leq
            \delta_{\mathrm{rhs}}^{(L)}.
        \end{equation}

        \item
        \begin{equation}\label{eq:app-Schur-red-comparison-bound}
            \norm{
                \widetilde T_{\mathfrak g,\mathrm{proj}}
                -
                \widetilde T_{\mathfrak g,3L}^{(L)}
            }_{L^2_{\mathrm e}(I)\to L^2_{\mathrm e}(I)}
            \leq
            \Delta_{\mathrm{proj}}^{(L)}.
        \end{equation}
    \end{enumerate}
\end{lemma}

\begin{proof}
    We first prove \eqref{eq:app-Tg-minus-TgL-full-bound}.  By the
    decomposition
    \[
        \widetilde T_{\mathfrak g}
        =
        \widetilde T_{\mathfrak g,\mathrm{sing}}
        +
        \widetilde T_{\mathfrak g,\mathrm{fr}},
        \qquad
        \widetilde T_{\mathfrak g}^{(L)}
        =
        \widetilde T_{\mathfrak g,\mathrm{sing}}^{(L)}
        +
        \widetilde T_{\mathfrak g,\mathrm{fr}}^{(L)},
    \]
    the singular part is estimated by
    \eqref{eq:part4-singular-model-difference-bound}:
    \[
        \begin{aligned}
            \norm{
                \widetilde T_{\mathfrak g,\mathrm{sing}}
                -
                \widetilde T_{\mathfrak g,\mathrm{sing}}^{(L)}
            }_{L^2_{\mathrm e}\to L^2_{\mathrm e}}
            &\leq
            \frac{\mathsf C_{\log\sin}^{L^2}}{\pi}
            \norm{
                a_{\mathfrak g}
                -
                a_{\mathfrak g}^{(L)}
            }_{L^2(I)}
            \left(
                \norm{a_{\mathfrak g}}_{L^\infty(I)}
                +
                \norm{a_{\mathfrak g}^{(L)}}_{L^\infty(I)}
            \right).
        \end{aligned}
    \]
    For the finite-rank part, using
    \(\norm{f_1\otimes f_2}_{L^2\to L^2}
    =
    \norm{f_1}_{L^2}\norm{f_2}_{L^2}\), we get
    \[
    \begin{aligned}
        &\norm{
            \widetilde T_{\mathfrak g,\mathrm{fr}}
            -
            \widetilde T_{\mathfrak g,\mathrm{fr}}^{(L)}
        }_{L^2_{\mathrm e}\to L^2_{\mathrm e}}
        \\
        &\qquad\leq
        \frac2\pi
        \left[
            \norm{
                a_{\mathfrak g,\log}
                -
                a_{\mathfrak g,\log}^{(L)}
            }_{L^2(I)}
            \norm{a_{\mathfrak g}}_{L^2(I)}
            +
            \norm{a_{\mathfrak g,\log}^{(L)}}_{L^2(I)}
            \norm{
                a_{\mathfrak g}
                -
                a_{\mathfrak g}^{(L)}
            }_{L^2(I)}
        \right].
    \end{aligned}
    \]
    The finite-rank estimate above, together with
    \eqref{eq:part4-singular-model-difference-bound} and
    \eqref{eq:definition of Delta_T^{(L)}}, gives
    \eqref{eq:app-Tg-minus-TgL-full-bound}. The same decomposition gives
    \[
        \norm{
            \widetilde T_{\mathfrak g,\mathrm{sing}}^{(L)}
        }_{L^2_{\mathrm e}\to L^2_{\mathrm e}}
        \leq
        \frac1\pi
        \norm{a_{\mathfrak g}^{(L)}}_{L^\infty(I)}^2
        \norm{\mathcal K_{\sin}}_{L^2_{\mathrm e}\to L^2_{\mathrm e}}
        =
        (\log2)
        \norm{a_{\mathfrak g}^{(L)}}_{L^\infty(I)}^2,
    \]
    while
    \[
        \norm{
            \widetilde T_{\mathfrak g,\mathrm{fr}}^{(L)}
        }_{L^2_{\mathrm e}\to L^2_{\mathrm e}}
        \leq
        \frac2\pi
        \norm{a_{\mathfrak g,\log}^{(L)}}_{L^2(I)}
        \norm{a_{\mathfrak g}^{(L)}}_{L^2(I)}.
    \]
    This proves \eqref{eq:app-TgL-full-norm-bound}.

    We next prove \eqref{eq:app-Schur-rhs-comparison-bound}.  Since
    \(a_{\mathfrak g}\) is proportional to
    \(\widetilde\Phi_{\mathfrak g}\), the finite-rank part gives
    \[
        P_{\widetilde\Phi_{\mathfrak g}^{\perp}}
        \widetilde T_{\mathfrak g,\mathrm{fr}}
        \widetilde\Phi_{\mathfrak g}
        =
        \frac{\norm{a_{\mathfrak g}}_{L^2(I)}^2}{\pi}
        P_{\widetilde\Phi_{\mathfrak g}^{\perp}}
        \left(
            \widetilde\Phi_{\mathfrak g}\log|\cos|
        \right),
    \]
    and similarly,
    \[
        P_{(\widetilde\Phi_{\mathfrak g}^{(L)})^\perp}
        \widetilde T_{\mathfrak g,\mathrm{fr}}^{(L)}
        \widetilde\Phi_{\mathfrak g}^{(L)}
        =
        \frac{\norm{a_{\mathfrak g}^{(L)}}_{L^2(I)}^2}{\pi}
        P_{(\widetilde\Phi_{\mathfrak g}^{(L)})^\perp}
        \left(
            \widetilde\Phi_{\mathfrak g}^{(L)}\log|\cos|
        \right).
    \]
    Moreover, since \(a_{\mathfrak g}^{(L)}\) and
    \(\widetilde\Phi_{\mathfrak g}^{(L)}\) belong to \(\mathcal S_L\), the singular
    part
    \[
        P_{(\widetilde\Phi_{\mathfrak g}^{(L)})^\perp}
        \widetilde T_{\mathfrak g,\mathrm{sing}}^{(L)}
        \widetilde\Phi_{\mathfrak g}^{(L)}
    \]
    belongs to \(\mathcal S_{3L}\).  Hence the projection
    \(\Pi_{\leq 3L}\) only creates a tail from the logarithmic finite-rank
    contribution, and this tail is bounded by the first term in
    the sum defining $\delta_{\mathrm{rhs}}^{(L)}$ in \eqref{eq:definition of Delta_{rhs}^{(L)}}.

    It remains to estimate the unprojected difference. For the singular
    contribution, the triangle inequality, Lemma~\ref{lem:app-rank-one-projection-difference},
    and \eqref{eq:part4-singular-model-difference-bound} give
    \[
    \begin{aligned}
        \norm{
            P_{\widetilde\Phi_{\mathfrak g}^{\perp}}
            \widetilde T_{\mathfrak g,\mathrm{sing}}
            \widetilde\Phi_{\mathfrak g}
            &-
            P_{(\widetilde\Phi_{\mathfrak g}^{(L)})^\perp}
            \widetilde T_{\mathfrak g,\mathrm{sing}}^{(L)}
            \widetilde\Phi_{\mathfrak g}^{(L)}
        }_{L^2(I)}
        \\
        &\leq
        \frac{\mathsf C_{\log\sin}^{L^2}}{\pi}
        \norm{
            a_{\mathfrak g}
            -
            a_{\mathfrak g}^{(L)}
        }_{L^2(I)}
        \left(
            \norm{a_{\mathfrak g}}_{L^\infty(I)}
            +
            \norm{a_{\mathfrak g}^{(L)}}_{L^\infty(I)}
        \right)\\
        &\quad+
        2(\log2)
        \norm{a_{\mathfrak g}^{(L)}}_{L^\infty(I)}^2
        \norm{
            \widetilde\Phi_{\mathfrak g}
            -
            \widetilde\Phi_{\mathfrak g}^{(L)}
        }_{L^2(I)}.
    \end{aligned}
    \]
    For the logarithmic finite-rank contribution, Lemma~\ref{lem:app-rank-one-projection-difference}
    gives
    \[
    \begin{aligned}
        &\norm{
            P_{\widetilde\Phi_{\mathfrak g}^{\perp}}
            \left(
                \widetilde\Phi_{\mathfrak g}\log|\cos|
            \right)
            -
            P_{(\widetilde\Phi_{\mathfrak g}^{(L)})^\perp}
            \left(
                \widetilde\Phi_{\mathfrak g}^{(L)}\log|\cos|
            \right)
        }_{L^2(I)}
        \\
        &\qquad\leq
        \norm{
            \left(
                \widetilde\Phi_{\mathfrak g}
                -
                \widetilde\Phi_{\mathfrak g}^{(L)}
            \right)
            \log|\cos|
        }_{L^2(I)}
        +
        \norm{
            \widetilde\Phi_{\mathfrak g}
            -
            \widetilde\Phi_{\mathfrak g}^{(L)}
        }_{L^2(I)}
        \norm{
            \widetilde\Phi_{\mathfrak g}^{(L)}
            \log|\cos|
        }_{L^2(I)}.
    \end{aligned}
    \]
    Finally,
    \[
        \left|
            \norm{a_{\mathfrak g}}_{L^2(I)}^2
            -
            \norm{a_{\mathfrak g}^{(L)}}_{L^2(I)}^2
        \right|
        \leq
        \norm{
            a_{\mathfrak g}
            -
            a_{\mathfrak g}^{(L)}
        }_{L^2(I)}
        \left(
            \norm{a_{\mathfrak g}}_{L^2(I)}
            +
            \norm{a_{\mathfrak g}^{(L)}}_{L^2(I)}
        \right).
    \]
    Substituting the projection-tail, singular, finite-rank, and coefficient
    estimates above into \eqref{eq:definition of Delta_{rhs}^{(L)}}, and using
    \eqref{eq:app-rank-one-proj-diff-lipschitz}, gives
    \eqref{eq:app-Schur-rhs-comparison-bound}.

    It remains to prove \eqref{eq:app-Schur-red-comparison-bound}.  We split
    \[
        \widetilde T_{\mathfrak g,\mathrm{proj}}
        -
        \widetilde T_{\mathfrak g,3L}^{(L)}
        =
        \left(
            \widetilde T_{\mathfrak g,\mathrm{proj}}
            -
            \widetilde T_{\mathfrak g,\mathrm{proj}}^{(L)}
        \right)
        +
        \left(
            \widetilde T_{\mathfrak g,\mathrm{proj}}^{(L)}
            -
            \widetilde T_{\mathfrak g,3L}^{(L)}
        \right).
    \]
    By \eqref{eq:BS-model-replacement-bound},
    \[
        \norm{
            \widetilde T_{\mathfrak g,\mathrm{proj}}
            -
            \widetilde T_{\mathfrak g,\mathrm{proj}}^{(L)}
        }_{L^2_{\mathrm e}(I)\to L^2_{\mathrm e}(I)}
        \leq
        \delta_{a_{\mathfrak g}\to a_{\mathfrak g}^{(L)}}^{\mathrm{BS}}.
    \]
    For the second term, using the tail bounds for \(\mathcal K_{\sin}\) gives
    \[
    \begin{aligned}
        \norm{
            \widetilde T_{\mathfrak g,\mathrm{proj}}^{(L)}
            -
            \widetilde T_{\mathfrak g,3L}^{(L)}
        }_{L^2_{\mathrm e}\to L^2_{\mathrm e}}
        &\leq
        \frac1\pi
        \norm{a_{\mathfrak g}^{(L)}}_{L^\infty(I)}^2
        \norm{
            \Pi_{>3L}\mathcal K_{\sin}\Pi_{>3L}
        }_{L^2_{\mathrm e}\to L^2_{\mathrm e}}
        \\
        &\quad+
        \frac2\pi
        \norm{a_{\mathfrak g}^{(L)}}_{L^\infty(I)}^2
        \norm{
            \Pi_{>2L}\mathcal K_{\sin}
        }_{L^2_{\mathrm e}\to L^2_{\mathrm e}}
        \\
        &\leq
        \frac{
            \norm{a_{\mathfrak g}^{(L)}}_{L^\infty(I)}^2
        }{
            2(3L+1)
        }
        +
        \frac{
            \norm{a_{\mathfrak g}^{(L)}}_{L^\infty(I)}^2
        }{
            2L+1
        }.
    \end{aligned}
    \]
    The preceding tail estimate,
    \eqref{eq:BS-model-replacement-bound}, and
    \eqref{eq:definition of Delta_{red}^{(L)}} give
    \eqref{eq:app-Schur-red-comparison-bound}.
\end{proof}

\section{DCT-I estimates for the limiting matrix entries}

This section records the cosine coefficient and DCT-I approximation
estimates used for the limiting matrix entries.

\subsection{Cosine coefficient estimates}

Recall \(a_{\mathfrak g}\) and \(a_{\mathfrak g,\log}\) from
\eqref{eq:definition of a_g, tidle T_g} and
\eqref{eq:definition of a_g,log}. Set
\begin{equation}\label{eq:definition-of-C-tail}
    \mathsf C_{a_{\mathfrak g},\mathrm{coef\text{-}tail}}
    \coloneqq
    \frac{
        |a_{\mathfrak g}'(0^+)|
        +
        |a_{\mathfrak g}'(\pi/2^-)|
        +
        \operatorname{TV}
        \left(
            a_{\mathfrak g}';
            \left(
                0,\frac{\pi}{2}
            \right)
        \right)
    }{\pi}.
\end{equation}
Since
\(a_{\mathfrak g}(0)=a_{\mathfrak g}(\pi/2)=0\), integration by parts
gives
\begin{equation}\label{eq:ag-cosine-tail-bound}
    |\mathcal C_k[a_{\mathfrak g}]|
    \leq
    \frac{
        \mathsf C_{a_{\mathfrak g},\mathrm{coef\text{-}tail}}
    }{
        k^2
    },
    \qquad
    k\geq1.
\end{equation}

\begin{definition}[Constants for the logarithmic coefficient estimates]
    \begin{equation}\label{eq:app-g-side-log-estimate-constants}
        \begin{aligned}
            \mathsf C_{\mathrm{coef},\mathrm{end}}
            &\coloneqq
            \frac{2}{\pi}
            \norm{a_{\mathfrak g}'}_{L^\infty(I)},
            \\
            \mathsf C_{\mathrm{coef},\log}
            &\coloneqq
            \frac{4}{\pi}
            \left[
                \left(
                    \frac74+\frac{\pi^2}{16}
                \right)
                \norm{a_{\mathfrak g}'}_{L^\infty(I)}
                +
                \frac14
                \int_0^{\pi/2}
                    |a_{\mathfrak g}''(\theta)|
                \,d\theta
            \right],
            \\
            \mathsf C_{\mathrm{coef},0}
            &\coloneqq
            \mathsf C_{\mathrm{coef},\log}
            \log\frac{\pi}{2}
            +
            \mathsf C_{\mathrm{coef},\mathrm{end}},
            \qquad
            \Sigma_{\log,2}
            \coloneqq
            \sum_{r=1}^{\infty}
                \frac{\log r}{r^2},
            \\
            \mathsf C_{\mathrm{alias},\log}
            &\coloneqq
            \frac{\pi^2}{6}
            \mathsf C_{\mathrm{coef},\log},
            \\
            \mathsf C_{\mathrm{alias},0}
            &\coloneqq
            \mathsf C_{\mathrm{coef},\log}
            \left(
                \frac{\pi^2}{6}\log\frac{\pi}{2}
                +
                \frac{5\pi^2}{24}\log2
                +
                \frac54\Sigma_{\log,2}
            \right)
            +
            \frac{\pi^2}{6}
            \mathsf C_{\mathrm{coef},\mathrm{end}}.
        \end{aligned}
    \end{equation}
\end{definition}

\begin{lemma}[Log-weighted coefficients and DCT-I aliasing]\label{lem:app-ag-logcos-coeff-alias}
    With the constants defined in
    \eqref{eq:app-g-side-log-estimate-constants}, for every \(k\geq2\),
    \begin{equation}\label{eq:app-ag-logcos-coeff-bound}
        |\mathcal C_k[a_{\mathfrak g,\log}]|
        \leq
        \frac{\mathsf C_{\mathrm{coef},\log}\log k+\mathsf C_{\mathrm{coef},0}}{k^2}.
    \end{equation}
    The function \(a_{\mathfrak g,\log}\) is extended continuously to \(\pi/2\) by setting
    \(a_{\mathfrak g,\log}(\pi/2)=0\). For \(N\geq2\) and \(0\leq k\leq N-1\),
    \begin{equation}\label{eq:app-ag-logcos-dct-alias-bound}
        \left|
            \mathcal D_{N,k}[a_{\mathfrak g,\log}]
            -
            \mathcal C_k[a_{\mathfrak g,\log}]
        \right|
        \leq
        \frac{
            \mathsf C_{\mathrm{alias},\log}\log N+\mathsf C_{\mathrm{alias},0}
        }{
            N^2
        }.
    \end{equation}
\end{lemma}

\begin{proof}
    Set
    \[
        f_1(\tau)
        \coloneqq
        a_{\mathfrak g}
        \left(
            \frac{\pi}{2}-\tau
        \right),
        \qquad
        0\leq\tau\leq\frac{\pi}{2}.
    \]
    Since \(a_{\mathfrak g}(\pi/2)=0\),
    \[
        f_1(0)=0,
        \qquad
        |f_1(\tau)|
        \leq
        \norm{a_{\mathfrak g}'}_{L^\infty(I)}\tau.
    \]
    For \(k\geq2\), the change of variables
    \(\tau=\pi/2-\theta\) gives
    \[
        \mathcal C_k[a_{\mathfrak g,\log}]
        =
        \frac4\pi(-1)^k
        \int_0^{\pi/2}
            f_1(\tau)
            \log(\sin\tau)
            \cos(2k\tau)
        \,d\tau.
    \]
    Set
    \[
        \tau_0^{(k)}
        \coloneqq
        \frac1k,
        \qquad
        L_k
        \coloneqq
        \log k+\log\frac{\pi}{2},
    \]
    and, on
    \([\tau_0^{(k)},\pi/2]\), set
    \[
        f_2(\tau)
        \coloneqq
        f_1(\tau)\log(\sin\tau).
    \]
    Since
    \(\sin\tau\geq2\tau/\pi\) and
    \(\cot\tau\leq1/\tau\) on this interval,
    \[
        |f_2(\tau_0^{(k)})|
        \leq
        \frac{
            \norm{a_{\mathfrak g}'}_{L^\infty(I)}
        }{k}
        L_k,
        \qquad
        |f_2'(\tau_0^{(k)})|
        \leq
        \norm{a_{\mathfrak g}'}_{L^\infty(I)}
        (L_k+1).
    \]
    In the sense of finite measures on
    \([\tau_0^{(k)},\pi/2]\),
    \[
        df_2'
        =
        \log(\sin\tau)\,df_1'
        +
        2\cot\tau\,f_1'(\tau)\,d\tau
        -
        \csc^2\tau\,f_1(\tau)\,d\tau.
    \]
    Therefore,
    \[
        \begin{aligned}
            \operatorname{TV}_{[\tau_0^{(k)},\pi/2]}(f_2')
            &\leq
            \int_{\tau_0^{(k)}}^{\pi/2}
                |\log(\sin\tau)|
            \,d|f_1'|
            +
            2
            \int_{\tau_0^{(k)}}^{\pi/2}
                \cot\tau\,|f_1'(\tau)|
            \,d\tau
            +
            \int_{\tau_0^{(k)}}^{\pi/2}
                \csc^2\tau\,|f_1(\tau)|
            \,d\tau
            \\
            &\leq
            L_k
            \int_0^{\pi/2}
                |a_{\mathfrak g}''(\theta)|
            \,d\theta
            +
            2
            \norm{a_{\mathfrak g}'}_{L^\infty(I)}
            \int_{\tau_0^{(k)}}^{\pi/2}
                \frac{d\tau}{\tau}
            +
            \frac{\pi^2}{4}
            \norm{a_{\mathfrak g}'}_{L^\infty(I)}
            \int_{\tau_0^{(k)}}^{\pi/2}
                \frac{d\tau}{\tau}
            \\
            &=
            \left[
                \int_0^{\pi/2}
                    |a_{\mathfrak g}''(\theta)|
                \,d\theta
                +
                \left(
                    2+\frac{\pi^2}{4}
                \right)
                \norm{a_{\mathfrak g}'}_{L^\infty(I)}
            \right]
            L_k.
        \end{aligned}
    \]
    Since
    \(f_2(\pi/2)=f_2'(\pi/2)=0\), two integrations by parts give
    \[
        \begin{aligned}
            \left|
                \int_{\tau_0^{(k)}}^{\pi/2}
                    f_2(\tau)\cos(2k\tau)
                \,d\tau
            \right|
            &\leq
            \frac{|f_2(\tau_0^{(k)})|}{2k}
            +
            \frac{
                |f_2'(\tau_0^{(k)})|
                +
                \operatorname{TV}_{[\tau_0^{(k)},\pi/2]}(f_2')
            }{
                4k^2
            }
            \\
            &\leq
            \left[
                \left(
                    \frac54+\frac{\pi^2}{16}
                \right)
                \norm{a_{\mathfrak g}'}_{L^\infty(I)}
                +
                \frac14
                \int_0^{\pi/2}
                    |a_{\mathfrak g}''(\theta)|
                \,d\theta
            \right]
            \frac{L_k}{k^2}
            \\
            &\quad+
            \frac{
                \norm{a_{\mathfrak g}'}_{L^\infty(I)}
            }{
                4k^2
            }.
        \end{aligned}
    \]
    On the endpoint interval,
    \[
        \begin{aligned}
            \left|
                \int_0^{\tau_0^{(k)}}
                    f_1(\tau)
                    \log(\sin\tau)
                    \cos(2k\tau)
                \,d\tau
            \right|
            &\leq
            \norm{a_{\mathfrak g}'}_{L^\infty(I)}
            \int_0^{1/k}
                \tau
                \left(
                    |\log\tau|
                    +
                    \log\frac{\pi}{2}
                \right)
            \,d\tau
            \\
            &=
            \frac{
                \norm{a_{\mathfrak g}'}_{L^\infty(I)}
            }{
                2k^2
            }
            \left(
                L_k+\frac12
            \right).
        \end{aligned}
    \]
    The estimates on \([0,\tau_0^{(k)}]\) and \([\tau_0^{(k)},\pi/2]\), together with
    \eqref{eq:app-g-side-log-estimate-constants}, give \eqref{eq:app-ag-logcos-coeff-bound}.
    To prove \eqref{eq:app-ag-logcos-dct-alias-bound}, set
    \[
        \Delta_{N,k}
        \coloneqq
        \mathcal D_{N,k}[a_{\mathfrak g,\log}]
        -
        \mathcal C_k[a_{\mathfrak g,\log}].
    \]
    For \(0\leq k\leq N-1\) and \(r\geq1\),
    \[
        2rN+k\geq2rN,
        \qquad
        2rN-k\geq(2r-1)N,
    \]
    and
    \[
        \log(2rN+k)
        \leq
        \log N+\log(4r),
        \qquad
        \log(2rN-k)
        \leq
        \log N+\log(2r).
    \]
    The aliasing identities \eqref{eq:app-dct1-aliasing-zero} and \eqref{eq:app-dct1-aliasing-pos}, together with the coefficient bound \eqref{eq:app-ag-logcos-coeff-bound}, give
    \[
        \begin{aligned}
            |\Delta_{N,k}|
            &\leq
            \frac1{N^2}
            \sum_{r=1}^{\infty}
            \frac{
                \mathsf C_{\mathrm{coef},\log}
                \bigl(
                    \log N+\log(4r)
                \bigr)
                +
                \mathsf C_{\mathrm{coef},0}
            }{
                (2r)^2
            }
            \\
            &\quad+
            \frac1{N^2}
            \sum_{r=1}^{\infty}
            \frac{
                \mathsf C_{\mathrm{coef},\log}
                \bigl(
                    \log N+\log(2r)
                \bigr)
                +
                \mathsf C_{\mathrm{coef},0}
            }{
                (2r-1)^2
            }.
        \end{aligned}
    \]
    Using
    \[
        \sum_{r\geq1}\frac1{(2r)^2}
        =
        \frac{\pi^2}{24},
        \quad
        \sum_{r\geq1}\frac1{(2r-1)^2}
        =
        \frac{\pi^2}{8},
        \quad
        \sum_{r\geq1}
            \frac{\log r}{(2r)^2}
        =
        \frac14\Sigma_{\log,2},
        \quad
        \sum_{r\geq1}
            \frac{\log r}{(2r-1)^2}
        \leq
        \Sigma_{\log,2},
    \]
    we obtain
    \[
        \begin{aligned}
            |\Delta_{N,k}|
            &\leq
            \frac1{N^2}
            \bigg[
                \frac{\pi^2}{6}
                \mathsf C_{\mathrm{coef},\log}
                \log N
                +
                \mathsf C_{\mathrm{coef},\log}
                \left(
                    \frac{5\pi^2}{24}\log2
                    +
                    \frac54\Sigma_{\log,2}
                \right)
                +
                \frac{\pi^2}{6}
                \mathsf C_{\mathrm{coef},0}
            \bigg]
            \\
            &=
            \frac{
                \mathsf C_{\mathrm{alias},\log}\log N
                +
                \mathsf C_{\mathrm{alias},0}
            }{
                N^2
            }.
        \end{aligned}
    \]
    The last equality follows from
    \eqref{eq:app-g-side-log-estimate-constants}, and hence
    \eqref{eq:app-ag-logcos-dct-alias-bound} follows.
\end{proof}

\subsection{DCT-I approximation estimates}

\begin{definition}[Constants for the DCT-I approximation bounds]
    Let \(N\geq\max\{2,3L_Z\}\). For \(m=0,1,2\) and \(j=1,2\), set
    \begin{equation}\label{eq:app-g-side-Smj-def}
        \mathsf Z_{\ell^1,m}^{(j)}
        \coloneqq
        \sum_{k=0}^{L_Z}
            |z_{m,k}^{(j)}|.
    \end{equation}
    Also set
    \begin{equation*}
        \begin{aligned}
            \mathsf P_{\log,\mathrm{tail}}(\tau)
            \coloneqq{}&
            \frac{
                \mathsf C_{\mathrm{coef},\log}^2
            }{3}
            \tau^2
            +
            \left(
                \frac{
                    2\mathsf C_{\mathrm{coef},\log}^2
                }{9}
                +
                \frac{
                    2\mathsf C_{\mathrm{coef},\log}
                    \mathsf C_{\mathrm{coef},0}
                }{3}
            \right)\tau
            \\
            &+
            \frac{
                2\mathsf C_{\mathrm{coef},\log}^2
            }{27}
            +
            \frac{
                2\mathsf C_{\mathrm{coef},\log}
                \mathsf C_{\mathrm{coef},0}
            }{9}
            +
            \frac{
                \mathsf C_{\mathrm{coef},0}^2
            }{3}.
        \end{aligned}
    \end{equation*}
    Define
    \begin{align}
        \delta_{a_{\mathfrak g},\mathrm{DCT}}(N)^2
        &\coloneqq
        \mathsf C_{a_{\mathfrak g},\mathrm{coef\text{-}tail}}^2
        \left[
            \frac{\pi^5}{72}
            \frac{N+1}{N^4}
            +
            \frac{\pi}{6}
            \frac1{(N-1)^3}
        \right],
        \label{eq:app-g-side-delta-ag-DCT-def}
        \\
        \delta_{\mathcal K_{\sin},\mathrm{DCT}}^{(j,m)}(N)^2
        &\coloneqq
        \left(
            \mathsf C_{a_{\mathfrak g},\mathrm{coef\text{-}tail}}
            \mathsf Z_{\ell^1,m}^{(j)}
        \right)^2
        \left[
            \frac{\pi^7}{36}
            \left(
                (\log2)^2
                +
                \frac{\pi^2}{48}
            \right)
            \frac1{(N-L_Z)^4}
            +
            \frac{2\pi^3}{5}
            \frac1{(N-1)^5}
        \right],
        \label{eq:app-g-side-delta-Ksin-DCT-def}
        \\
        \delta_{a_{\mathfrak g,\log},\mathrm{DCT}}(N)^2
        &\coloneqq
        \frac{\pi}{2}(N+1)
        \frac{
            \left(
                \mathsf C_{\mathrm{alias},\log}\log N
                +
                \mathsf C_{\mathrm{alias},0}
            \right)^2
        }{
            N^4
        }
        +
        \frac{\pi}{2(N-1)^3}
        \mathsf P_{\log,\mathrm{tail}}
        \bigl(
            \log(N-1)
        \bigr),
        \label{eq:app-g-side-delta-ag-log-DCT-def}
        \\
        \delta_{\mathcal I_0,\mathrm{DCT}}^{(j,m)}(N)
        &\coloneqq
        \frac{\pi^3}{6}
        \frac{
            \mathsf C_{a_{\mathfrak g},\mathrm{coef\text{-}tail}}
            \mathsf Z_{\ell^1,m}^{(j)}
        }{
            (N-L_Z)^2
        },
        \nonumber
        \\
        \delta_{\mathcal I_{\cos},\mathrm{DCT}}^{(j,m)}(N)
        &\coloneqq
        \left(
            \pi\log2
            +
            \frac{\pi}{2}
            \sum_{k=1}^{N-1}\frac1k
        \right)
        \frac{\pi^2}{6}
        \frac{
            \mathsf C_{a_{\mathfrak g},\mathrm{coef\text{-}tail}}
            \mathsf Z_{\ell^1,m}^{(j)}
        }{
            (N-L_Z)^2
        }
        +
        \frac{\pi}{2}
        \frac{
            \mathsf C_{a_{\mathfrak g},\mathrm{coef\text{-}tail}}
            \mathsf Z_{\ell^1,m}^{(j)}
        }{
            N(N-L_Z-1)
        }.
        \label{eq:app-g-side-delta-Icos-DCT-def}
    \end{align}
\end{definition}

\begin{lemma}[DCT-I approximation bounds]
    Let \(N\geq\max\{2,3L_Z\}\). Then, for \(m=0,1,2\) and \(j=1,2\),
    \begin{align}
        \norm{
            a_{\mathfrak g}
            -
            \mathcal R_Na_{\mathfrak g}
        }_{L^2(I)}
        &\leq
        \delta_{a_{\mathfrak g},\mathrm{DCT}}(N),
        \label{eq:app-g-side-ag-DCT-bound}
        \\
        \norm{
            a_{\mathfrak g,\log}
            -
            \mathcal R_N(a_{\mathfrak g,\log})
        }_{L^2(I)}
        &\leq
        \delta_{a_{\mathfrak g,\log},\mathrm{DCT}}(N),
        \label{eq:app-g-side-ag-log-DCT-bound}
        \\
        \norm{
            \mathcal K_{\sin}[a_{\mathfrak g}Z_m^{(j)}]
            -
            \mathcal K_{\sin}
            \left[
                \mathcal R_N(a_{\mathfrak g}Z_m^{(j)})
            \right]
        }_{L^2(I)}
        &\leq
        \delta_{\mathcal K_{\sin},\mathrm{DCT}}^{(j,m)}(N),
        \label{eq:app-g-side-Ksin-DCT-bound}
        \\
        \left|
            \mathcal I_0[a_{\mathfrak g}Z_m^{(j)}]
            -
            \mathcal I_0
            \left[
                \mathcal R_N(a_{\mathfrak g}Z_m^{(j)})
            \right]
        \right|
        &\leq
        \delta_{\mathcal I_0,\mathrm{DCT}}^{(j,m)}(N),
        \label{eq:app-g-side-I0-DCT-bound}
        \\
        \left|
            \mathcal I_{\cos}[a_{\mathfrak g}Z_m^{(j)}]
            -
            \mathcal I_{\cos}
            \left[
                \mathcal R_N(a_{\mathfrak g}Z_m^{(j)})
            \right]
        \right|
        &\leq
        \delta_{\mathcal I_{\cos},\mathrm{DCT}}^{(j,m)}(N).
        \label{eq:app-g-side-Icos-DCT-bound}
    \end{align}
\end{lemma}

\begin{proof}
    For \(0\leq k\leq N-1\), set
    \[
        \delta_{N,k}[f]
        \coloneqq
        \mathcal D_{N,k}[f]-\mathcal C_k[f],
        \qquad
        0\leq k\leq N-1.
    \]
    By \eqref{eq:ag-cosine-tail-bound},
    \begin{equation}\label{eq:app-g-side-ag-cosine-tail-used}
        |\mathcal C_n[a_{\mathfrak g}]|
        \leq
        \frac{\mathsf C_{a_{\mathfrak g},\mathrm{coef\text{-}tail}}}{n^2},
        \qquad n\geq1.
    \end{equation}
    Set \(h_m^{(j)}\coloneqq a_{\mathfrak g}Z_m^{(j)}\). The product formula for cosine coefficients gives, for \(n\geq L_Z+1\),
    \[
    \begin{aligned}
        |\mathcal C_n[h_m^{(j)}]|
        &\leq
        |z_{m,0}^{(j)}|
        |\mathcal C_n[a_{\mathfrak g}]|
        +
        \frac12
        \sum_{q=1}^{L_Z}
        |z_{m,q}^{(j)}|
        \left(
            |\mathcal C_{n-q}[a_{\mathfrak g}]|
            +
            |\mathcal C_{n+q}[a_{\mathfrak g}]|
        \right)                                      \\
        &\leq
        \frac{
            \mathsf C_{a_{\mathfrak g},\mathrm{coef\text{-}tail}}
        }{
            (n-L_Z)^2
        }
        \sum_{q=0}^{L_Z}|z_{m,q}^{(j)}|
        =
        \frac{
            \mathsf C_{a_{\mathfrak g},\mathrm{coef\text{-}tail}}\mathsf Z_{\ell^1,m}^{(j)}
        }{
            (n-L_Z)^2
        }.
    \end{aligned}
    \]
    Hence, by Lemma~\ref{lem:app-alias-bound-n2}, for \(0\leq k\leq N-1\),
    \begin{equation}\label{eq:app-g-side-proof-alias-a-hjm}
        |\delta_{N,k}[a_{\mathfrak g}]|
        \leq
        \frac{\pi^2}{6}
        \frac{\mathsf C_{a_{\mathfrak g},\mathrm{coef\text{-}tail}}}{N^2},
        \qquad
        |\delta_{N,k}[h_m^{(j)}]|
        \leq
        \frac{\pi^2}{6}
        \frac{
            \mathsf C_{a_{\mathfrak g},\mathrm{coef\text{-}tail}}\mathsf Z_{\ell^1,m}^{(j)}
        }{
            (N-L_Z)^2
        }.
    \end{equation}
    By cosine orthogonality on \(I\),
    \[
    \begin{aligned}
        \norm{
            a_{\mathfrak g}
            -
            \mathcal R_Na_{\mathfrak g}
        }_{L^2(I)}^2
        &=
        \pi|\delta_{N,0}[a_{\mathfrak g}]|^2
        +
        \frac{\pi}{2}
        \sum_{k=1}^{N-1}
        |\delta_{N,k}[a_{\mathfrak g}]|^2
        +
        \frac{\pi}{2}
        \sum_{k=N}^{\infty}
        |\mathcal C_k[a_{\mathfrak g}]|^2.
    \end{aligned}
    \]
    Using \eqref{eq:app-g-side-proof-alias-a-hjm} for the first two terms and
    \eqref{eq:app-g-side-ag-cosine-tail-used} for the tail gives
    \[
    \begin{aligned}
        \norm{
            a_{\mathfrak g}
            -
            \mathcal R_Na_{\mathfrak g}
        }_{L^2(I)}^2
        &\leq
        \frac{\pi}{2}(N+1)
        \left(
            \frac{\pi^2}{6}
            \frac{\mathsf C_{a_{\mathfrak g},\mathrm{coef\text{-}tail}}}{N^2}
        \right)^2
        +
        \frac{\pi}{2}
        \mathsf C_{a_{\mathfrak g},\mathrm{coef\text{-}tail}}^2
        \sum_{k=N}^{\infty}\frac1{k^4} \\
        &\leq
        \mathsf C_{a_{\mathfrak g},\mathrm{coef\text{-}tail}}^2
        \left[
            \frac{\pi^5}{72}
            \frac{N+1}{N^4}
            +
            \frac{\pi}{6}
            \frac1{(N-1)^3}
        \right]
        =
        \delta_{a_{\mathfrak g},\mathrm{DCT}}(N)^2.
    \end{aligned}
    \]
    This proves \eqref{eq:app-g-side-ag-DCT-bound}.

    The operator \(\mathcal K_{\sin}\) is diagonal in the even cosine basis,
    with eigenvalues
    \[
        \kappa_0=\pi\log2,
        \qquad
        \kappa_k=\frac{\pi}{2k},
        \qquad k\geq1.
    \]
    Therefore, with \(h=h_m^{(j)}\), cosine orthogonality gives
    \[
    \begin{aligned}
        \norm{
            \mathcal K_{\sin}[h]
            -
            \mathcal K_{\sin}[\mathcal R_Nh]
        }_{L^2(I)}^2
        &=
        \pi|\kappa_0\delta_{N,0}[h]|^2
        +
        \frac{\pi}{2}
        \sum_{k=1}^{N-1}
        |\kappa_k\delta_{N,k}[h]|^2  
        +
        \frac{\pi}{2}
        \sum_{k=N}^{\infty}
        |\kappa_k\mathcal C_k[h]|^2.
    \end{aligned}
    \]
    The low-mode aliasing contribution is bounded by
    \[
    \begin{aligned}
        \pi|\kappa_0\delta_{N,0}[h]|^2
        +
        \frac{\pi}{2}
        \sum_{k=1}^{N-1}
        |\kappa_k\delta_{N,k}[h]|^2
        &\leq
        \pi^3
        \left(
            (\log2)^2+\frac{\pi^2}{48}
        \right)
        \left(
            \frac{\pi^2}{6}
            \frac{
                \mathsf C_{a_{\mathfrak g},\mathrm{coef\text{-}tail}}\mathsf Z_{\ell^1,m}^{(j)}
            }{
                (N-L_Z)^2
            }
        \right)^2                                      \\
        &=
        \left(
            \mathsf C_{a_{\mathfrak g},\mathrm{coef\text{-}tail}}\mathsf Z_{\ell^1,m}^{(j)}
        \right)^2
        \frac{\pi^7}{36}
        \left(
            (\log2)^2+\frac{\pi^2}{48}
        \right)
        \frac1{(N-L_Z)^4}.
    \end{aligned}
    \]
    Here we used \(\sum_{k=1}^{N-1}k^{-2}\leq \pi^2/6\). For the high-mode tail, since \(N\geq 3L_Z\), we have \(k-L_Z\geq k/2\) for \(k\geq N\), and hence
    \[
        |\kappa_k\mathcal C_k[h]|
        \leq
        \kappa_k\frac{
            \mathsf C_{a_{\mathfrak g},\mathrm{coef\text{-}tail}}\mathsf Z_{\ell^1,m}^{(j)}
        }{
            (k-L_Z)^2
        }
        \leq
        \frac{
            2\pi \mathsf C_{a_{\mathfrak g},\mathrm{coef\text{-}tail}}\mathsf Z_{\ell^1,m}^{(j)}
        }{
            k^3
        }.
    \]
    Thus
    \[
    \begin{aligned}
        \frac{\pi}{2}
        \sum_{k=N}^{\infty}
        |\kappa_k\mathcal C_k[h]|^2
        &\leq
        2\pi^3
        \left(
            \mathsf C_{a_{\mathfrak g},\mathrm{coef\text{-}tail}}\mathsf Z_{\ell^1,m}^{(j)}
        \right)^2
        \sum_{k=N}^{\infty}\frac1{k^6} 
        \leq
        \frac{2\pi^3}{5}
        \left(
            \mathsf C_{a_{\mathfrak g},\mathrm{coef\text{-}tail}}\mathsf Z_{\ell^1,m}^{(j)}
        \right)^2
        \frac1{(N-1)^5}.
    \end{aligned}
    \]
    The low-mode estimate obtained from
    \eqref{eq:app-g-side-proof-alias-a-hjm}, the high-mode estimate above,
    and \eqref{eq:app-g-side-delta-Ksin-DCT-def} give
    \[
        \norm{
            \mathcal K_{\sin}[h_m^{(j)}]
            -
            \mathcal K_{\sin}
            \left[
                \mathcal R_Nh_m^{(j)}
            \right]
        }_{L^2(I)}^2
        \leq
        \delta_{\mathcal K_{\sin},\mathrm{DCT}}^{(j,m)}(N)^2.
    \]
    Taking square roots gives \eqref{eq:app-g-side-Ksin-DCT-bound}.

    For \(a_{\mathfrak g,\log}\), Lemma~\ref{lem:app-ag-logcos-coeff-alias} gives 
    \[
        |\mathcal C_k[a_{\mathfrak g,\log}]|
        \leq
        \frac{\mathsf C_{\mathrm{coef},\log}\log k+\mathsf C_{\mathrm{coef},0}}{k^2},\quad 
        |\delta_{N,k}[a_{\mathfrak g,\log}]|
        \leq
        \frac{
            \mathsf C_{\mathrm{alias},\log}\log N+\mathsf C_{\mathrm{alias},0}
        }{
            N^2
        },
    \]
    for \(k\geq2\) and \(0\leq k\leq N-1\) respectively. Thus cosine orthogonality yields
    \[
        \begin{aligned}
            \norm{
                a_{\mathfrak g,\log}
                -
                \mathcal R_N(a_{\mathfrak g,\log})
            }_{L^2(I)}^2
            &\leq
            \frac{\pi}{2}(N+1)
            \frac{
                (\mathsf C_{\mathrm{alias},\log}\log N+\mathsf C_{\mathrm{alias},0})^2
            }{
                N^4
            }
            \\
            &\quad+
            \frac{\pi}{2}
            \sum_{k=N}^{\infty}
            \frac{
                (\mathsf C_{\mathrm{coef},\log}\log k+\mathsf C_{\mathrm{coef},0})^2
            }{
                k^4
            }.
        \end{aligned}
    \]
    The remaining tail is estimated by
    \[
        \begin{aligned}
            \sum_{k=N}^{\infty}
            \frac{
                (\mathsf C_{\mathrm{coef},\log}\log k+\mathsf C_{\mathrm{coef},0})^2
            }{
                k^4
            }
            &\leq
            \int_{N-1}^{\infty}
            \frac{
                (\mathsf C_{\mathrm{coef},\log}\log x+\mathsf C_{\mathrm{coef},0})^2
            }{
                x^4
            }dx
            \\
            &=
            \frac1{(N-1)^3}
            \mathsf P_{\log,\mathrm{tail}}(\log(N-1)).
        \end{aligned}
    \]
    The cosine-orthogonality estimate above, the tail estimate above, and
    \eqref{eq:app-g-side-delta-ag-log-DCT-def} give
    \[
        \norm{
            a_{\mathfrak g,\log}
            -
            \mathcal R_N(a_{\mathfrak g,\log})
        }_{L^2(I)}^2
        \leq
        \delta_{a_{\mathfrak g,\log},\mathrm{DCT}}(N)^2.
    \]
    Taking square roots gives
    \eqref{eq:app-g-side-ag-log-DCT-bound}.

    The estimate for \(\mathcal I_0\) follows from the zero-mode identity
    \[
        \mathcal I_0[h_m^{(j)}]
        =
        \pi\mathcal C_0[h_m^{(j)}],
        \qquad
        \mathcal I_0[\mathcal R_Nh_m^{(j)}]
        =
        \pi\mathcal D_{N,0}[h_m^{(j)}].
    \]
    Hence, by \eqref{eq:app-g-side-proof-alias-a-hjm},
    \[
        \left|
            \mathcal I_0[h_m^{(j)}]
            -
            \mathcal I_0[\mathcal R_Nh_m^{(j)}]
        \right|
        =
        \pi|\delta_{N,0}[h_m^{(j)}]|
        \leq
        \frac{\pi^3}{6}
        \frac{
            \mathsf C_{a_{\mathfrak g},\mathrm{coef\text{-}tail}}\mathsf Z_{\ell^1,m}^{(j)}
        }{
            (N-L_Z)^2
        }
        =
        \delta_{\mathcal I_0,\mathrm{DCT}}^{(j,m)}(N).
    \]
    This proves \eqref{eq:app-g-side-I0-DCT-bound}.

    Finally, the Fourier series
    \[
        \log|\cos\theta|
        =
        -\log2
        -
        \sum_{k=1}^{\infty}
            \frac{(-1)^k}{k}
            \cos(2k\theta)
    \]
    gives, for any even \(h\),
    \[
        \mathcal I_{\cos}[h]
        =
        -\pi\log2\,\mathcal C_0[h]
        -
        \frac{\pi}{2}
        \sum_{k=1}^{\infty}
        \frac{(-1)^k}{k}\mathcal C_k[h].
    \]
    Applying this identity to \(h=h_m^{(j)}\) and to \(\mathcal R_Nh\), we obtain
    \[
    \begin{aligned}
        \left|
            \mathcal I_{\cos}[h]
            -
            \mathcal I_{\cos}[\mathcal R_Nh]
        \right|                                                  
        \leq
        \pi\log2\,|\delta_{N,0}[h]|
        +
        \frac{\pi}{2}
        \sum_{k=1}^{N-1}
        \frac{
            |\delta_{N,k}[h]|
        }{
            k
        }
        +
        \frac{\pi}{2}
        \sum_{k=N}^{\infty}
        \frac{
            |\mathcal C_k[h]|
        }{
            k
        }.
    \end{aligned}
    \]
    By \eqref{eq:app-g-side-proof-alias-a-hjm}, the first two terms are bounded by
    \[
        \left(
            \pi\log2
            +
            \frac{\pi}{2}
            \sum_{k=1}^{N-1}\frac1k
        \right)
        \frac{\pi^2}{6}
        \frac{
            \mathsf C_{a_{\mathfrak g},\mathrm{coef\text{-}tail}}\mathsf Z_{\ell^1,m}^{(j)}
        }{
            (N-L_Z)^2
        }.
    \]
    For the remaining tail,
    \[
    \begin{aligned}
        \frac{\pi}{2}
        \sum_{k=N}^{\infty}
        \frac{
            |\mathcal C_k[h]|
        }{
            k
        }
        &\leq
        \frac{\pi}{2}
        \mathsf C_{a_{\mathfrak g},\mathrm{coef\text{-}tail}}\mathsf Z_{\ell^1,m}^{(j)}
        \sum_{k=N}^{\infty}
        \frac1{k(k-L_Z)^2}                            
        \leq
        \frac{\pi}{2}
        \frac{
            \mathsf C_{a_{\mathfrak g},\mathrm{coef\text{-}tail}}\mathsf Z_{\ell^1,m}^{(j)}
        }{
            N(N-L_Z-1)
        }.
    \end{aligned}
    \]
    The preceding two estimates and
    \eqref{eq:app-g-side-delta-Icos-DCT-def} give
    \eqref{eq:app-g-side-Icos-DCT-bound}.
\end{proof}

\subsection{Comparison of the right-hand sides and residuals}

\begin{definition}[Constants for the right-hand-side and residual estimates]
    Let \(N\geq\max\{2,3L_Z\}\). For \(m=0\), set
    \begin{equation}\label{eq:app-g-side-eta-rhs0-def}
        \delta_{\mathrm{rhs},0}^{(j)}(N)
        \coloneqq
        \delta_{a_{\mathfrak g},\mathrm{DCT}}(N),
        \qquad
        j=1,2.
    \end{equation}
    For \(m=1,2\), set
    \begin{equation}\label{eq:app-g-side-eta-rhsm-def}
        \begin{aligned}
            \delta_{\mathrm{rhs},m}^{(j)}(N)
            \coloneqq
            \frac1\pi
            \bigg[
                &
                \delta_{a_{\mathfrak g},\mathrm{DCT}}(N)
                \norm{
                    \mathcal K_{\sin}[\Omega_m^{(j)}]
                }_{L^\infty(I)}
                +
                \delta_{a_{\mathfrak g,\log},\mathrm{DCT}}(N)
                \left|
                    \mathcal I_0[\Omega_m^{(j)}]
                \right|
                \\
                &+
                \delta_{a_{\mathfrak g},\mathrm{DCT}}(N)
                \left|
                    \mathcal I_{\cos}[\Omega_m^{(j)}]
                \right|
            \bigg].
        \end{aligned}
    \end{equation}
    For \(m=0,1,2\), set
    \begin{equation}\label{eq:app-g-side-eta-res-def}
        \begin{aligned}
            \delta_{\mathrm{res\text{-}diff},m}^{(j)}(N)
            \coloneqq{}&
            \frac{c_j}{\pi}
            \bigg[
                \delta_{a_{\mathfrak g},\mathrm{DCT}}(N)
                \norm{
                    \mathcal K_{\sin}
                    \left[
                        \mathcal R_N
                        (a_{\mathfrak g}Z_m^{(j)})
                    \right]
                }_{L^\infty(I)}
                \\
                &\quad+
                \norm{a_{\mathfrak g}}_{L^\infty(I)}
                \delta_{\mathcal K_{\sin},\mathrm{DCT}}^{(j,m)}(N)
                \\
                &\quad+
                \delta_{a_{\mathfrak g,\log},\mathrm{DCT}}(N)
                \left|
                    \mathcal I_0
                    \left[
                        \mathcal R_N
                        (a_{\mathfrak g}Z_m^{(j)})
                    \right]
                \right|
                \\
                &\quad+
                \norm{a_{\mathfrak g,\log}}_{L^2(I)}
                \delta_{\mathcal I_0,\mathrm{DCT}}^{(j,m)}(N)
                \\
                &\quad+
                \delta_{a_{\mathfrak g},\mathrm{DCT}}(N)
                \left|
                    \mathcal I_{\cos}
                    \left[
                        \mathcal R_N
                        (a_{\mathfrak g}Z_m^{(j)})
                    \right]
                \right|
                \\
                &\quad+
                \norm{a_{\mathfrak g}}_{L^2(I)}
                \delta_{\mathcal I_{\cos},\mathrm{DCT}}^{(j,m)}(N)
            \bigg]
            +
            \delta_{\mathrm{rhs},m}^{(j)}(N).
        \end{aligned}
    \end{equation}
\end{definition}

\begin{lemma}[DCT-I comparison of the right-hand sides and residuals]
\label{lem:app-g-side-residual-comparison}
    Let \(N\geq\max\{2,3L_Z\}\). Then, for \(m=0,1,2\) and \(j=1,2\),
    \begin{align}
        \norm{
            B_m^{(j)}-B_{m,N}^{(j)}
        }_{L^2(I)}
        &\leq
        \delta_{\mathrm{rhs},m}^{(j)}(N),
        \label{eq:app-g-side-B-minus-BN-bound}
        \\
        \norm{
            \operatorname{Res}_m^{(j)}
            -
            \operatorname{Res}_{m,N}^{(j)}
        }_{L^2(I)}
        &\leq
        \delta_{\mathrm{res\text{-}diff},m}^{(j)}(N).
        \label{eq:app-g-side-Res-minus-ResN-bound}
    \end{align}
\end{lemma}

\begin{proof}
    For \(m=0\), the definitions \eqref{eq:g-side-B-def} and \eqref{eq:g-side-discrete-rhs-def} give
    \[
        B_0^{(j)}-B_{0,N}^{(j)}
        =
        a_{\mathfrak g}
        -
        \mathcal R_Na_{\mathfrak g}.
    \]
    Hence
    \eqref{eq:app-g-side-ag-DCT-bound} and
    \eqref{eq:app-g-side-eta-rhs0-def} give
    \eqref{eq:app-g-side-B-minus-BN-bound} for \(m=0\).

    For \(m=1,2\), the definitions \eqref{eq:K_log-definition}, \eqref{eq:g-side-B-def}, and \eqref{eq:g-side-discrete-rhs-def} give
    \[
        \begin{aligned}
            B_m^{(j)}-B_{m,N}^{(j)}
            =
            \frac1\pi
            \bigg[
                &
                \left(
                    a_{\mathfrak g}
                    -
                    \mathcal R_Na_{\mathfrak g}
                \right)
                \mathcal K_{\sin}[\Omega_m^{(j)}]
                +
                \left(
                    a_{\mathfrak g,\log}
                    -
                    \mathcal R_Na_{\mathfrak g,\log}
                \right)
                \mathcal I_0[\Omega_m^{(j)}]
                \\
                &+
                \left(
                    a_{\mathfrak g}
                    -
                    \mathcal R_Na_{\mathfrak g}
                \right)
                \mathcal I_{\cos}[\Omega_m^{(j)}]
            \bigg].
        \end{aligned}
    \]
    The Cauchy--Schwarz inequality,
    \eqref{eq:app-g-side-ag-DCT-bound},
    \eqref{eq:app-g-side-ag-log-DCT-bound}, and
    \eqref{eq:app-g-side-eta-rhsm-def} give
    \eqref{eq:app-g-side-B-minus-BN-bound} for \(m=1,2\).

    Set
    \(
        f_m^{(j)}
        \coloneqq
        a_{\mathfrak g}Z_m^{(j)}.
    \)
    The definitions \eqref{eq:g-side-continuous-residual-def} and \eqref{eq:g-side-discrete-residual-def} give
    \[
        \operatorname{Res}_m^{(j)}
        -
        \operatorname{Res}_{m,N}^{(j)}
        =
        B_m^{(j)}
        -
        B_{m,N}^{(j)}
        +
        \frac{c_j}{\pi}
        \left(
            a_{\mathfrak g}
            \mathcal K_{\log}[f_m^{(j)}]
            -
            \mathscr A_{m,N}^{(j)}
        \right).
    \]
    By
    \eqref{eq:K_log-definition} and
    \eqref{eq:g-side-discrete-action-def},
    \[
        \begin{aligned}
            a_{\mathfrak g}
            \mathcal K_{\log}[f_m^{(j)}]
            -
            \mathscr A_{m,N}^{(j)}
            &=
            \left(
                a_{\mathfrak g}
                -
                \mathcal R_Na_{\mathfrak g}
            \right)
            \mathcal K_{\sin}
            \left[
                \mathcal R_Nf_m^{(j)}
            \right]
            +
            a_{\mathfrak g}
            \left\{
                \mathcal K_{\sin}[f_m^{(j)}]
                -
                \mathcal K_{\sin}
                \left[
                    \mathcal R_Nf_m^{(j)}
                \right]
            \right\}
            \\
            &\quad+
            \left(
                a_{\mathfrak g,\log}
                -
                \mathcal R_Na_{\mathfrak g,\log}
            \right)
            \mathcal I_0
            \left[
                \mathcal R_Nf_m^{(j)}
            \right]
            +
            a_{\mathfrak g,\log}
            \left\{
                \mathcal I_0[f_m^{(j)}]
                -
                \mathcal I_0
                \left[
                    \mathcal R_Nf_m^{(j)}
                \right]
            \right\}
            \\
            &\quad+
            \left(
                a_{\mathfrak g}
                -
                \mathcal R_Na_{\mathfrak g}
            \right)
            \mathcal I_{\cos}
            \left[
                \mathcal R_Nf_m^{(j)}
            \right]
            +
            a_{\mathfrak g}
            \left\{
                \mathcal I_{\cos}[f_m^{(j)}]
                -
                \mathcal I_{\cos}
                \left[
                    \mathcal R_Nf_m^{(j)}
                \right]
            \right\}.
        \end{aligned}
    \]
    Applying
    \eqref{eq:app-g-side-ag-DCT-bound},
    \eqref{eq:app-g-side-ag-log-DCT-bound},
    \eqref{eq:app-g-side-Ksin-DCT-bound},
    \eqref{eq:app-g-side-I0-DCT-bound}, and
    \eqref{eq:app-g-side-Icos-DCT-bound} to these six terms, and then
    using
    \eqref{eq:app-g-side-B-minus-BN-bound} and
    \eqref{eq:app-g-side-eta-res-def}, gives
    \eqref{eq:app-g-side-Res-minus-ResN-bound}.
\end{proof}

 \bibliographystyle{abbrv}
\bibliography{reference}

\end{document}